\documentclass[12pt]{article}
\usepackage{amsmath,graphicx,amssymb,amsfonts,amsthm}
\usepackage{verbatim,enumitem,lipsum,float}
\usepackage{color,epstopdf,epsfig,lscape,mathabx,natbib,longtable,setspace}
\usepackage{xr}

\newcommand{\indep}{\;\, \rule[0em]{.03em}{.67em} \hspace{-.25em}
\rule[0em]{.65em}{.03em}
\hspace{-.25em}\rule[0em]{.03em}{.67em}\;\,}

\newtheorem{Th}{\underline{\bf Theorem}}

\newtheorem{Pro}{Proposition}
\newtheorem{Lem}{\underline{\bf Lemma}}

\def\bse{\begin{eqnarray*}}
\def\ese{\end{eqnarray*}}
\def\be{\begin{eqnarray}}
\def\ee{\end{eqnarray}}
\def\bsq{\begin{equation*}}
\def\esq{\end{equation*}}
\def\bq{\begin{equation}}
\def\eq{\end{equation}}

\def\var{\hbox{var}}

\def\wh{\widehat}
\def\wt{\widetilde}
\def\eff{_{\rm eff}}

\def\n{\nonumber}

\def\vecl{\mbox{vecl}}

\def\sumi{\sum_{i=1}^n}
\def\sumj{\sum_{j=1}^n}
\def\trans{^{\rm T}}

\def\bb{\boldsymbol\beta}

\def\bg{\boldsymbol\gamma}

\def\0{{\bf 0}}
\def\A{{\bf A}}
\def\a{{\bf a}}
\def\B{{\bf B}}

\def\c{{\bf c}}

\def\g{{\bf g}}

\def\h{{\bf h}}

\def\I{{\bf I}}

\def\K{{\bf K}}

\def\T{{\bf T}}

\def\bS{{\bf S}}

\def\m{{\bf m}}
\def\v{{\bf v}}

\def\T{{\bf T}}

\def\X{{\bf X}}
\def\x{{\bf x}}
\def\I{{\bf I}}

\def\z{{\bf z}}

\def\bLam{{\boldsymbol \Lambda}}
\def\blam{{\boldsymbol \lambda}}
\def\bkappa{\boldsymbol\kappa}

\def\brho{\boldsymbol\rho}

\def\bq{\begin{equation}}
\def\eq{\end{equation}}
\def\pr{\hbox{pr}}
\def\wh{\widehat}
\def\wt{\widetilde}
\def\trans{^{\rm T}}

\def\log{\hbox{log}}

\def\squarebox#1{\hbox to #1{\hfill\vbox to #1{\vfill}}}

\def\var{\hbox{var}}

\def\bse{\begin{eqnarray*}}
\def\ese{\end{eqnarray*}}
\def\be{\begin{eqnarray}}
\def\ee{\end{eqnarray}}
\def\bsq{\begin{equation*}}
\def\esq{\end{equation*}}
\def\bq{\begin{equation}}
\def\eq{\end{equation}}
\def\pr{\hbox{pr}}
\def\wh{\widehat}
\def\wt{\widetilde}

\def\log{\hbox{log}}

\def\trans{^{\rm T}}

\def\boxit#1{\vbox{\hrule\hbox{\vrule\kern6pt\vbox{\kern6pt#1\kern6pt}\kern6pt\vrule}\hrule}}

\def\spacingset#1{\renewcommand{\baselinestretch}%
	{#1}\small\normalsize} \spacingset{1}

\begin{document}
\thispagestyle{empty}
\baselineskip 17pt
\renewcommand {\thepage}{}
\pagenumbering{arabic}
\begin{center}
{\Large\bf Semiparametric regression of mean residual life  with censoring and covariate dimension reduction}

\vspace{.3cm}
\centerline{\bf Ge Zhao}
\vspace{.1cm}
\centerline{\it Department of Mathematics and Statistics}
\centerline{\it Portland State University, Portland, OR 97201}
\vspace{.3cm}
\centerline{\bf Yanyuan Ma}
\vspace{.1cm}
\centerline{\it Department of Statistics, Pennsylvania State University, State College, PA 16802}
\vspace{.3cm}
\centerline{\bf Huazhen Lin}
\vspace{.1cm}
\centerline{\it Center of Statistical Research,  School of Statistics}
\centerline{\it Southwestern University of
Finance and Economics, Chengdu, China 611130}
\vspace{.3cm}
\centerline{and {\bf Yi Li}}
\vspace{.1cm}
\centerline{\it Department of Biostatistics, University of Michigan,
Ann Arbor, MI 48109}\vspace{.55cm}
\fontsize{9}{12pt plus.6pt minus .6pt}\selectfont

\end{center}

\begin{abstract}
	We propose a new class of semiparametric regression models of  mean residual life for censored outcome data. The models, which enable us to  estimate the expected remaining survival time and generalize commonly used mean residual life models, also conduct covariate dimension reduction. Using the geometric approaches in semiparametrics literature and the martingale properties with survival data, we propose a  flexible inference procedure that relaxes the parametric assumptions on the dependence of mean residual life on covariates and
	how long a patient has lived. We show that the estimators for the covariate effects are  root-$n$ consistent, asymptotically normal, and  semiparametrically efficient. With the unspecified mean residual life function, we  provide a nonparametric estimator for predicting the  residual life of a given subject, and  establish the root-$n$ consistency and asymptotic normality for this estimator.
	Numerical experiments are conducted to illustrate the feasibility of the proposed estimators. We apply the method to analyze a national kidney transplantation dataset to further demonstrate the utility of the work.
\end{abstract}

{\it Keywords:} Mean residual life, kidney transplant,  nonparametric estimation, semiparametric efficiency.

\doublespacing

\section{Introduction}

End stage renal disease (ESRD) is  one of the most lethal diseases  
in the US \citep{ferri2017}. An ESRD  patient 
loses kidney functions permanently and has to rely on
 renal replacement   to survive. The  common renal 
replacement therapies are dialysis and kidney transplant, where the latter is associated with better outcomes and qualify of life.
Due to severe shortages in kidney supplies,  there are far more renal failure patients in need of kidney
transplants than donors available, and 
the decision for organ allocation  is often based on 
the Estimated Post-Transplant Survival (EPTS) score (\url{https://srtr.transplant.hrsa.gov/}),
which assesses a patient's overall survival 
post transplantation.

EPTS is calculated from a Cox model which includes
age, diabetes status, prior solid organ transplant and time on dialysis as covariates. With this model, however, candidates with  diabetes or a prior solid organ transplant  or being on dialysis
treatment  for many years will have a lower priority for  transplantation. On the other hand,  younger patients tend to
have a higher  priority for  transplantation and age has effectively become the most dominant
factor.  A more efficient system should allocate organs to those who will benefit most  among all the eligible patients \citep{bertsimas2013fairness,israni2014new}. One way to quantify the transplant efficiency is to compare the improvement of  patients'
expected residual life with and without transplantation \citep{wolfe2008calculating}, where the expected residual life characterizes the remaining survival time
given that a patient has survived up to a certain time.  The residual life expectancy is advantageous to the overall life expectancy as  the former  takes into
account the most updated information
\citep{lin2016semiparametrically} once an organ is available.

Since the seminal work of \cite{oakesdasu1990},
much research has been conducted on mean residual life models. For example,
\cite{maguluri1994estimation}  proposed a univariate
proportional mean residual life model and provided estimation
and
inference tools;
\cite{oakes2003inference} established the theoretical properties
of
the methods in \cite{oakesdasu1990};
\cite{chencheng2005}  studied the proportional mean
residual
life model and proposed an alternative estimator
for the covariate effect through a
partial-score construction, similar to the partial likelihood
approach;
\cite{chenetal2005} employed the inverse probability weighting
approach to estimate the covariate effects;
\cite{muller2005time} extended the
mean residual life model to incorporate time-varying covariates;
\cite{chen2006linear} proposed  a linear residual life model, and
\cite{chen2007} proposed an additive mean residual life
model.
These models
  inspired  quantile residual life
model research \citep{jeong2008nonparametric,ma2010semiparametric}.
However, all of 
these works assumed that the
mean or quantile residual life depends on the covariates as well as  how long a patient has lived through parametric functions. Violations of the model assumptions may lead to biased results \citep{mclain2011nonparametric,liu2019nonparametric}, which motivates our work to relax these parametric  assumptions. 

Our work is also motivated by the kidney transplant data from
the U.S. Scientific Registry of Transplant Recipients \citep{leppke2013scientific,israni2014new}, which feature,
in addition to survival outcome data, rich information, such as treatment history, comorbidity conditions and demographic variables. To strike a balance
between interpretation and flexibility when quantifying the  covariates' impacts on the expected residual life so as to identify patients
who may benefit most from transplantation, we propose a semiparametric 
mean residual life regression model with
covariate dimension reduction.  Our model relaxes the parametric assumptions on the dependence of mean residual life on covariates and	how long a patient has lived, as opposed to many
existing residual life models, and in the meantime conducts covariate dimension reduction  via a data driven fashion \citep{ma2013efficient}.
Using the geometric approaches in semiparametrics literature and the martingale properties with survival data, 
we obtain semiparametrically efficient estimators \citep{bickeletal1994,tsiatis2006semiparametric} 
for the covariate effects. With the unspecified mean residual life function, we  provide a nonparametric estimator for predicting the  residual life of a given subject, and  establish the root-$n$ consistency and asymptotic normality for this estimator.
We  analyze the kidney transplant data and quantify
the  benefit with residual life expectancy.

The remainder is organized as follows. In Section \ref{sec:model}, we propose the mean residual life model and introduce the
notation. In Section \ref{sec:estimation}, we derive an
efficient estimator and discuss its
properties.  We also provide the estimation of the mean
residual life function. The large sample properties
of the estimator and the estimated mean residual life function
are established in Section \ref{sec:asymp}. We assess the finite
sample properties
of the methods using simulation studies
in Section \ref{sec:sim} and apply it to analyze the kidney
transplant data in
Section \ref{sec:app}. We conclude the paper with a
discussion in Section \ref{sec:conclusion}.

\section{Semiparametric regression of mean residual life}\label{sec:model}

Denote the failure time and the covariates by $T\in{\cal R}^+$ and
$\X\in{\cal R}^p$, respectively. For any given $t>0$ (e.g. the time length a patient has lived so far), we propose
a semiparametric regression model of mean residual life
\be\label{eq:model}
E(T-t\mid T\ge t,\X)=m(t,\bb\trans\X),
\ee
where $\bb\in{\cal R}^{p\times
d}$ is the coefficient matrix with $d < p$,  $m$ is an
unspecified positive function of $t$ and
$\bb\trans\X$. 
Here $d$ is the number of indices. When $d=1$, the 
model reduces to the single index model; when $ 
1<d<p$, it corresponds to a  dimension reduction
	structure; when $d=p$, the model is completely
	nonparametric. Our analysis first focuses on a fixed $d$, 
	followed by  selecting $d$ with a data driven fashion as 
	discussed in Section~\ref{sec:app}.
Assume further that $T$ is subject to random
right censoring, i.e. $Z=\min(T,C)$ and $\Delta=I(T\le
C)$, where  the censoring time $C$ satisfies $C\indep T\mid
\X$. Model (\ref{eq:model}) generalizes  many existing mean
residual life models. For example, it includes the proportional mean
residual life
model \citep{oakesdasu1990,
	chenetal2005, chencheng2005}
as a special case, by specifying
$m(t,\bb\trans\X)=m_0(t)\exp(\bb\trans\X)$ and $d=1$;
it reduces to the additive model \citep{chen2007}  by specifying $m(t,\bb\trans\X)=m_0(t)+\bb\trans\X$ and setting $d=1$.

We assume that the observed data  $(\X_i,
Z_i,\Delta_i)$, $i=1, \dots, n$, are  independently and identically distributed  realizations of $(\X,
Z,\Delta)$. To make  (\ref{eq:model})  identifiable, we fix
the
upper $d\times d$ block of
$\bb$
to be $\I_d$. We aim to estimate the
unspecified  function $m(\cdot,\cdot)$, as well as
the column space of $\bb$, which is
equivalent to estimating the lower $(p-d)\times d$ block of
$\bb$. To associate the
covariates with the upper and lower parts of $\bb$, we write $\X=(\X_u\trans,
\X_l\trans)\trans$, where $\X_u\in{\cal R}^d$ and $\X_l\in{\cal R}^{p-d}$.

Model (\ref{eq:model})  leads to 
\be
S(t\mid\X)&=&S(t,\bb\trans\X)=\frac{m(0,\bb\trans\X)}{m(t,\bb\trans\X)}\exp\left\{-\int_0^t\frac{1}{m(u,\bb\trans\X)}du\right\},\nonumber\\
\lambda(t,\bb\trans\X)&=&\frac{m_1(t,\bb\trans\X)+1}{m(t,\bb\trans\X)},\label{eq:hazard}
\ee
where
$\lambda(t,\bb\trans\X)$ is the hazard function of $T$
conditional on $\X$, and  $m_1(t,\bb\trans\X)=\partial
m(t,\bb\trans\X)/\partial t$. Likewise, we denote by
$\m_2(t,\bb\trans\X)=\partial
m(t,\bb\trans\X)/\partial (\bb\trans\X)$.
We can re-express (\ref{eq:hazard}) as
\be\label{eq:cumuhazard}
m(t,\bb\trans\X) = e^{\Lambda(t,\bb\trans\X)}\int_t^\infty
e^{-\Lambda(s,\bb\trans\X)}ds
\ee
where $\Lambda(\cdot,\bb\trans\X)$ is the cumulative hazard
function  \citep{maguluri1994estimation}. 



\section{A semiparametrically efficient estimator}\label{sec:estimation}

\subsection{Nuisance tangent spaces}\label{sec:space}
Denote the conditional survival function, cumulative hazard
function, hazard function and  probability density function
(pdf) of the censoring time $C$ by
$S_c(z,\X)=\pr (C\ge z\mid \X)$,
$\Lambda_c(z,\X)=-\log S_c(z,\X)$,
$\lambda_c(z,\X)=\partial\Lambda_c(z,\X)/\partial z$ and
$f_c(z,\X)=-\partial S_c(z,\X)/\partial z$ with
$z<\tau$, where
$\tau <\infty$ is the upper bound of the follow-up time.
Let
	$p(\X)\equiv\pr(C=\tau\mid\X)$,
	$S_c(\tau,\X)=f_c(\tau,\X)=p(\X)$, and	
$\lambda_c(\tau, \X)=1$.
	Here, 
	$\lambda_c(z,\X)$ and $f_c(z,\X)$ are absolutely
	continuous on $(0, \tau)$, but with 
	a discontinuity point at $\tau$.
For any $\bb$,  in addition to the survival
function
$S(z,\bb\trans\X)$ and the hazard function
$\lambda(z,\bb\trans\X)$  in (\ref{eq:hazard}), we also
define the pdf
$f(z,\bb\trans\X)=-\partial S(z,\bb\trans\X)/\partial z$ and the
cumulative hazard function
$\Lambda(z,\bb\trans\X)=-\log S(z,\bb\trans\X)$. We
  write
$\blam_2(s,\bb\trans\X)=\partial\lambda(s,\bb\trans\X)/
\partial(\bb\trans\X)$
and
 $\blam_{20}(s,\bb\trans\X)=\partial\lambda_0(s,\bb\trans\X)/
\partial(\bb\trans\X)$
as the derivatives of the hazard functions with respect to
$\bb\trans\X$, where $\lambda_0(z,\bb\trans\X)$ is given
in (\ref{eq:hazard}) when the mean
residual function $m(\cdot,\cdot)=m_0(\cdot,\cdot)$,
where the subscript ``0" indicates the truth.

The pdf of $(\X,Z,\Delta)$  is
\be\label{eq:pdf}
&&f_{\X,Z,\Delta}(\x,z,\delta) 
=f_\X(\x)
\lambda(z,\bb\trans\x)^\delta
e^{-\int_0^z\lambda(s,\bb\trans\x)ds}
\lambda_c(z,\x)^{1-\delta}
e^{-\int_0^z\lambda_c(s,\x)ds},
\ee
where $f_\X(\x)$ is the pdf of $\X$. 
We view the pdf in (\ref{eq:pdf}) as a semiparametric model
where
all unknown components, except for $\bb$, are infinite dimensional 
nuisance parameters.  The parameters $\bb$ are  
parameters of interest with a finite
dimension. We will  estimate
$\bb$ by using
a geometric approach, which  
avoids modeling the hazard function
$\lambda(z,\bb\trans\X)$
to be the product of a baseline function of  $z$  and
a specific covariate function such as $\exp(\bb\trans\X)$, as in a proportional hazards model. 
This entails more flexibility for the model.

Let $Y(t) = I(Z\ge
t)$ be the at risk process and
$N(t) = I(Z\le
t)\Delta$ be the counting process.
Define the filtration ${\cal F}_t$ to be $\sigma\{
N(u), Y(u), \X, 0\le u < t\}$,  and let $M(t)$ be
the martingale with respect to ${\cal F}_t$, i.e.
$M(t,\bb\trans\X)=N(t)-\int_0^t
Y(s)\lambda(s,\bb\trans\X)ds$.
The nuisance tangent space, which will be utilized for deriving our
estimator, is obtained as follows. 

\begin{Pro}\label{pro:nuisance}
	The nuisance tangent space is
	${\cal T}={\cal T}_1\oplus{\cal T}_2\oplus{\cal T}_3$, where each component corresponds to $f_\X$, $m(\cdot, \cdot)$ and $\lambda_c$, 
respectively. Specifically,	\bse
	{\cal T}_1&=&[\a(\X): E\{\a(\X)\}=\0, \a(\X)\in{\cal
	R}^{(p-d)d}, \var\{\a(\X)\}<\infty],\\
	{\cal T}_2&=&[\int_0^\infty
	\left\{\frac{\h_1(s,\bb\trans\X)}{m_1(s,\bb\trans\X)+1}
	-\frac{\h(s,\bb\trans\X)}{m(s,\bb\trans\X)}\right\}
	dM(s,\bb\trans\X):\\
&&\forall\h(z,\bb\trans\X)\in{\cal
	R}^{(p-d)d}, \var\{\h(z,\bb\trans\X)\}<\infty],\\
	{\cal T}_3&=&[\int_0^\infty\h(s,\X)dM_c(s,\X):\forall\h(z,\X)\in{\cal
	 R}^{(p-d)d},\var\{\h(z,\X)\}<\infty].
	\ese
\end{Pro}
The derivation of Proposition \ref{pro:nuisance} is provided in
Supplement \ref{app:nuisance}.

\subsection{Derivation of an efficient score function}\label{sec:score}

Taking the derivative of the logarithm of (\ref{eq:pdf}) 
with
respect to  $\bb$, we obtain
the score function
\bse
\bS_{\bb}(\Delta,Z,\X)
=\int_0^\infty
\left\{\frac{\m_{12}(s,\bb\trans\X)}{m_1(s,\bb\trans\X)+1}
-\frac{\m_2(s,\bb\trans\X)
}{m(s,\bb\trans\X)}\right\}\otimes\X_l
dM(s,\bb\trans\X),
\ese
where
$\m_{12}(t,\bb\trans\X)\equiv
\partial\m_2(t,\bb\trans\X)/\partial t$ and  $\X_l$
is the lower $p-d$ components in $\X$.
We can verify that, at $\bb_0$, $\bS_{\bb}(\Delta,Z,\X)\perp{\cal T}_1$ and
$\bS_{\bb}(\Delta,Z,\X)\perp{\cal T}_3$ due to the martingale
properties.
To look for an efficient score by projecting
$\bS_{\bb}(\Delta, Z,\X)$ at $\bb_0$ to ${\cal T}_2$, we search for
$\h^*(s,\bb_0\trans\X)$ such that
\bse
&&\bS\eff(\Delta,Z,\X)\\
&=&\bS_{\bb}(\Delta,Z,\X)
-\int_0^\infty
\left\{\frac{\h_1^*(s,\bb_0\trans\X)}{m_1(s,\bb_0\trans\X)+1}
-\frac{\h^*(s,\bb_0\trans\X)}{m(s,\bb_0\trans\X)}\right\}
dM(s,\bb_0\trans\X)\\
&=&\int_0^\infty
\left\{\frac{\m_{12}(s,\bb_0\trans\X)\otimes\X_l-\h_1^*(s,\bb_0\trans\X)}{m_1(s,\bb_0\trans\X)+1}
-\frac{\m_2(s,\bb_0\trans\X)
\otimes\X_l-\h^*(s,\bb_0\trans\X)}{m(s,\bb_0\trans\X)}\right\}
dM(s,\bb_0\trans\X)
\ese
is orthogonal to ${\cal T}_2$, which implies that, for any $\h(s,\bb_0\trans\X)$, it must hold that
\be\label{eq:must}
0
=E\left[\int_0^\infty \a(s,\bb_0\trans\X)\trans
\left\{\frac{\h_1(s,\bb_0\trans\X)}{m_1(s,\bb_0\trans\X)+1}
-\frac{\h(s,\bb_0\trans\X)}{m(s,\bb_0\trans\X)}\right\}
ds\right],
\ee
where  
\be\label{eq:a}
&&\hspace{-0.6cm}\a(s,\bb_0\trans\X)\equiv
E\left[\left\{\frac{\m_{12}(s,\bb_0\trans\X)\otimes\X_l-\h_1^*(s,\bb_0\trans\X)}{m_1(s,\bb_0\trans\X)+1}
-\frac{\m_2(s,\bb_0\trans\X)
	\otimes\X_l-\h^*(s,\bb_0\trans\X)}{m(s,\bb_0\trans\X)}\right\} \right.\n\\
&& \hspace{2.6cm}\times	S_c(s,\X)\mid\bb_0\trans\X\bigg]
S(s,\bb_0\trans\X)\frac{m_1(s,\bb_0\trans\X)+1}{m(s,\bb_0\trans\X)}.
\ee
We can choose any $\h(s,\bb_0\trans\X)$ function. Specifically, by letting $\h(s,\bb_0\trans\X)=0$ for $s<t$ and
$\h(s,\bb_0\trans\X)=\c(\bb_0\trans\X)$ for $s\ge t$ with an
arbitrary
function $\c(\bb_0\trans\X)$, we obtain
$
\a(t,\bb_0\trans\X)/\{m_1(t,\bb_0\trans\X)+1\}-\int_t^\infty\a(s,\bb_0\trans\X)/m(s,\bb_0\trans\X)ds=\0.
$
Solving this integral equation leads to
\bse
\a(t,\bb_0\trans\X)=\{m_1(t,\bb_0\trans\X)+1\}
\exp\left\{-\int_0^t\frac{m_1(s,\bb_0\trans\X)+1}{m(s,\bb_0\trans\X)}ds\right\}\c(\bb_0\trans\X),
\ese
for  function $\c(\cdot)$. Thus, reusing (\ref{eq:must}), we require that for all
$\h(t,\bb_0\trans\X)$,
\bse
0&=&E\left[\int_0^\infty \{m_1(t,\bb_0\trans\X)+1\}
\exp\left\{-\int_0^t\frac{m_1(s,\bb_0\trans\X)+1}{m(s,\bb_0\trans\X)}ds\right\}\c(\bb_0\trans\X)\trans\right.\\
&&\left.\hspace{2.6cm}\times\left\{\frac{\h_1(t,\bb_0\trans\X)}{m_1(t,\bb_0\trans\X)+1}
-\frac{\h(t,\bb_0\trans\X)}{m(t,\bb_0\trans\X)}\right\}dt\right]\\
&=&-E\left\{\c(\bb_0\trans\X)\trans
\h(0,\bb_0\trans\X)\right\}.
\ese
Letting $\h(0,\bb_0\trans\X)=\c(\bb_0\trans\X)$ yields the only
possibility of $\c(\bb_0\trans\X)=\0$, hence
$\a(t,\bb_0\trans\X)=\0$.
Inserting the expression of
$\a(t,\bb_0\trans\X)$ into (\ref{eq:a}), we have 
\bse
\frac{\h_1^*(t,\bb_0\trans\X)}{m_1(t,\bb_0\trans\X)+1}
-\frac{\h^*(t,\bb_0\trans\X)}{m(t,\bb_0\trans\X)}
=\left\{\frac{\m_{12}(t,\bb_0\trans\X)}{m_1(t,\bb_0\trans\X)+1}
-\frac{\m_2(t,\bb_0\trans\X)
}{m(t,\bb_0\trans\X)}\right\} \otimes \frac{E\left\{\X_l
S_c(t,\X)\mid\bb_0\trans\X\right\}}{E\left\{S_c(t,\X)\mid\bb_0\trans\X\right\}}.
\ese
Thus an efficient score is
\be\label{eq:effscore}
&&\bS\eff(\Delta,Z,\X)\\
&=&\int_0^\infty
\left\{\frac{\m_{12}(s,\bb_0\trans\X)}{m_1(s,\bb_0\trans\X)+1}
-\frac{\m_2(s,\bb_0\trans\X) }{m(s,\bb_0\trans\X)}\right\}\otimes
\left[\X_l-
\frac{E\left\{\X_l
S_c(s,\X)\mid\bb_0\trans\X\right\}}{E\left\{S_c(s,\X)\mid\bb_0\trans\X\right\}}\right]
dM(s,\bb_0\trans\X).\n
\ee

\subsection{Construction of a semiparametrically efficient estimator}\label{sec:estimator}

Based on  (\ref{eq:effscore}), we
construct a semiparametrically efficient estimator of $\bb$.
First, a consistent estimating equation can be obtained 
from
$E\{\bS\eff(\Delta,Z,\X)\mid\X\}=\0$, due to
$E\{dM(t,\bb_0\trans\X)\mid\X\}=0$. Hence, to preserve the mean zero
property and to simplify the computation, we can replace
$
{\m_{12}(s,\bb_0\trans\X)}/\{m_1(s,\bb_0\trans\X)+1\}
-{\m_2(s,\bb_0\trans\X) }/{m(s,\bb_0\trans\X)}
$
by an arbitrary function of $s$ and $\bb_0\trans\X$, say $\g(s,\bb_0\trans\X)$, and
still obtain
\bse
E\left(\int_0^\infty\g(s,\bb_0\trans\X)
\otimes\left[\X_l-
\frac{E\left\{\X_l
	S_c(s,\X)\mid\bb_0\trans\X\right\}}
{E\left\{S_c(s,\X)\mid\bb_0\trans\X\right\}}\right]dM(s,\bb_0\trans\X)\right)=\0.
\ese
This provides a richer class of estimators than the estimator based on
$\bS\eff$ alone.

Second, we can obtain
\be\label{eq:useful}
\frac{E\left\{ \X_lY(t)\mid\bb_0\trans\X\right\}}
{E\left\{Y(t)\mid\bb_0\trans\X\right\}}
=\frac{E\left\{\X_l
	S_c(t,\X)\mid\bb_0\trans\X\right\}}
{E\left\{S_c(t,\X)\mid\bb_0\trans\X\right\}},
\ee
where we define (\ref{eq:useful})
to be 
${E\left\{\X_l
	p(\X)\mid\bb_0\trans\X\right\}}/
{E\left\{p(\X)\mid\bb_0\trans\X\right\}}$ when $t>\tau$, with $p(\X)$ defined in Section \ref{sec:space}.
We then verify that 
\be\label{eq:expecteff}
E\left(\int_0^\infty\g(s,\bb_0\trans\X)
\otimes\left[\X_l-
\frac{E\left\{\X_l
	S_c(s,\X)\mid\bb_0\trans\X\right\}}
{E\left\{S_c(s,\X)\mid\bb_0\trans\X\right\}}\right]dN(s)\right)=\0.
\ee
The proof of (\ref{eq:useful}) and
(\ref{eq:expecteff}) is given
in Supplement \ref{app:eff}.
This implies that we can construct estimating equations of the
form
\be\label{eq:general}
\sumi \Delta_i\g(Z_i,\bb\trans\X_i)
\otimes\left[\X_{li}-
\frac{\wh E\left\{\X_{li}
	Y_i(Z_i)\mid\bb\trans\X_i\right\}}
{\wh E\left\{Y_i(Z_i)\mid\bb\trans\X_i\right\}}\right]=\0
\ee
for any $\g(\cdot,\cdot)$, with
\be
\wh E\left\{Y_i(Z_i)\mid\bb\trans\X_i\right\}
&=& \frac{\sumj K_h(\bb\trans\X_j-\bb\trans\X_i)I(Z_j\ge Z_i)}
{\sumj K_h(\bb\trans\X_j-\bb\trans\X_i)},\label{eq:expectY}\\
\wh E\left\{\X_{li}Y_i(Z_i)\mid\bb\trans\X_i\right\}
&=& \frac{\sumj K_h(\bb\trans\X_j-\bb\trans\X_i)\X_{lj}I(Z_j\ge
	Z_i)}
{\sumj K_h(\bb\trans\X_j-\bb\trans\X_i)}.\label{eq:expectXY}
\ee
Here, $\hat E\{Y_i(Z_i)\mid \bb\trans\X_i\}\equiv
\hat E\{Y_i(t)\mid\bb\trans\X_i\}|_{t=Z_i}$
and similarly for other terms.

Third, when we choose to estimate the unknown components
$m(\cdot,\cdot)$, $m_1(\cdot,\cdot)$, $\m_2(\cdot,\cdot)$ and
$\m_{12}(\cdot,\cdot)$ nonparametrically,
we then obtain the
efficient estimator of $\bb$ by solving
\be\label{eq:eff}
\sumi \Delta_i
\left\{\frac{\wh\m_{12}(Z_i,\bb\trans\X_i)}{\wh m_1(Z_i,\bb\trans\X_i)+1}
-\frac{\wh\m_2(Z_i,\bb\trans\X_i) }{\wh m(Z_i,\bb\trans\X_i)}
\right\}
\otimes\left[\X_{li}-
\frac{\wh E\left\{\X_{li}
	Y_i(Z_i)\mid\bb\trans\X_i\right\}}
{\wh E\left\{Y_i(Z_i)\mid\bb\trans\X_i\right\}}\right]=\0.
\ee
We estimate $m(t,\bb\trans\X)$ nonparametrically  via estimating $\Lambda(t,\bb\trans\X)$ by
\be
\wh\Lambda(t,\bb\trans\X)&=&\sumi\int_{0}^{t}
\frac{K_h(\bb\trans\X_i-\bb\trans\X)}{\sumj
Y_j(s)K_h(\bb\trans\X_j-\bb\trans\X)}
dN_i(s),\label{eq:Lambda}
\ee
at any $t$,
and using (\ref{eq:cumuhazard}) to obtain
\be\label{eq:m}
\wh m(t,\bb\trans\X)
=e^{\wh\Lambda(t,\bb\trans\X)}
\int_t^\infty
e^{-\wh\Lambda(s,\bb\trans\X)}ds,
\ee
where $\wh\Lambda(t,\bb\trans\X)$ is a kernel
smoothed version of the Nelson-Aalen estimator
\citep{ramlau1983choice, andersenborgan}.
The estimators of the derivatives of $\Lambda(t,\bb\trans\X)$ and $m(t,\bb\trans\X)$ are given in Supplement \ref{sec:nonpara}.
To bypass the zero-denominator issue when applying
the nonparametric estimator to a finite sample case, we propose a trimmed version in Supplement \ref{sec:trimmed}. and show that it retains the same asymptotic properties proved in Section 
  \ref{sec:asymp}.


	
	
	

\section{Asymptotic properties and semiparametric efficiency}\label{sec:asymp}

We list the  regularity conditions for the results of
root-$n$ consistency and  asymptotic normality of the estimators proposed
in Section \ref{sec:estimation}.
We  also establish semiparametric
efficiency of the
estimator obtained by solving (\ref{eq:eff}).  

\begin{enumerate}[label=C\arabic*]
	
	\item\label{assum:kernel} ({\itshape kernel function})
	The kernel function $K_h(\cdot)=h^{-d}K(\cdot/h)$ where
	$K(\a)=\prod_{j=1}^{d}K(a_j)$ for
	$\a=(a_1,...,a_d)\trans$ is symmetric on each individual
	entry and $K(a_j)$ is differentiable,  
	decreasing when $x\ge0$, and $\int K(x)dx=1$,
	$\int x^j K(x)dx=0$, for $1\le j<\nu$, $0<\int x^\nu
	K(x)dx<\infty$, and
	$\int K^2(x)dx$,
	$\int x^2 K^2(x)dx$,
	$\int K'^2(x)dx$,
	$\int x^2K'^2(x)dx$,
	$\int K''^2(x)dx$,
	$\int x^2K''^2(x)dx$ are all bounded. When there is no confusion, 
	we use the same $K$ for both
	univariate and multivariate kernel functions for simplicity.


	\item\label{assum:bandwidth}({\itshape bandwidths})
	The bandwidths $h$ and $b$ satisfy $h\to0$,	$nh^{2\nu}\to0$, $b\to0$ and
	$nh^{d+2}b\to\infty$, where $2\nu> d+1$.

	\item\label{assum:fbeta}({\itshape density functions of
	covariates})
	For all $\bb\in{\cal B}$, the 
parameter space,
	the probability density function of $\bb\trans\X$,
	$f_{\bb\trans\X}(\bb\trans\x)$, has a compact support and is
	bounded away from zero and $\infty$. Further
	$f_{\bb\trans\X}(\bb\trans\x)$ has a 
	derivative, up to the fourth order, that is bounded uniformly on the support.
	
	\item\label{assum:exi}({\itshape smoothness})
	For all $\bb\in{\cal B}$ and $t>0$,
	$E\{\X_jI(Z_j\ge
	t)\mid\bb\trans\X_j=\bb\trans\x\}$,
	and its 
	derivatives, up to the fourth order, are bounded uniformly as functions of $\bb\trans\x$;
	 $E\{\X_j\X_j\trans I(Z_j\ge
	t)\mid\bb\trans\X_j=\bb\trans\x\}$
	and its first and second order
	derivatives are bounded uniformly as functions of $\bb\trans\x$.

	\item\label{assum:survivalfunction}({\itshape survival
	function})
	For all $\bb\in{\cal B}$, it holds that,  at any  $t$,
	$\partial^{i+j}S(t,\bb\trans\x)/\partial
	t^i\partial (\bb\trans\x)^j$, $\partial^{i+j}E\{S_c(t,\X)\mid
	\bb\trans\x\}/\partial t^i\partial (\bb\trans\x)^j$ and
	$\partial^{i+j}f(t,\bb\trans\x)/\partial t^i\partial
	(\bb\trans\x)^j$ exist and  are bounded and bounded away
	from zero, for all $i\ge0, j\ge0, i+j\le 4$. In
          addition, $S_c(\tau,\X)$ is bounded
          way from zero.

	\item\label{assum:bounded}({\itshape boundedness})
The true parameter $\bb_0$ is an interior point in  $\cal B$
and ${\cal B}$ is bounded.

	\item\label{assum:unique}({\itshape uniqueness})
	The equation
	\bse
	E\left( \Delta\left\{\frac{\m_{12}(Z,\bb\trans\X)}{
	m_1(Z,\bb\trans\X)+1}
	-\frac{\m_2(Z,\bb\trans\X) }{m(Z,\bb\trans\X)}
	\right\}
	\otimes\left[\X_{l}-
	\frac{ E\left\{\X_{l}
		Y(Z)\mid\bb\trans\X\right\}}
	{E\left\{Y(Z)\mid\bb\trans\X\right\}}\right]\right)
	=\0
	\ese
	has a unique solution in $\cal B$.

\end{enumerate}

 Conditions \ref{assum:kernel} and
\ref{assum:bandwidth} are commonly assumed in kernel
regression analysis \citep{silverman1986density,ma2013efficient}. Conditions
\ref{assum:fbeta}-\ref{assum:survivalfunction} assume boundedness
of event time, censoring time, covariates and their
expectations, which hold for real datasets. The smoothness of several functions is imposed by constraining
their derivatives, which are  common conditions \citep{silverman1978weak}. It is natural to make a boundedness assumption on
the parameter  space
$\cal B$ as in Condition \ref{assum:bounded} in practical
problems \citep{hardle1997semiparametric}. 
Condition \ref{assum:unique} precludes 
that the estimating equation 
is degenerate.

Condition \ref{assum:fbeta} can be slightly modified  for  the trimmed estimators.
\begin{enumerate}[label=C\arabic*$'$, start=3]
	\item\label{assum:fbetarelax}({\itshape The density of 
		index, relaxed.})
	Uniformly for any $\bb$ in a local neighborhood of
	$\bb_0$, the density function 
	of $\bb\trans\X$, i.e. $f_{\bb\trans\X}(\v)$, is bounded, and
	there exists a constant $\epsilon>0$ such that 
	$\int_{\{\v:f_{\bb\trans\X}(\v)\le
		d_n\}}f_{\bb\trans\X}(\v)d\v<n^{-\epsilon}$
	for sufficiently large $n$. Here $d_n\rightarrow 0$ as
	$n\rightarrow\infty$, 
	and $n^{-\epsilon}=O(h^2+n^{-1/2}h^{-1/2})$,
	where $h$ satisfies Condition
	\ref{assum:bandwidth}.  In addition, 
	the  
	derivatives of $f_{\bb\trans\X}(\cdot)$, up to the fourth order, are bounded.
\end{enumerate}
Condition \ref{assum:fbetarelax}, 
 weaker than Condition \ref{assum:fbeta}, requires the
tail of $f_{\bb\trans\X}$ to be sufficiently thin to ensure the
near zero values of $f_{\bb\trans\X}(\cdot)$  not affect the
overall performance of our estimator. It guarantees that the trimmed  nonparametric estimators will retain the same asymptotic properties discussed below.

Theorems \ref{th:consistency} and \ref{th:eff} demonstrate
the root-$n$ consistency and asymptotical normality of the
profile parameter estimator $\wh\bb$.  The proofs  are given in Supplement \ref{app:consistofb} and
\ref{app:asympofb}.

\begin{Th}\label{th:consistency}
	Under Conditions
	\ref{assum:kernel}-\ref{assum:unique},
	the estimator, $\wh\bb$, obtained by solving (\ref{eq:general}) or
(\ref{eq:eff}) 
is
	consistent,
	i.e.
	$\wh\bb-\bb\to\0$ in probability when $n\to\infty$.
\end{Th}

\begin{Th}\label{th:eff}
	Under Conditions
	\ref{assum:kernel}-\ref{assum:unique},
	the estimator, $\wh\bb$, obtained by solving (\ref{eq:general}) or
(\ref{eq:eff}) 
satisfies
$
	\sqrt{n}(\wh\bb-\bb)\to N(\0, \A^{-1}\B{\A^{-1}}\trans)
$
	in distribution when $n\to\infty$, where
\bse
\A&=&E\left\{\frac{\partial}{\partial\vecl(\bb)\trans}
\vecl\left(\Delta\g(Z,\bb\trans\X)
\otimes\left[\a(\X_{l})-
\frac{E\left\{\a(\X_{l})
	Y(Z)\mid\bb\trans\X\right\}}
{E\left\{Y(Z)\mid\bb\trans\X\right\}}\right]\right)\right\},\\
\B&=&E\left\{\vecl\left(\Delta\g(Z,\bb\trans\X)
\otimes\left[\a(\X_{l})-
\frac{E\left\{\a(\X_{l})
	Y(Z)\mid\bb\trans\X\right\}}
{E\left\{Y(Z)\mid\bb\trans\X\right\}}\right]\right)^{\otimes2}\right\}.
\ese
Here $\vecl(\A)$ represents the vectorization of the lower
$(p-d)\times d$ block of a generic matrix $\A$ and
$\A^{\otimes2}=\A\A\trans$ for any matrix or vector $\A$.
Note that in (\ref{eq:general}), $\a(\X_l)=\X_l$ and in
(\ref{eq:eff}), $\a(\X_l)=\X_l$,
$g(Z,\bb\trans\X)=\wh\m_{12}(Z,\bb\trans\X)/\{\wh
m_1(Z,\bb\trans\X)+1\}
-\wh\m_{2}(Z,\bb\trans\X)/\wh
m(Z,\bb\trans\X)$.
Further, the estimator, $\wh\bb$, obtained from solving
(\ref{eq:eff}) is semiparametrically efficient and  satisfies
	\bse
	\sqrt{n}(\wh\bb-\bb)\to N(\0, [E\{\bS\eff ^{\otimes2}
	(\Delta,
	Z,\X)\}]^{-1})
	\ese
	in distribution, where $\bS\eff(\Delta, Z,
	\X)$ is given
	in (\ref{eq:effscore}).
\end{Th}

With $\bS\eff$ being a martingale, 
\bse
&&E\{\bS\eff ^{\otimes2} (\Delta,
Z,\X)\}\\
&=&E\left\{\int_0^\infty
\left(\left\{\frac{\m_{12}(s,\bb\trans\X)}{
	m_1(s,\bb\trans\X)+1}
-\frac{\m_2(s,\bb\trans\X) }{m(s,\bb\trans\X)}\right\}
\otimes\left[\X_l-
\frac{E\left\{\X_l
	S_c(s,\X)\mid\bb\trans\X\right\}}
{E\left\{S_c(s,\X)\mid\bb\trans\X\right\}}\right]\right)^{\otimes2}
dN(s)\right\},
\ese
which leads to a consistent estimator of $E\{\bS\eff ^{\otimes2} (\Delta,
	Z,\X)\}$ as follows
	\bse
	\frac{1}{n}\sumi\delta_i
	\left(\frac{\wh\blam_{2}(z_i,\wh\bb\trans\x_i)}{\wh\lambda(z_i,\wh\bb\trans\x_i)}
	\otimes\left[\x_{il}-
	\frac{\wh E\left\{\X_{l}
		Y(z_i)\mid\wh\bb\trans\x_i\right\}}
	{\wh
	E\left\{Y(z_i)\mid\wh\bb\trans\x_i\right\}}\right]\right)^{\otimes2}.
	\ese
	Here $\wh
	E\left\{Y(z_i)\mid\wh\bb\trans\x_i\right\}$,
	$\wh E\left\{\X_{l}
	Y(z_i)\mid\wh\bb\trans\x_i\right\}$,
	$\wh\lambda(z_i,\wh\bb\trans\x_i)$
	and $\wh\blam_{2}(z_i,\wh\bb\trans\x_i)$
	are
	given in (\ref{eq:expectY}), (\ref{eq:expectXY}),
	(\ref{eq:lambda}) and (\ref{eq:lambda1}) respectively.

\begin{Th}\label{th:m}
	Under Conditions
	\ref{assum:kernel}-\ref{assum:unique}, the
	nonparametric estimator $\wh m(t,\wh\bb\trans\x)$ satisfies
	\bse
	\sqrt{nh}\left\{\wh{m}(t,\wh\bb\trans\x)-
	{m}(t,\bb\trans\x)\right\}\to
	N\{0,\sigma^2_m(t,\bb\trans\x)\}
	\ese
in distribution, where
	\bse
	\sigma^2_m(t,\bb\trans\x)
	&=&e^{2\Lambda(t,\bb\trans\x)}\frac{\int
		K^2(u)du}{f_{\bb\trans\X}(\bb\trans\x)}
	\int_{0}^{\infty}\frac{\lambda(r,\bb\trans\x)}
	{E\{I(Z\ge r)\mid\bb\trans\x\}}\left\{I(r<t)\int_{t}^{\infty}
	e^{-\Lambda(s,\bb\trans\x)}ds\right.\\
&&\left.+\int_{\max(r,t)}^{\infty}e^{-\Lambda(s,\bb\trans\x)}ds
	\right\}^2dr.
	\ese
\end{Th}
The proof  is provided in
Supplement \ref{app:m}. Theorem
\ref{th:m} implies that   
we can even
estimate the variance $\sigma^2_m(t,\bb\trans\x)$, without estimating
$\lambda$ or $f_{\bb\trans\X}(\bb\trans\x)$, by using 
	\bse
	&&	\hspace{-1cm}\wh\sigma_m^2(t,\wh\bb\trans\x)\\
	&=&e^{2\wh\Lambda(t,\wh\bb\trans\x)}\int K^2(u)du
	\sum_{i=1}^{n_T}\left(\frac{\wh\Lambda(\tau_{(i)},	
	\wh\bb\trans\x)-\wh\Lambda(\tau_{(i-1)},\wh\bb\trans\x)}	
	{1/n\sumj
	Y_j(\tau_{(i-1)})K_h(\wh\bb\trans\x_j-\wh\bb\trans\x)}\right.\\
	&&\left[I(\tau_{(i-1)}<t)\sum_{j=1}^{n_t}I(\tau_{(j)}>
	t)e^{-\wh\Lambda(\tau_{(j-1)},\wh\bb\trans\x)}\{\tau_{(j)}-\max(t,\tau_{(j-1)})\}
	\right.\\
	&&\left.\left.+\sum_{j=1}^{n_t}I\{\tau_{(j)}>\max(t,\tau_{(i-1)})\}
	e^{-\wh\Lambda(\tau_{(j-1)},\wh\bb\trans\x)}\{\tau_{(j)}-\max(r,\tau_{(i-1)},\tau_{(j-1)})\}\right]^2
	\right),
	\ese
where $n_T$ is the total number of the observed events and
$n_t$ is the number of the observed events  up to time $t$.

\section{Simulation}\label{sec:sim}
The section features three simulation studies for evaluating the 
finite sample performance of our
method  in estimating  $\bb$ and
$m(t,\bb\trans\X)$.  For
	comparisons, we also
	implement three estimators,   two  for the proportional mean
	residual life model,  named as ``PM1''
	and
	``PM2'' for \cite{chencheng2005} and \cite{chenetal2005}  respectively, and one for  the additive
	mean residual life model \citep{chen2007}, named as
	``additive''.
All of these competing methods implicitly require $d=1$. 

\noindent
{\bf Study 1:} We generate event times with a  hazard function
$
\lambda(t,\bb\trans\X) =
{2te^{\bb\trans\X}}/\{1+(1+t^2)e^{\bb\trans\X}\}
$
such that the true mean residual life is
\bse
m(t,\bb\trans\X)=\frac{1+e^{\bb\trans\X}+t^2e^{\bb\trans\X}}{(e^{\bb\trans\X}+e^{2\bb\trans\X})^{1/2}}\left[\frac{\pi}{2}-\tan^{-1}\left\{t\sqrt{\frac{e^{\bb\trans\X}}{1+e^{\bb\trans\X}}}\right\}\right].
\ese
Each component in $\X$ is generated independently from the
standard normal
distribution. We consider $d=1, p=10$ and
set the true parameters
to be $\bb=(1,-0.75,0.4,-0.2,0.15,0,$
$-0.15,0.2,-0.4,0.75)\trans$.

\noindent
{\bf Study 2:}  We generate the event times from a
log-logistic distribution with shape parameter 8.0 and scale 
parameter $\exp(\bb\trans\X)$. The corresponding mean residual life
function is
\bse
m(t,\bb\trans\X) = 
\frac{\exp(\bb\trans\X)}{8}\left[1+\left\{\frac{t}{\exp(\bb\trans\X)}\right\}^8\right]
 \int_{z(t)}^{1}(1-u)^{7/8}u^{1/8}du,
\ese
where $z(t) = 1/[1+\{t\exp(-\bb\trans\X)\}^8]$.
Each component in $\X$ is generated independently from the uniform
distribution over $[0,1]$. We consider $d=1, p=7$ and
set the true parameter
$\bb=(1,1.3, -1.3,1,-0.5,0.5,-0.5)\trans$.

\noindent
{\bf Study 3:}
We consider $d=2$ and $p=6$ and  generate the event times with  a hazard function
$
\lambda(t,\bb\trans\X)=t\{e^{(\bb\trans\X)_1}+e^{(\bb\trans\X)_2}\}
$
where $(\bb\trans\X)_k$ denotes the $k$th entry of $\bb\trans\X$,
and the true parameters are
$\bb=\{(1,0,2.75,-0.75,-1,2)\trans,(0,1,-3.125,-1.125,1,-2)\trans\}$.
The mean residual life time function is
\bse
m(t,\bb\trans\X) =
e^{\frac{1}{2}t^2\left\{e^{(\bb\trans\X)_1}+e^{(\bb\trans\X)_2}\right\}}\int_{t}^{\infty}e^{-\frac{1}{2}s^2\left\{e^{(\bb\trans\X)_1}+e^{(\bb\trans\X)_2}\right\}}ds.
\ese

For  each simulation configuration, we generate 1,000 data
sets with $n=500$, and generate
the censoring times from a gamma
	distribution with different parameter values to achieve various
	censoring rates. 

\begin{figure}[H]
	\centering
	\includegraphics[width=5cm]{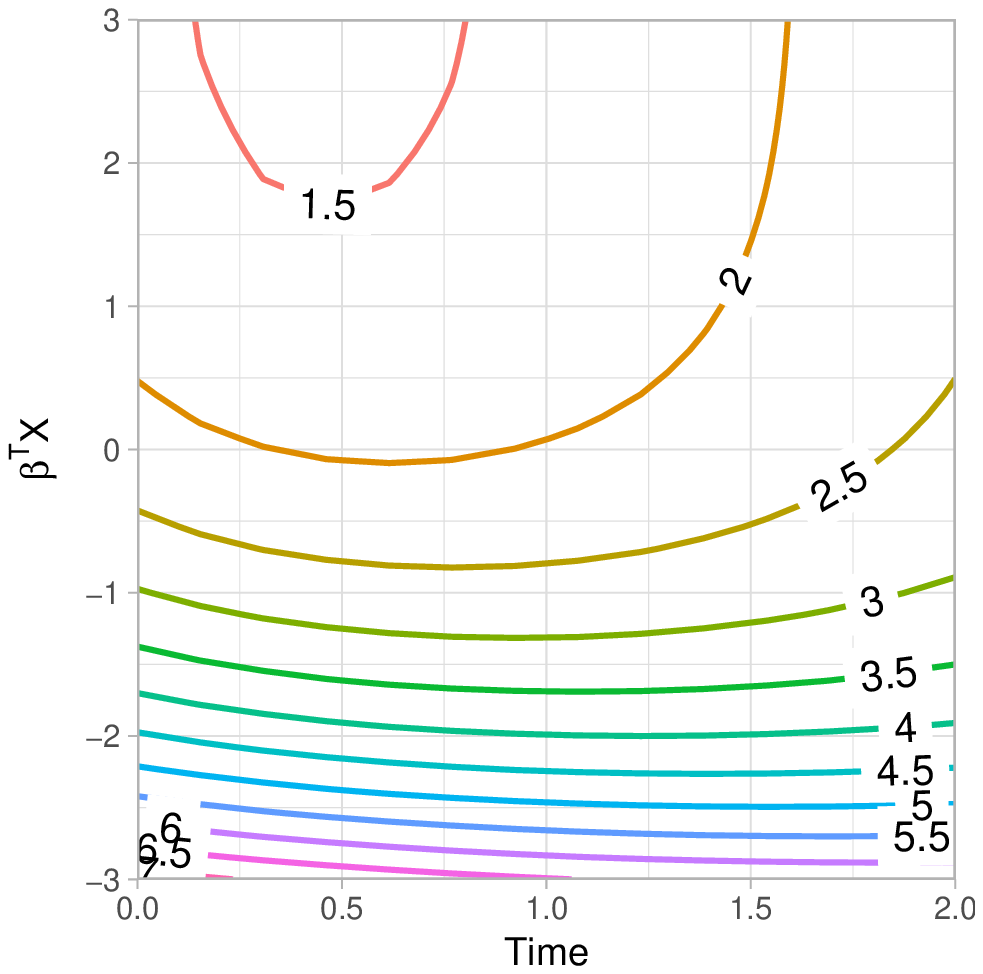}
	\includegraphics[width=5cm]{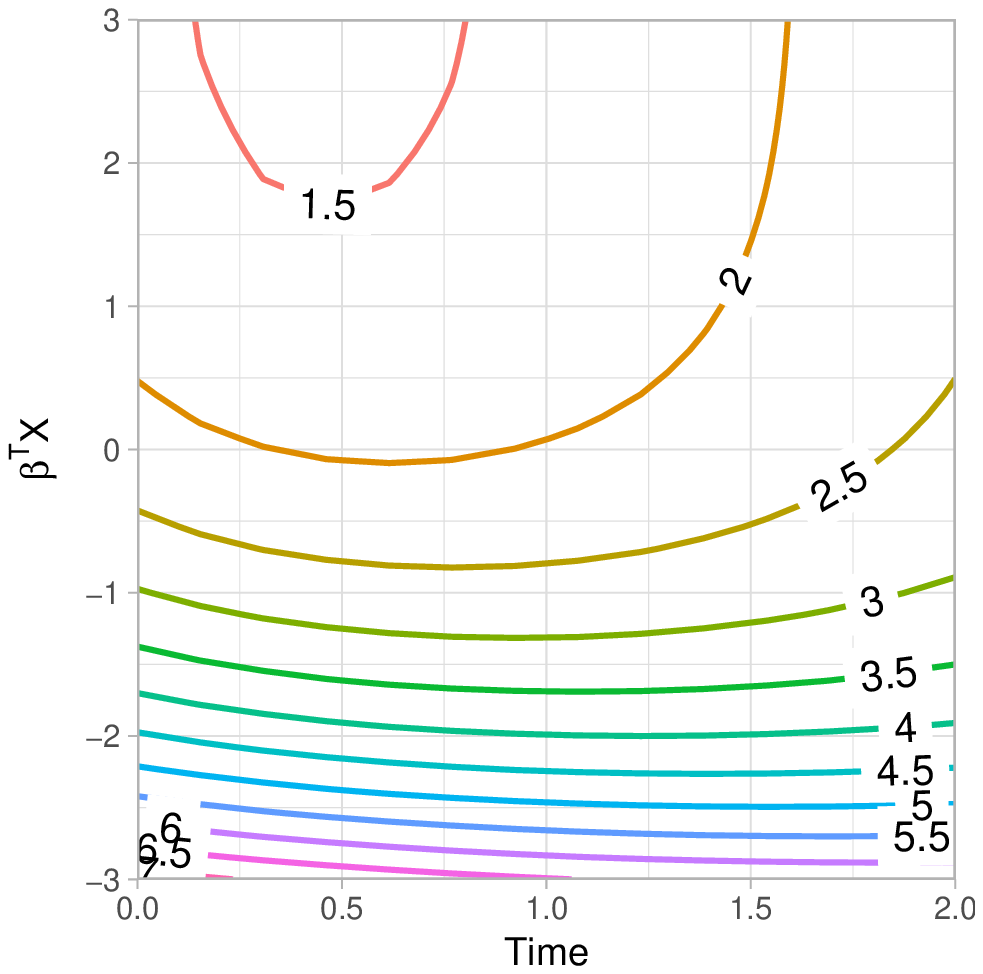}
	\includegraphics[width=5cm]{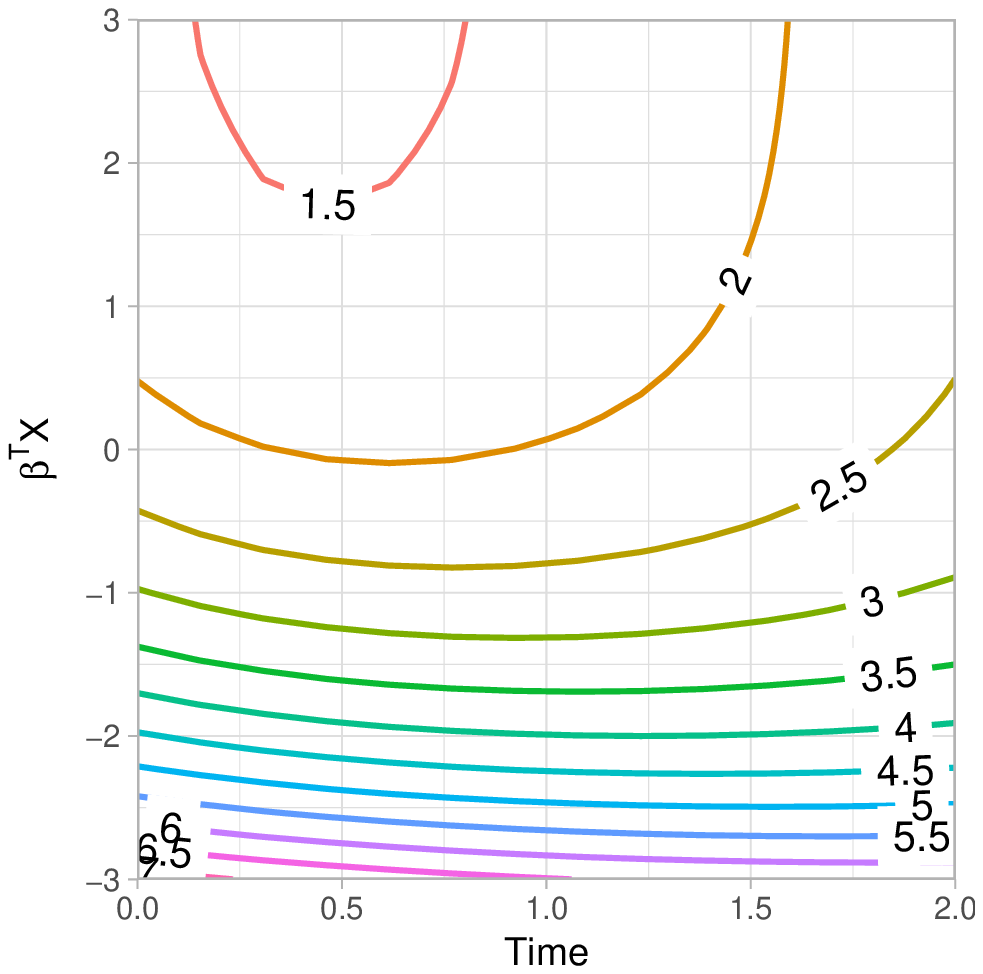}\\
	\includegraphics[width=5cm]{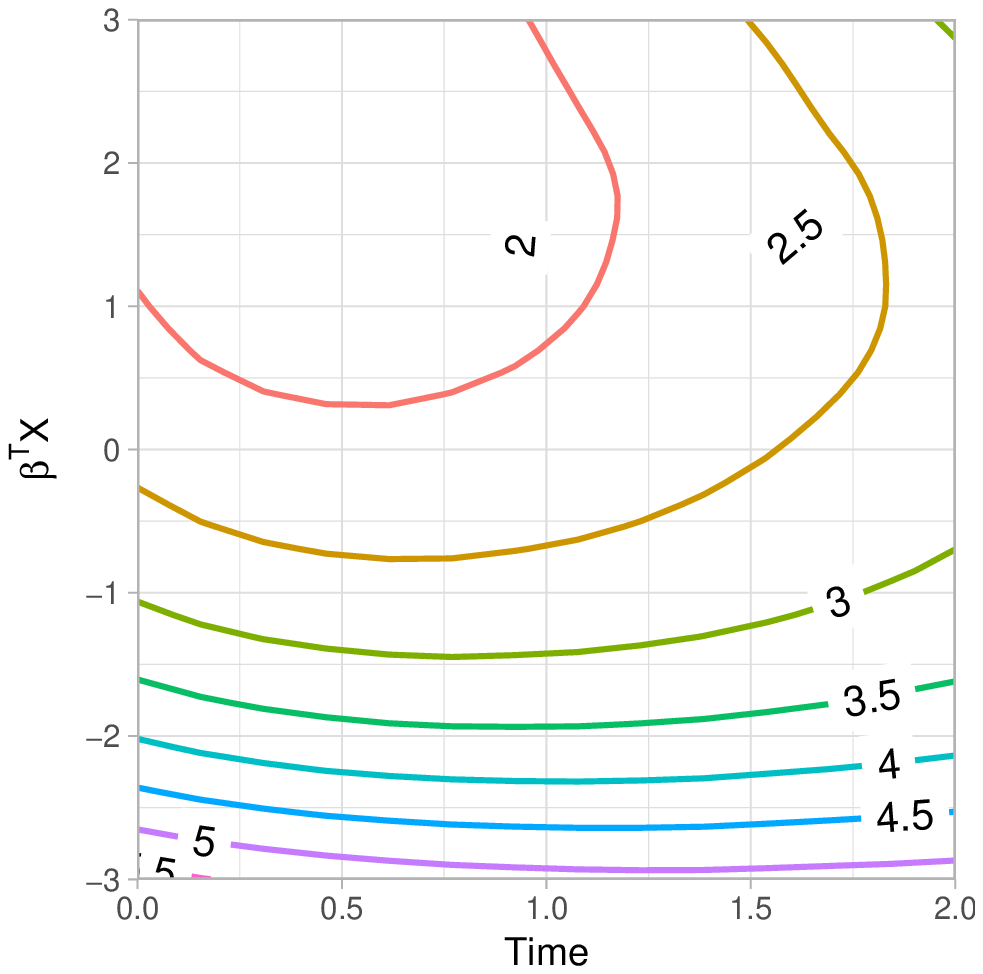}
	\includegraphics[width=5cm]{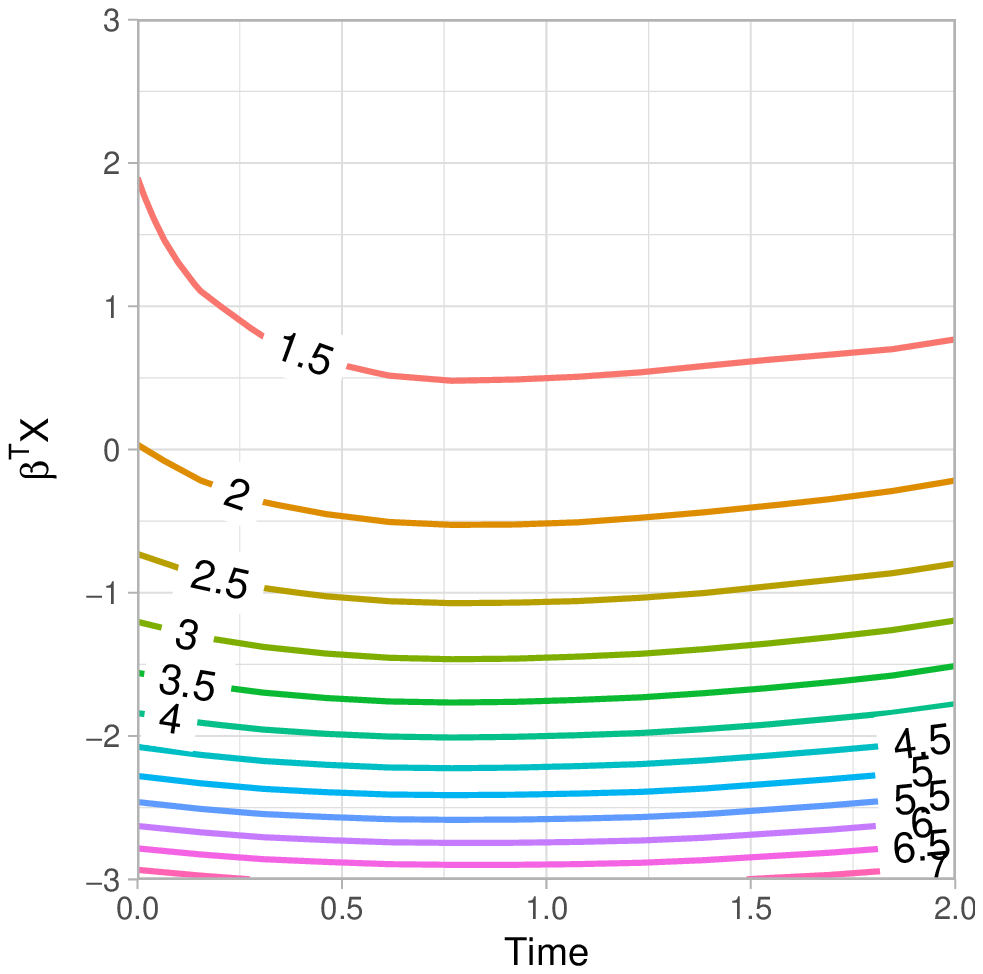}
	\includegraphics[width=5cm]{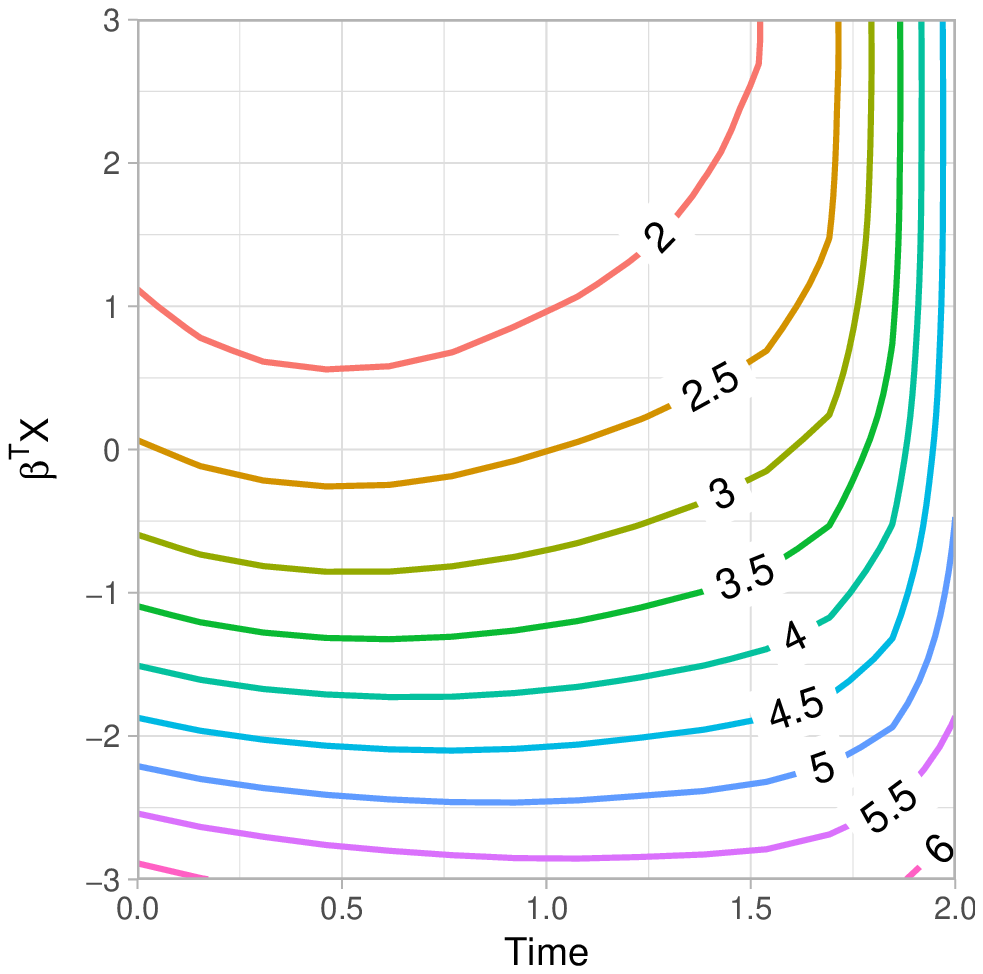}\\
	\includegraphics[width=5cm]{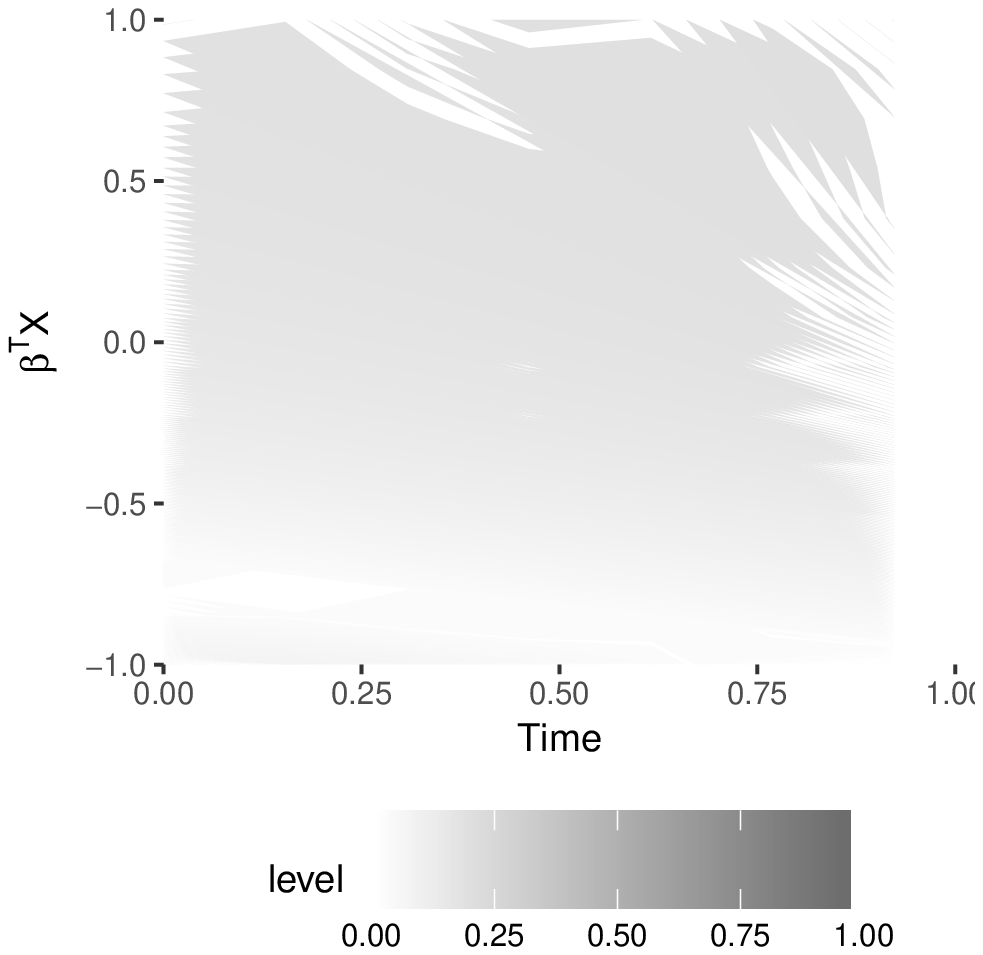}
	\includegraphics[width=5cm]{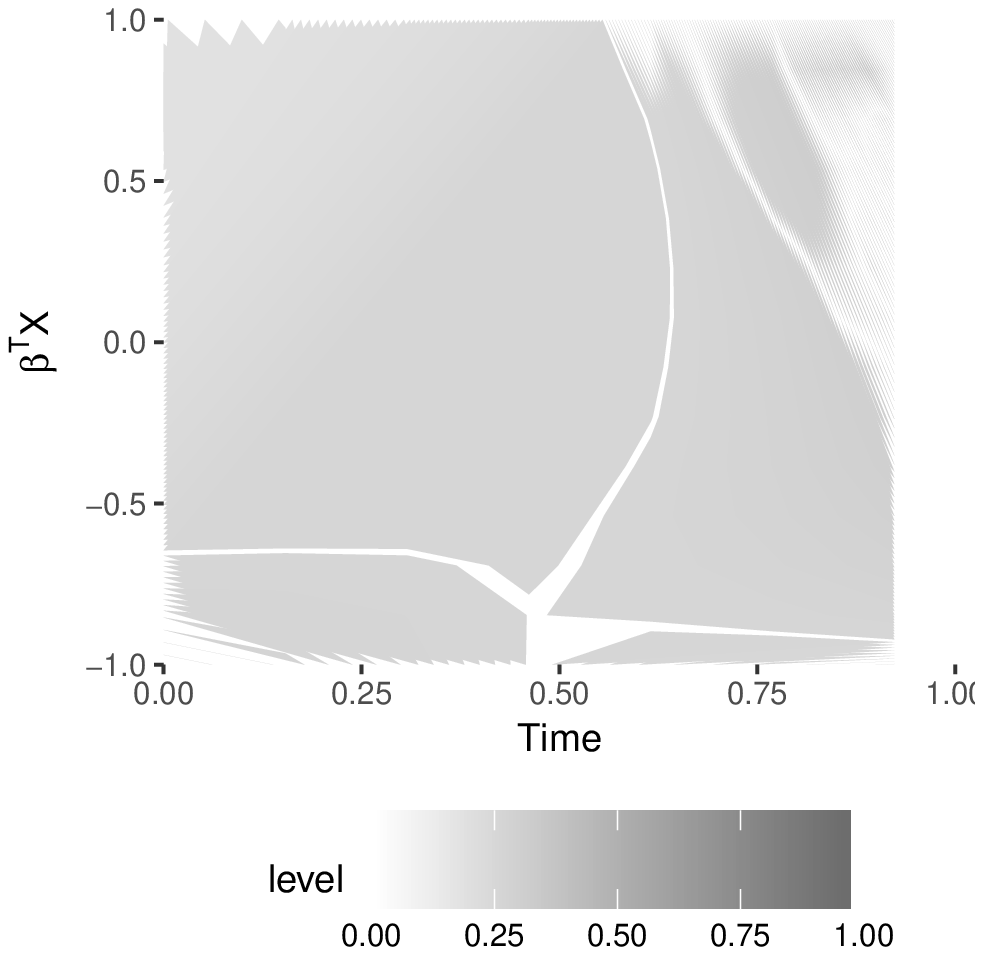}
	\includegraphics[width=5cm]{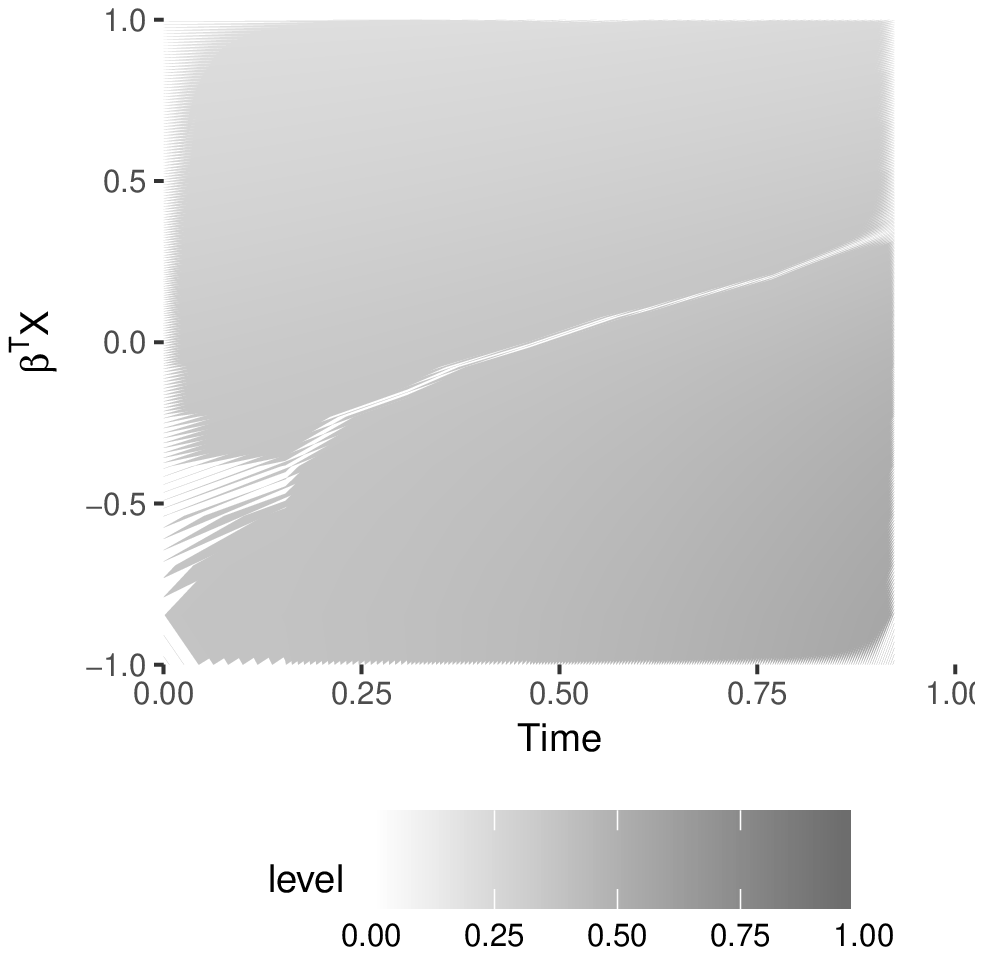}	
	\caption{Performance of the semiparametric method on mean
		residual life
		function of study 1.
		First row: contour plot of true $m(t,\bb\trans\X)$;
		Second
		row: contour plot of $\wh m(t,\bb\trans\X)$;
		Third row: contour plot of $|\wh
		m(t,\bb\trans\X)-m(t,\bb\trans\X)|$. Left to right
		columns: no censoring; 20\% censoring rate; 40\%
		censoring rate.
	}
	\label{fig:simu1contour}
\end{figure}

The results for the estimation of
$\bb$ under Study 1 are given in Table
\ref{tab:simu1}, with three
censoring rates, 0\%, 20\% and 40\%. 
 The proposed 
method has much smaller biases and standard deviations, whereas
``PM1",``PM2" 
and ``additive"  are 
biased with larger standard deviations. 
The performances of all of the estimators deteriorate when the censoring rate increases, though our method still outperforms the others. 
We also 
demonstrate the true and estimated mean residual life functions 
in Figure \ref{fig:simu1contour}.
The error plots, shown in the last row of Figure
\ref{fig:simu1contour}, demonstrate
that our method fare well for
estimating $m(t,\bb\trans\x)$ when $t$ is not too large.
The contour plots reveal that bias
increases as censoring rate increases and the estimation 
deteriorates when $t$ is large.
We further illustrate the dependence of $\wh
m(t,\bb\trans\x)$ on
$t$ at fixed $\bb\trans\x$ and the dependence of $\wh
m(t,\bb\trans\x)$ on
$\bb\trans\x$ at fixed $t$ values in
Figure \ref{fig:simu1curve}.
These results show an overall satisfactory  performance of our semiparametric method. 
As  Figure \ref{fig:simu1contour} reveals, the performance of our method is better when $t$ is in the
interior of the range because more observations are
available for
the local estimation, as opposed to a larger $t$ with fewer observations still available. In contrast, 
 regardless of the magnitude of $t$, the mean   residual 
 life function estimated by ``PM1'', ``PM2'' and ``additive'' is severely biased, as shown 
in the last three rows from Figure \ref{fig:simu1curve}. This is because 
these models assume a pre-determined
functional
form of the mean residual life, which in this case is misspecified.

Tables \ref{tab:simu2}
and \ref{tab:simu3} report the results of Studies 2 and 3  in relate to $\wh\bb$,
respectively. 
We also provide the average
estimation result
of $m(t,\bb\trans\x)$ using a contour plot
in Figure \ref{fig:simu2contour} and Figure 
\ref{fig:simu3contour}. 
We further plot $\wh m(t,\bb\trans\x)$
as a function of $\bb\trans\x$ at
fixed $t$
and as a function of $t$ at fixed $\bb\trans\x$ respectively in
Figure \ref{fig:simu2curve} and Figure \ref{fig:simu3curve}.
Similar to the conclusion in the first simulation
study, the performance of estimating $\bb$ by our proposed estimator 
is satisfactory. The performance of the mean residual 
life estimation is
better when $t$ is smaller,  deteriorates when
$t$ becomes larger, and   is better for smaller censoring rates. 
When $t$ is fixed, $\wh m(t,\bb\trans\x)$ from the proposed 
method has a stable
performance in the entire range of $\bb\trans\x$ for all of the 
censoring rates. We also show
$\wh m(t,\bb\trans\x)$ as a function of
$t$ when $\bb\trans\x$ is fixed, in which the proposed estimator is 
consistent, while ``PM1'', ``PM2'' and ``additive'' deviate 
from the true curve.

To recap,  for estimating $\bb$, 
 the proposed method yields much smaller biases 
than ``PM1'', ``PM2'' and ``additive''; for estimating the mean 
residual life function, ``PM1'', ``PM2'' and ``additive'' 
estimators deviate much from the truth when the model is 
misspecified, but the proposed estimator recovers the truth well.

\section{Analysis of the Kidney Transplant Data}\label{sec:app}

We apply the proposed method to analyze a kidney transplant data 
set from the
U.S. Scientific Registry of Transplant Recipients. The registry is
maintained by the United Network for Organ Sharing and Organ Procurement and
Transplantation Network (UNOS/OPTN) and includes all waitlisted
kidney transplant candidates and transplant recipients in the United
States (\url{https://unos.org/}). To evaluate the possible benefit of  transplantation, we 
use the residual life to estimate how much longer a patient can 
survive if he receives a transplant than otherwise.

To avoid the confounding cohort effect,  we consider the patients who were waitlisted in the same year of 2011. There were 42,217 patients in this cohort with an average followup  of 907 days after waitlisting. During the followup, a  total of  22,295
patients received kidney transplants. The
response variable is the survival time in days
($T_i$) starting from waitlisting.
In the transplant group, 5.82\% of the observations were
censored, while in the non-transplant group, the censoring rate was
26.27\%. The covariates $\X_{app}$ included in our analysis are
gender ($X_{1}$), race ($X_{2}$), cold ischemia time
($X_{3}$),
insurance coverage ($X_{4}$), body mass index ($X_{5}$),
diagnosis
type ($X_{6}$), peak PRA/CPRA ($X_{7}$), previous malignancy
status
($X_{8}$) and diabete indicator ($X_{9}$).
In the transplant group of patients, we further
included the waiting time  $X_{w}$, calculated as the difference
between operation  and waitlisting. The hypothesis to test is that a patient with more 
prompt transplant operation may get more benefit.
The goal of this
study is to
quantify the potential residual life increment if a patient
receives a kidney transplant given the covariate profile.

As the survival trajectories for patients who received transplants might differ from 
those for patients who did not, we  analyze the associations between the residual survival time
and the covariates using model (\ref{eq:model}) for the transplant and non-transplant groups, separately. 
To proceed,  we first 
determine the number of indices $d$
through Validated
Information Criterion (VIC) \citep{ma2015validated},
where the smallest VIC value corresponds to the selected $d$ value.
In our analysis,  $d=1$ is chosen based on this criterion for both models and  we fix $d=1$ for these two models.  The results are reported in  Table
\ref{tab:app1}.
The index vector  is normalized by
fixing the first component (gender) at 1, hence only 9 coefficient estimates are
reported in the transplant group and 8  in the non-transplant group. 

We first focus on the model for the transplant group.
Our
	analysis finds that several covariates, such as the
	body mass index ($X_{5}$), have no significant effects
	on the mean residual life,  which agrees with the
	previous  studies 
	\citep{friedman2003demographics}. The waiting time $X_{w}$ turns out not
  significant in the model for the transplant group.
To confirm this, we investigate the possible  confounding 
between $X_{w}$ and other covariates. Specifically, we exclude
$X_{w}$ and perform the same analysis for the operation group 
using
the remaining 9 covariates. The VIC results in dimension $d=1$ as 
well and we report the results in Table
\ref{tab:app1}. We then examine the dependence between $X_w$ and 
the
index
$\wh\bb\trans\X_{app}$, where $\wh\bb$ contains the values reported 
in the
middle part of Table
\ref{tab:app1}, via the distance correlation test
	\citep{szekely2014partial}. The $p$-value is less than
	$0.001$, indicating a significant confounding effect. That is, the correlation between  $X_w$ and the index $\wh\bb\trans\X_{app}$
 masks the direct link between 
the waiting time and expected residual life.
Indeed, the 
	waiting time affects survival  
	benefits in a very complex way 
	\citep{meier2000effect,gill2005impact} and hence might infer  
	complicated strategies for
	organ allocation \citep{meier2002}.
	
We also estimate the mean residual
function, i.e. $\wh m_{\rm treat}(t, \wh\bb_{\rm 
treat}\trans\X_{app})$,
where the subscript ``treat'' indicates
the transplant group.
We estimate $\wh m_{\rm treat}(t, \wh\bb_{\rm
treat}\trans\X_{app})$  with and without including $X_{w}$  
as a covariate. As the numerical results
are almost identical,  we only focus on the result without
including $X_{w}$ as a covariate in our later discussion.

We also present the model results for the patients that did not receive kidney transplantation
during the followup in Table \ref{tab:app1},
which reveals that 
 the sets of significant variables differ across the  transplant and non-transplant groups.  For example, BMI ($X_{5}$), which is not 
	significant in the transplant group, has a positive effect 
	on 
	the mean residual life in the non-transplant group, which 
	has been noted in the literature
	\citep{kalantar2005survival}.
	 Another
	example is diagnosis type ($X_{6}$), though not
	significant in the
	non-transplant group, 
	has a significantly negative effect in the transplant group. 
	This can be
	easily understood
	because hypertension 
	\citep{frei1995pre}, polycystic kidney disease 
	\citep{kasiske2003diabetes}
	and vascular disease 
	\citep{grimm1997proteinuria} are high-risk 
	factors for the post-transplant survival benefits.
	The last example is the previous malignancy status 
	($X_{8}$),  which has opposite effects in 
	the non-transplant and transplant groups. We
	believe this is caused by the correlation between covariates.
Specifically, 	the type of malignancies 
	\citep{Brattstroem2013},
	the waiting time between malignancy 
	diagnosed and transplantation 
	\citep{pennI1997} and the
	malignancy treatment before transplantation
	\citep{oechslin1996pretransplant}, which 
	all affect patients' survival,  are highly 
	correlated 	with $X_8$.
	
\begin{figure}[H]
	\centering
	\includegraphics[trim={0 0 0
		0},width=6cm]{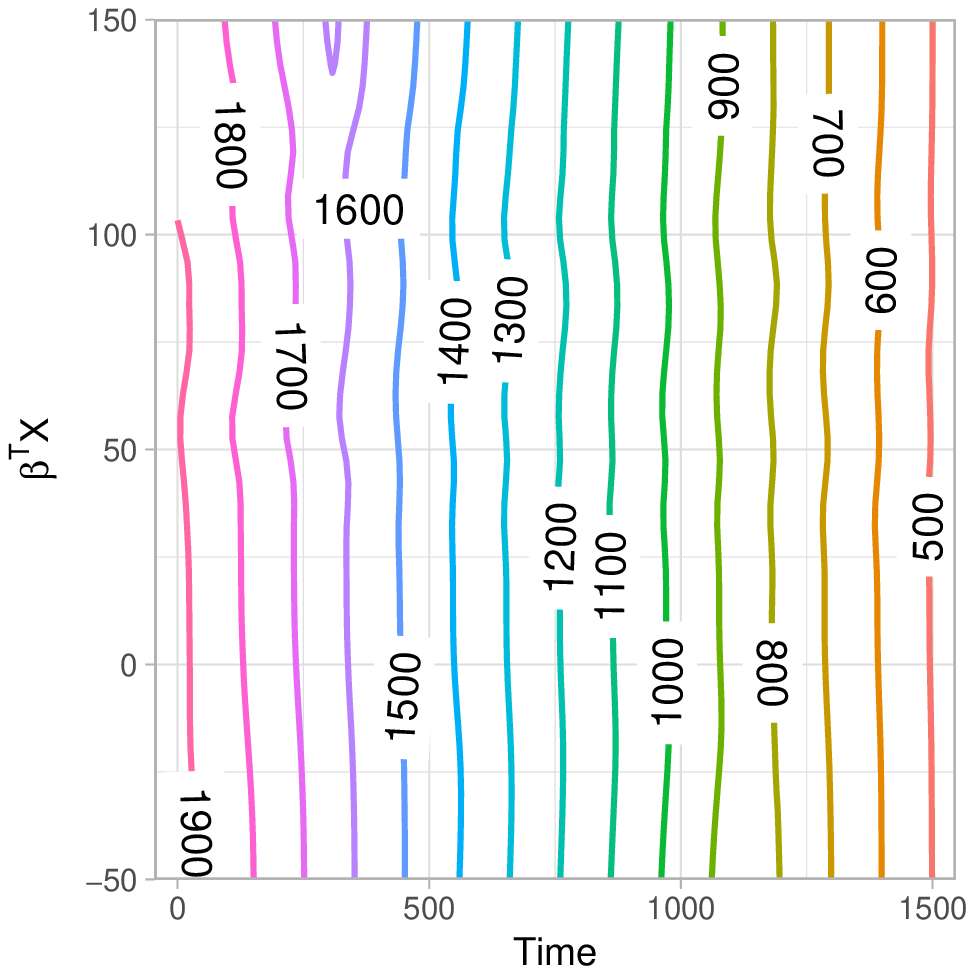}
	\includegraphics[trim={0 0 0
		0},width=6cm]{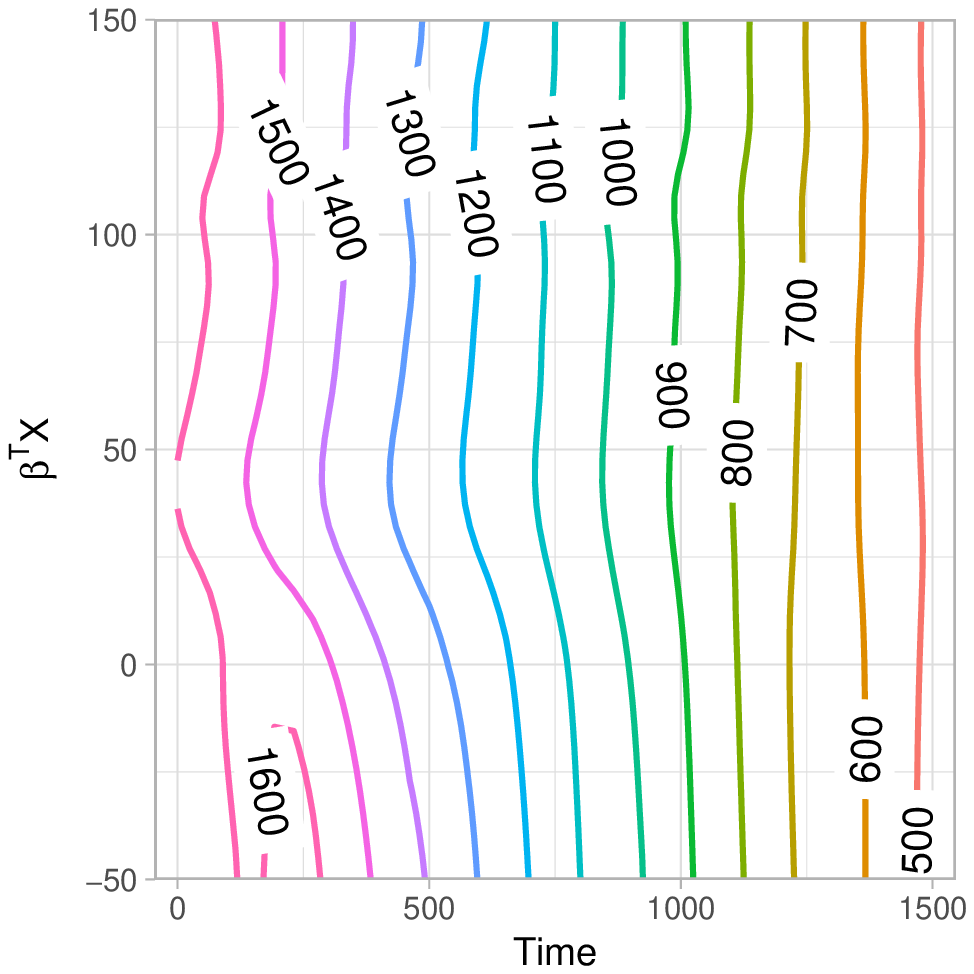}
	\caption{Mean residual life function in kidney transplant 
		application study. Left panel: Patient receive kidney
		transplant. Right panel: Patients did not received kidney
		transplant.}
	\label{fig:appcontour}
\end{figure}

Denote the  estimated mean residual life time by
$\wh m_{\rm wait}(t, \wh\bb_{\rm wait}\trans\X_{app})$, where
the
subscript
``wait'' indicates
the waiting (or the non-transplant) group.
The contour plot of both $\wh m_{\rm treat}(t,\wh\bb_{\rm
  treat}\trans\X_{app})$ and $\wh m_{\rm wait}(t, \wh\bb_{\rm
  wait}\trans\X_{app})$ are given in Figure \ref{fig:appcontour}.

\begin{figure}[H]
	\centering
	\includegraphics[trim={0 0 0
		0},width=4.5cm]{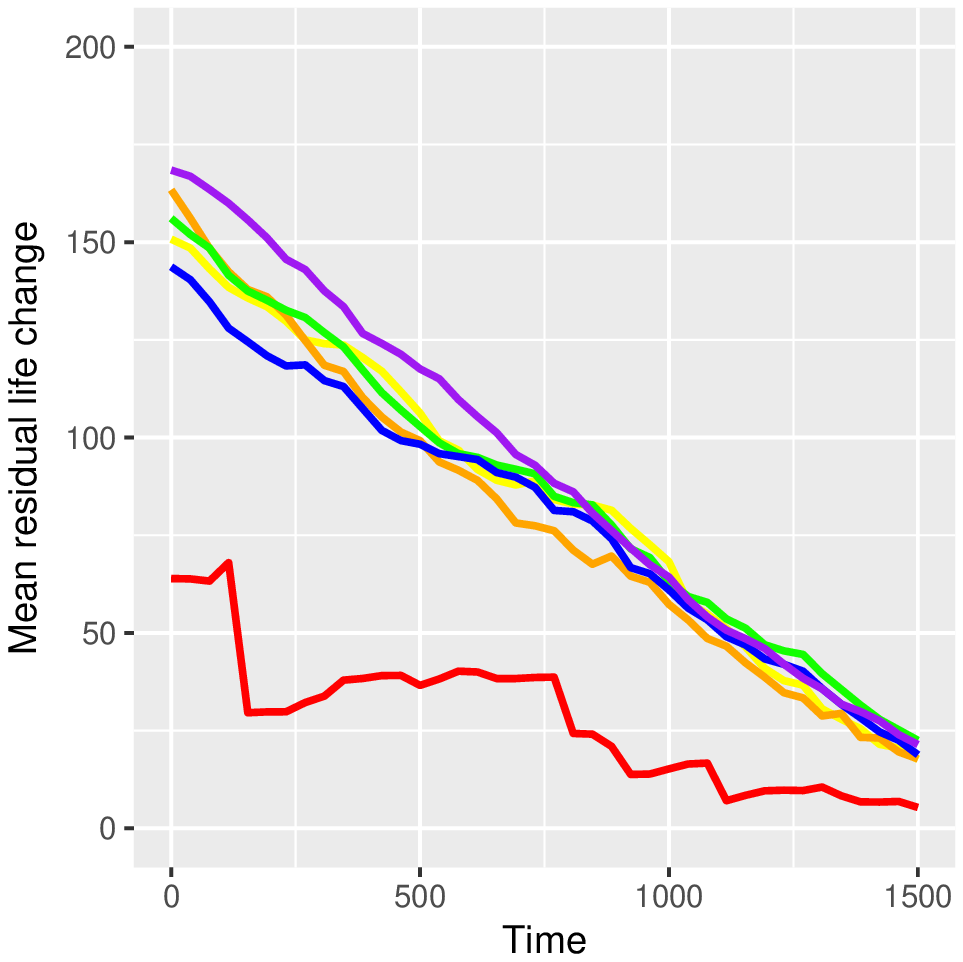}
	\includegraphics[trim={0 0 0
		0},width=4.5cm]{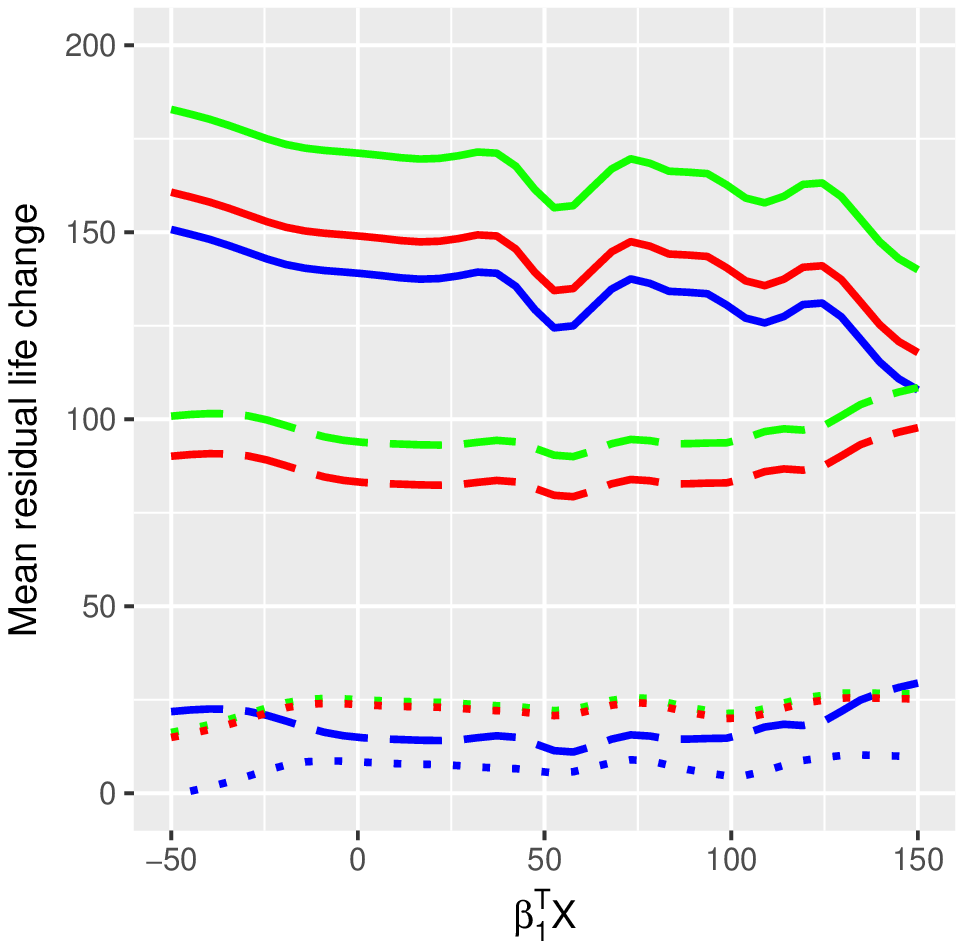}
	\includegraphics[trim={0 0 0
		0},width=4.5cm]{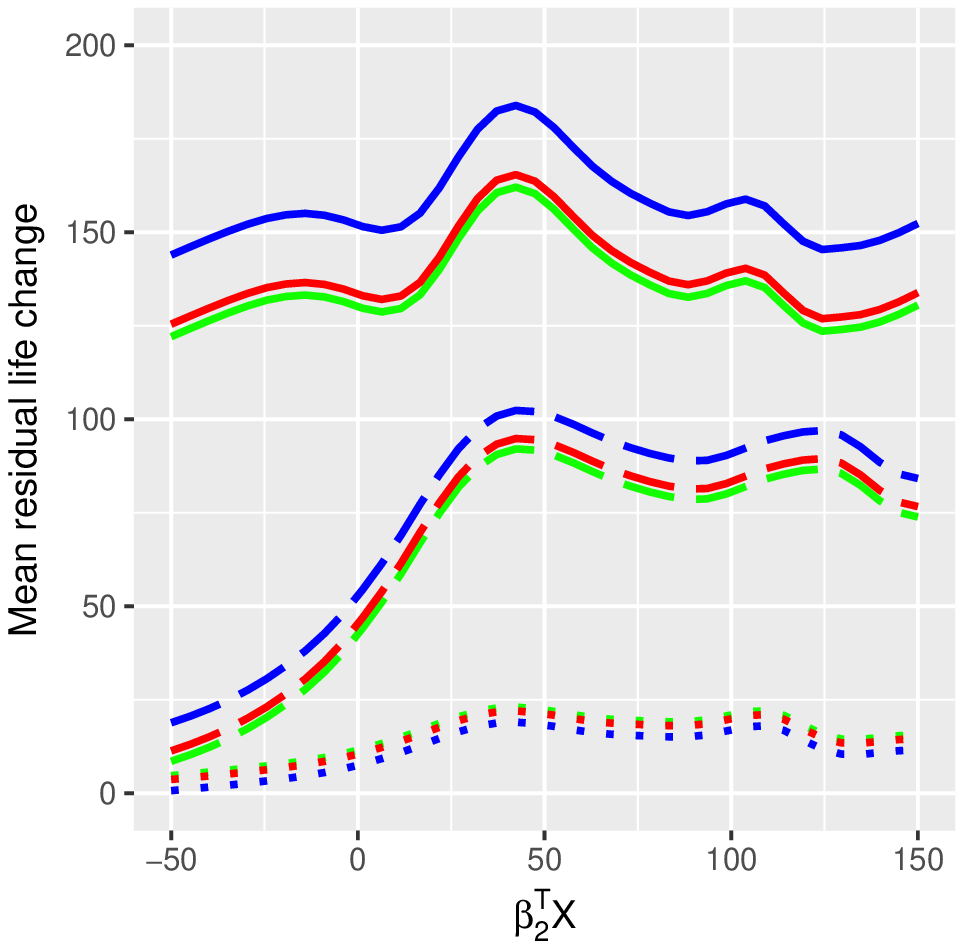}
	\caption{Mean residual life difference between receiving and
		not receiving transplant. Left panel: time is the
		variable
		at the following setting of $[X_3,X_5,X_6,X_8]$, red:
		$[100,50,2,3]$, blue: $[20,45,9,1]$, orange:
		$[80,30,4,1]$,
		yellow: $[50,30,2,1]$, green: $[50,15,7,1]$. The rest
		variables remain the same at
		$[X_1,X_2,X_4,X_7,X_9]=[1,1,1,50,0]$. Middle panel:
		$\bb_{\rm treat}\trans\X$ is the variable with the
		following
		setting,
		blue: $\bb_{\rm wait}\trans\X=-40$, green: $\bb_{\rm
			wait}\trans\X=50$ ,
		red:
		$\bb_{\rm wait}\trans\X=100$, solid: $t=100$, dashed:
		$t=500$,
		dotted: $t=1000$. Right panel: $\bb_{\rm wait}\trans\X$
		is the
		variable with the following setting, blue:
		$\bb_{\rm treat}\trans\X=-40$, green: $\bb_{\rm
			treat}\trans\X=50$ , red:
		$\bb_{\rm treat}\trans\X=100$, solid: $t=100$, dashed:
		$t=500$,
		dotted: $t=1000$.}
	\label{fig:appcurve1}
\end{figure}

Given a patient with 
$\X$ and alive at time $t$,
 $\wh m_{\rm treat}(t, \wh\bb_{\rm 
treat}\trans\X_{app})-\wh
m_{\rm wait}(t,
\wh\bb_{\rm wait}\trans\X_{app})$
provides an estimate of his/her mean residual life time
difference between
receiving and not receiving kidney transplant. Because the
difference is a function of three variables, namely $t$,
$\wh\bb_{\rm treat}\trans\X_{app}$ and $\wh\bb_{\rm 
wait}\trans\X_{app}$, we present
the
difference using various
plots. Figure \ref{fig:appcurve1}  plots curves that
change with $t$ at several fixed
$(\wh\bb_{\rm treat}\trans\X_{app}, \wh\bb_{\rm 
wait}\trans\X_{app})$ values, curves
that change
with $\wh\bb_{\rm treat}\trans\X_{app}$
at several fixed
$(t, \wh\bb_{\rm wait}\trans\X_{app})$ values, and
curves that change
with $\wh\bb_{\rm wait}\trans\X_{app}$
at several fixed
$(t, \wh\bb_{\rm treat}\trans\X_{app})$ values.
Further, Figure \ref{fig:appcontour2} plots the contour
of
the mean residual life time as a function of
$(\wh\bb_{\rm treat}\trans\X_{app}, \wh\bb_{\rm 
wait}\trans\X_{app})$ at a fixed
time
$t$,
as a function of
$(t, \wh\bb_{\rm wait}\trans\X_{app})$ at a fixed index value
$\wh\bb_{\rm treat}\trans\X_{app}$ and
as a function of
$(t, \wh\bb_{\rm treat}\trans\X_{app})$ at a fixed index value
$\wh\bb_{\rm wait}\trans\X_{app}$.

In summary, we find that when the waiting
time is less than 500 days, regardless of the
index values $\wh\bb_{\rm treat}\trans\X_{app}$ and
$\wh\bb_{\rm wait}\trans\X_{app}$,
150 to 300 more days on average  can be gained if a patient
receives transplantation; with the
waiting time between 500
and 1000 days, kidney transplants can still lead to a reasonably
large
improvement if the patient's
index value $\bb_{\rm wait}\trans\X_{app}$ is positive (indicating overall good health), 
regardless of
his/her
$\bb_{\rm treat}\trans\X_{app}$
value. However, if the  waiting time is more  than 1000 days or the waiting time is more than 500 days but with a negative 
index $\bb_{\rm wait}\trans\X_{app}$  (indicating  poor health conditions), the benefit of kidney transplant is the least.

\section{Discussion}\label{sec:conclusion}
The work is to address 
a severe shortage of organs that are needed to sustain ESRD patients' life and aims to design a feasible strategy to increase
the potential efficiency brought by each available kidney. Instead of
evaluating the patients' expected survival time,  as is done
in current strategy, we
propose to consider
the potential residual life prolonged by kidney transplant. We
compare the patients' expected residual life with and without
transplant and use the difference to gauge the potential benefit 
gained from the transplant. Patients with larger
differences may  have a higher
priority for organ allocation than those
 with smaller incremental values.

 We have proposed a  flexible semiparametric regression model of mean residual life, which  relaxes the parametric assumptions on the dependence of mean residual life on covariates and
	how long a patient has lived.
	To strike a balance
between interpretation and flexibility,  our procedure also enables reduce the covariate dimension $p$ to $d$:
when $d=1$, the  model falls to the  single index model, while $d=p$ corresponds to a completely
nonparametric model. We suggest to use the Validated
Information Criterion  \citep{ma2015validated} to choose $d$,
 which seems to fare well in practice.


\newpage
\LTcapwidth=\textwidth
\spacing{1.5}
\begin{longtable}{cc|ccccccccc}
	\caption{Results of  study 1, based on 1000
		simulations with sample size 500. ``Prop.'' is the
		semiparametric method, ``PM1'' and ``PM2'' are the
		proportional mean residual life
		methods, ``additive'' is the additive
		method. ``bias'' is the average absolute bias of
		each component in $\wh\bb$, ``sd''
		is the
		sample standard deviation 
		of the
		corresponding estimators, ``MSE" is the mean squared 
		error. }
	\label{tab:simu1}\\\hline
	\vspace{.2cm}
		&& $\beta_2$ &$\beta_3$ &$\beta_4$ &$\beta_5$
		&$\beta_6$      &$\beta_7$  &$\beta_8$  &$\beta_9$
		&$\beta_{10}$\\
		&true& $-0.75$ &$0.4$ &$-0.2$ &$0.15$ &$0.0$
		&$-0.15$ &$0.2$
		&$-0.4$ &$0.75$ \\\hline
		&&\multicolumn{9}{c}{No censoring}\\
		Prop.	& bias&0.002 &0.004 &0.001 &0.000  &0.002 &0.005 
		&0.003 &0.002 &0.003\\
		& sd&0.091 &0.090 &0.092 &0.091 &0.088 &0.083 &0.098 
		&0.081 &0.092\\
		&MSE&0.008 &0.008 &0.008 &0.008 &0.008 &0.007 &0.01  
		&0.006 &0.008\\
		PM1	& bias &0.057 &0.027 &0.014 &0.018 &0.003 &0.025 
		&0.019 &0.004 &0.059\\
		& sd   &0.383 &0.319 &0.299 &0.300 &0.307 &0.307 &0.296 
		&0.323 &0.382\\
		&MSE&0.152 &0.105 &0.091 &0.114 &0.103 &0.097 &0.088 &0.105 
		&0.151\\
		PM2	& bias&0.129 &0.069 &0.004 &0.043 &0.017 &0.029 
		&0.014 &0.046 &0.153\\
		& sd &0.824 &0.736 &0.695 &0.665 &0.649 &0.626 &0.631 
		&0.687 &0.761\\
		&MSE&0.532 &0.392 &0.431 &0.377 &0.422 &0.358 &0.311 &0.397 
		&0.513\\
		additive& bias&0.066 &0.023 &0.017 &0.017 &0.004 &0.026 
		&0.021 &0.012 &0.060 \\
		& sd&0.396 &0.339 &0.322 &0.316 &0.346 &0.337 &0.335 
		&0.340  &0.395\\
		&MSE &0.184 &0.127 &0.114 &0.124 &0.137 &0.114 &0.115 &0.119 
		&0.184\\
		&&\multicolumn{9}{c}{20\% censoring}\\		
		Prop.	& bias&0.002 &0.016 &0.000  &0.001 &0.003 &0.003 
		&0.001 &0.006 &0.001\\
		& sd&0.122 &0.127 &0.124 &0.123 &0.116 &0.118 &0.116 
		&0.124 &0.120 \\
		&MSE&0.015 &0.016 &0.015 &0.015 &0.013 &0.014 &0.013 
		&0.015 &0.014\\
		PM1	& bias&0.098 &0.054 &0.029 &0.020  &0.006 &0.022 
		&0.065 &0.032 &0.087\\
		& sd&0.569 &0.444 &0.413 &0.419 &0.411 &0.430 &0.435 
		&0.442 &0.500\\
		&MSE&0.333 &0.200 &0.172 &0.176 &0.168 &0.185 &0.193 
		&0.197 &0.257\\
		PM2	& bias&0.136 &0.071 &0.009 &0.076 &0.025 &0.058 
		&0.054 &0.097 &0.131\\
		& sd  &1.548 &1.132 &1.113 &1.147 &1.171 &1.194 &1.110 
		&1.333 &1.545\\
		&MSE  &2.413 &1.284 &1.237 &1.320 &1.371 &1.427 &1.233 
		&1.784 &2.401\\
		additive& bias&0.094 &0.052 &0.033 &0.011 &0.007 &0.024 
		&0.056 &0.038 &0.070\\
		& sd &0.512 &0.403 &0.413 &0.372 &0.388 &0.409 &0.386 
		&0.435 &0.461\\
		&MSE  &0.271 &0.165 &0.171 &0.139 &0.151 &0.168 &0.152 
		&0.191 &0.217\\
		&&\multicolumn{9}{c}{40\% censoring}\\		
		Prop.	& bias&0.003 &0.016 &0.005 &0.004 &0.010 &0.010 
		&0.007 &0.000  &0.005\\
		& sd   &0.175 &0.190 &0.155 &0.169 &0.204 &0.181 &0.161 
		&0.140 &0.188\\
		&MSE&0.030 &0.036 &0.024 &0.029 &0.042 &0.033 &0.026 
		&0.02  &0.035\\
		PM1	& bias&0.121 &0.067 &0.049 &0.031 &0.000  &0.026 
		&0.046 &0.037 &0.106\\
		& sd  &0.734 &0.482 &0.541 &0.394 &0.410 &0.429 &0.456 
		&0.447 &0.606\\
		&MSE&0.553 &0.237 &0.295 &0.156 &0.168 &0.184 &0.209 
		&0.201 &0.378\\
		PM2	& bias&0.232 &0.068 &0.063 &0.129 &0.066 &0.242 
		&0.047 &0.200 &0.186\\
		& sd   &5.141 &3.178 &4.173 &3.268 &4.487 &4.783 &3.618 
		&4.341 &5.914\\
		&MSE&26.45 &10.09 &17.39 &10.68 &20.11 &22.91 &13.08  
		&18.86 &34.97\\
		additive& bias&0.107 &0.063 &0.042 &0.036 &0.017 &0.027 
		&0.033 &0.041 &0.101\\
		& sd  &0.566 &0.462 &0.475 &0.398 &0.376 &0.427 &0.420 
		&0.420 &0.625\\
		&MSE  &0.332 &0.217 &0.227 &0.160 &0.141 &0.183 &0.178 
		&0.178 &0.401\\
		\hline
\end{longtable}

\newpage

\begin{longtable}{c|ccccccccc}
	\caption{Analysis of the kidney transplant data, where ``SE''
		is the estimated standard error of $\wh\bb$}
	\label{tab:app1}\\\hline
		& $\wh\bb_{2}$& $\wh\bb_{3}$&
		$\wh\bb_{4}$&$\wh\bb_{5}$& $\wh\bb_{6}$& $\wh\bb_{7}$&
		$\wh\bb_{8}$& $\wh\bb_{9}$& $\wh\bb_{w}$\\\hline
		&\multicolumn{9}{c}{Transplant including waiting time
			$X_w$}\\
		Estimate&
		-0.928& 0.488&1.045& 0.016&-1.057&1.006&0.475&
		-1.041&-0.053 \\
		SE& 0.293&0.254&0.315&0.231& 0.277& 0.294& 0.258&
		0.338&0.194 \\
		p-value & $<0.001$ &0.027& $<0.001$& 0.471& $<0.001$&
		$<0.001$& 0.033& 0.001&0.392\\		\hline
		&\multicolumn{9}{c}{Transplant excluding waiting time
			$X_w$}\\
		Estimate&
		-0.875& 0.477& 0.912& 0.061&-0.779& 0.966& 0.379&-0.958\\
		SE& 0.280&0.109& 0.218& 0.155&0.360&0.263& 0.162& 0.269\\
		p-value & 0.002& $<0.001$& $<0.001$&0.694&0.030&
		$<0.001$&0.019& $<0.001$\\		\hline
		&\multicolumn{9}{c}{Non-transplant}\\
		Estimate &
		-0.875&  1.016&  1.116& 1.061&-0.033&
		0.995&-1.106&-1.123&\\
		SE& 0.221& 0.263& 0.363& 0.338& 0.158& 0.218& 0.432&
		0.332&\\
		p-value &  $<0.001$&  $<0.001$& 0.002&0.002&0.834&
		$<0.001$&0.011& $<0.001$&\\	\hline
\end{longtable}

\newpage

\section*{Supplement to
``Semiparametric regression of mean residual life  with censoring and covariate dimension reduction"}
\setcounter{equation}{0}\renewcommand{\theequation}{S.\arabic{equation}}
\setcounter{subsection}{0}\renewcommand{\thesubsection}{S.\arabic{subsection}}
\setcounter{figure}{0}\renewcommand{\thefigure}{S.\arabic{figure}}
\setcounter{table}{0}\renewcommand{\thetable}{S.\arabic{table}}

\baselineskip=18pt

\subsection{Proof of Proposition \ref{pro:nuisance}  \label{app:nuisance}}
\noindent Proof: 
Let ${\cal T}_1$, ${\cal T}_2$ and ${\cal T}_3$ be the nuisance 
tangent 
spaces  corresponding to $f_\X$, $m(\cdot, \cdot)$ and $\lambda_c$ 
respectively. 
The result of ${\cal T}_1$ follows obviously. 
To obtain ${\cal T}_2$, let 
$m(z,\bb\trans\X)+\bg\trans\h(z,\bb\trans\X)$ be a sub model of 
$m(z,\bb\trans\X)$, where $\h(z,\bb\trans\X)\in{\cal
	R}^{(p-d)d}$ with $\var\{\h(z,\bb\trans\X)\}<\infty,$ 
differentiate the log of  (\ref{eq:pdf}) with 
respect  to $\bg$ and evaluate it at
$\bg=0$. Then, ${\cal T}_2$ is 
\bse 
&&\Delta \left\{\frac{\h_1(z,\bb\trans\X)}{m_1(z,\bb\trans\X)+1}
-\frac{\h(z,\bb\trans\X)}{m(z,\bb\trans\X)}\right\}\\
&&\quad 
-\int_0^z 
\frac{\h_1(s,\bb\trans\X)m(s,\bb\trans\X) 
	-\{m_1(s,\bb\trans\X)+1\}\h(s,\bb\trans\X) 
}{m^2(s,\bb\trans\X)}dt\\
&=&\Delta \left\{\frac{\h_1(z,\bb\trans\X)}{m_1(z,\bb\trans\X)+1}
-\frac{\h(z,\bb\trans\X)}{m(z,\bb\trans\X)}\right\}-\int_0^z 
\left\{\frac{\h_1(s,\bb\trans\X)}{m_1(s,\bb\trans\X)+1}-\frac{\h(s,\bb\trans\X)
 }{m(s,\bb\trans\X)}\right\}\lambda(s,\bb\trans\X)ds\\
&=&\int_0^\infty 
\left\{\frac{\h_1(s,\bb\trans\X)}{m_1(s,\bb\trans\X)+1}-\frac{\h(s,\bb\trans\X)
 }{m(s,\bb\trans\X)}\right\}dM(s,\bb\trans\X) 
\ese 
where $\h_1(z,\bb\trans\X)=\partial\h(z,\bb\trans\X)/\partial 
z$. 

	To obtain ${\cal T}_3$,
	let
	$\lambda_c(t,\X)\{1+\bg\trans\h(t,\X)\}$ be 
 a submodel of $\lambda_{c}(t,\X)$, 
 where $\h(z,\X)\in{\cal
	 R}^{(p-d)d},$ with $\var\{\h(z,\X)\}<\infty$.
	 We then obtain ${\cal T}_3$ as follows 
	\bse
	\frac{\partial\log f(\X,Z,\Delta)}{\partial \bg} | _{\bg=0}
	&=&(1-\Delta)\h(Z,\X)-\int_0^Z\h(s,\X)\lambda_c(s,\X)ds\\
	&=&\int_0^\infty\h(s,\X)dM_c(s,\X),
	\ese
	where $M_c(t,\X)=N_c(t)-\int_0^tI(Z\ge
        s)\lambda_c(s,\X)ds$
        is a martingale process (See Theorem 1.3.2 in \cite{fh1991}). A
        similar result was also established by
        \cite{prentice2003mixed} for a mixed discrete and continuous
        Cox regression model. Because
	$\lambda_c(t,\X)$ can be any positive function,
	$\h(s,\X)$ can be any function. This leads to the form of ${\cal T}_3$.

By taking  conditional expectations given $\X$, 
it follows that ${\cal T}_1\perp{\cal T}_2$ and 
${\cal T}_1\perp{\cal T}_3$. 
Further, ${\cal T}_2\perp{\cal T}_3$ because the martingale 
integrals  associated with $M(s,\bb\trans\X)$ and $M_c(s,\X)$ are 
independent 
conditional on $\X$ due to $T\indep C\mid \X$.
\qed

\subsection{Proof of Equations (\ref{eq:useful}) and
	(\ref{eq:expecteff})}\label{app:eff}
Proof: 
 First, we note that, when $t \le \tau$, 
\bse
E\left\{ \X_lY(t)\mid\bb\trans\X\right\}
&=&E[E\left\{\X_l
	I(T \ge t)I(C \ge t)\mid\bb\trans\X, \X \right\} | \bb\trans\X]
\\
&=&E[E\left\{\X_l
	I(T \ge t)I(C \ge t)\mid\X \right\} | \bb\trans\X]
\\
&=&E\left\{\X_l S(t, \bb\trans\X)
	S_c(t, \X) \mid\bb\trans\X\right\}
\\
&=& S(t, \bb\trans\X) E\left\{\X_l 
	S_c(t, \X) \mid\bb\trans\X\right\},
\ese
where 
the second to last equality holds because of $T\indep C\mid \X$
and that $\pr(T \ge t|\X)$ is a function of $\bb\trans\X$ only.
Similarly, for $t \le \tau$, we obtain that 
$$ E\left\{ Y(t)\mid\bb\trans\X\right\}
= S(t, \bb\trans\X) E\left\{ 
	S_c(t, \X) \mid\bb\trans\X\right\}.
$$
Hence, when $t\le \tau$,
\be \label{junk}
\frac{E\left\{ \X_lY(t)\mid\bb\trans\X\right\}}
{E\left\{Y(t)\mid\bb\trans\X\right\}}
&=&\frac{E\left\{\X_lS_c(t,\X)\mid\bb\trans\X\right\}}
{E\left\{S_c(t,\X)\mid\bb\trans\X\right\}}.
\ee

 Second, when $t>\tau$, $Y(t)=0$ and $S_c( t,\X)=0$. Hence, we have a $0/0$ scenario in (\ref{junk}), in which case,
we define (\ref{junk}) to be \bse
  \frac{E\left\{\X_lS_c(\tau,\X)\mid\bb\trans\X\right\}}
	{E\left\{S_c(\tau, \X)\mid\bb\trans\X\right\}}
   = \frac{E\left\{\X_l p(\X)\mid\bb\trans\X\right\}}
	{E\left\{p(\X)\mid\bb\trans\X\right\}},
	\ese
	a time-invariant constant.
	Here, $p(\X)$ is defined in Section \ref{sec:space}.
Hence, (\ref{eq:useful}) holds  over $[0, \infty)$, with the truth of $\bb$ being $\bb_0$. 

	
With a generic $\bb$, (\ref{junk}) leads to
\bse
&&E\left(\int_0^\infty\g(s,\bb\trans\X)
\otimes\left[\X_l-
\frac{E\left\{\X_l
	S_c(s,\X)\mid\bb\trans\X\right\}}
{E\left\{S_c(s,\X)\mid\bb\trans\X\right\}}\right]
Y(s)\lambda_0(s,\bb\trans\X)ds\right)\\
&=&E\left(\int_0^\infty\g(s,\bb\trans\X)
\otimes\left[E\{\X_lY(s)\mid\bb\trans\X\}-
\frac{E\left\{\X_l
	S_c(s,\X)\mid\bb\trans\X\right\}}
{E\left\{S_c(s,\X)\mid\bb\trans\X\right\}}
E\{Y(s)\mid\bb\trans\X\}\right]\lambda_0(s,\bb\trans\X)ds\right)\\
&=&\0,
\ese
as the quantity inside the square bracket is zero.
In addition, 
\bse
E\left(\int_0^\infty\g(s,\bb\trans\X)
\otimes\left[\X_l-
\frac{E\left\{\X_l
	S_c(s,\X)\mid\bb\trans\X\right\}}
{E\left\{S_c(s,\X)\mid\bb\trans\X\right\}}\right]dM(s,\bb\trans\X)\right)=\0,
\ese
because  $dM(s,\bb\trans\X)=dN(s)-Y(s)\lambda_0(s,\bb\trans\X)ds$ is a martingale.
Therefore, we have 
\bse
E\left(\int_0^\infty\g(s,\bb\trans\X)
\otimes\left[\X_l-
\frac{E\left\{\X_l
	S_c(s,\X)\mid\bb\trans\X\right\}}
{E\left\{S_c(s,\X)\mid\bb\trans\X\right\}}\right]dN(s)\right)=\0.
\ese
Hence,  (\ref{eq:expecteff}) holds with the truth of $\bb$ being $\bb_0$.
\qed

\subsection{Nonparametric Estimators}

\subsubsection{Nonparametric Estimators of Hazard and Mean Residual Life Functions and Their Derivatives}\label{sec:nonpara}

The nonparametric estimators of
$\bLam_2(t,\bb\trans\X)$,
$\lambda(t,\bb\trans\X)$, $\blam_2(t,\bb\trans\X)$, $m_1(t,\bb\trans\X)$, $\m_2(t,\bb\trans\X)$, 
and $\m_{12}(t,\bb\trans\X)$ are
\be
\wh\bLam_{2}(t,\bb\trans\X)&=&
-\sum_{i=1}^n \frac{I(Z_i\le
	t)\Delta_i \K_h'(\bb\trans\X_i-\bb\trans\X)}{\sumj
	I(Z_j\ge
	Z_i)K_h(\bb\trans\X_j-\bb\trans\X)}\n\\
&&+\sum_{i=1}^nI(Z_i\le t)\Delta_iK_h(\bb\trans\X_i-\bb\trans\X)
\frac{\sumj
	I(Z_j\ge Z_i) \K_h'(\bb\trans\X_j-\bb\trans\X)}{\{\sumj
	I(Z_j\ge
	Z_i)K_h(\bb\trans\X_j-\bb\trans\X)\}^2},\n\\\label{eq:Lambda2}\\
\wh{\lambda}(t,\bb\trans\X) &=& \sum_{i=1}^n
\frac{K_b(Z_i-t)\Delta_iK_h(\bb\trans\X_i-\bb\trans\X)}{\sumj I(Z_j\ge
	Z_i)K_h(\bb\trans\X_j-\bb\trans\X)},\label{eq:lambda}\\
\wh\blam_2(t,\bb\trans\X)&=&-\sum_{i=1}^n
\frac{K_b(Z_i-t)\Delta_i\K_h'(\bb\trans\X_i-\bb\trans\X)}{\sumj I(Z_j\ge
	Z_i)K_h(\bb\trans\X_j-\bb\trans\X)}\n\\
&&+\sum_{i=1}^nK_b(Z_i-t)\Delta_iK_h(\bb\trans\X_i-\bb\trans\X)
\frac{\sumj
	I(Z_j\ge
	Z_i)\K_h'(\bb\trans\X_j-\bb\trans\X)}{\{\sumj I(Z_j\ge
	Z_i)K_h(\bb\trans\X_j-\bb\trans\X)\}^2},\n\\
\label{eq:lambda1}
\wh
m_1(t,\bb\trans\X)&=&\wh\lambda(t,\bb\trans\X)e^{\wh\Lambda(t,\bb\trans\X)}
\int_t^\infty
e^{-\wh\Lambda(s,\bb\trans\X)}ds-1,\label{eq:m1}\\
\wh\m_2(t,\bb\trans\X)&=&\wh\bLam_2(t,\bb\trans\X)e^{\wh\Lambda(t,\bb\trans\X)}
\int_t^\infty
e^{-\wh\Lambda(s,\bb\trans\X)}ds\n\\
&&-e^{\wh\Lambda(t,\bb\trans\X)}
\int_t^\infty
\wh\bLam_2(s,\bb\trans\X)e^{-\wh\Lambda(s,\bb\trans\X)}ds,\label{eq:m2}\\
\wh\m_{12}(t,\bb\trans\X)&=&\wh\blam_2(t,\bb\trans\X)e^{\wh\Lambda(t,\bb\trans\X)}
\int_t^\infty
e^{-\wh\Lambda(s,\bb\trans\X)}ds\n\\
&&+\wh\lambda(t,\bb\trans\X)\wh\bLam_2(t,\bb\trans\X)e^{\wh\Lambda(t,\bb\trans\X)}
\int_t^\infty
e^{-\wh\Lambda(s,\bb\trans\X)}ds\n\\
&&-\wh\lambda(t,\bb\trans\X)e^{\wh\Lambda(t,\bb\trans\X)}
\int_t^\infty\wh\bLam_2(s,\bb\trans\X)
e^{-\wh\Lambda(s,\bb\trans\X)}ds.\label{eq:m12}
\ee

\subsubsection{Trimmed Nonparametric Estimators}\label{sec:trimmed}
The trimmed estimators of 
(\ref{eq:expectY}), 
(\ref{eq:expectXY}), (\ref{eq:Lambda}), and (\ref{eq:Lambda2})-(\ref{eq:lambda1}) are 
\be
&&\wh E\left\{Y_i(Z_i)\mid\bb\trans\X_i\right\}\n\\
&=& \frac{\sumj K_h(\bb\trans\X_j-\bb\trans\X_i)I(Z_j\ge Z_i)}
{\sumj K_h(\bb\trans\X_j-\bb\trans\X_i)}I\left\{\frac{1}{n}\sum_{k=1}^n 
K_h(\bb\trans\X_k-\bb\trans\X_i)>d_n\right\},\label{eq:expectYtrim}\\
&&\wh E\left\{\X_{li}Y_i(Z_i)\mid\bb\trans\X_i\right\}\n\\
&=& \frac{\sumj K_h(\bb\trans\X_j-\bb\trans\X_i)\X_{lj}I(Z_j\ge 
	Z_i)}
{\sumj K_h(\bb\trans\X_j-\bb\trans\X_i)}I\left\{\frac{1}{n}\sum_{k=1}^n 
K_h(\bb\trans\X_k-\bb\trans\X_i)>d_n\right\}.\label{eq:expectXYtrim}\\
&&\wh\Lambda(t,\bb\trans\X)\n\\
&=&\sumi\int_{0}^{t}
\frac{K_h(\bb\trans\X_i-\bb\trans\X)}{\sumj
	Y_j(s)K_h(\bb\trans\X_j-\bb\trans\X)}I\left\{\frac{1}{n}\sum_{k=1}^n 
Y_j(s)K_h(\bb\trans\X_k-\bb\trans\X)>d_n\right\}
dN_i(s),\label{eq:Lambdatrim}\\
&&\wh\bLam_2(t,\bb\trans\X)\n\\
&=&-\sumj 
\frac{K_b(Z_j-t) \Delta_j\K_h'(\bb\trans\X_j-\bb\trans\X)}{\sum_{k=1}^n 
	I(Z_k\ge
	Z_j)K_h(\bb\trans\X_k-\bb\trans\X)}I\left\{\frac{1}{n}\sum_{k=1}^n 
I(Z_k\ge
Z_j)K_h(\bb\trans\X_k-\bb\trans\X)>d_n\right\}\nonumber\\
& & + 
\sumj K_b(Z_j-t)\Delta_jK_h(\bb\trans\X_j-\bb\trans\X)\n\\
&&\quad\times\frac{\sum_{k=1}^n I(Z_k\ge
	Z_j)\K_h'(\bb\trans\X_k-\bb\trans\X)}{\{\sum_{k=1}^n  
	I(Z_k\ge 
	Z_j)K_h(\bb\trans\X_k-\bb\trans\X)\}^2}I\left\{\frac{1}{n}\sum_{k=1}^n 
I(Z_k\ge
Z_j)K_h(\bb\trans\X_k-\bb\trans\X)>d_n\right\},\n\\\label{eq:Lambda1trim}\\
&&\widehat{\lambda}(t,\bb\trans\X)\n\\
&=& \sumj 
\frac{K_b(Z_j-t) \Delta_jK_h(\bb\trans\X_j-\bb\trans\X)}{\sum_{k=1}^n 
	I(Z_k\ge
	Z_j)K_h(\bb\trans\X_k-\bb\trans\X)}I\left\{\frac{1}{n}\sum_{k=1}^n 
I(Z_k\ge
Z_j)K_h(\bb\trans\X_k-\bb\trans\X)>d_n\right\},\n\\\label{eq:lambdatrim}\\
&&\wh\blam_2(t,\bb\trans\X)\n\\
&=&-\sumj 
\frac{K_b(Z_j-t) \Delta_j\K_h'(\bb\trans\X_j-\bb\trans\X)}{\sum_{k=1}^n 
	I(Z_k\ge
	Z_j)K_h(\bb\trans\X_k-\bb\trans\X)}I\left\{\frac{1}{n}\sum_{k=1}^n 
I(Z_k\ge
Z_j)K_h(\bb\trans\X_k-\bb\trans\X)>d_n\right\}\nonumber\\
& & + 
\sumj K_b(Z_j-t)\Delta_jK_h(\bb\trans\X_j-\bb\trans\X)\n\\
&&\quad\times\frac{\sum_{k=1}^n I(Z_k\ge
	Z_j)\K_h'(\bb\trans\X_k-\bb\trans\X)}{\{\sum_{k=1}^n  
	I(Z_k\ge 
	Z_j)K_h(\bb\trans\X_k-\bb\trans\X)\}^2}I\left\{\frac{1}{n}\sum_{k=1}^n 
I(Z_k\ge
Z_j)K_h(\bb\trans\X_k-\bb\trans\X)>d_n\right\},\n\\\label{eq:lambda1trim}
\ee 
where $d_n$ satisfies \ref{assum:fbetarelax}.

\subsection{Two useful lemmas}\label{sec:prepare}

\subsubsection{Lemma \ref{lem:pre}}\label{sec:prooflem}
\begin{Lem}\label{lem:pre}
	Under the regularity conditions 
	\ref{assum:kernel}-\ref{assum:survivalfunction}  listed 
	above, 
	\be
	\wh E\left\{Y(Z)\mid\bb\trans\X\right\}
	&=&E\{Y(Z)\mid\bb\trans\X\}+O_p\{(nh)^{-1/2}+h^2\},\label{eq:lemeq1}\\
	\wh E\left\{\X Y(Z)\mid\bb\trans\X\right\}
	&=&E\{\X Y(Z)\mid\bb\trans\X\}+O_p\{(nh)^{-1/2}+h^2\}, 
	\label{eq:lemeq2}\\
	\wh{\lambda}(z,\bb\trans\X)&=&\lambda(z,\bb\trans\X)+O_p\{(nhb)^{-1/2}+h^2+b^2\}\label{eq:lemeq5}\\
	\wh\blam_2(z,\bb\trans\X)&=&\blam_2(z,\bb\trans\X)+O_p\{(nbh^3)^{-1/2}+h^2+b^2\}\label{eq:lemeq6}\\
	\wh{\Lambda}(z,\bb\trans\X)&=&\Lambda(z,\bb\trans\X)+O_p\{(nh)^{-1/2}+h^2\}\label{eq:lemeq7}\\
	\wh\bLam_2(z,\bb\trans\X)&=&\bLam_2(z,\bb\trans\X)+O_p\{(nh^3)^{-1/2}+h^2\}\label{eq:lemeq8}
	\ee
	uniformly for all $z,\bb\trans\X$.
\end{Lem}

\noindent Proof: For notation convenience, we prove the results for $d=1$. We prove 
\bse
	\wh E\left\{\X Y(Z)\mid\bb\trans\X\right\}
&=&E\{\X Y(Z)\mid\bb\trans\X\}+O_p\{(nh)^{-1/2}+h^2\}
\ese
and
\bse
	\wh\bLam_2(z,\bb\trans\X)&=&\bLam_2(z,\bb\trans\X)+O_p\{(nh^3)^{-1/2}+h^2\}.
\ese
and skip the remaining results because of the similar arguments.

First, for any $\X$ and $\bb$ in a local neighborhood of $\bb_0$,
	\be
\frac{1}{n}\sumj 
K_h(\bb\trans\X_j-\bb\trans\X)=f_{\bb\trans\X}(\bb\trans\X)+
O_p(n^{-1/2}h^{-1/2}+h^2),\label{eq:fbeta}
\ee
	To see this, the absolute bias of the left hand 
side of
(\ref{eq:fbeta}) is
\bse
&&\left|E\left\{\frac{1}{n}\sumj 
K_h(\bb\trans\X_j-\bb\trans\X)\right\}-f_{\bb\trans\X}(\bb\trans\X)\right|\\
&=&\left|\int\frac{1}{h}K\left(\frac{\bb\trans\x_j-\bb\trans\X}{h}\right)
f_{\bb\trans\X}(\bb\trans\x_j)d\bb\trans\x_j
-f_{\bb\trans\X}(\bb\trans\X)\right|\\
&=&\left|\int K(u)f_{\bb\trans\X}(\bb\trans\X+hu)du 
-f_{\bb\trans\X}(\bb\trans\X)\right|\\
&=&\left|\int 
K(u)\left\{f_{\bb\trans\X}(\bb\trans\X)+f_{\bb\trans\X}'
(\bb\trans\X)hu+\frac{1}{2}f''_{\bb\trans\X}(\xi)h^2u^2\right\}du
-f_{\bb\trans\X}(\bb\trans\X)\right|\\
&\le&\frac{h^2}{2}\sup_{\bb\trans\x} 
|f''_{\bb\trans\X}(\bb\trans\x)|\int u^2K(u)du,
\ese
where throughout the text, $\xi$ is between $\bb\trans\X$ and $\bb\trans\X+hu$.
The variance is
\bse
&&\var\left\{\frac{1}{n}\sumj 
K_h(\bb\trans\X_j-\bb\trans\X)\right\}\\
&=&\frac{1}{n}\var K_h(\bb\trans\X_j-\bb\trans\X)\\
&=&\frac{1}{n}\left[EK_h^2(\bb\trans\X_j-\bb\trans\X)- 
\left\{EK_h(\bb\trans\X_j-\bb\trans\X)\right\}^2\right]\\
&=&\frac{1}{n}\left[\int\frac{1}{h^2} 
K^2\{(\bb\trans\x_j-\bb\trans\X)/h\}f_{\bb\trans\X} 
(\bb\trans\x_j)d\bb\trans\x_j-f_{\bb\trans\X}^2(\bb\trans\X)+O(h^2)\right]\\
&=&\frac{1}{nh}\int
K^2(u)f_{\bb\trans\X}(\bb\trans\X+hu)du-\frac{1}{n}f_{\bb\trans\X}^2
(\bb\trans\X)+O(h^2/n)
\\
&\le&\frac{1}{nh}f_{\bb\trans\X}(\bb\trans\X)\int 
K^2(u)du
+\frac{h}{2n}\sup_{\bb\trans\X}|f_{\bb\trans\X}''(\bb\trans\X)|\int
u^2K^2(u)du
+\frac{1}{n}|f_{\bb\trans\X}^2(\bb\trans\X)|+O(h^2/n).
\ese
Therefore, applying the central limit theorem, we have that 
\bse
\frac{1}{n}\sumj 
K_h(\bb\trans\X_j-\bb\trans\X)=f_{\bb\trans\X}(\bb\trans\X)+O_p(n^{-1/2}h^{-1/2}+h^2)
\ese
for all
$\bb$ under conditions \ref{assum:kernel}-\ref{assum:fbeta}. Condition \ref{assum:fbeta}
also holds for any $\bb$ in a local neighborhood of $\bb_0$ due to the continuity.
Similarly, We have
\be
\label{eq:fbeta1}
-\frac{1}{n}\sumj 
K_h'(\bb\trans\X_j-\bb\trans\X)=f_{\bb\trans\X}'(\bb\trans\X)+
O_p(n^{-1/2}h^{-3/2}+h^2).
\ee


To show (\ref{eq:lemeq2}), the absolute bias and variance are
\bse
&&\left|
E\left\{\frac{1}{n}\sumj \X_j I(Z_j\ge z)
K_h(\bb\trans\X_j-\bb\trans\X)\right\}
-f_{\bb\trans\X}(\bb\trans\X)E\{\X_j(Z_j\ge
z)\mid\bb\trans\X\}\right|\\
&=&
\left|\frac{h^2}{2}\int
\frac{\partial^2}{\partial(\bb\trans\X)^2}
f_{\bb\trans\X}(\xi) E
\left\{\X_jI(Z_j\ge z)\mid\xi\right\}
u^2 K(u)du\right|\\
&\le&\frac{h^2}{2}\sup_{\bb\trans\X}\left|
\frac{\partial^2}{\partial(\bb\trans\X)^2}
f_{\bb\trans\X}(\bb\trans\X) E
\left\{\X_jI(Z_j\ge z)\mid\bb\trans\X\right\}\right|
\left\{\int u^2 K(u)du\right\},\\
\ese
and 
\bse
&&\var\left\{-\frac{1}{n}\sumj \X_j I(Z_j\ge z)  
K_h(\bb\trans\X_j-\bb\trans\X)\right\}\\
&\le&\frac{1}{nh} 
\sup_{\bb\trans\X}\left|f_{\bb\trans\X}(\bb\trans\X)E\{\X_j\X_j\trans
I(Z_j\ge z)\mid\bb\trans\X\}\right|\int K^2(u)du\\
&&+\frac{h}{2n} 
\sup_{\xi}\left|\frac{\partial^2}{\partial(\bb\trans\X)^2}
f_{\bb\trans\X}(\xi)E\{\X_j\X_j\trans
I(Z_j\ge z)\mid\xi\}\right|\int u^2K^2(u)du+O(1/n)
\ese
under conditions \ref{assum:kernel}-\ref{assum:fbeta}.
So
\be\label{eq:fbetaxi}
\frac{1}{n}\sumj \X_j I(Z_j\ge z)  K_h(\bb\trans\X_j-\bb\trans\X)
&=&f_{\bb\trans\X}(\bb\trans\X)E\{\X_j I(Z_j\ge 
z)\mid\bb\trans\X\}\n\\
&&+O_p(n^{-1/2}h^{-1/2}+h^2)
\ee
under conditions \ref{assum:kernel}-\ref{assum:fbeta}.

To show (\ref{eq:lemeq8}), let 
\bse
\wh\bLam_{21}(z,\bb\trans\X)&=&
-\sum_{i=1}^n \frac{I(Z_i\le 
	z)\Delta_iK_h'(\bb\trans\X_i-\bb\trans\X)}{\sumj 
	I(Z_j\ge
	Z_i)K_h(\bb\trans\X_j-\bb\trans\X)}\\
\wh\bLam_{22}(z,\bb\trans\X)&=&\sum_{i=1}^nI(Z_i\le 
z)\Delta_iK_h(\bb\trans\X_i-\bb\trans\X) \frac{\sumj I(Z_j\ge 
	Z_i)K_h'(\bb\trans\X_j-\bb\trans\X)}{\{\sumj I(Z_j\ge
	Z_i)K_h(\bb\trans\X_j-\bb\trans\X)\}^2}.
\ese
Then 
$\wh\bLam_2(z,\bb\trans\X)=\wh\bLam_{21}(z,\bb\trans\X)+\wh\bLam_{22}(z,\bb\trans\X)$.
To analyze $\wh\bLam_{21}$,
\bse
&&E\wh\bLam_{21}(z,\bb\trans\X)\\
&=&E\left[ \frac{-I(Z_i\le 
	z)\Delta_iK_h'(\bb\trans\X_i-\bb\trans\X)}{f_{\bb\trans\X}(\bb\trans\X)S(Z_i,\bb\trans\X)E\{S_c(Z_i,\X_j)\mid\bb\trans\X_j=\bb\trans\X,Z_i\}}\right]\\
&&+E\left[ \frac{1}{n}\sumi\frac{-I(Z_i\le 
	z)\Delta_iK_h'(\bb\trans\X_i-\bb\trans\X)}{f_{\bb\trans\X}(\bb\trans\X)S(Z_i,\bb\trans\X)E\{S_c(Z_i,\X_j)\mid\bb\trans\X_j=\bb\trans\X,Z_i\}}O_p(A)\right].
\ese
The first term is
\be
&&E\left[ \frac{-I(Z_i\le 
	z)\Delta_iK_h'(\bb\trans\X_i-\bb\trans\X)}{f_{\bb\trans\X}(\bb\trans\X)S(Z_i,\bb\trans\X)E\{S_c(Z_i,\X_j)\mid\bb\trans\X_j=\bb\trans\X,Z_i\}}\right]\n\\
&=&\int I(z_i\le 
z)\frac{\partial\left[f(z_i,\bb\trans\X)E\{S_c(z_i,\X_i)\mid\bb\trans\X,z_i\}f_{\bb\trans\X}(\bb\trans\X)\right]/\partial\bb\trans\X}{f_{\bb\trans\X}(\bb\trans\X)
	S(z_i,\bb\trans\X)E\{S_c( 
	z_i,\X_j)\mid\bb\trans\X,z_i\}}dz_i\n\\
&&-\frac{h^2\partial^3}{6\partial(\bb\trans\X)^3}\iint I(z_i\le 
z)
\frac{f_{\bb\trans\X}(\xi)f(z_i,\xi)E\{S_c(z_i,\X_i)\mid\xi,z_i\}}{f_{\bb\trans\X}(\bb\trans\X)
	S(z_i,\bb\trans\X)E\{S_c( 
	z_i,\X_j)\mid\bb\trans\X,z_i\}}u^3K'(u)dudz_i.\n\\
\label{eq:expectationOf1}
\ee
Hence 
\bse
&&\left|E\left[ \frac{-I(Z_i\le 
	z)\Delta_iK_h'(\bb\trans\X_i-\bb\trans\X)} 
{f_{\bb\trans\X}(\bb\trans\X)S(Z_i,\bb\trans\X)E\{S_c(Z_i,\X_j) 
	\mid\bb\trans\X_j=\bb\trans\X,Z_i\}}\right]\right.\\
&&\left.-\int I(z_i\le z) 
\frac{\partial\left[f(z_i,\bb\trans\X)E\{S_c(z_i,\X_i) 
	\mid\bb\trans\X,z_i\}f_{\bb\trans\X}(\bb\trans\X)\right] 
	/\partial\bb\trans\X}{f_{\bb\trans\X}(\bb\trans\X)S(z_i,\bb\trans\X)
	E\{S_c(z_i,\X_j)\mid\bb\trans\X,z_i\}}dz_i\right|\\
&\le&h^2\sup_{z,\bb\trans\X,\xi}\left| 
\frac{\partial^3}{6\partial(\bb\trans\X)^3}\int\frac{I(z_i\le z) 
	f_{\bb\trans\X}(\xi)f(z_i,\xi)E\{S_c(z_i,\X_i)
	\mid\xi,z_i\}}{f_{\bb\trans\X}(\bb\trans\X) 
	S(z_i,\bb\trans\X)E\{S_c(z_i,\X_j)\mid\bb\trans\X,z_i\}}dz_i\right|\\
&&\times\int |u^3K'(u)|du\\
&=&O(h^2)
\ese
under Condition 
\ref{assum:kernel}-\ref{assum:survivalfunction}. Similarly, we conclude that
\bse
\frac{1}{n}\sumi\frac{-I(Z_i\le 
	z)\Delta_iK_h'(\bb\trans\X_i-\bb\trans\X)}{f_{\bb\trans\X}(\bb\trans\X)S(Z_i,\bb\trans\X)E\{S_c(Z_i,\X_j)\mid\bb\trans\X_j=\bb\trans\X,Z_i\}}O_p(A)=O_p\{h^2+(nh)^{-1/2}\}
\ese
under conditions 
\ref{assum:kernel}-\ref{assum:survivalfunction} due 
to $A=O_p\{h^2+(nh)^{-1/2}\}$. Therefore
\bse
E\wh\bLam_{21}(z,\bb\trans\X)&=&\int I(z_i\le 
z)\frac{\partial\left[f(z_i,\bb\trans\X)E\{S_c(z_i,\X_i)\mid\bb\trans\X,z_i\}f_{\bb\trans\X}(\bb\trans\X)\right]/\partial\bb\trans\X}{f_{\bb\trans\X}(\bb\trans\X)
	S(z_i,\bb\trans\X)E\{S_c( 
	z_i,\X_j)\mid\bb\trans\X,z_i\}}dz_i\\
&&+O\{(nh)^{-1/2}+b^2+h^2\}.
\ese
For $\wh\bLam_{22}$, let 
$B=	-1/n\sumj I(Z_j\ge Z_i)  K_h'(\bb\trans\X_j-\bb\trans\X)
-\partial
f_{\bb\trans\X}(\bb\trans\X)E\{I(Z_j\ge 
Z_i)\mid\bb\trans\X\}/\partial 
\bb\trans\X$,
then
\bse
&&\wh\bLam_{22}(z,\bb\trans\X)\\
&=&-\frac{1}{n}\sumi I(Z_i\le 
z)\Delta_iK_h(\bb\trans\X_i-\bb\trans\X) 
\frac{\partial\left[
	f_{\bb\trans\X}(\bb\trans\X)S(z_i,\bb\trans\X)E\{S_c(z_i,\X_j)\mid\bb\trans\X,z_i\}\right]/
	\partial 
	\bb\trans\X}{f_{\bb\trans\X}^2(\bb\trans\X)S^2(z_i,\bb\trans\X)E^2\{S_c(z_i,\X_j)\mid\bb\trans\X,z_i\}}\\
&&\times\left\{1+O_p(B)+O_p(A)\right\}.
\ese
We have
\be
&&E\left[-\frac{1}{n}\sumi I(Z_i\le 
z)\Delta_iK_h(\bb\trans\X_i-\bb\trans\X) 
\frac{\partial\left[
	f_{\bb\trans\X}(\bb\trans\X)E\{I(Z_j\ge 
	Z_i)\mid\bb\trans\X,Z_i\}\right]/ 
	\partial 
	\bb\trans\X}{f_{\bb\trans\X}^2(\bb\trans\X)E^2\{I(Z_j\ge 
	Z_i)\mid\bb\trans\X,Z_i\}}\right]\n\\
&=&-\iint I(z_i\le z)K_h(\bb\trans\X_i-\bb\trans\X) 
\frac{\partial
	\left[f_{\bb\trans\X}(\bb\trans\X)S(z_i,\bb\trans\X)E\{S_c(z_i,\X_j)\mid\bb\trans\X,z_i\}\right]/
	\partial 
	\bb\trans\X}{f_{\bb\trans\X}^2(\bb\trans\X)S^2(z_i,\bb\trans\X)E^2\{S_c(
	z_i,\X_j)\mid\bb\trans\X,z_i\}}\n\\
&&\times 
f(z_i,\bb\trans\X_i)E\{S_c(z_i,\X_i)\mid\bb\trans\X_i,z_i\}f_{\bb\trans\X}(\bb\trans\X_i)dz_id\bb\trans\X_i\n\\
&=&-\int I(z_i\le z)\frac{\partial
	\left[f_{\bb\trans\X}(\bb\trans\X)S(z_i,\bb\trans\X)E\{S_c( 
	z_i,\X_j)\mid\bb\trans\X,z_i\}\right]/\partial\bb\trans\X} 
{f_{\bb\trans\X}(\bb\trans\X) 
	S^2(z_i,\bb\trans\X)E\{S_c(z_i,\X_j)\mid\bb\trans\X,z_i\}} 
f(z_i,\bb\trans\X)dz_i\n\\
&&-h^2\iint  I(z_i\le z)\frac{\partial
	\left[f_{\bb\trans\X}(\bb\trans\X)S(z_i,\bb\trans\X)E\{S_c( 
	z_i,\X_j)\mid\bb\trans\X,z_i\}\right]/ \partial 
	\bb\trans\X}{f_{\bb\trans\X}^2(\bb\trans\X)S^2(z_i,\bb\trans\X)E^2\{S_c(
	z_i,\X_j)\mid\bb\trans\X,z_i\}}\n\\
&&\times \frac{\partial^2}{2\partial (\bb\trans\X)^2} 
E\{S_c(z_i,\X_i)\mid\xi\}f(z_i,\xi)f_{\bb\trans\X}(\xi)u^2K(u)dz_idu,\label{eq:expectationOf2}
\ee
therefore
\bse
&&\left|E\left[-\frac{1}{n}\sumi I(Z_i\le 
z)\Delta_iK_h(\bb\trans\X_i-\bb\trans\X) \frac{\partial\left[
	f_{\bb\trans\X}(\bb\trans\z)E\{I(Z_j\ge 
	Z_i)\mid\bb\trans\X,Z_i\}\right]/ 
	\partial 
	\bb\trans\X}{f_{\bb\trans\X}^2(\bb\trans\X)E^2\{I(Z_j\ge 
	Z_i)\mid\bb\trans\X,Z_i\}}\right]\right.\\
&&\left.+\int I(z_i\le z)\frac{\partial
	\left[f_{\bb\trans\X}(\bb\trans\X)S(z_i,\bb\trans\X)E\{S_c( 
	z_i,\X_j)\mid\bb\trans\X,z_i\}\right]/ \partial 
	\bb\trans\X}{f_{\bb\trans\X}(\bb\trans\X)S^2(z_i,\bb\trans\X)E\{S_c(
	z_i,\X_j)\mid\bb\trans\X,z_i\}}f(z,\bb\trans\X)dz_i\right|\\
&\le&h^2\sup_{z,z_i,\bb\trans\X,\xi}\left| 
\frac{\partial
	\left[f_{\bb\trans\X}(\bb\trans\X)S(z_i,\bb\trans\X)E\{S_c( 
	z_i,\X_j)\mid\bb\trans\X,z_i\}\right]/ \partial 
	\bb\trans\X}{f_{\bb\trans\X}^2(\bb\trans\X)S^2(z_i,\bb\trans\X)E^2\{S_c(
	z_i,\X_j)\mid\bb\trans\X,z_i\}}\right.\\
&&\times \frac{\partial^2}{2\partial (\bb\trans\X)^2}
E\{S_c(z_i,\X_i)\mid\xi,z_i\}f_{\bb\trans\X}(\xi)\Bigg|\left|\int
u^2K(u)du\right|\\
&=&O(h^2)
\ese
under conditions 
\ref{assum:kernel}-\ref{assum:survivalfunction}. 
Recall $B=O_p(n^{-1/2}h^{-3/2}+h^2)$, then
\bse
&&-\frac{1}{n}\sumi I(Z_i\le 
z)\Delta_iK_h(\bb\trans\X_i-\bb\trans\X) 
\frac{\partial\left[
	f_{\bb\trans\X}(\bb\trans\X)E\{I(Z_j\ge
	Z_i)\mid\bb\trans\X,Z_i\}\right]/ \partial
	\bb\trans\X}{f_{\bb\trans\X}^2(\bb\trans\X)E^2\{I(Z_j\ge
	Z_i)\mid\bb\trans\X,Z_i\}}O_p(B)\\
&=&O_p(n^{-1/2}h^{-3/2}+h^2),\\
&&-\frac{1}{n}\sumi I(Z_i\le 
z)\Delta_iK_h(\bb\trans\X_i-\bb\trans\X) 
\frac{\partial\left[
	f_{\bb\trans\X}(\bb\trans\X)E\{I(Z_j\ge
	Z_i)\mid\bb\trans\X,Z_i\}\right]/ \partial
	\bb\trans\X}{f_{\bb\trans\X}^2(\bb\trans\X)E^2\{I(Z_j\ge
	Z_i)\mid\bb\trans\X,Z_i\}}O_p(A)\\
&=&O_p(n^{-1/2}h^{-1/2}+h^2)
\ese
under condition 
\ref{assum:kernel}-\ref{assum:survivalfunction}. 
Therefore
\bse
E\wh\bLam_{22}&=&-\int I(z_i\le z)\frac{\partial
	\left[f_{\bb\trans\X}(\bb\trans\X)S(z_i,\bb\trans\X)E\{S_c( 
	z_i,\X_j)\mid\bb\trans\X,z_i\}\right]/ \partial 
	\bb\trans\X}{f_{\bb\trans\X}(\bb\trans\X)S^2(z_i,\bb\trans\X)E\{S_c(
	z_i,\X_j)\mid\bb\trans\X,z_i\}}f(z_i,\bb\trans\X)dz_i\\
&&+O(n^{-1/2}h^{-3/2}+h^2)
\ese 
In addition
\bse
&&\int \left(I(z_i\le 
z)\frac{\partial\left[f(z_i,\bb\trans\X)E\{S_c(z_i,\X_j)\mid\bb\trans\X,z_i\}f_{\bb\trans\X}(\bb\trans\X)\right]/\partial\bb\trans\X}{f_{\bb\trans\X}(\bb\trans\X)
	S(z_i,\bb\trans\X)E\{S_c( 
	z_i,\X_j)\mid\bb\trans\X,z_i\}}\right.\\
&&-\left.I(z_i\le z)\frac{\partial
	\left[f_{\bb\trans\X}(\bb\trans\X)S(z_i,\bb\trans\X)E\{S_c( 
	z_i,\X_j)\mid\bb\trans\X,z_i\}\right]/ \partial 
	\bb\trans\X}{f_{\bb\trans\X}(\bb\trans\X)S^2(z_i,\bb\trans\X)E\{S_c(
	z_i,\X_j)\mid\bb\trans\X,z_i\}}f(z_i,\bb\trans\X)\right)dz_i\\
&=&\bLam_2(z,\bb\trans\X).
\ese
Combining $E\wh\bLam_{21}(z,\bb\trans\X)$ and 
$E\wh\bLam_{22}(z,\bb\trans\X)$ gives
\bse
\left|E\wh\bLam_2(z,\bb\trans\X)-\bLam_2(z,\bb\trans\X)\right|
=O(n^{-1/2}h^{-3/2}+h^2)
\ese
under conditions 
\ref{assum:kernel}-\ref{assum:survivalfunction}.

The variance of $\wh\bLam_2(z,\bb\trans\X)$ is
\bse
\var\left\{\wh\bLam_2(z,\bb\trans\X)\right\}
&=&\var\left\{\wh\bLam_{21}(z,\bb\trans\X)+\wh\bLam_{22}(z,\bb\trans\X)\right\}\\
&\le&2\var\{\wh\bLam_{21}(z,\bb\trans\X)\}+2\var\{\wh\bLam_{22}(z,\bb\trans\X)\}.
\ese 
The first term
\bse
&&\var\{\wh\bLam_{21}(z,\bb\trans\X)\}\\
&\le&\frac{2}{n}\var\left[\frac{-I(Z_i\le 
	z)\Delta_iK_h'(\bb\trans\X_i-\bb\trans\X)}{f_{\bb\trans\X}(\bb\trans\X)S(Z_i,\bb\trans\X)E\{S_c(Z_i,\X_j)\mid\bb\trans\X_i=\bb\trans\X,Z_i\}}\right]\\
&&+2\var\left[\frac{1}{n}\sumi\frac{-I(Z_i\le 
	z)\Delta_iK_h'(\bb\trans\X_i-\bb\trans\X)}{f_{\bb\trans\X}(\bb\trans\X)S(Z_i,\bb\trans\X)E\{S_c(Z_i,\X_j)\mid\bb\trans\X_i=\bb\trans\X,Z_i\}}O_p(A)\right]\\
\ese
The first part is
\be
&&\frac{2}{n}\var\left[\frac{-I(Z_i\le 
	z)\Delta_iK_h'(\bb\trans\X_i-\bb\trans\X)}{f_{\bb\trans\X}(\bb\trans\X)S(Z_i,\bb\trans\X)E\{S_c(Z_i,\X_j)\mid\bb\trans\X_i=\bb\trans\X,Z_i\}}\right]\n\\
&=&\frac{2}{n}E\left[\frac{-I(Z_i\le 
	z)\Delta_iK_h'(\bb\trans\X_i-\bb\trans\X)}{f_{\bb\trans\X}(\bb\trans\X)S(Z_i,\bb\trans\X)E\{S_c(Z_i,\X_j)\mid\bb\trans\X_i=\bb\trans\X,Z_i\}}\right]^2\n\\
&&-\frac{2}{n}\left(E\left[\frac{-I(Z_i\le 
	z)\Delta_iK_h'(\bb\trans\X_i-\bb\trans\X)}{f_{\bb\trans\X}(\bb\trans\X)S(Z_i,\bb\trans\X)E\{S_c(Z_i,\X_j)\mid\bb\trans\X_i=\bb\trans\X,Z_i\}}\right]\right)^2\n\\
&\le&\frac{2}{nh^3}\int \frac{I(z_i\le z)f(z_i,\bb\trans\X)}
{f_{\bb\trans\X}(\bb\trans\X)S^2(z_i,\bb\trans\X)E\{S_c(z_i,\X_j)\mid\bb\trans\X,z_i\}}dz_i\left\{\int
K'^2(u)du\right\}\n\\
&&+\frac{1}{nh}\sup_{z,z_i,\bb\trans\X,\xi}\left| 
\frac{\partial^2}{\partial(\xi)^2} 
\frac{f_{\bb\trans\X}(\xi)E\{S_c(z_i,\X_i) 
	\mid\xi,z_i\}}{f_{\bb\trans\X}^2(\bb\trans\X) 
	S^2(z_i,\bb\trans\X)E^2\{S_c(z_i,\X_j)\mid\bb\trans\X,z_i\}}\right|\n\\
&&\times\left\{\int u^2K'^2(u)du\right\}+O(1/n)\n\\
&=&O\{1/(nh^3)\}\label{eq:varianceOf1}
\ee
under condition
\ref{assum:kernel}-\ref{assum:survivalfunction}.
The second part is
\bse
&&2\var\left[\frac{1}{n}\sumi\frac{-I(Z_i\le 
	z)\Delta_iK_h'(\bb\trans\X_i-\bb\trans\X)}{f_{\bb\trans\X}(\bb\trans\X)S(Z_i,\bb\trans\X)E\{S_c(Z_i,\X_j)\mid\bb\trans\X_i=\bb\trans\X,Z_i\}}O_p(A)\right]\\
&\le&2E\left[\frac{1}{n}\sumi\frac{-I(Z_i\le 
	z)\Delta_iK_h'(\bb\trans\X_i-\bb\trans\X)}{f_{\bb\trans\X}(\bb\trans\X)S(Z_i,\bb\trans\X)E\{S_c(Z_i,\X_j)\mid\bb\trans\X_i=\bb\trans\X,Z_i\}}O_p(A)\right]^2\\
&=&\Bigg\{\left(E\left[\frac{-I(Z_i\le 
	z)\Delta_iK_h'(\bb\trans\X_i-\bb\trans\X)}{f_{\bb\trans\X}(\bb\trans\X)S(Z_i,\bb\trans\X)E\{S_c(Z_i,\X_j)\mid\bb\trans\X_i=\bb\trans\X,Z_i\}}\right]\right)^2\\
&&+\frac{1}{n}\var\left[\frac{-I(Z_i\le 
	z)\Delta_iK_h'(\bb\trans\X_i-\bb\trans\X)}{f_{\bb\trans\X}(\bb\trans\X)S(Z_i,\bb\trans\X)E\{S_c(Z_i,\X_j)\mid\bb\trans\X_i=\bb\trans\X,Z_i\}}\right]\Bigg\}O\{1/(nh)+h^4\}\\
&=& O\{1/(nh)+h^4\}
\ese
under conditions 
\ref{assum:kernel}-\ref{assum:survivalfunction}, where the second last equation is because of 
(\ref{eq:expectationOf1}) and 
(\ref{eq:varianceOf1}).
Therefore 
$
\var\{\wh\bLam_{21}(z,\bb\trans\X)\}=O\{1/(nh^3)\}
$
under conditions 
\ref{assum:kernel}-\ref{assum:survivalfunction}.

For $\wh\bLam_{22}(z,\bb\trans\X)$,
\bse
&&\var\{\wh\bLam_{22}(z,\bb\trans\X)\}\\
&\le&\frac{2}{n}\var\left[I(Z_i\le 
z)\Delta_iK_h(\bb\trans\X_i-\bb\trans\X) 
\frac{\partial\left[
	f_{\bb\trans\X}(\bb\trans\X)E\{I(Z_j\ge 
	Z_i)\mid\bb\trans\X\}\right]/ 
	\partial 
	\bb\trans\X}{f_{\bb\trans\X}^2(\bb\trans\X)E^2\{I(Z_j\ge 
	Z_i)\mid\bb\trans\X\}}\right]\\
&&+4\var\left[\frac{1}{n}\sumi I(Z_i\le 
z)\Delta_iK_h(\bb\trans\X_i-\bb\trans\X) \frac{\partial\left[
	f_{\bb\trans\X}(\bb\trans\X)E\{I(Z_j\ge 
	Z_i)\mid\bb\trans\X\}\right]/ 
	\partial 
	\bb\trans\X}{f_{\bb\trans\X}^2(\bb\trans\X)E^2\{I(Z_j\ge 
	Z_i)\mid\bb\trans\X\}}O_p(B)\right]\\
&&+4\var\left[\frac{1}{n}\sumi I(Z_i\le 
z)\Delta_iK_h(\bb\trans\X_i-\bb\trans\X) \frac{\partial\left[
	f_{\bb\trans\X}(\bb\trans\X)E\{I(Z_j\ge 
	Z_i)\mid\bb\trans\X\}\right]/ 
	\partial 
	\bb\trans\X}{f_{\bb\trans\X}^2(\bb\trans\X)E^2\{I(Z_j\ge 
	Z_i)\mid\bb\trans\X\}}O_p(A)\right].
\ese 
The first part is
\be
&&\frac{2}{n}\var\left[I(Z_i\le 
z)\Delta_iK_h(\bb\trans\X_i-\bb\trans\X) 
\frac{\partial\left[
	f_{\bb\trans\X}(\bb\trans\X)E\{I(Z_j\ge 
	Z_i)\mid\bb\trans\X,Z_i\}\right]/ 
	\partial 
	\bb\trans\X}{f_{\bb\trans\X}^2(\bb\trans\X)E^2\{I(Z_j\ge 
	Z_i)\mid\bb\trans\X,Z_i\}}\right]\n\\
&\le&\frac{2}{nh}\int I(z_i\le 
z)f(z_i,\bb\trans\X)\frac{\left(\partial\left[
	f_{\bb\trans\X}(\bb\trans\X)S(z_i,\bb\trans\X)E\{S_c(z_i,\X_j)\mid\bb\trans\X\}\right]/
	\partial 
	\bb\trans\X\right)^2}{f_{\bb\trans\X}^3(\bb\trans\X)S^4(z_i,\bb\trans\X)E^3\{S_c(z_i,\X_j)\mid\bb\trans\X,z_i\}}dz_i\n\\
&&\times \left\{\int K^2(u) du\right\}\n\\
&&+\frac{h}{n}\sup_{z,z_i,\bb\trans\X,\xi} 
\left|\frac{\left(\partial\left[
	f_{\bb\trans\X}(\bb\trans\X)S(z_i,\bb\trans\X)E\{S_c(z_i,\X_i)\mid\bb\trans\X\}\right]/
	\partial 
	\bb\trans\X\right)^2}{f_{\bb\trans\X}^4(\bb\trans\X)S^4(z_i,\bb\trans\X)E^4\{S_c(z_i,\X_i)\mid\bb\trans\X,z_i\}}\right.\n\\
&&\times \frac{\partial^2}{\partial 
	(\xi)^2}E\{S_c(z_i,X_i)\mid\xi\}f_{\bb\trans\X}(\xi)\Bigg|\left\{\int
u^2K^2(u)du\right\}+O(1/n)\n\\
&=&O\{1/(nh)\}\label{eq:varianceOf2}
\ee
under conditions 
\ref{assum:kernel}-\ref{assum:survivalfunction}. 
The second part is 
\bse
&&4\var\left[\frac{1}{n}\sumi I(Z_i\le 
z)\Delta_iK_h(\bb\trans\X_i-\bb\trans\X) 
\frac{\partial\left[
	f_{\bb\trans\X}(\bb\trans\X)E\{I(Z_j\ge 
	Z_i)\mid\bb\trans\X\}\right]/ 
	\partial 
	\bb\trans\X}{f_{\bb\trans\X}^2(\bb\trans\X)E^2\{I(Z_j\ge 
	Z_i)\mid\bb\trans\X,Z_i\}}O_p(B)\right]\\
&\le&4E\left(\left[\frac{1}{n}\sumi I(Z_i\le 
z)\Delta_iK_h(\bb\trans\X_i-\bb\trans\X) \frac{\partial\left[
	f_{\bb\trans\X}(\bb\trans\X)E\{I(Z_j\ge 
	Z_i)\mid\bb\trans\X\}\right]/ 
	\partial 
	\bb\trans\X}{f_{\bb\trans\X}^2(\bb\trans\X)E^2\{I(Z_j\ge 
	Z_i)\mid\bb\trans\X,Z_i\}}\right]^2\right)\\
&&\times O\{1/(nh^3)+h^4\}\\
&=&O\{1/(nh^3)+h^4\}
\ese
under conditions \ref{assum:kernel}-\ref{assum:survivalfunction}, where the second last equation is because of 
(\ref{eq:expectationOf2}) and 
(\ref{eq:varianceOf2}). The last part is
\bse
&&4\var\left[\frac{1}{n}\sumi I(Z_i\le 
z)\Delta_iK_h(\bb\trans\X_i-\bb\trans\X) 
\frac{\partial\left[
	f_{\bb\trans\X}(\bb\trans\X)E\{I(Z_j\ge 
	Z_i)\mid\bb\trans\X\}\right]/ 
	\partial 
	\bb\trans\X}{f_{\bb\trans\X}^2(\bb\trans\X)E^2\{I(Z_j\ge 
	Z_i)\mid\bb\trans\X,Z_i\}}O_p(A)\right]\\
&\le&4E\left(\left[\frac{1}{n}\sumi I(Z_i\le 
z)\Delta_iK_h(\bb\trans\X_i-\bb\trans\X) \frac{\partial\left[
	f_{\bb\trans\X}(\bb\trans\X)E\{I(Z_j\ge 
	Z_i)\mid\bb\trans\X\}\right]/ 
	\partial 
	\bb\trans\X}{f_{\bb\trans\X}^2(\bb\trans\X)E^2\{I(Z_j\ge 
	Z_i)\mid\bb\trans\X,Z_i\}}\right]^2\right)\\
&&\times O\{(nh)^{-1}+h^4\}\\
&=&O\{(nh)^{-1}+h^4\}
\ese
under conditions \ref{assum:kernel}-\ref{assum:survivalfunction}  
where the second last equation is because of 
(\ref{eq:expectationOf2}) and 
(\ref{eq:varianceOf2}).
Therefore
$
\var\{\wh\bLam_{22}(z,\bb\trans\X)\}=O\{1/(nh^3)\}
$
under conditions 
\ref{assum:kernel}-\ref{assum:survivalfunction}.

Summarizing the results above,
$\var\{\wh\bLam_2(z,\bb\trans\X)\}=O\{1/(nh^3)\}$. 
Hence
the estimator $\wh\bLam_2(z,\bb\trans\X)$ satisfies
\be
\wh\bLam_2(z,\bb\trans\X)=\bLam_2(z,\bb\trans\X)+O_p\{(nh^3)^{-1/2}+h^2\}\n
\ee
under conditions 
\ref{assum:kernel}-\ref{assum:survivalfunction}.

To show that the trimmed estimators have the same asymptotic results, we prove (\ref{eq:expectYtrim}) and skip the others. For further reading about the trimmed
kernel estimators, please see Appendix A.2 of
\cite{hardle1989investigating}. For notational
simplicity, let $\wh f(\bb\trans\X)\equiv n^{-1}\sumj
K_h(\bb\trans\X_j-\bb\trans\X)$. The absolute bias of the trimmed
estimator is
\bse
&&\left|E\left[\frac{\sumj K_h(\bb\trans\X_j-\bb\trans\X)I(Z_j\ge Z)}{\sumj 
	K_h(\bb\trans\X_j-\bb\trans\X)}I\left\{\wh f(\bb\trans\X)> 
d_n\right\}\right]-E\{I(Z_j\ge Z)\mid\bb\trans\X\}\right|\\
&\le&\left|E\left[\frac{\sumj K_h(\bb\trans\X_j-\bb\trans\X)I(Z_j\ge Z)}{\sumj K_h(\bb\trans\X_j-\bb\trans\X)}I\left\{\wh f(\bb\trans\X)> d_n\right\}\right]\right.\\
&&\quad\left.-E\left[\frac{\sumj K_h(\bb\trans\X_j-\bb\trans\X)I(Z_j\ge Z)}{\sumj K_h(\bb\trans\X_j-\bb\trans\X)}I\left\{ f_{\bb\trans\X}(\bb\trans\X)> d_n\right\}\right]\right|\\
&&+\left|E\left[\frac{\sumj K_h(\bb\trans\X_j-\bb\trans\X)I(Z_j\ge Z)}{\sumj K_h(\bb\trans\X_j-\bb\trans\X)}I\left\{f_{\bb\trans\X}(\bb\trans\X)> d_n\right\}\right]\right.\\
&&\quad\left.-E\left[\frac{n^{-1}\sumj K_h(\bb\trans\X_j-\bb\trans\X)I(Z_j\ge Z)}{f_{\bb\trans\X}(\bb\trans\X)}I\left\{f_{\bb\trans\X}(\bb\trans\X)> d_n\right\}\right]\right|\\
&&+\left|E\left[\frac{n^{-1}\sumj K_h(\bb\trans\X_j-\bb\trans\X)I(Z_j\ge 
	Z)}{f_{\bb\trans\X}(\bb\trans\X)}I\left\{f_{\bb\trans\X}(\bb\trans\X)> 
d_n\right\}\right]-E\{I(Z_j\ge Z)\mid\bb\trans\X\}\right|.
\ese
When Condition \ref{assum:fbeta} is replaced by Condition \ref{assum:fbetarelax}, the first term satisfies
\bse
&&\left|E\left[\frac{\sumj K_h(\bb\trans\X_j-\bb\trans\X)I(Z_j\ge Z)}{\sumj K_h(\bb\trans\X_j-\bb\trans\X)}I\left\{\wh f(\bb\trans\X)> d_n\right\}\right]\right.\\
&&-\left.E\left[\frac{\sumj K_h(\bb\trans\X_j-\bb\trans\X)I(Z_j\ge Z)}{\sumj K_h(\bb\trans\X_j-\bb\trans\X)}I\left\{ f_{\bb\trans\X}(\bb\trans\X)> d_n\right\}\right]\right|\\
&\le& \left|E\left[\frac{\sumj K_h(\bb\trans\X_j-\bb\trans\X)I(Z_j\ge
	Z)}{\sumj K_h(\bb\trans\X_j-\bb\trans\X)}I\left\{\wh
f(\bb\trans\X)> d_n,f_{\bb\trans\X}(\bb\trans\X)\le
d_n\right\}\right]\right|\\
&&{+
	\left|E\left[\frac{\sumj K_h(\bb\trans\X_j-\bb\trans\X)I(Z_j\ge
		Z)}{\sumj K_h(\bb\trans\X_j-\bb\trans\X)}I\left\{\wh
	f(\bb\trans\X)\le d_n,f_{\bb\trans\X}(\bb\trans\X)>
	d_n\right\}\right]\right|}\\
&\le& \left|E\left[I\left\{\wh
f(\bb\trans\X)> d_n,f_{\bb\trans\X}(\bb\trans\X)\le
d_n\right\}\right]\right|+
\left|E\left[I\left\{\wh
f(\bb\trans\X)\le d_n,f_{\bb\trans\X}(\bb\trans\X)>
d_n\right\}\right]\right|\\
&\le&\left|E\left[I\left\{f_{\bb\trans\X}(\bb\trans\X)\le
d_n\right\}\right]\right|
+\left|E\left[I\left\{\wh f_{\bb\trans\X}(\bb\trans\X)\le
d_n\right\}\right]\right|\\
&=&O_p\{n^{-\epsilon}+h^2+(nh)^{-1/2}\}\\
&=&O_p\{h^2+(nh)^{-1/2}\}.
\ese
The second term is
\bse
&&\left|E\left[\frac{\sumj K_h(\bb\trans\X_j-\bb\trans\X)I(Z_j\ge Z)}{\sumj K_h(\bb\trans\X_j-\bb\trans\X)}I\left\{f_{\bb\trans\X}(\bb\trans\X)> d_n\right\}\right]\right.\\
&&\quad\left.-E\left[\frac{1/n\sumj K_h(\bb\trans\X_j-\bb\trans\X)I(Z_j\ge Z)}{f_{\bb\trans\X}(\bb\trans\X)}I\left\{f_{\bb\trans\X}(\bb\trans\X)> d_n\right\}\right]\right|\\
&\le& \left|E\left[\frac{1/n\sumj K_h(\bb\trans\X_j-\bb\trans\X)I(Z_j\ge Z)}{f_{\bb\trans\X}(\bb\trans\X)+O_p\{h^2+(nh)^{-1/2}\}}I\left\{f_{\bb\trans\X}(\bb\trans\X)> d_n\right\}\right]\right.\\
&&\quad\left.-E\left[\frac{1/n\sumj K_h(\bb\trans\X_j-\bb\trans\X)I(Z_j\ge Z)}{f_{\bb\trans\X}(\bb\trans\X)}I\left\{f_{\bb\trans\X}(\bb\trans\X)> d_n\right\}\right]\right|\\
&\le& E\left[\frac{1/n\sumj K_h(\bb\trans\X_j-\bb\trans\X)I(Z_j\ge Z)}{f_{\bb\trans\X}(\bb\trans\X)}I\left\{f_{\bb\trans\X}(\bb\trans\X)> d_n\right\}\right]O_p\{h^2+(nh)^{-1/2}\}\\
&=&O_p\{h^2+(nh)^{-1/2}\}.
\ese
The third term is
\bse
&&\left|E\left[\frac{1/n\sumj K_h(\bb\trans\X_j-\bb\trans\X)I(Z_j\ge 
	Z)}{f_{\bb\trans\X}(\bb\trans\X)}I\left\{f_{\bb\trans\X}(\bb\trans\X)> 
d_n\right\}\right]-E\{I(Z_j\ge Z)\mid\bb\trans\X\}\right|\\
&=&\left|E\left[\frac{ K_h(\bb\trans\X_j-\bb\trans\X)I(Z_j\ge 
	Z)}{f_{\bb\trans\X}(\bb\trans\X)}I\left\{f_{\bb\trans\X}(\bb\trans\X)> 
d_n\right\}\right]-E\{I(Z_j\ge Z)\mid\bb\trans\X\}\right|\\
&=&E\left|E\left[I(Z_j\ge Z)I\left\{f_{\bb\trans\X}(\bb\trans\X)> 
d_n\right\}\mid\bb\trans\X\right]-E\{I(Z_j\ge 
Z)\mid\bb\trans\X\}+O_p(h^2)\right|\\
&=&E\left[I(Z_j\ge Z)I\left\{f_{\bb\trans\X}(\bb\trans\X)\le d_n\right\}\right]+O_p(h^2)\\
&\le&E\left[I\left\{f_{\bb\trans\X}(\bb\trans\X)\le d_n\right\}\right]+O_p(h^2)\\
&=&O_p\{n^{-\epsilon}+h^2+(nh)^{-1/2}\}\\
&=&O_p\{h^2+(nh)^{-1/2}\}.
\ese
It is easy to see the variance of this trimmed estimator,
\bse
\var\left[\frac{\sumj K_h(\bb\trans\X_j-\bb\trans\X)I(Z_j\ge Z)}{\sumj K_h(\bb\trans\X_j-\bb\trans\X)}I\left\{\wh f(\bb\trans\X)> d_n\right\}\right]=O_p\{(nh)^{-1/2}\}.
\ese
Therefore
\bse
\frac{\sumj K_h(\bb\trans\X_j-\bb\trans\X)I(Z_j\ge Z)}{\sumj K_h(\bb\trans\X_j-\bb\trans\X)}I\left\{\wh f(\bb\trans\X)> d_n\right\}=E\{I(Z_j\ge Z)\mid\bb\trans\X\}+O_p\{h^2+(nh)^{-1/2}\}.
\ese

	We give detailed proof for uniform result of (\ref{eq:lemeq2}) only. Because
	the domain of $\bb\trans\X$ is compact, we divide it into rectangular
	regions. In each region, the distance between a point $\bb\trans\x$ in
	this region and the nearest grid point is less than $n^{-2}$. We need only $N\le Cn^2$ grid points, where $C$ is a constant. Let the grid points be
	$\bkappa_1, \dots, \bkappa_N$. Let
	$\wh\brho(\bb\trans\X)=\wh E\{\X
	Y(Z)\mid\bb\trans\X\}$ and
	$\brho(\bb\trans\X)=E\{\X
	Y(Z)\mid\bb\trans\X\}$. Then for any
	$(\bb\trans\X)$, there exists a $\bkappa_i$, $1 \le i \le N$, such that \bse
	|\wh\brho(\bb\trans\X) -\brho(\bb\trans\X)| &\le&|\wh\brho(\bkappa_i)
	-\brho(\bkappa_i)| +|\wh\brho(\bb\trans\X) -\wh\brho(\bkappa_i)|+
	|\brho(\bb\trans\X) -\brho(\bkappa_i)|\\ &\le&|\wh\brho(\bkappa_i)
	-\brho(\bkappa_i)| +D_1n^{-2}, 
\ese 
for an absolute constant $D_1$ under
	Conditions \ref{assum:kernel} and \ref{assum:survivalfunction}. Thus, for
	any $D\ge D_1$,
 \bse &&\pr(\sup_{\bb\trans\x}|\wh\brho(\bb\trans\X)-\brho(\bb\trans\X)|>2D[h^2+\{\log n (nh)^{-1}\}^{1/2}])\\
	&\le&\pr(\sup_{\bkappa_i}|\wh\brho(\bkappa_i )-\brho(\bkappa_i)|>2D[h^2+\{\log n
	(nh)^{-1}\}^{1/2}] -D_1n^{-2})\\
	&\le&\pr(\sup_{\bkappa_i}|\wh\brho(\bkappa_i )-\brho(\bkappa_i)|>D[h^2+\{\log n
	(nh)^{-1}\}^{1/2}]) 
\ese 
under Condition
        \ref{assum:bandwidth}. Using
	Bernstein's inequality on $\wh\brho(\bkappa_i)$,  under Conditions
        \ref{assum:kernel}-\ref{assum:survivalfunction}, we have
\bse 
\pr
	[|\wh\brho(\bkappa_i )-\brho(\bkappa_i)|\ge A\{\log
        n/(nh)\}^{1/2}] 
&\le&
	2\exp\left\{\frac{-n A^2\log n/(nh)}{2D_2h^{-1}+2/3AD_3(\log
		n)^{1/2}(nh)^{-1/2}}\right\}\\ 
&=&2\exp\left\{\frac{- A^2\log
		n}{2D_2+2/3AD_3(\log
		nh/n)^{1/2}}\right\}\\
&\le&2\exp\left(\frac{- A^2\log
		n}{2D_2+AD_3}\right)
\ese 
for all $A>D_3+\sqrt{D_3^2+4D_2}$, where $D_2$ and $D_3$ are constants satisfying
\bse
\var\{\wh\brho(\bb\trans\X) -\brho(\bb\trans\X)\}&\le& \frac{D_2}{nh},\\
\left|\frac{K_h(\bb\trans\X_i-\bb\trans\X)\X_iI(Z_i\ge Z)}{1/n\sumj K_h(\bb\trans\X_j-\bb\trans\X)}-\brho(\bb\trans\X)\right|&\le& D_3 \text{ with probability 1}.
\ese
This leads to 
\bse \pr
	[\sup_{\bkappa_i}|\wh\brho(\bkappa_i )-E\wh\brho(\bkappa_i)|\ge A\{\log
	n/(nh)\}^{1/2}] &\le& 2Cn^{2}\exp\left(\frac{- A^2\log n}{2D_2+AD_3}\right)\\
	&=&2C\exp\left[\{2- A^2/(2D_2+AD_3)\}\log n\right]\to0 
\ese
 because $A>D_3+\sqrt{D_3^2+4D_2}$.
	Combining the above results, for $A_1=\max(A,D)$, \bse
	&&\pr(\sup_{\bb\trans\x}|\wh\brho(\bb\trans\X,
	\bb)-\brho(\bb\trans\X)|>2A_1[h^2+\{\log n (nh)^{-1}\}^{1/2}])\\ &\le&
	\pr(\sup_{\bkappa_i}|\wh\brho(\bkappa_i )-\brho(\bkappa_i)|>A_1[h^2+\{\log n
	(nh)^{-1}\}^{1/2}])\\ &\le&\pr\{\sup_{\bkappa_i}|\wh\brho(\bkappa_i
	)-\brho(\bkappa_i)|>A_1h^2\} +\pr([\sup_{\bkappa_i}|\wh\brho(\bkappa_i
	)-\brho(\bkappa_i)|\ge A_1\{\log n (nh)^{-1}\}^{1/2}])\\ &\to&0. \ese 
	The uniform convergence results concerning (\ref{eq:lemeq5})-(\ref{eq:lemeq8})
	are slightly different because these functions contain the additional component
	$Z$. Nevertheless, under Condition \ref{assum:survivalfunction}, the support of
	$(\bb\trans\X_i,Z_i)$ or $(\bb\trans\X_i, Z_j)$  is also bounded so we
	can similarly divide the region using $N\le
        Cn^{2+2}$ grid points while the 
	distance of a point to the nearest grid point is less than $n^{-2}$. The rest of
	the analysis can then be similarly carried out  as above, then the uniform convergence is established. \qed

\subsubsection{Lemma \ref{th:Lambda}}
\begin{Lem}\label{th:Lambda}
	The estimator $\wh{\Lambda}(t,\wh\bb\trans\X)$ has the 
	expansion
	\bse
	&&\sqrt{nh}\left\{\wh{\Lambda}(t,\wh\bb\trans\X)- 
	{\Lambda}(t,\bb\trans\X)\right\}\\
	&=&
	\sqrt{\frac{h}{n}}\sumi\int_{0}^{t}\frac{I\left\{
		\sumj Y_j(s)K_h(\bb\trans\X_j-\bb\trans\X)>0\right\}}
	{f_{\bb\trans\X}(\bb\trans\X)E\{I(z\ge s)\mid\bb\trans\X\}} 
	K_h(\bb\trans\X_i-\bb\trans\X)dM_i(s,\bb\trans\X)+o_p(1),
	\ese
	and
	satisfies 
	\bse
	\sqrt{nh}\left\{\wh{\Lambda}(t,\wh\bb\trans\X)-{\Lambda}(t,\bb\trans\X)\right\}\to
	N\{0,\sigma^2(t,\bb\trans\X)\}
	\ese
	in distribution when $n\to\infty$ for all $t,\bb\trans\X$ 
	under Conditions 
	\ref{assum:kernel}-\ref{assum:survivalfunction}, where
	\bse
	\sigma^2(t,\bb\trans\X)=\int 
	K^2(u)du\int_{0}^{t}\frac{\lambda(s,\bb\trans\X)} 
	{f_{\bb\trans\X}(\bb\trans\X)E\{I(Z\ge 
		s)\mid\bb\trans\X\}}ds.
	\ese
\end{Lem}
\noindent Proof:
For notational convenience, let $d=1$ and 
$\nu=2$. 
For any 
$t$ and 
$\bb\trans\X$, define 
\bse 
\phi_n(s,\bb\trans\X)&\equiv&\sumj 
Y_j(s)K_h(\bb\trans\X_j-\bb\trans\X),\\
Q_n(t,\bb\trans\X) 
&\equiv&\wh\bLam(t,\bb\trans\X)-\int_{0}^{t}I\left\{ 
\phi_n(s,\bb\trans\X)>0\right\}\lambda(s,\bb\trans\X)ds,\\
D_n(t,\bb\trans\X)&\equiv&\int_{0}^{t}\lambda(s,\bb\trans\X) 
\left[1-I\left\{\phi_n(s,\bb\trans\X)>0\right\}\right]ds.
\ese 
Then 
\bse 
\sqrt{nh}\left\{\wh{\Lambda}(t,\bb\trans\X)- 
{\Lambda}(t,\bb\trans\X)\right\}=\sqrt{nh}Q_n(t,\bb\trans\X) 
-\sqrt{nh}D_n(t,\bb\trans\X). 
\ese 
We first show that $\sqrt{nh}D_n(t,\bb\trans\X)\to 0$ in 
probability uniformly. It suffices to show that 
\bse 
\sqrt{nh}\left[1-I\left\{\phi_n(s,\bb\trans\X)>0\right\}\right]\overset{p}{\to}
0, 
\ese 
which is equivalent to show that for any $\epsilon>0$, 
\bse 
\Pr\left(\sqrt{nh}\left[1-I\left\{\phi_n(s,\bb\trans\X)>0\right\}\right]>\epsilon\right){\to}
0. 
\ese 
Now for $\epsilon\ge\sqrt{nh}$, the above automatically holds. 
For $\epsilon<\sqrt{nh}$, this is equivalent to show 
\be 
\Pr\left\{n^{-1}\phi_n(s,\bb\trans\X)\le 0\right\}{\to} 
0.\label{eq:sumto0}
\ee 
Because
$n^{-1}\phi_n(s,\bb\trans\X)=f_{\bb\trans\X}(\bb\trans\X) 
E\{I(Z_j\ge s)\mid\bb\trans\X\}+O_p\{h^2+(nh)^{-1/2}\}$, 
(\ref{eq:sumto0}) is equivalent to 
\bse 
\Pr\left[f_{\bb\trans\X}(\bb\trans\X)E\{I(Z_j\ge 
s)\mid\bb\trans\X\}+O_p\{h^2+(nh)^{-1/2}\}
\le 0\right]{\to} 0, 
\ese 
which automatically holds under Condition \ref{assum:fbeta} and 
\ref{assum:survivalfunction}. 
Hence $\sqrt{nh}D_n(t,\bb\trans\X)\to 0$ in 
probability uniformly. 
Second we inspect the asymptotic property of 
$\sqrt{nh}Q_n(t,\bb\trans\X) 
$. 
Recall $M_i(s,\bb\trans\X)$ is the martingale corresponding 
to the 
counting process $N_i(s)$ and satisfies 
$dM_i(s,\bb\trans\X)=dN_i(s)-Y_i(s)\lambda(s,\bb\trans\X)ds$. 
\be 
&&\sqrt{nh}Q_n(t,\bb\trans\X)\n\\
&=&\int_{0}^{t}\sqrt{nh}\frac{1}{\sumj 	
	Y_j(s)K_h(\bb\trans\X_j-\bb\trans\X)}\sumi 
K_h(\bb\trans\X_i-\bb\trans\X)dN_i(s)\n\\
&&-\int_{0}^{t}\sqrt{nh}I\left\{\phi_n(s,\bb\trans\X)>0\right\}\lambda(s,\bb\trans\X)ds\n\\
&=&\int_{0}^{t}\sqrt{nh}\frac{I\{\phi_n(s,\bb\trans\X)\le 
	0\}}{\sumj 
	Y_j(s)K_h(\bb\trans\X_j-\bb\trans\X)}\sumi 
K_h(\bb\trans\X_i-\bb\trans\X)dN_i(s)\n\\
&&+\int_{0}^{t}\sqrt{nh}\frac{I\{\phi_n(s,\bb\trans\X)> 
	0\}}{\sumj Y_j(s)K_h(\bb\trans\X_j-\bb\trans\X)}\sumi 
K_h(\bb\trans\X_i-\bb\trans\X)dN_i(s)\n\\
&&-\int_{0}^{t}\sqrt{nh}\frac{I\{\phi_n(s,\bb\trans\X)> 
	0\}}{\sumj Y_j(s)K_h(\bb\trans\X_j-\bb\trans\X)}\sumi 
K_h(\bb\trans\X_i-\bb\trans\X)Y_i(s)\lambda(s,\bb\trans\X)ds\n\\
&=&\int_{0}^{t}\sqrt{nh}\frac{I\{\phi_n(s,\bb\trans\X)>0\}} 
{\sumj Y_j(s)K_h(\bb\trans\X_j-\bb\trans\X)}\sumi 
K_h(\bb\trans\X_i-\bb\trans\X)dM_i(s,\bb\trans\X)\label{eq:Qn1}\\
&&+\int_{0}^{t}\sqrt{nh}\frac{I\{\phi_n(s,\bb\trans\X)\le 
	0\}}{\sumj Y_j(s)K_h(\bb\trans\X_j-\bb\trans\X)}\sumi 
K_h(\bb\trans\X_i-\bb\trans\X)dN_i(s).\label{eq:Qn2}
\ee 
We decompose (\ref{eq:Qn1}) as 
\be 
&&\int_{0}^{t}\sqrt{nh}\frac{I\{\phi_n(s,\bb\trans\X)>0\}} 
{\sumj Y_j(s)K_h(\bb\trans\X_j-\bb\trans\X)}\sumi 
K_h(\bb\trans\X_i-\bb\trans\X)dM_i(s,\bb\trans\X)\n\\
&=&\sqrt{\frac{h}{n}}\sumi\int_{0}^{t} 
\frac{I\{\phi_n(s,\bb\trans\X)>0\}}{1/n\sumj 
	Y_j(s)K_h(\bb\trans\X_j-\bb\trans\X)} 
K_h(\bb\trans\X_i-\bb\trans\X)dM_i(s,\bb\trans\X)\n\\
&=&\sqrt{\frac{h}{n}}\sumi\int_{0}^{t}I\left\{\phi_n(s,\bb\trans\X)>0\right\}\left[\frac{1}
{f_{\bb\trans\X}(\bb\trans\X)E\{I(Z_j\ge 
	s)\mid\bb\trans\X\}}\right.\n\\
&&\left.-\frac{1/n\sumj Y_j(s)K_h(\bb\trans\X_j-\bb\trans\X)} 
{f^2_{\bb\trans\X}(\bb\trans\X)E^2\{I(Z_j\ge 
	s)\mid\bb\trans\X\}}  
+\frac{1}{f_{\bb\trans\X}(\bb\trans\X)E\{I(Z_j\ge 	
	s)\mid\bb\trans\X\}} +O_p\{h^4+(nh)^{-1}\}\right]\n\\
&&\times K_h(\bb\trans\X_i-\bb\trans\X)dM_i(s,\bb\trans\X)\n\\
&=&Q_{n1}-Q_{n2}+o_p(1),\label{eq:Qn1decompose}
\ee 
where 
\bse 
Q_{n1}&=&\sqrt{\frac{h}{n}}\sumi\int_{0}^{t}\frac{I\left\{\phi_n(s,\bb\trans\X)>0\right\}}
{f_{\bb\trans\X}(\bb\trans\X)E\{I(Z\ge 
	s)\mid\bb\trans\X\}}K_h(\bb\trans\X_i-\bb\trans\X)dM_i(s,\bb\trans\X),\\
Q_{n2}&=&\sqrt{\frac{h}{n}}\sumi\int_{0}^{t}
\frac{I\left\{\phi_n(s,\bb\trans\X)>0\right\}} 
{f^2_{\bb\trans\X}(\bb\trans\X)E^2\{I(Z_j\ge 
	s)\mid\bb\trans\X\}}\n\\
&&\times\left[\frac{1}{n}\sumj 
Y_j(s)K_h(\bb\trans\X_j-\bb\trans\X) 
-f_{\bb\trans\X}(\bb\trans\X)
E\{I(Z\ge s)\mid\bb\trans\X\}\right]\n\\
&&	\times K_h(\bb\trans\X_i-\bb\trans\X)dM_i(s,\bb\trans\X)\n\\
&=&\sqrt{\frac{h}{n^3}}\sumi\sumj\int_{0}^{t} 
\frac{I\left\{\phi_n(s,\bb\trans\X)>0\right\}} 
{f^2_{\bb\trans\X}(\bb\trans\X)E^2\{I(Z_j\ge 
	s)\mid\bb\trans\X\}}\n\\
&&\times\left[Y_j(s)K_h(\bb\trans\X_j-\bb\trans\X)- 
f_{\bb\trans\X}(\bb\trans\X)E\{I(Z\ge s)\mid\bb\trans\X\}\right] 
\\
&&\times K_h(\bb\trans\X_i-\bb\trans\X)dM_i(s,\bb\trans\X)\n, 
\ese 
and the remaining term in (\ref{eq:Qn1decompose}) is $o_p(1)$ 
because 
$\sqrt{n/h} O_p\{h^4+(nh)^{-1}\} 
=O_p\{n^{1/2}h^{7/2}+(nh^3)^{-1/2}\}=o_p(1)$ by Condition 
\ref{assum:bandwidth}.

Using the U-statistic property, $Q_{n2}$ has leading order terms 
$Q_{n21}+Q_{n22}-Q_{n23}$, where 
\bse 
Q_{n21}&=&\sqrt{\frac{h}{n}}E\Bigg(\sumi\int_{0}^{t} 
\frac{I\left\{\phi_n(s,\bb\trans\X)>0\right\}} 
{f^2_{\bb\trans\X}(\bb\trans\X)E^2\{I(Z_j\ge 
	s)\mid\bb\trans\X\}}\n\\
&&\times\left[Y_j(s)K_h(\bb\trans\X_j-\bb\trans\X)- 
f_{\bb\trans\X}(\bb\trans\X)E\{I(Z\ge s)\mid\bb\trans\X\}\right] 
\\
&&\times 
K_h(\bb\trans\X_i-\bb\trans\X)dM_i(s,\bb\trans\X) 
\mid\Delta_i,\bb\trans\X_i,Z_i\Bigg),\\
Q_{n22}&=&\sqrt{\frac{h}{n}}E\Bigg(\sumj 
\int_{0}^{t}\frac{I\left\{\phi_n(s,\bb\trans\X)>0\right\}} 
{f^2_{\bb\trans\X}(\bb\trans\X)E^2\{I(Z_j\ge 
	s)\mid\bb\trans\X\}}\n\\
&&\times\left[Y_j(s)K_h(\bb\trans\X_j-\bb\trans\X)- 
f_{\bb\trans\X}(\bb\trans\X)E\{I(Z\ge s)\mid\bb\trans\X\}\right] 
\\
&&\times 
K_h(\bb\trans\X_i-\bb\trans\X)dM_i(s,\bb\trans\X) 
\mid\Delta_j,\bb\trans\X_j,Z_j\Bigg),\\
Q_{n23}&=&\sqrt{nh}E\Bigg(\int_{0}^{t} 
\frac{I\left\{\phi_n(s,\bb\trans\X)>0\right\}} 
{f^2_{\bb\trans\X}(\bb\trans\X)E^2\{I(Z_j\ge 
	s)\mid\bb\trans\X\}}\n\\
&&\times\left[Y_j(s)K_h(\bb\trans\X_j-\bb\trans\X)- 
f_{\bb\trans\X}(\bb\trans\X)E\{I(Z\ge 
s)\mid\bb\trans\X\}\right]\\
&&\times 
E\left\{K_h(\bb\trans\X_i-\bb\trans\X)dM_i(s,\bb\trans\X)\right\}\Bigg).
\ese 
$I\left\{\phi_n(s,\bb\trans\X)>0\right\}
=I\left[f_{\bb\trans\X}(\bb\trans\X)E\{I(Z_j\ge 
s)\mid\bb\trans\X\}+O_p\{h^2+(nh)^{-1}\}>0\right]=1 
$
almost surely. Thus, almost surely, 
\bse 
Q_{n21}
&=&\sqrt{\frac{h}{n}}\sumi\Bigg(\int_{0}^{t} 
\frac{1}{f_{\bb\trans\X}(\bb\trans\X)E\{I(Z_j\ge 
	s)\mid\bb\trans\X\}}\n\\
&&\times E\left[Y_j(s)K_h(\bb\trans\X_j-\bb\trans\X)- 
f_{\bb\trans\X}(\bb\trans\X)E\{I(Z\ge s)\mid\bb\trans\X\}\right] 
\\
&&\times 
K_h(\bb\trans\X_i-\bb\trans\X)dM_i(s,\bb\trans\X)\Bigg)\\
&=&\sqrt{\frac{h}{n}}\sumi 
\int_{0}^{t}\frac{O(h^2)}{f_{\bb\trans\X}(\bb\trans\X)E\{I(Z_j\ge
	s)\mid\bb\trans\X\}}
K_h(\bb\trans\X_i-\bb\trans\X)dM_i(s,\bb\trans\X)\\
&\to& 0. 
\ese 
uniformly as $h\to 0$. 
Similarly, almost surely, 
\bse 
Q_{n22}
%
&=&\sqrt{\frac{h}{n}}\sumj 
\int_{0}^{t}\frac{1}{f^2_{\bb\trans\X}(\bb\trans\X)E^2\{I(Z_j\ge 
	s)\mid\bb\trans\X\}}\n\\
&&\times\left[Y_j(s)K_h(\bb\trans\X_j-\bb\trans\X)- 
f_{\bb\trans\X}(\bb\trans\X)E\{I(Z\ge s)\mid\bb\trans\X\}\right] 
\\
&&\times 
E\{K_h(\bb\trans\X_i-\bb\trans\X)dM_i(s,\bb\trans\X)\}\\
&=&0.
\ese 
Obviously, $Q_{n23}=E(Q_{n22})=0$, 
hence $Q_{n2}\to0$ in probability as $n\to\infty$. 

For (\ref{eq:Qn2}) 
\bse 
&&\int_{0}^{t}\sqrt{nh}\frac{I\{\phi_n(s,\bb\trans\X)\le 
	0\}}{\sumj Y_j(s)K_h(\bb\trans\X_j-\bb\trans\X)}\sumi 
K_h(\bb\trans\X_i-\bb\trans\X)dN_i(s) 
\to 0 
\ese 
in probability uniformly. 
We have obtained 
\be 
&&\sqrt{nh}Q_n(t,\bb\trans\X)\n\\
&=&\sqrt{\frac{h}{n}}\sumi\int_{0}^{t}\frac{I\left\{\phi_n(s,\bb\trans\X)>0\right\}}
{f_{\bb\trans\X}(\bb\trans\X)E\{I(Z\ge 
	s)\mid\bb\trans\X\}}K_h(\bb\trans\X_i-\bb\trans\X)dM_i(s,\bb\trans\X)
\label{eq:Qn1a}\\
&&+o_p(1).\n 
\ee 
Applying martingale 
central limit theorem on (\ref{eq:Qn1a}),
we have 
\be 
&&\frac{h}{n}\sumi\int_{0}^{t}\frac{\lambda(s,\bb\trans\X) 
	I\{\phi_n(s,\bb\trans\X)>0\}}{[f_{\bb\trans\X}(\bb\trans\X)E\{I(Z\ge
	s)\mid\bb\trans\X\}]^2} 
K^2_h(\bb\trans\X_i-\bb\trans\X)Y_i(s)ds\n\\
&=&\int_{0}^{t}\frac{\lambda(s,\bb\trans\X)I\{\phi_n(s,\bb\trans\X)>0\}}
{[f_{\bb\trans\X}(\bb\trans\X)E\{I(Z\ge 
	s)\mid\bb\trans\X\}]^2}\frac{1}{n}\sumi 
hK^2_h(\bb\trans\X_i-\bb\trans\X)Y_i(s)ds\n\\
&=&\int_{0}^{t}\frac{\lambda(s,\bb\trans\X) 
	I\{\phi_n(s,\bb\trans\X)>0\}}{[f_{\bb\trans\X}(\bb\trans\X) 
	E\{I(Z\ge s)\mid\bb\trans\X\}]^2}\n\\
&&\times \left[f_{\bb\trans\X}(\bb\trans\X)E\{I(Z_i\ge 
s)\mid\bb\trans\X\}\int 
K^2(u)du+O_p(n^{-1/2}h^{-1/2}+h^2)\right]ds\n\\
&\overset{p}{\to}&\int 
K^2(u)du\int_{0}^{t}\frac{\lambda(s,\bb\trans\X)} 
{f_{\bb\trans\X}(\bb\trans\X)E\{I(Z_j\ge 
	s)\mid\bb\trans\X\}}ds\n\\
&=&\sigma^2(t,\bb\trans\X)\label{eq:MCLTa}.
\ee 
Next we inspect the following integration for any $\epsilon>0$. 
\bse 
&&\sumi\int_{0}^{t}\frac{h}{n}\frac{I\{\phi_n(s,\bb\trans\X)>0\}}
{[f_{\bb\trans\X}(\bb\trans\X)E\{I(Z\ge s)\mid\bb\trans\X\}]^2}
K^2_h(\bb\trans\X_i-\bb\trans\X)Y_i(s)\\
&&\times 
I\left[\sqrt{\frac{h}{n}}\left|\frac{I\{\phi_n(s,\bb\trans\X)>0\}}
{f_{\bb\trans\X}(\bb\trans\X)E\{I(Z\ge	s)\mid\bb\trans\X\}} 
K_h(\bb\trans\X_i-\bb\trans\X)\right|>\epsilon\right] 
\lambda(s,\bb\trans\X)ds\\
&=&\int_{0}^{t}\frac{I\{\phi_n(s,\bb\trans\X)>0\}} 
{[f_{\bb\trans\X}(\bb\trans\X)E\{I(Z\ge s)\mid\bb\trans\X\}]^2}
\frac{1}{n}\sumi \Bigg(hK^2_h(\bb\trans\X_i-\bb\trans\X)Y_i(s)\\
&&\times 
I\left[\sqrt{\frac{h}{n}}\left|\frac{I\{\phi_n(s,\bb\trans\X)>0\}}
{f_{\bb\trans\X}(\bb\trans\X)E\{I(Z\ge s)\mid\bb\trans\X\}} 
K_h(\bb\trans\X_i-\bb\trans\X)\right|>\epsilon\right]\Bigg)\lambda(s,\bb\trans\X)ds\\
&\le&\int_{0}^{t}\lambda(s,\bb\trans\X) 
\frac{I\{\phi_n(s,\bb\trans\X)>0\}}{[f_{\bb\trans\X}(\bb\trans\X)
	E\{I(Z\ge s)\mid\bb\trans\X\}]^2}
\left\{\frac{1}{n}\sumi 
hK^2_h(\bb\trans\X_i-\bb\trans\X)Y_i(s)\right\}\\
&&\times \sup_{1\le i\le n} 
I\left[\left|f_{\bb\trans\X}(\bb\trans\X)E\{I(Z\ge 
s)\mid\bb\trans\X\}\right|<\frac{m_i}{\epsilon\sqrt{nh}}\right]ds,
\ese 
where 
$m_i=|I\{\phi_n(s,\bb\trans\X)>0\}K\{(\bb\trans\X_i-\bb\trans\X)/h\}|$,
which is bounded following Condition \ref{assum:kernel}. 
In the above display, 
\bse 
\sup_{1\le i\le 
	n}I\left[\left|f_{\bb\trans\X}(\bb\trans\X)E\{I(Z_j\ge 
s)\mid\bb\trans\X\}\right|<\frac{m_i}{\epsilon\sqrt{nh}}\right]=0
\ese 
as long as $n$ is large enough because the right hand side 
converges to 0 by 
Condition \ref{assum:bandwidth} but the left hand side will be 
always 
larger than 0 by conditions \ref{assum:fbeta} and 
\ref{assum:survivalfunction}. 
On the other hand, 
\bse 
&&\frac{1}{n}\sumi hK^2_h(\bb\trans\X_i-\bb\trans\X)Y_i(s)\\
&=&f_{\bb\trans\X}(\bb\trans\X)E\{I(Z_i\ge 
s)\mid\bb\trans\X\}\int 
K^2(u)du+O_p(n^{-1/2}h^{-1/2}+h^2)\\
&\to&0 
\ese 
in probability uniformly. 
Hence 
\be 
&&\lim_{n\to\infty}\sumi\int_{0}^{t}\frac{h}{n} 
\frac{I\{\phi_n(s,\bb\trans\X)>0\}}{[f_{\bb\trans\X}(\bb\trans\X)
	E\{I(Z\ge s)\mid\bb\trans\X\}]^2}
K^2_h(\bb\trans\X_i-\bb\trans\X)Y_i(s)\n\\
&&\times 
I\left[\sqrt{\frac{h}{n}}\left|\frac{I\{\phi_n(s,\bb\trans\X)>0\}}
{f_{\bb\trans\X}(\bb\trans\X)E\{I(Z\ge s)\mid\bb\trans\X\}} 
K_h(\bb\trans\X_i-\bb\trans\X)\right|>\epsilon\right]\lambda(s,\bb\trans\X)ds=0\n\\
\label{eq:MCLTb}
\ee 
with probability 1 uniformly for any $\epsilon>0$. 

In summary 
\bse 
\sqrt{nh}\left\{\wh{\Lambda}(t,\bb\trans\X)-{\Lambda}(t,\bb\trans\X)\right\}\to
N\{0,\sigma^2(t,\bb\trans\X)\}
\ese 
uniformly. 
\qed

\subsection{Proof of Theorem 
	\ref{th:consistency}}\label{app:consistofb}
Because the result regarding (\ref{eq:eff}) is the most 
difficult to
establish, we provide only the proof concerning (\ref{eq:eff}),
the result concerning (\ref{eq:general}) 
is based on a similar proof.

For each $n$, let $\wh\bb_n$ satisfy
\bse
\frac{1}{n}\sumi 
\Delta_i\frac{\wh\blam_1(Z_i,\wh\bb_n\trans\X_i)}{\wh\lambda(Z_i,\wh\bb_n\trans\X_i)}
\otimes\left[\X_{li}-
\frac{\wh E\left\{\X_{li} 
	Y_i(Z_i)\mid\wh\bb_n\trans\X_i\right\}}
{\wh E\left\{Y_i(Z_i)\mid\wh\bb_n\trans\X_i\right\}}\right]=\0.
\ese
Under condition \ref{assum:bounded}, there exists a subsequence 
of
$\wh\bb_n, n=1, 2,\dots$, 
that converges. For notational simplicity,  we still write 
$\wh\bb_n,
n=1, 2, \dots, $
as the subsequence that  
converges and let the limit be $\bb^*$. 

From the uniform convergence in (\ref{eq:lemeq1}), 
(\ref{eq:lemeq2}),
(\ref{eq:lemeq5}), (\ref{eq:lemeq6}) given in
Lemma \ref{lem:pre},
\bse
&&\frac{1}{n}\sumi 
\Delta_i\frac{\wh\blam_1(Z_i,\wh\bb_n\trans\X_i)}{\wh\lambda(Z_i,\wh\bb_n\trans\X_i)}
\otimes\left[\X_{li}-
\frac{\wh E\left\{\X_{li} 
	Y_i(Z_i)\mid\wh\bb_n\trans\X_i\right\}}
{\wh E\left\{Y_i(Z_i)\mid\wh\bb_n\trans\X_i\right\}}\right]\\
&=&\frac{1}{n}\sumi 
\Delta_i\frac{\blam_1(Z_i,\wh\bb_n\trans\X_i)+O_p\{(nbh^3)^{-1/2}+h^2+b^2\}}{\lambda(Z_i,\wh\bb_n\trans\X_i)+O_p\{(nbh)^{-1/2}+h^2+b^2\}}\\
&&
\otimes\left[\X_{li}-
\frac{E\left\{\X_{li} 
	Y_i(Z_i)\mid\wh\bb_n\trans\X_i\right\}+O_p\{(nh)^{-1/2}+h^2\}}
{E\left\{Y_i(Z_i)\mid\wh\bb_n\trans\X_i\right\}+O_p\{(nh)^{-1/2}+h^2\}}\right]\\
&=&\frac{1}{n}\sumi 
\Delta_i\left[\frac{\blam_1(Z_i,\wh\bb_n\trans\X_i)}{\lambda(Z_i,\wh\bb_n\trans\X_i)}+O_p\{(nbh^3)^{-1/2}+h^2+b^2\}\right]\\
&&
\otimes\left[\X_{li}-
\frac{E\left\{\X_{li} 
	Y_i(Z_i)\mid\wh\bb_n\trans\X_i\right\}}
{E\left\{Y_i(Z_i)\mid\wh\bb_n\trans\X_i\right\}}+O_p\{(nh)^{-1/2}+h^2\}\right]\\
&=&\frac{1}{n}\sumi 
\Delta_i\frac{\blam_1(Z_i,\wh\bb_n\trans\X_i)}{\lambda(Z_i,\wh\bb_n\trans\X_i)}
\otimes\left[\X_{li}-
\frac{E\left\{\X_{li} 
	Y_i(Z_i)\mid\wh\bb_n\trans\X_i\right\}}
{E\left\{Y_i(Z_i)\mid\wh\bb_n\trans\X_i\right\}}\right]+o_p(1).
\ese
Thus,  for sufficiently large $n$,
\bse
&&\frac{1}{n}\sumi 
\Delta_i\frac{\blam_1(Z_i,\wh\bb_n\trans\X_i)}{\lambda(Z_i,\wh\bb_n\trans\X_i)}
\otimes\left[\X_{li}-
\frac{E\left\{\X_{li} 
	Y_i(Z_i)\mid\wh\bb_n\trans\X_i\right\}}
{E\left\{Y_i(Z_i)\mid\wh\bb_n\trans\X_i\right\}}\right]\\
&=&\frac{1}{n}\sumi 
\Delta_i\frac{\blam_1(Z_i,{\bb^*}\trans\X_i)}{\lambda(Z_i,{\bb^*}\trans\X_i)}
\otimes\left[\X_{li}-
\frac{E\left\{\X_{li} 
	Y_i(Z_i)\mid{\bb^*}\trans\X_i\right\}}
{E\left\{Y_i(Z_i)\mid{\bb^*}\trans\X_i\right\}}\right]+O_p(\wh\bb_n-\bb^*)\\
&=&\frac{1}{n}\sumi 
\Delta_i\frac{\blam_1(Z_i,{\bb^*}\trans\X_i)}{\lambda(Z_i,{\bb^*}\trans\X_i)}
\otimes\left[\X_{li}-
\frac{E\left\{\X_{li} 
	Y_i(Z_i)\mid{\bb^*}\trans\X_i\right\}}
{E\left\{Y_i(Z_i)\mid{\bb^*}\trans\X_i\right\}}\right]+o_p(1),
\ese
under Condition 
\ref{assum:kernel}-\ref{assum:bandwidth}, where
the last equality is because $\wh\bb_n$ converges to $\bb^*$.
In addition,
\bse
&&\frac{1}{n}\sumi 
\Delta_i\frac{\blam_1(Z_i,{\bb^*}\trans\X_i)}{\lambda(Z_i,{\bb^*}\trans\X_i)}
\otimes\left[\X_{li}-
\frac{E\left\{\X_{li} 
	Y_i(Z_i)\mid{\bb^*}\trans\X_i\right\}}
{E\left\{Y_i(Z_i)\mid{\bb^*}\trans\X_i\right\}}\right]\\
&=&E\left(\Delta\frac{\blam_1(Z,{\bb^*}\trans\X)}{\lambda(Z,{\bb^*}\trans\X)}
\otimes\left[\X_{l}-
\frac{E\left\{\X_{l} 
	Y(Z)\mid{\bb^*}\trans\X\right\}}
{E\left\{Y(Z)\mid{\bb^*}\trans\X_i\right\}}\right]\right)+o_p(1)
\ese
under Condition \ref{assum:kernel}-\ref{assum:bandwidth}. Thus, 
for sufficient 
large $n$
\bse
\0&=&\frac{1}{n}\sumi 
\Delta_i\frac{\wh\blam_1(Z_i,\wh\bb_n\trans\X_i)}{\wh\lambda(Z_i,\wh\bb_n\trans\X_i)}
\otimes\left[\X_{li}-
\frac{\wh E\left\{\X_{li} 
	Y_i(Z_i)\mid\wh\bb_n\trans\X_i\right\}}
{\wh E\left\{Y_i(Z_i)\mid\wh\bb_n\trans\X_i\right\}}\right]\\
&=&\frac{1}{n}\sumi 
\Delta_i\frac{\blam_1(Z_i,\wh\bb_n\trans\X_i)}{\lambda(Z_i,\wh\bb_n\trans\X_i)}
\otimes\left[\X_{li}-
\frac{E\left\{\X_{li} 
	Y_i(Z_i)\mid\wh\bb_n\trans\X_i\right\}}
{E\left\{Y_i(Z_i)\mid\wh\bb_n\trans\X_i\right\}}\right]+o_p(1)\\
&=&\frac{1}{n}\sumi 
\Delta_i\frac{\blam_1(Z_i,{\bb^*}\trans\X_i)}{\lambda(Z_i,{\bb^*}\trans\X_i)}
\otimes\left[\X_{li}-
\frac{E\left\{\X_{li} 
	Y_i(Z_i)\mid{\bb^*}\trans\X_i\right\}}
{E\left\{Y_i(Z_i)\mid{\bb^*}\trans\X_i\right\}}\right]+o_p(1)\\
&=&E\left(\Delta\frac{\blam_1(Z,{\bb^*}\trans\X)}{\lambda(Z,{\bb^*}\trans\X)}
\otimes\left[\X_{l}-
\frac{E\left\{\X_{l} 
	Y(Z)\mid{\bb^*}\trans\X\right\}}
{E\left\{Y(Z)\mid{\bb^*}\trans\X_i\right\}}\right]\right)+o_p(1)
\ese
under conditions \ref{assum:kernel}-\ref{assum:bandwidth} and 
\ref{assum:bounded}. Note that 
\bse
E\left(\Delta\frac{\blam_1(Z,{\bb^*}\trans\X)}{\lambda(Z,{\bb^*}\trans\X)}
\otimes\left[\X_{l}-
\frac{E\left\{\X_{l} 
	Y(Z)\mid{\bb^*}\trans\X\right\}}
{E\left\{Y(Z)\mid{\bb^*}\trans\X_i\right\}}\right]\right)
\ese
is a nonrandom quantity that does not depend on $n$, hence it is zero.
Thus the uniqueness requirement in Condition \ref{assum:unique} 
ensures that
$\bb^*=\bb_0$.

We now show that the subsequence that converges includes all but a
finite number of $n$'s. Assume this is not the case, then we can
obtain an infinite sequence of $\wh\bb_n$'s that do not converge 
to
$\bb^*$. As an infinite sequence in a compact set $\cal B$, we 
can
thus obtain another subsequence that converges, say to
$\bb^{**}\ne\bb^*$. Identical derivation as before then leads to
$\bb^{**}=\bb_0$, which is a contradiction to 
$\bb^{**}\ne\bb^*$.
Thus we conclude
$
\wh\bb-\bb_0\to\0
$
in probability when $n\to\infty$ under condition
\ref{assum:kernel}-\ref{assum:bounded}.
\qed

\subsection{Proof of Theorem \ref{th:eff}}\label{app:asympofb}

We only provide the proof concerning (\ref{eq:eff}); 
the result concerning (\ref{eq:general}) 
follows by using a similar and simpler proof.

We first expand (\ref{eq:eff}) as
\be
\0
&=&n^{-1/2}\sumi 
\Delta_i\frac{\wh\blam_1(Z_i,\wh\bb\trans\X_i)}{\wh\lambda(Z_i,\wh\bb\trans\X_i)}
\otimes\left[\X_{li}-
\frac{\wh E\left\{\X_{li} 
	Y_i(Z_i)\mid\wh\bb\trans\X_i\right\}}
{\wh 
E\left\{Y_i(Z_i)\mid\wh\bb\trans\X_i\right\}}\right]\nonumber\\
&=&n^{-1/2}\sumi 
\Delta_i\frac{\wh\blam_1(Z_i,\bb_0\trans\X_i)}{\wh\lambda(Z_i,\bb_0\trans\X_i)}
\otimes\left[\X_{li}-
\frac{\wh E\left\{\X_{li} 
	Y_i(Z_i)\mid\bb_0\trans\X_i\right\}}
{\wh 
E\left\{Y_i(Z_i)\mid\bb_0\trans\X_i\right\}}\right]\label{eq:main}\\
&&+\frac{1}{n}\sumi 
\left\{\frac{\partial}{\partial(\X_i\trans\bb)}
\left(\Delta_i\frac{\wh\blam_1(Z_i,\bb\trans\X_i)}{\wh\lambda(Z_i,\bb\trans\X_i)}
\otimes\left[\X_{li}-
\frac{\wh E\left\{\X_{li} 
	Y_i(Z_i)\mid\bb\trans\X_i\right\}}
{\wh
	E\left\{Y_i(Z_i)\mid\bb\trans\X_i\right\}}\right]\right)\otimes\X_{li}\trans\right\}\Bigg|_{\bb=\wt\bb}\nonumber\\
	\label{eq:easy}\\
&&\times\sqrt{n}(\wh\bb-\bb_0),\nonumber
\ee
where $\wt\bb$ is on the line connecting $\bb_0$ and $\wh\bb$.

We first consider (\ref{eq:easy}). Because of Theorem
\ref{th:consistency} and Lemma \ref{lem:pre}, 
\be
&&\frac{1}{n}\sumi 
\left\{\frac{\partial}{\partial(\X_i\trans\bb)}
\left(\Delta_i\frac{\wh\blam_1(Z_i,\bb\trans\X_i)}{\wh\lambda(Z_i,\bb\trans\X_i)}
\otimes\left[\X_{li}-
\frac{\wh E\left\{\X_{li} 
	Y_i(Z_i)\mid\bb\trans\X_i\right\}}
{\wh
	E\left\{Y_i(Z_i)\mid\bb\trans\X_i\right\}}\right]\right)\otimes\X_{li}\trans\right\}\Bigg|_{\bb=\wt\bb}\n\\
&=&\frac{1}{n}\sumi 
\left\{\frac{\partial}{\partial(\X_i\trans\bb_0)}
\left(\Delta_i\frac{\wh\blam_1(Z_i,\bb_0\trans\X_i)}{\wh\lambda(Z_i,\bb_0\trans\X_i)}
\otimes\left[\X_{li}-
\frac{\wh E\left\{\X_{li} 
	Y_i(Z_i)\mid\bb_0\trans\X_i\right\}}
{\wh
	E\left\{Y_i(Z_i)\mid\bb_0\trans\X_i\right\}}\right]\right)\otimes\X_{li}\trans\right\}\n\\
&&+o_p(1)\n\\
&=&-\frac{1}{n}\sumi \left(\Delta_i
\frac{\wh\blam_1^{\otimes2}(Z_i,\bb_0\trans\X_i)}{\wh\lambda^2(Z_i,\bb_0\trans\X_i)}
\otimes\left[\X_{li}-
\frac{\wh E\left\{\X_{li} 
	Y_i(Z_i)\mid\bb_0\trans\X_i\right\}}
{\wh
	E\left\{Y_i(Z_i)\mid\bb_0\trans\X_i\right\}}\right]\otimes\X_{li}\trans\right)\label{eq:easy1}\\
&&+\frac{1}{n}\sumi 
\frac{\Delta_i}{\wh\lambda(Z_i,\bb_0\trans\X_i)}
\frac{\partial}{\partial(\X_i\trans\bb_0)}
\left(\wh\blam_1(Z_i,\bb_0\trans\X_i)
\otimes\left[\X_{li}-
\frac{\wh E\left\{\X_{li} 
	Y_i(Z_i)\mid\bb_0\trans\X_i\right\}}
{\wh
	E\left\{Y_i(Z_i)\mid\bb_0\trans\X_i\right\}}\right]\right)\otimes\X_{li}\trans\n\\
\label{eq:easy2}\\
&&+o_p(1)\nonumber.
\ee
Because of Lemma \ref{lem:pre}, (\ref{eq:easy1}) converges 
uniformly
in probability to
\bse
&&-E \left(\int_0^\infty 
\frac{\blam_1^{\otimes2}(s,\bb_0\trans\X)}{\lambda^2(s,\bb_0\trans\X)}
\otimes\left[\X_l-
\frac{ E\left\{\X_l 
	Y(s)\mid\bb_0\trans\X\right\}}
{E\left\{Y(s)\mid\bb_0\trans\X\right\}}\right]\otimes\X_l\trans
dN(s)\right)\\
&=&-E \left(\int_0^\infty 
\frac{\blam_1^{\otimes2}(s,\bb_0\trans\X)}{\lambda^2(s,\bb_0\trans\X)}
\otimes\left[\X_l-
\frac{ E\left\{\X_l 
	Y(s)\mid\bb_0\trans\X\right\}}
{E\left\{Y(s)\mid\bb_0\trans\X\right\}}\right]\otimes\X_l\trans 
Y(s)\lambda(s,\bb_0\trans\X)ds\right)\\
&=&-E \left(\int_0^\infty 
\frac{\blam_1^{\otimes2}(s,\bb_0\trans\X)}{\lambda(s,\bb_0\trans\X)}
\otimes\left[\X_l-
\frac{ E\left\{\X_l 
	Y(s)\mid\bb_0\trans\X\right\}}
{E\left\{Y(s)\mid\bb_0\trans\X\right\}}\right]\otimes
\left[\X_l
-\frac{ E\left\{\X_l 
	Y(s)\mid\bb_0\trans\X\right\}}
{E\left\{Y(s)\mid\bb_0\trans\X\right\}}\right]
\trans 
Y(s)ds\right)\\
&&-E \left(\int_0^\infty 
\frac{\blam_1^{\otimes2}(s,\bb_0\trans\X)}{\lambda(s,\bb_0\trans\X)}
\otimes\left[\X_l-
\frac{ E\left\{\X_l 
	Y(s)\mid\bb_0\trans\X\right\}}
{E\left\{Y(s)\mid\bb_0\trans\X\right\}}\right]\otimes
\frac{ E\left\{\X_l 
	Y(s)\mid\bb_0\trans\X\right\}}
{E\left\{Y(s)\mid\bb_0\trans\X\right\}}\trans 
Y(s)ds\right)\\
&=&-E\{\bS\eff(\Delta,Z,\X)^{\otimes2}\},
\ese
where the last equality is because the second term above is zero 
by 
first taking expectation conditional on $\bb_0\trans\X$.

Similarly, from Lemma \ref{lem:pre},  the term in 
(\ref{eq:easy2})
converges uniformly in probability to the
limit of
\bse
E \left\{\frac{\Delta_i}{\lambda(Z_i,\bb_0\trans\X_i)}
\frac{\partial}{\partial(\X_i\trans\bb_0)}
\left(\wh\blam_1(Z_i,\bb_0\trans\X_i)
\otimes \left[\X_{li}-
\frac{ E\left\{\X_{li} Y_i(Z_i)\mid\bb_0\trans\X_i\right\}}
{E\left\{Y_i(Z_i)\mid\bb_0\trans\X_i\right\}}
\right]
\right)\otimes\X_{li}\trans\right\}.
\ese
Now let $\wh\blam_{1,-i}(Z, \bb_0\trans\X)$ be the
leave-one-out version of $\wh\blam_{1}(Z, \bb_0\trans\X)$, i.e. 
it
is constructed the same as $\wh\blam_{1}(Z, \bb_0\trans\X)$ 
except
that the $i$th observation is not used.
Obviously, 
\bse
&&\frac{\Delta_i}{\lambda(Z_i,\bb_0\trans\X_i)}
\frac{\partial}{\partial(\X_i\trans\bb_0)}
\left(\wh\blam_1(Z_i,\bb_0\trans\X_i)
\otimes \left[\X_{li}-
\frac{ E\left\{\X_{li} 
	Y_i(Z_i)\mid\bb_0\trans\X_i\right\}}
{E\left\{Y_i(Z_i)\mid\bb_0\trans\X_i\right\}}\right]\right)\otimes\X_{li}\trans\\
&-&\frac{\Delta_i}{\lambda(Z_i,\bb_0\trans\X_i)}
\frac{\partial}{\partial(\X_i\trans\bb_0)}
\left(\wh\blam_{1,-i}(Z_i,\bb_0\trans\X_i)
\otimes \left[\X_{li}-
\frac{ E\left\{\X_{li} 
	Y_i(Z_i)\mid\bb_0\trans\X_i\right\}}
{E\left\{Y_i(Z_i)\mid\bb_0\trans\X_i\right\}}\right]\right)\otimes\X_{li}\trans\\
&=&o_p(1).
\ese
Let $E_i$ mean taking expectation with respect to the $i$th
observation conditional on all other observations, then
\bse
&&
E_i\left\{\frac{\Delta_i}{\lambda(Z_i,\bb_0\trans\X_i)}
\frac{\partial}{\partial(\X_i\trans\bb_0)}
\left(\wh\blam_{1,-i}(Z_i,\bb_0\trans\X_i)
\otimes \left[\X_{li}-
\frac{ E\left\{\X_{li} 
	Y_i(Z_i)\mid\bb_0\trans\X_i\right\}}
{E\left\{Y_i(Z_i)\mid\bb_0\trans\X_i\right\}}\right]\right)\otimes\X_{li}\trans\right\}\\
&=&E_i \left\{
\frac{\partial}{\partial\bb_0}\int 
\wh\blam_{1,-i}(s,\bb_0\trans\X_i) 
\otimes \left[\X_{li}-
\frac{ E\left\{\X_{li} 
	Y_i(s)\mid\bb_0\trans\X_i\right\}}
{E\left\{Y_i(s)\mid\bb_0\trans\X_i\right\}}\right]
E\{Y_i(s)\mid\X_i\}ds\right\}\\
&=&\frac{\partial}{\partial\bb_0} E_i \left\{
\int 
\wh\blam_{1,-i}(s,\bb_0\trans\X_i) 
\otimes \left[\X_{li}-
\frac{ E\left\{\X_{li} 
	Y_i(s)\mid\bb_0\trans\X_i\right\}}
{E\left\{Y_i(s)\mid\bb_0\trans\X_i\right\}}\right]
Y_i(s)ds\right\}\\
&=&\0.
\ese
Here, the last equality is because the integrand has expectation zero
conditional on $\bb_0\trans\X_i$ and all other observations, and the third to last 
equality is because
the expectation is with respect to $\X_i$ and does not involve 
$\bb_0$. 
Therefore, the term in (\ref{eq:easy2}) converges in probability
uniformly to
\bse
E\left\{\frac{\Delta_i}{\lambda(Z_i,\bb_0\trans\X_i)}
\frac{\partial}{\partial(\X_i\trans\bb_0)}
\left(\wh\blam_{1,-i}(Z_i,\bb_0\trans\X_i)
\otimes \left[\X_{li}-
\frac{ E\left\{\X_{li} 
	Y_i(Z_i)\mid\bb_0\trans\X_i\right\}}
{E\left\{Y_i(Z_i)\mid\bb_0\trans\X_i\right\}}\right]\right)\otimes\X_{li}\trans\right\}
=0
\ese
Combining the results concerning (\ref{eq:easy1}) and
(\ref{eq:easy2}), thus the expression in 
(\ref{eq:easy}) 
is 
$-E\{\bS\eff(\Delta,Z,\X)^{\otimes2}\}+o_p(1)$.

Next we decompose (\ref{eq:main}) into 
\be\label{eq:Ts}
n^{-1/2}\sumi 
\Delta_i\frac{\wh\blam_1(Z_i,\bb_0\trans\X_i)}{\wh\lambda(Z_i,\bb_0\trans\X_i)}
\otimes\left[\X_{li}-
\frac{\wh E\left\{\X_{li} 
	Y_i(Z_i)\mid\bb_0\trans\X_i\right\}}
{\wh E\left\{Y_i(Z_i)\mid\bb_0\trans\X_i\right\}}\right]
=\T_1+\T_2+\T_3+\T_4,
\ee
where
\bse
\T_1
&=&n^{-1/2}\sumi 
\Delta_i\frac{\blam_1(Z_i,\bb_0\trans\X_i)}{\lambda(Z_i,\bb_0\trans\X_i)}
\otimes\left[\X_{li}-
\frac{E\left\{\X_{li} 
	Y_i(Z_i)\mid\bb_0\trans\X_i\right\}}
{E\left\{Y_i(Z_i)\mid\bb_0\trans\X_i\right\}}\right],\\
\T_2&=&n^{-1/2}\sumi 
\Delta_i\left\{\frac{\wh\blam_1(Z_i,\bb_0\trans\X_i)}{\wh\lambda(Z_i,\bb_0\trans\X_i)}
-\frac{\blam_1(Z_i,\bb_0\trans\X_i)}{\lambda(Z_i,\bb_0\trans\X_i)}\right\}\otimes\left[\X_{li}-
\frac{E\left\{\X_{li} 
	Y_i(Z_i)\mid\bb_0\trans\X_i\right\}}
{E\left\{Y_i(Z_i)\mid\bb_0\trans\X_i\right\}}\right],\\
\T_3&=&n^{-1/2}\sumi 
\Delta_i\frac{\blam_1(Z_i,\bb_0\trans\X_i)}{\lambda(Z_i,\bb_0\trans\X_i)}
\otimes\left[\frac{E\left\{\X_{li} 
	Y_i(Z_i)\mid\bb_0\trans\X_i\right\}}
{E\left\{Y_i(Z_i)\mid\bb_0\trans\X_i\right\}}-
\frac{\wh E\left\{\X_{li} 
	Y_i(Z_i)\mid\bb_0\trans\X_i\right\}}
{\wh
	E\left\{Y_i(Z_i)\mid\bb_0\trans\X_i\right\}}\right],\\
\T_4&=&n^{-1/2}\sumi \Delta_i
\left\{\frac{\wh\blam_1(Z_i,\bb_0\trans\X_i)}{\wh\lambda(Z_i,\bb_0\trans\X_i)}-
\frac{\blam_1(Z_i,\bb_0\trans\X_i)}{\lambda(Z_i,\bb_0\trans\X_i)}\right\}\\
&&\otimes\left[\frac{E\left\{\X_{li} 
	Y_i(Z_i)\mid\bb_0\trans\X_i\right\}}
{E\left\{Y_i(Z_i)\mid\bb_0\trans\X_i\right\}}-
\frac{\wh E\left\{\X_{li} 
	Y_i(Z_i)\mid\bb_0\trans\X_i\right\}}
{\wh
	E\left\{Y_i(Z_i)\mid\bb_0\trans\X_i\right\}}\right].
\ese

First, 
\bse
\T_2&=&n^{-1/2}\sumi 
\int\left\{\frac{\wh\blam_1(s,\bb_0\trans\X_i)}{\wh\lambda(s,\bb_0\trans\X_i)}
-\frac{\blam_1(s,\bb_0\trans\X_i)}{\lambda(s,\bb_0\trans\X_i)}\right\}\otimes\left[\X_{li}-
\frac{E\left\{\X_{li} 
	Y_i(s)\mid\bb_0\trans\X_i\right\}}
{E\left\{Y_i(s)\mid\bb_0\trans\X_i\right\}}\right]dN_i(s)\\
&=&o_p\left(n^{-1/2}\sumi \int\left[\X_{li}-
\frac{E\left\{\X_{li} 
	Y_i(s)\mid\bb_0\trans\X_i\right\}}
{E\left\{Y_i(s)\mid\bb_0\trans\X_i\right\}}\right]Y_i(s)\lambda(s,\bb_0\trans\X_{li})ds\right)\\
&=&o_p(1),
\ese
where the last equality above is because the quantity inside the parentheses is 
a mean zero normal random
quantity of order $O_p(1)$. Further,
\bse
\T_3&=&n^{-1/2}\sumi 
\Delta_i\frac{\blam_1(Z_i,\bb_0\trans\X_i)}{\lambda(Z_i,\bb_0\trans\X_i)}
\otimes\left(-
\frac{\wh E\left\{\X_{li} 
	Y_i(Z_i)\mid\bb_0\trans\X_i\right\}}{E\left\{Y_i(Z_i)\mid\bb_0\trans\X_i\right\}}\right.\n\\
&&\left. +\frac{\wh 
	E\left\{Y_i(Z_i)\mid\bb_0\trans\X_i\right\}E\left\{\X_{li}Y_i(Z_i)\mid\bb_0\trans\X_i\right\}}
{[E\left\{Y_i(Z_i)\mid\bb_0\trans\X_i\right\}]^2}\right) 
+o_p(1)\n\\
&=&
n^{-1/2}\sumi 
\Delta_i\frac{\blam_1(Z_i,\bb_0\trans\X_i)}{\lambda(Z_i,\bb_0\trans\X_i)}
\otimes\left(-
\frac{
	n^{-1}\sumj 
	K_h(\bb_0\trans\X_j-\bb_0\trans\X_i)\X_{lj}I(Z_j\ge
	Z_i)}{f_{\bb_0\trans\X}(\bb_0\trans\X_i) 
	E\left\{Y_i(Z_i)\mid\bb_0\trans\X_i\right\}}\right.\n\\
&&\left.+\frac{E\left\{\X_{li}Y_i(Z_i)\mid\bb_0\trans\X_i\right\}
	\{n^{-1}\sumj 
	K_h(\bb_0\trans\X_j-\bb_0\trans\X_i)-f_{\bb_0\trans\X}(\bb_0\trans\X_i)
	\}}{f_{\bb_0\trans\X}(\bb_0\trans\X_i) 
	E\left\{Y_i(Z_i)\mid\bb_0\trans\X_i\right\}}\right.\n\\
&&\left. 
+\frac{E\left\{\X_{li}Y_i(Z_i)\mid\bb_0\trans\X_i\right\}}
{[E\left\{Y_i(Z_i)\mid\bb_0\trans\X_i\right\}]^2}
\left[\frac{
	n^{-1}\sumj K_h(\bb_0\trans\X_j-\bb_0\trans\X_i)I(Z_j\ge
	Z_i)}{f_{\bb_0\trans\X}(\bb_0\trans\X_i) }\right.\right.\n\\
&&\left.\left.-\frac{E\left\{Y_i(Z_i)\mid\bb_0\trans\X_i\right\}
	\{n^{-1}\sumj 
	K_h(\bb_0\trans\X_j-\bb_0\trans\X_i)-f_{\bb_0\trans\X}(\bb_0\trans\X_i)
	\}}{f_{\bb_0\trans\X}(\bb_0\trans\X_i) }\right]
\right) +o_p(1)\n\\
&=&n^{-3/2}\sumi\sumj 
\Delta_i\frac{\blam_1(Z_i,\bb_0\trans\X_i)}{\lambda(Z_i,\bb_0\trans\X_i)}
\otimes\left[
-\frac{K_h(\bb_0\trans\X_j-\bb_0\trans\X_i)\X_{lj}I(Z_j\ge 
	Z_i)}{f_{\bb_0\trans\X}(\bb_0\trans\X_i) 
	E\left\{Y_i(Z_i)\mid\bb_0\trans\X_i\right\}}\right.\n\\
&&\left. +
\frac{E\left\{\X_{li}Y_i(Z_i)\mid\bb_0\trans\X_i\right\} 
	K_h(\bb_0\trans\X_j-\bb_0\trans\X_i)I(Z_j\ge 
	Z_i)}{f_{\bb_0\trans\X}(\bb_0\trans\X_i) 
	[E\left\{Y_i(Z_i)\mid\bb_0\trans\X_i\right\}]^2
}\right] +o_p(1)\n\\
&=&\T_{31}+\T_{32}+\T_{33}+o_p(1),
\ese
where
\bse
\T_{31}&=&n^{-1/2}\sumi\Delta_i
\frac{\blam_1(Z_i,\bb_0\trans\X_i)}{\lambda(Z_i,\bb_0\trans\X_i)}
\otimes E\left[
-\frac{K_h(\bb_0\trans\X_j-\bb_0\trans\X_i)\X_{lj}I(Z_j\ge 
	Z_i)}{f_{\bb_0\trans\X}(\bb_0\trans\X_i) 
	E\left\{Y_i(Z_i)\mid\bb_0\trans\X_i\right\}}\right.\n\\
&&\left. +
\frac{E\left\{\X_{li}Y_i(Z_i)\mid\bb_0\trans\X_i\right\} 
	K_h(\bb_0\trans\X_j-\bb_0\trans\X_i)I(Z_j\ge 
	Z_i)}{f_{\bb_0\trans\X}(\bb_0\trans\X_i) 
	[E\left\{Y_i(Z_i)\mid\bb_0\trans\X_i\right\}]^2
}\mid \Delta_i, Z_i, \X_i\right]\\
\T_{32}&=&n^{-1/2}\sumj E 
\left(\Delta_i\frac{\blam_1(Z_i,\bb_0\trans\X_i)}{\lambda(Z_i,\bb_0\trans\X_i)}
\otimes\left[
-\frac{K_h(\bb_0\trans\X_j-\bb_0\trans\X_i)\X_{lj}I(Z_j\ge 
	Z_i)}{f_{\bb_0\trans\X}(\bb_0\trans\X_i) 
	E\left\{Y_i(Z_i)\mid\bb_0\trans\X_i\right\}}\right.\right.\n\\
&&\left. \left.+
\frac{E\left\{\X_{li}Y_i(Z_i)\mid\bb_0\trans\X_i\right\} 
	K_h(\bb_0\trans\X_j-\bb_0\trans\X_i)I(Z_j\ge 
	Z_i)}{f_{\bb_0\trans\X}(\bb_0\trans\X_i) 
	[E\left\{Y_i(Z_i)\mid\bb_0\trans\X_i\right\}]^2
}\right]\mid \Delta_j, Z_j,\X_j\right)\\
\T_{33}&=&-n^{1/2}E\left(\Delta_i
\frac{\blam_1(Z_i,\bb_0\trans\X_i)}{\lambda(Z_i,\bb_0\trans\X_i)}
\otimes E\left[
-\frac{K_h(\bb_0\trans\X_j-\bb_0\trans\X_i)\X_{lj}I(Z_j\ge 
	Z_i)}{f_{\bb_0\trans\X}(\bb_0\trans\X_i) 
	E\left\{Y_i(Z_i)\mid\bb_0\trans\X_i\right\}}\right.\right.\n\\
&&\left.\left. +
\frac{E\left\{\X_{li}Y_i(Z_i)\mid\bb_0\trans\X_i\right\} 
	K_h(\bb_0\trans\X_j-\bb_0\trans\X_i)I(Z_j\ge 
	Z_i)}{f_{\bb_0\trans\X}(\bb_0\trans\X_i) 
	[E\left\{Y_i(Z_i)\mid\bb_0\trans\X_i\right\}]^2
}\right]\right).
\ese
Here we used U-statistic property in the last equality above. 
Now when 
$nh^4\to0$,
\bse
\T_{31}&=&n^{-1/2}\sumi\Delta_i
\frac{\blam_1(Z_i,\bb_0\trans\X_i)}{\lambda(Z_i,\bb_0\trans\X_i)}
\otimes \left[
-\frac{E\{\X_{li}Y_i(Z_i)\mid 
	\bb_0\trans\X_i\}}{E\left\{Y_i(Z_i)\mid\bb_0\trans\X_i\right\}}\right.\\
&&\left.+
\frac{E\left\{\X_{li}Y_i(Z_i)\mid\bb_0\trans\X_i\right\}E\left\{Y_i(Z_i)\mid\bb_0\trans\X_i\right\}}{
	[E\left\{Y_i(Z_i)\mid\bb_0\trans\X_i\right\}]^2
}\right]+O(n^{1/2}h^2)\\
&=&o_p(1).
\ese
Thus, $\T_{33}=o_p(1)$ as well.
To analyze $\T_{32}$, 
\bse
\T_{32}&=&n^{-1/2}\sumj E 
\left(\Delta_i\frac{\blam_1(Z_i,\bb_0\trans\X_i)}{\lambda(Z_i,\bb_0\trans\X_i)}
\otimes\left[
-\frac{K_h(\bb_0\trans\X_j-\bb_0\trans\X_i)\X_{lj}I(Z_j\ge 
	Z_i)}{f_{\bb_0\trans\X}(\bb_0\trans\X_i) E\left\{I(Z\ge
	Z_i)\mid\bb_0\trans\X=\bb_0\trans\X_i, 
	Z_i\right\}}\right.\right.\n\\
&&\left. \left.+
\frac{E\left\{\X_{l}I(Z\ge 
Z_i)\mid\bb_0\trans\X=\bb_0\trans\X_i, Z_i\right\} 
	K_h(\bb_0\trans\X_j-\bb_0\trans\X_i)I(Z_j\ge 
	Z_i)}{f_{\bb_0\trans\X}(\bb_0\trans\X_i) [E\left\{I(Z\ge
	Z_i)\mid\bb_0\trans\X=\bb_0\trans\X_i, Z_i\right\}]^2
}\right]\mid \Delta_j, Z_j,\X_j\right)\\
&=&n^{-1/2}\sumj E\left\{E 
\left(\Delta_i\frac{\blam_1(Z_i,\bb_0\trans\X_i)}{\lambda(Z_i,\bb_0\trans\X_i)}
\otimes\left[
-\frac{\x_{lj}I(z_j\ge 
	Z_i)}{f_{\bb_0\trans\X}(\bb_0\trans\X_i) E\left\{I(Z\ge
	Z_i)\mid\bb_0\trans\X=\bb_0\trans\X_i, 
	Z_i\right\}}\right.\right.\right.\n\\
&&\left. \left.\left.+
\frac{E\left\{\X_{l}I(Z\ge 
Z_i)\mid\bb_0\trans\X=\bb_0\trans\X_i, Z_i\right\} 
	I(z_j\ge 
	Z_i)}{f_{\bb_0\trans\X}(\bb_0\trans\X_i) [E\left\{I(Z\ge
	Z_i)\mid\bb_0\trans\X=\bb_0\trans\X_i, Z_i\right\}]^2
}\right]\mid \bb_0\trans\X_i\right)
K_h(\bb_0\trans\x_j-\bb_0\trans\X_i)\right\}\\
&=&n^{-1/2}\sumj 
E\left(\int_0^{z_j}\frac{\blam_1(s,\bb_0\trans\x_j)}{E\left\{S_c(s,\X)\mid\bb_0\trans\X=\bb_0\trans\x_j\right\}}\right.\\
&&\left.\otimes\left[
\frac{E 
\left\{\X_lS_c(s,\X)\mid\bb_0\trans\X=\bb_0\trans\x_j\right\}
}{E\left\{S_c(s,\X)\mid\bb_0\trans\X=\bb_0\trans\x_j\right\}
}-\x_{lj}\right]
S_c(s,\X_i)ds
\mid\bb_0\trans\X_i=\bb_0\trans\x_j\right) +O_p(n^{1/2}h^2)\\
&=&n^{-1/2}\sumj \int 
Y_j(s)\lambda(s,\bb_0\trans\x_j)\frac{\blam_1(s,\bb_0\trans\x_j)}{\lambda(s,\bb_0\trans\x_j)}\otimes\left[
\frac{E\left\{\X_{lj}Y_j(s)\mid\bb_0\trans\x_j\right\}
}{E\left\{Y_j(s)\mid\bb_0\trans\x_j\right\}
}-\x_{lj}\right]ds +O_p(n^{1/2}h^2).
\ese
When $nh^4\to0$, plugging the results of $\T_1$ and $\T_{32}$ to
(\ref{eq:Ts}), the expression in (\ref{eq:main}) 
is 
\bse
&&n^{-1/2}\sumi 
\Delta_i\frac{\wh\blam_1(Z_i,\bb_0\trans\X_i)}{\wh\lambda(Z_i,\bb_0\trans\X_i)}
\otimes\left[\X_{li}-
\frac{\wh E\left\{\X_{li} 
	Y_i(Z_i)\mid\bb_0\trans\X_i\right\}}
{\wh E\left\{Y_i(Z_i)\mid\bb_0\trans\X_i\right\}}\right]\\
&=&n^{-1/2}\sumi 
\int\frac{\blam_1(t,\bb_0\trans\X_i)}{\lambda(t,\bb_0\trans\X_i)}
\otimes\left[\X_{li}-
\frac{E\left\{\X_{li} 
	Y_i(t)\mid\bb_0\trans\X_i\right\}}
{E\left\{Y_i(t)\mid\bb_0\trans\X_i\right\}}\right]dM_i(t)+o_p(1)\\
&=&n^{-1/2}\sumi\bS\eff(\Delta_i,Z_i,\X_i)+o_p(1).
\ese

Finally, 
\bse
\T_4&=&n^{-1/2}\sumi \Delta_i
\left\{\frac{\wh\blam_1(Z_i,\bb_0\trans\X_i)}{\wh\lambda(Z_i,\bb_0\trans\X_i)}-
\frac{\blam_1(Z_i,\bb_0\trans\X_i)}{\lambda(Z_i,\bb_0\trans\X_i)}\right\}\\
&&\times\left[\frac{E\left\{\X_{li} 
	Y_i(Z_i)\mid\bb_0\trans\X_i\right\}}
{E\left\{Y_i(Z_i)\mid\bb_0\trans\X_i\right\}}-
\frac{\wh E\left\{\X_{li} 
	Y_i(Z_i)\mid\bb_0\trans\X_i\right\}}
{\wh
	E\left\{Y_i(Z_i)\mid\bb_0\trans\X_i\right\}}\right]\\
&=&o_p\left(
n^{-1/2}\sumi \Delta_i
\left[\frac{E\left\{\X_{li} 
	Y_i(Z_i)\mid\bb_0\trans\X_i\right\}}
{E\left\{Y_i(Z_i)\mid\bb_0\trans\X_i\right\}}-
\frac{\wh E\left\{\X_{li} 
	Y_i(Z_i)\mid\bb_0\trans\X_i\right\}}
{\wh
	E\left\{Y_i(Z_i)\mid\bb_0\trans\X_i\right\}}\right]\right)\\
&=&o_p\left(n^{-1/2}\sumi \int Y_i(s)\lambda(s,\bb_0\trans\x_i)
\left[
\frac{E\left\{\X_{li}Y_i(s)\mid\bb_0\trans\x_i\right\}
}{E\left\{Y_i(s)\mid\bb_0\trans\x_i\right\}
}-\x_{li}\right]ds\right)+o_p(n^{1/2}h^2)\\
&=&o_p(1),
\ese
where the last equality is because the integrand has mean zero
conditional on $\bb_0\trans\X$, and the second to last equality is 
obtained
following the same derivation of $\T_3$.
Using these results in (\ref{eq:main}), combined with the 
results on
(\ref{eq:easy}), it is clear that the theorem holds. 
\qed

\subsection{Proof of Theorem \ref{th:m}}\label{app:m}

We expand 
$\sqrt{nh}\left\{\wh{m}(t,\bb\trans\x)-{m}(t,\bb\trans\x)\right\}$
 as 
\be 
&&\sqrt{nh}\left\{\wh{m}(t,\bb\trans\x)-{m}(t,\bb\trans\x)\right\}\n\\
&=&\sqrt{nh}\left\{e^{{\wh\Lambda}(t,\bb\trans\x)}- 
e^{{\Lambda}(t,\bb\trans\x)}\right\}\int_{t}^{\infty} 
e^{-{\Lambda}(s,\bb\trans\x)}ds\label{eq:mnormal1}\\
&&+\sqrt{nh}e^{{\Lambda}(t,\bb\trans\x)}\int_{t}^{\infty} 
\left\{e^{-{\wh\Lambda}(s,\bb\trans\x)}-e^{-{\Lambda}(s,\bb\trans\x)}\right\}ds.\label{eq:mnormal2}\\
&&+\sqrt{nh}\left\{e^{{\wh\Lambda}(t,\bb\trans\x)}- 
e^{{\Lambda}(t,\bb\trans\x)}\right\}\int_{t}^{\infty}
\left\{e^{-{\wh\Lambda}(s,\bb\trans\x)}-e^{-{\Lambda}(s,\bb\trans\x)}\right\}ds\label{eq:mnormal3}.
\ee 
It is easy to see that the term in (\ref{eq:mnormal3}) satisfies 
\bse 
&&\sqrt{nh}\left\{e^{{\wh\Lambda}(t,\bb\trans\x)}- 
e^{{\Lambda}(t,\bb\trans\x)}\right\}\int_{t}^{\infty}
\left\{e^{-{\wh\Lambda}(s,\bb\trans\x)}- 
e^{-{\Lambda}(s,\bb\trans\x)}\right\}ds\\
&=&\sqrt{nh}O_p\{\wh\Lambda(t,\bb\trans\x)-\Lambda(t,\bb\trans\x)\}
\int_{t}^{\infty}e^{-\Lambda(s,\bb\trans\x)} 
O_p\{\wh\Lambda(s,\bb\trans\x)-\Lambda(s,\bb\trans\x)\}ds\\
&=&O_p(\sqrt{nh})O_p\{h^4+(nh)^{-1}\}\\
&=&o_p(1) 
\ese 
by Condition \ref{assum:bandwidth}.

We inspect the terms in (\ref{eq:mnormal1}) and 
(\ref{eq:mnormal2}). 
For (\ref{eq:mnormal1}), based on Lemma 
\ref{lem:pre}, 
\bse 
&&\sqrt{nh}\left\{e^{{\wh\Lambda}(t,\bb\trans\x)}- 
e^{\Lambda(t,\bb\trans\x)}\right\}\int_{t}^{\infty} 
e^{-{\Lambda}(s,\bb\trans\x)}ds\\
&=&\sqrt{nh}e^{\Lambda(t,\bb\trans\x)} 
\left(\wh\Lambda(t,\bb\trans\x)-\Lambda(t,\bb\trans\x)+ 
O_p[\{\wh\Lambda(t,\bb\trans\x)-\Lambda(t,\bb\trans\x)\}^2]\right)
\int_{t}^{\infty}
e^{-{\Lambda}(s,\bb\trans\x)}ds\\
&=&\sqrt{nh}e^{\Lambda(t,\bb\trans\x)} 
\left\{\wh\Lambda(t,\bb\trans\x)-\Lambda(t,\bb\trans\x)\right\} 
\int_{t}^{\infty}e^{-{\Lambda}(s,\bb\trans\x)}ds+o_p(1), 
\ese 
where the last step uses Condition 
\ref{assum:bandwidth}.

For (\ref{eq:mnormal3}), using Condition \ref{assum:bandwidth} 
as well, we get 
\bse 
&&\sqrt{nh}e^{{\Lambda}(t,\bb\trans\x)}\int_{t}^{\infty} 
\left\{e^{-{\wh\Lambda}(s,\bb\trans\x)}-e^{-{\Lambda}(s,\bb\trans\x)}\right\}ds\\
&=&\sqrt{nh}e^{{\Lambda}(t,\bb\trans\x)}\int_{t}^{\infty} 
\left[\wh\Lambda(s,\bb\trans\x)-\Lambda(s,\bb\trans\x)+O_p\{(nh)^{-1}+h^4\}\right]e^{-\Lambda(s,\bb\trans\x)}ds\\
&=&\sqrt{nh}e^{{\Lambda}(t,\bb\trans\x)}\int_{t}^{\infty} 
\left\{\wh\Lambda(s,\bb\trans\x)-\Lambda(s,\bb\trans\x)\right\} 
e^{-\Lambda(s,\bb\trans\x)}ds+o_p(1). 
\ese

Now combine the leading terms in (\ref{eq:mnormal1}) and 
(\ref{eq:mnormal2}) and use the expansion of 
$\wh\Lambda(t,\bb\trans\x)-\Lambda(t,\bb\trans\x)$ in Lemma 
\ref{th:Lambda}, 
\be 
&&\sqrt{nh}e^{\Lambda(t,\bb\trans\x)}\left\{\wh\Lambda(t,\bb\trans\x)-\Lambda(t,\bb\trans\x)\right\}\int_{t}^{\infty}
e^{-{\Lambda}(s,\bb\trans\x)}ds\n\\
&&+\sqrt{nh}e^{{\Lambda}(t,\bb\trans\x)}\int_{t}^{\infty}
\left\{\wh\Lambda(s,\bb\trans\x)-\Lambda(s,\bb\trans\x)\right\} 
e^{-\Lambda(s,\bb\trans\x)}ds\n\\
&=&e^{\Lambda(t,\bb\trans\x)}\sumi\int_{0}^{\infty}\sqrt{\frac{h}{n}}
\frac{I\left\{\phi_n(r,\bb\trans\x)>0\right\}K_h(\bb\trans\X_i-\bb\trans\x)}
{f_{\bb\trans\X}(\bb\trans\x)E\{I(Z\ge 
	r)\mid\bb\trans\x\}}\n\\
&&\times\left\{I(r<t)\int_{t}^{\infty}
e^{-\Lambda(s,\bb\trans\x)}ds+\int_{\max(r,t)}^{\infty} 
e^{-\Lambda(s,\bb\trans\x)}ds\right\}dM_i(r,\bb\trans\x) 
+o_p(1).\label{eq:mnormallong}
\ee 

Note that $I\left\{\phi_n(r,\bb\trans\x)>0\right\}=1$ almost 
surely and 
according to Lemma \ref{lem:pre}
\bse 
\frac{1}{n}\sumi hK_h^2(\bb\trans\X_i-\bb\trans\x)Y_i(r)=
f_{\bb\trans\X}(\bb\trans\x)E\{I(Z\ge 
r)\mid\bb\trans\x\}\int K^2(u)du +o_p(1). 
\ese 
The leading term in 
(\ref{eq:mnormallong}) converges to 
$N\{0,\sigma^2_m(t,\bb\trans\x)\}$ uniformly by martingale 
central limit 
theorem, 
where 
\bse 
&&\sigma^2_m(t,\bb\trans\x)\\
&=&e^{2\Lambda(t,\bb\trans\x)}\frac{\int 
	K^2(u)du}{f_{\bb\trans\X}(\bb\trans\x)}
\int_{0}^{\infty}\frac{\lambda(r,\bb\trans\x)}
{E\{I(Z\ge r)\mid\bb\trans\x\}}\\
&&\times\left\{I(r<t)\int_{t}^{\infty}
e^{-\Lambda(s,\bb\trans\x)}ds+\int_{\max(r,t)}^{\infty} 
e^{-\Lambda(s,\bb\trans\x)}ds\right\}^2dr. 
\ese 
Therefore 
$\sqrt{nh}\left\{\wh{m}(t,\bb\trans\x)-{m}(t,\bb\trans\x)\right\}\to
 N\{0,\sigma^2_m(t,\bb\trans\x)\}$ uniformly for all $t$ and 
$\bb\trans\x$. 
\qed 
\newpage

\LTcapwidth=\textwidth

\spacing{1.5}
\begin{longtable}{cc|cccccc}
	\caption{Results of  study 2, based on 1000
		simulations with sample size 500. ``Prop.'' is the
		semiparametric method, ``PM1'' and ``PM2" are the
		proportional mean residual life
		methods, ``additive'' is the additive
		method. ``bias'' is the average absolute bias of
		each component in $\wh\bb$, ``sd'' 
		is the
		sample standard deviation 
		of the
		corresponding estimators. ``MSE" is the mean squared 
		error.}
	\label{tab:simu2}\\\hline
	\vspace{.2cm}
		&& $\beta_2$ &$\beta_3$ &$\beta_4$ &$\beta_5$
		&$\beta_6$      &$\beta_7$ \\
		&true& $1.3$ &$-1.3$ &$1$ &$-0.5$ &$0.5$      &$-0.5$
		\\\hline
		&&\multicolumn{6}{c}{No censoring}\\
		Prop.& bias&0.001 &0.002 &0.000 &0.000 &0.001 &0.000 \\
		& sd &0.030 &0.030 &0.026 &0.026 &0.022 &0.023\\
		& MSE&0.001 &0.001 &0.001 &0.001 &0.000 &0.001\\
		PM1	& bias&0.003 &0.003 &0.003 &0.002 &0.002 &0.001\\
		& sd  &0.061 &0.061 &0.051 &0.040 &0.040 &0.042\\
		& MSE&0.004 &0.004 &0.003 &0.002 &0.002 &0.002\\
		PM2	& bias&0.005 &0.005 &0.004 &0.003 &0.002 &0.001\\
		& sd  &0.073 &0.071 &0.059 &0.046 &0.046 &0.049\\
		& MSE&0.005 &0.005 &0.004 &0.002 &0.002 &0.002\\
		additive& bias&0.011 &0.011 &0.004 &0.008 &0.009 &0.010\\
		& sd    &0.100 &0.101 &0.086 &0.077 &0.080 &0.079\\
		& MSE&0.010 &0.010 &0.007 &0.006 &0.006 &0.006\\
		&&\multicolumn{6}{c}{20\% censoring}\\		
		Prop.	& bias&0.001 &0.001 &0.001 &0.001 &0.002 &0.002\\
		& sd &0.032 &0.040 &0.035 &0.031 &0.040 &0.032\\
		& MSE&0.002 &0.002 &0.001 &0.001 &0.002 &0.001\\
		PM1	& bias&0.008 &0.014 &0.003 &0.005 &0.000  &0.001\\
		& sd&0.065 &0.066 &0.054 &0.047 &0.046 &0.047\\
		& MSE&0.004 &0.004 &0.003 &0.002 &0.002 &0.002\\
		PM2	& bias&0.034 &0.274 &0.003 &0.187 &0.064 &0.186\\
		& sd&0.108 &0.192 &0.087 &0.133 &0.084 &0.134\\
		& MSE&0.012 &0.112 &0.007 &0.053 &0.011 &0.052\\
		additive& bias&0.001 &0.014 &0.005 &0.006 &0.005 &0.007\\
		& sd &0.127 &0.136 &0.115 &0.089 &0.099 &0.096\\
		& MSE&0.016 &0.018 &0.013 &0.008 &0.009 &0.009\\
		&&\multicolumn{6}{c}{40\% censoring}\\
		Prop.	& bias&0.002 &0.001 &0.001 &0.001 &0.001 &0.002\\
		& sd &0.043 &0.054 &0.044 &0.038 &0.039 &0.044\\
		& MSE&0.002 &0.002 &0.002 &0.001 &0.001 &0.001\\
		PM1	& bias&0.012 &0.017 &0.005 &0.006 &0.001 &0.002\\
		& sd &0.077 &0.078 &0.069 &0.056 &0.056 &0.059\\
		& MSE&0.006 &0.006 &0.005 &0.003 &0.003 &0.003\\
		PM2	& bias&0.078 &0.625 &0.013 &0.423 &0.156 &0.428\\
		& sd&0.160 &0.307 &0.140 &0.214 &0.134 &0.218\\
		& MSE&0.030 &0.486 &0.019 &0.226 &0.042 &0.231\\
		additive& bias&0.021 &0.032 &0.006 &0.001 &0.004 &0.003\\
		& sd    &0.158 &0.152 &0.132 &0.107 &0.109 &0.103\\
		& MSE&0.023 &0.023 &0.017 &0.011 &0.012 &0.011\\
		\hline
\end{longtable}

\newpage
\begin{longtable}{cc|cccccccc}
	\caption{Results of  study 3, based on 1000
		simulations with sample size 500. ``mean'' is the average
		absolute bias of $(\wh\bb)$ of each component
		in $\bb$, ``sd'' and ``MSE" are the
		sample standard deviation and mean squared error of the
		corresponding estimations.
	}
	\label{tab:simu3}\\\hline
	\vspace{.2cm}
		&& $\beta_{31}$ &$\beta_{41}$ &$\beta_{51}$ &$\beta_{61}$
		&$\beta_{32}$      &$\beta_{42}$  &$\beta_{52}$
		&$\beta_{62}$\\
		&true& $2.75$ &$-0.75$ &$-1.0$ &$2.0$
		&$-3.125$
		&$-1.125$
		&$1.0$ &$-2.0$\\\hline
		&&\multicolumn{8}{c}{No censoring}\\
		Prop.	& mean&0.067 &0.047 &0.025 &0.056 &0.083 &0.052
		&0.023 &0.051\\
		& sd  &0.441 &0.300   &0.304 &0.383 &0.553 &0.363 &0.368
		&0.423\\
		&MSE&0.199 &0.092 &0.093 &0.149 &0.312 &0.134 &0.136
		&0.182\\
		\hline
		& & \multicolumn{8}{c}{20\% censoring}\\		
		Prop.	& mean&0.085 &0.036 &0.032 &0.056 &0.063 &0.044
		&0.033 &0.052\\
		& sd  &0.511 &0.404 &0.453 &0.490  &0.719 &0.472 &0.457
		&0.503\\
		&MSE&0.268 &0.164 &0.206 &0.243 &0.520  &0.224 &0.210
		&0.255\\
		\hline
		&&\multicolumn{8}{c}{40\% censoring}\\		
		Prop.	& mean&0.083 &0.034 &0.008 &0.075 &0.109 &0.050
		&0.019 &0.040 \\
		& sd  &0.583 &0.469 &0.474 &0.519 &0.676 &0.525 &0.540
		&0.565\\
		&MSE&0.347 &0.221 &0.225 &0.275 &0.468 &0.278 &0.292
		&0.321\\
		\hline
\end{longtable}

\newpage
 
\begin{figure}[H]
	\centering
	\includegraphics[width=4cm]{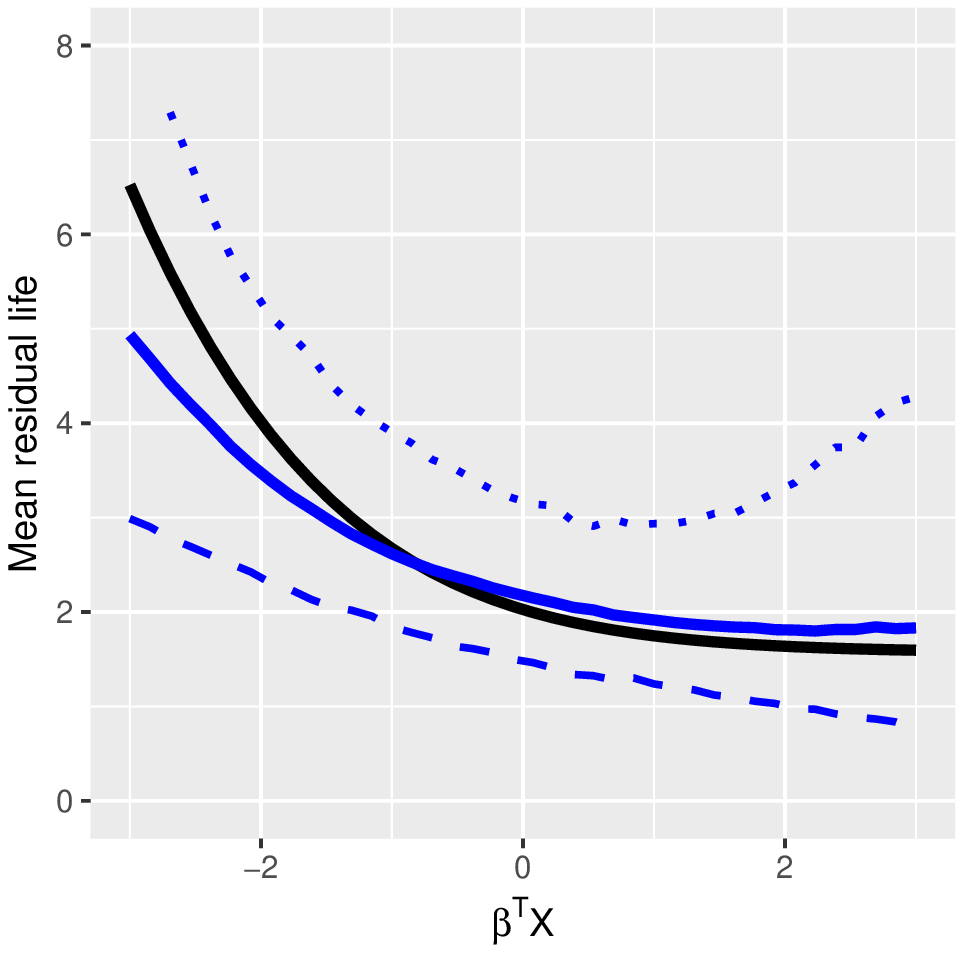}
	\includegraphics[width=4cm]{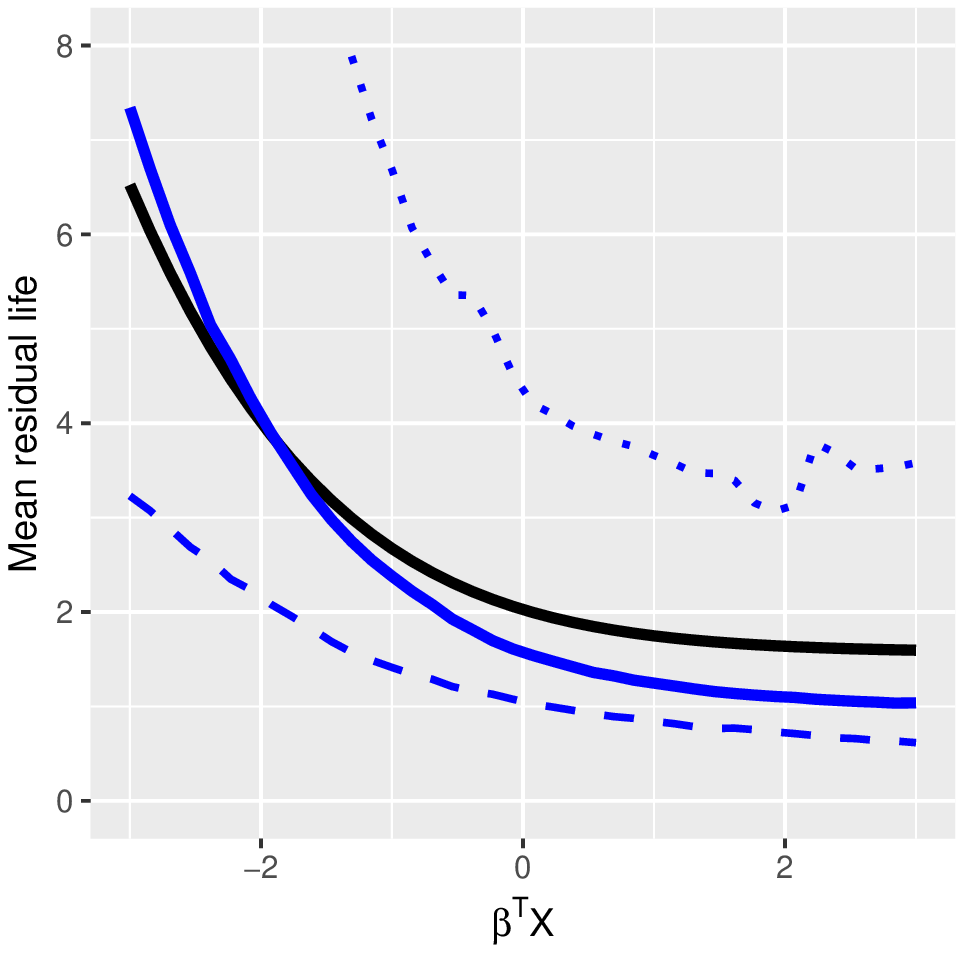}
	\includegraphics[width=4cm]{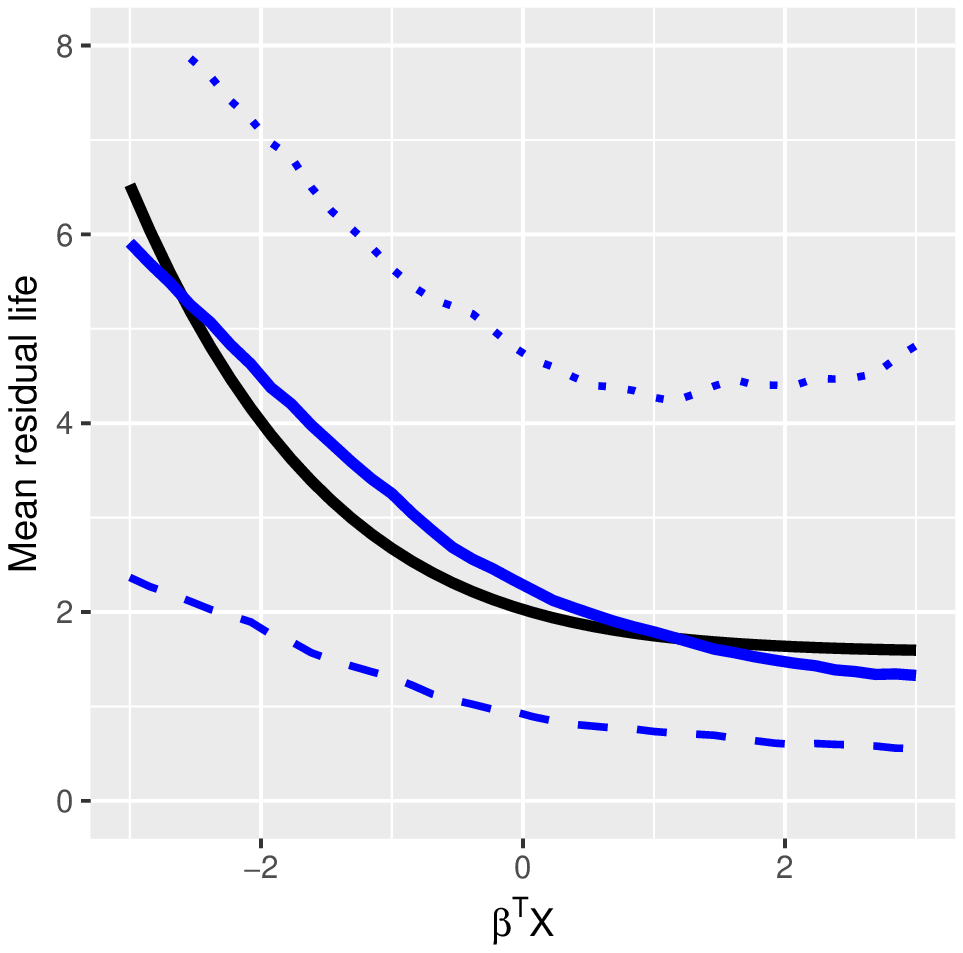}\\
	\includegraphics[width=4cm]{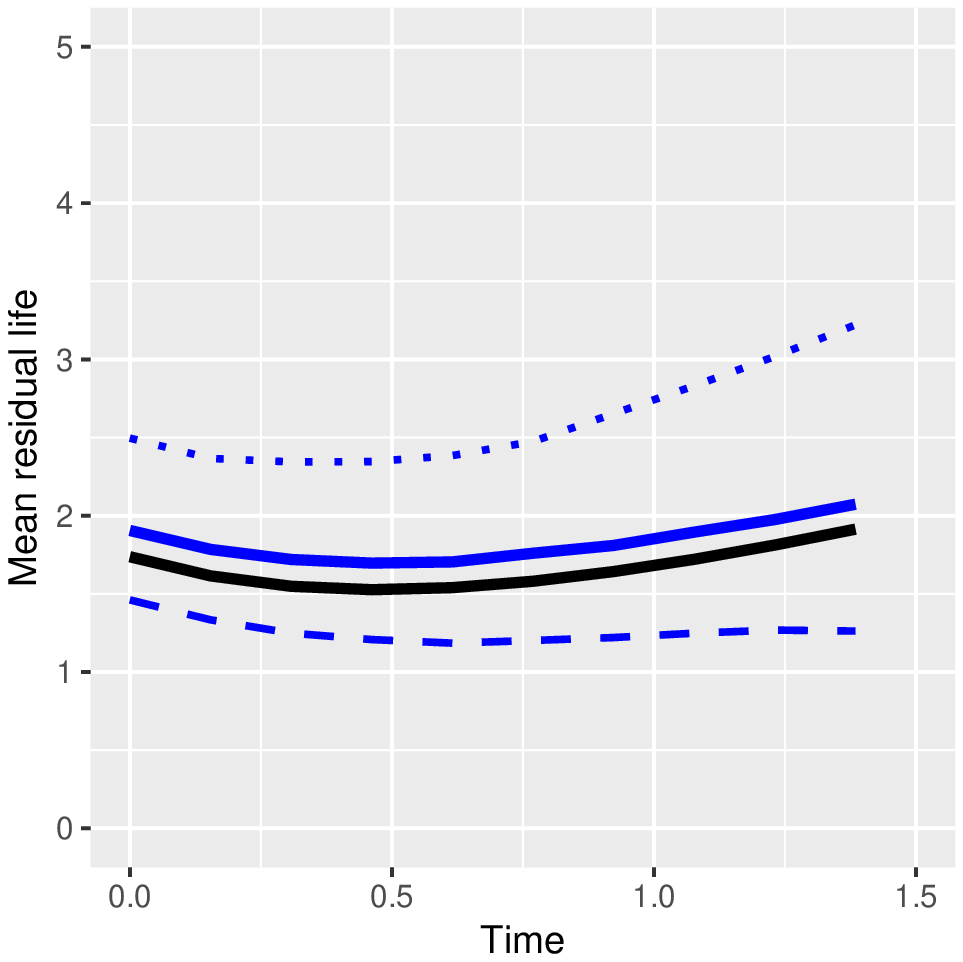}
	\includegraphics[width=4cm]{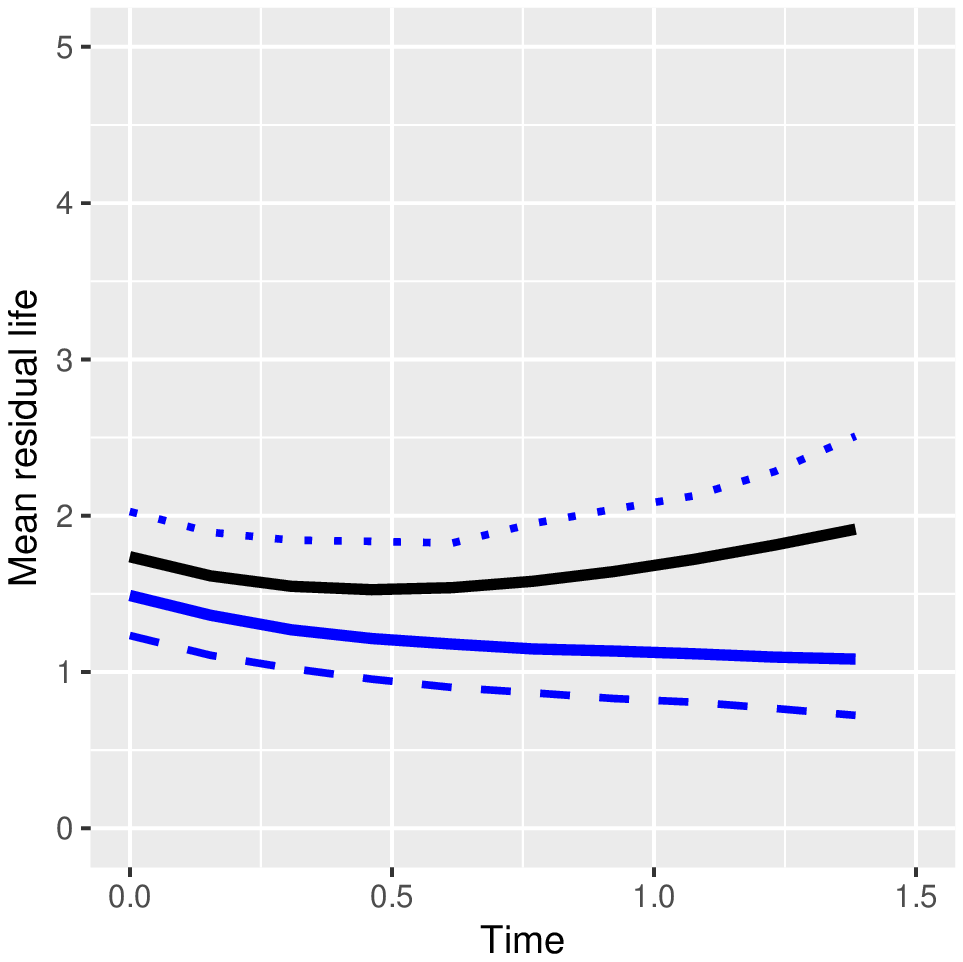}
	\includegraphics[width=4cm]{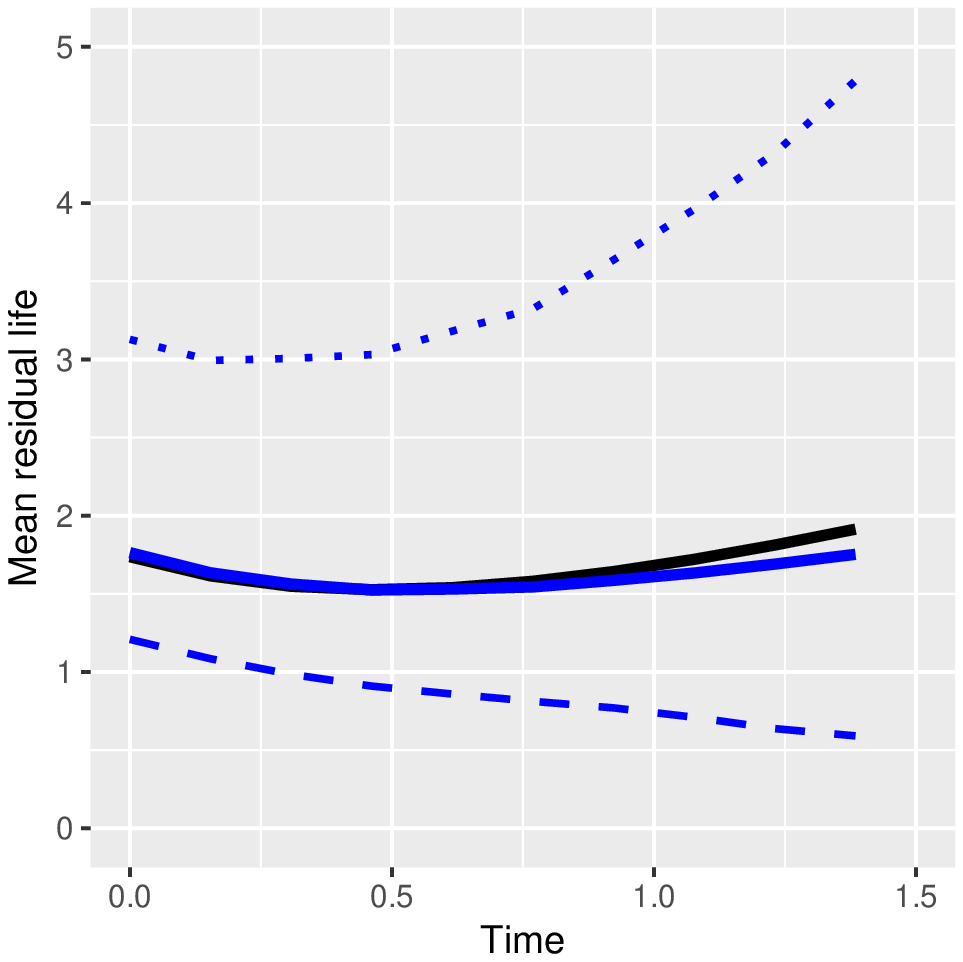}\\
	\includegraphics[width=4cm]{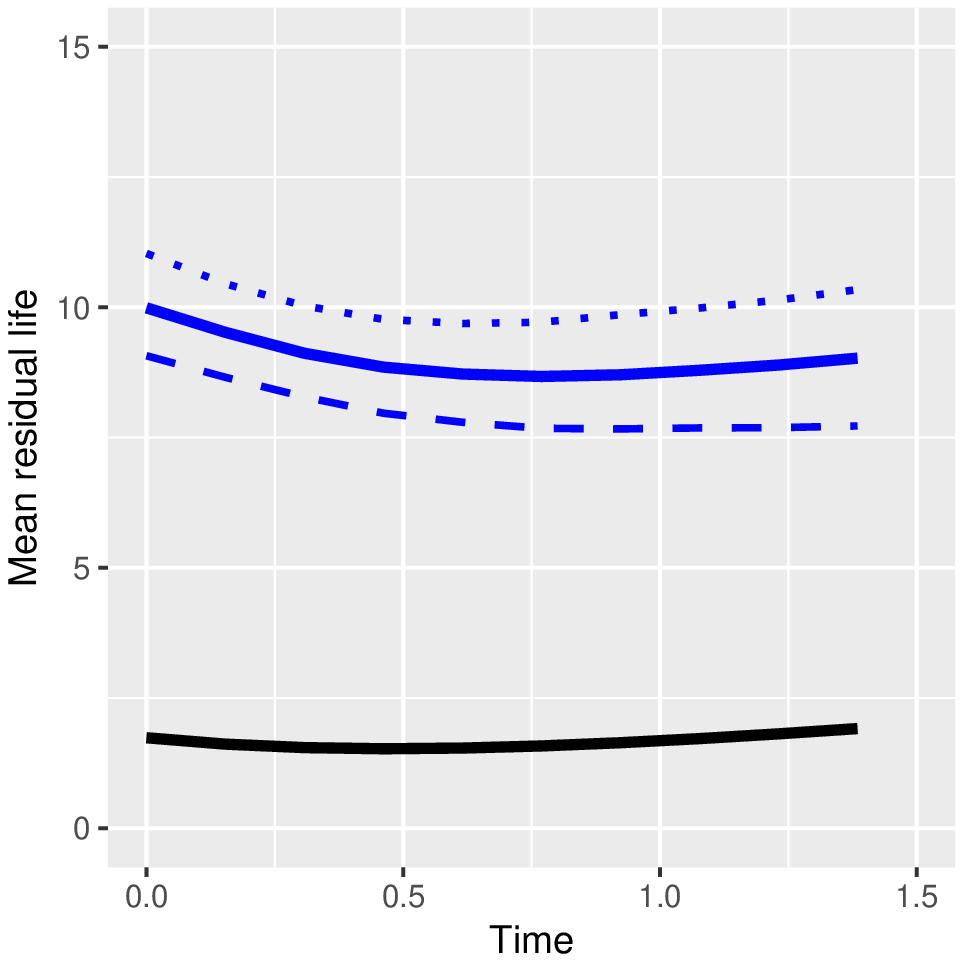}
	\includegraphics[width=4cm]{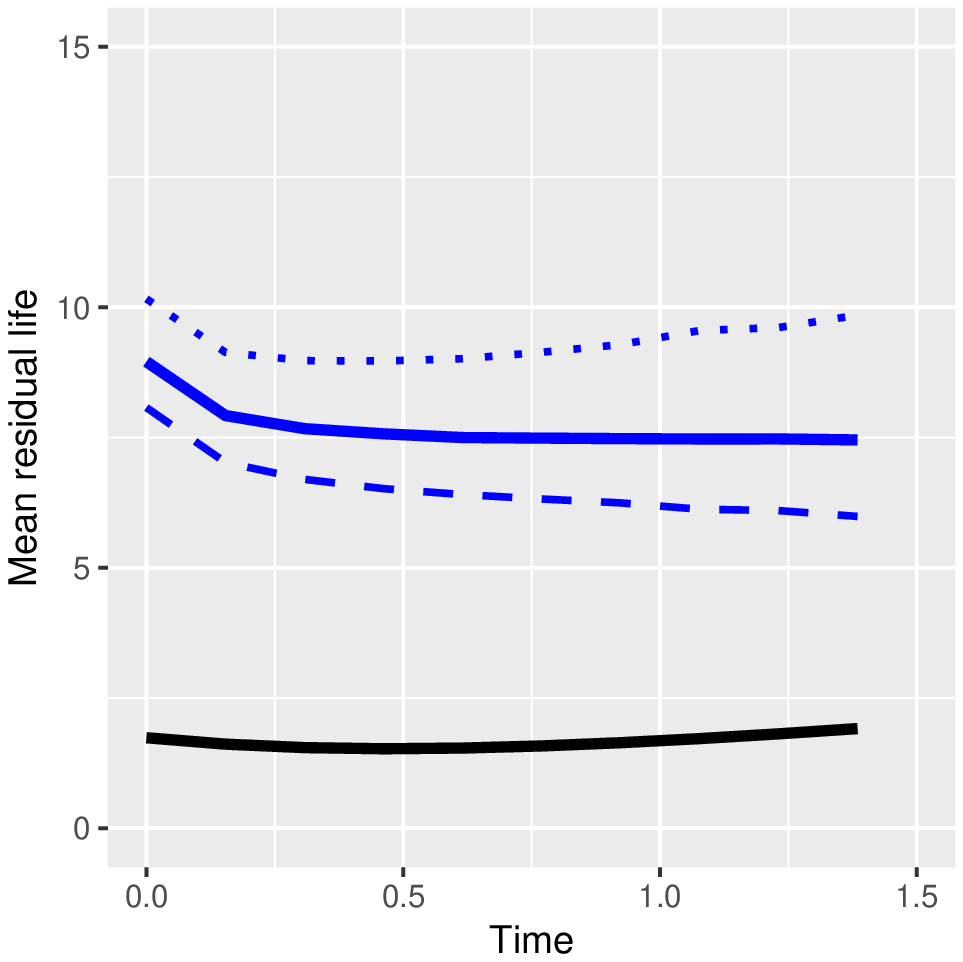}
	\includegraphics[width=4cm]{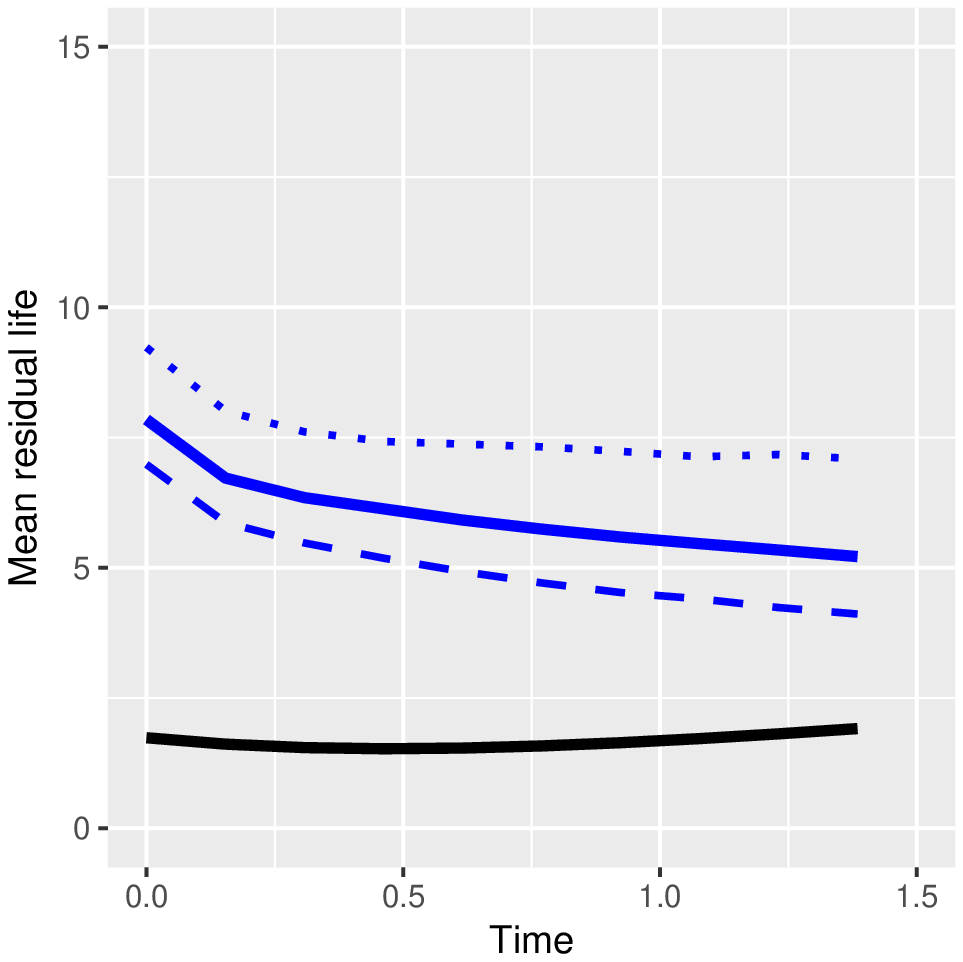}\\
	\includegraphics[width=4cm]{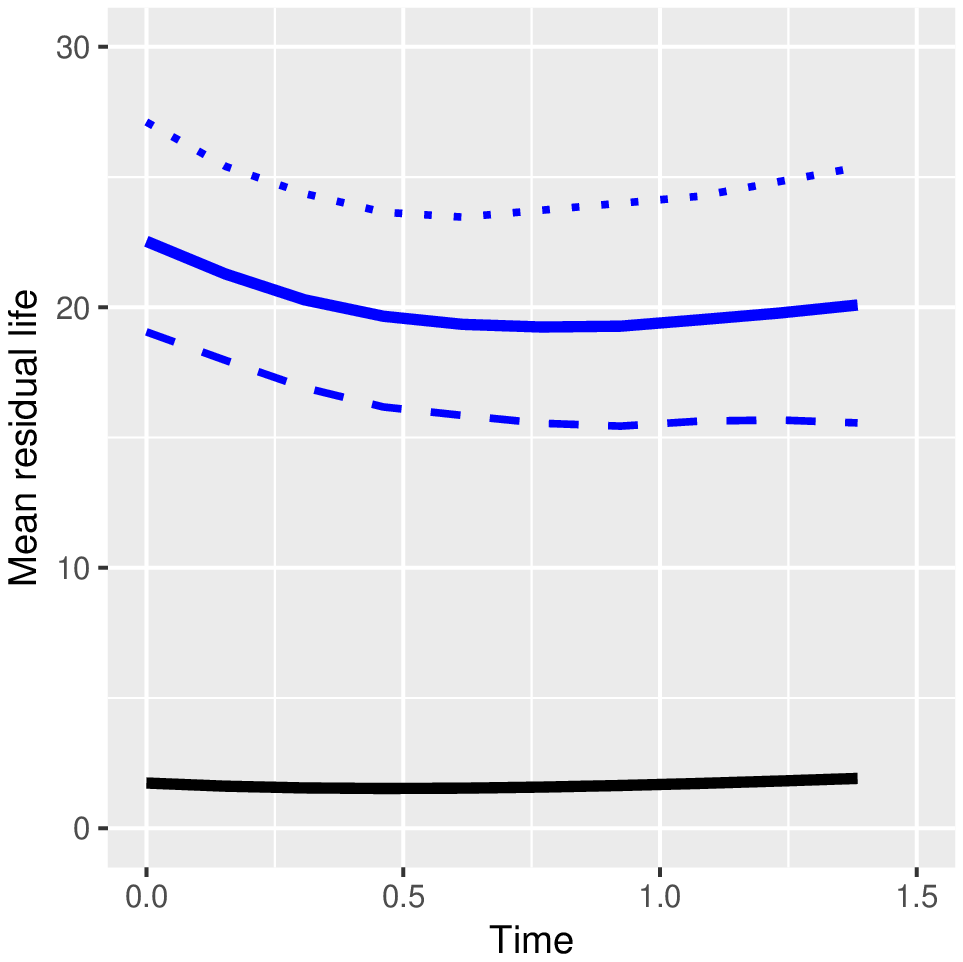}
	\includegraphics[width=4cm]{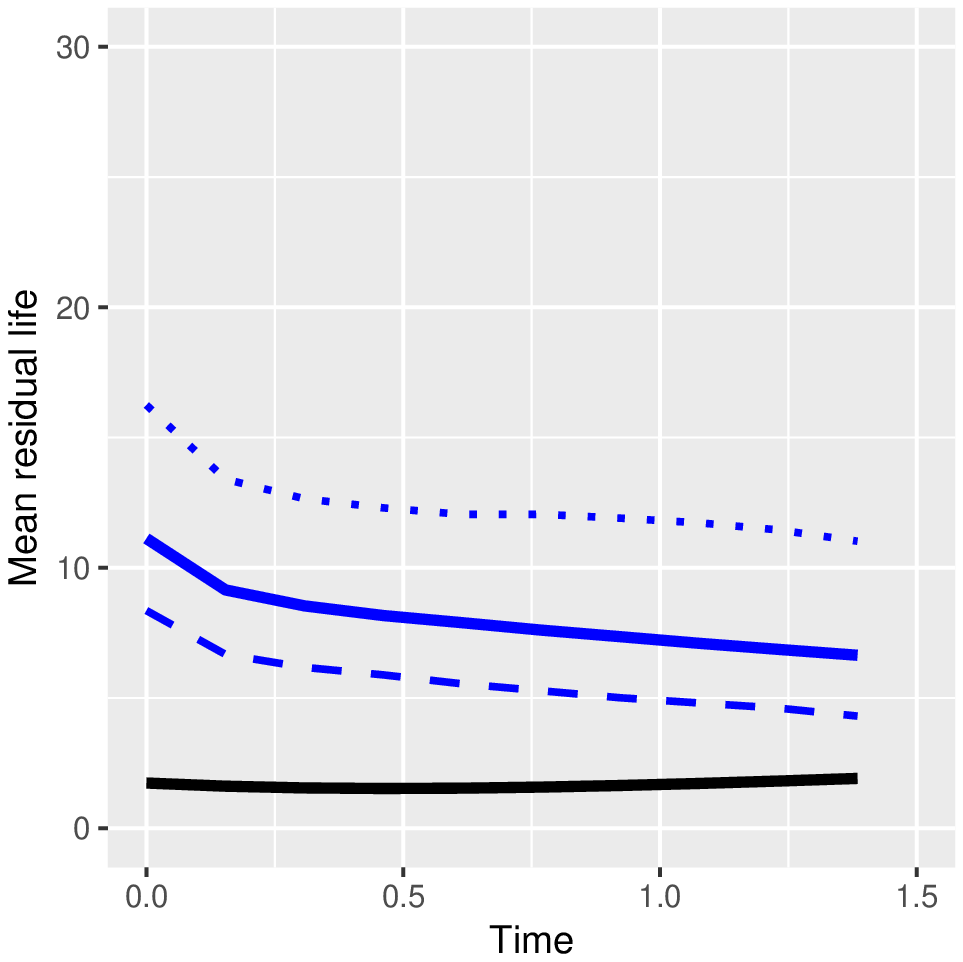}
	\includegraphics[width=4cm]{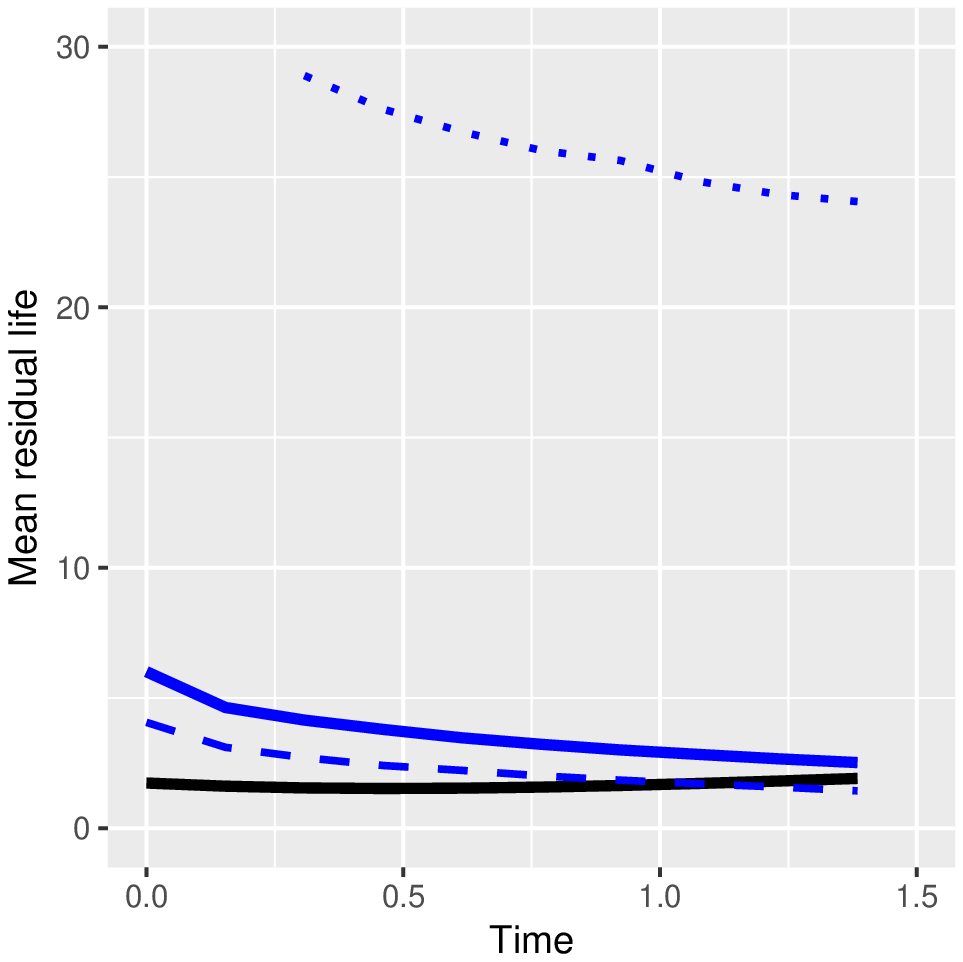}\\
	\includegraphics[width=4cm]{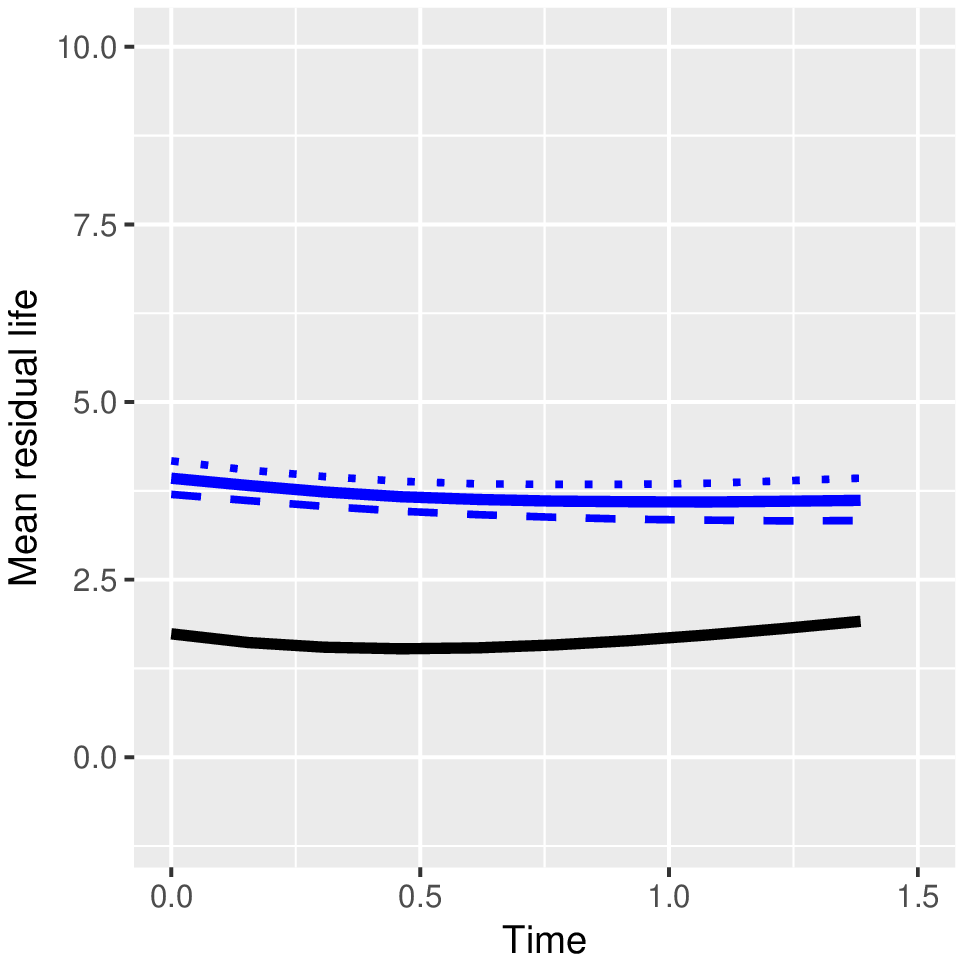}
	\includegraphics[width=4cm]{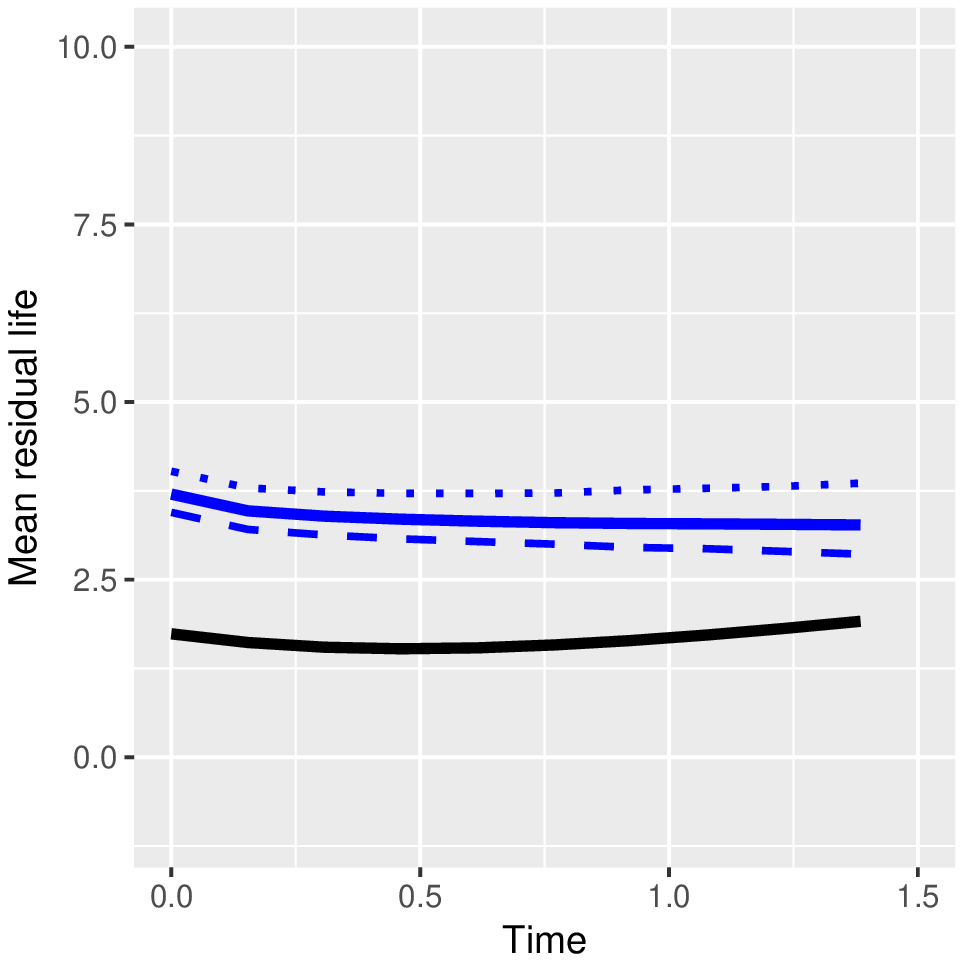}
	\includegraphics[width=4cm]{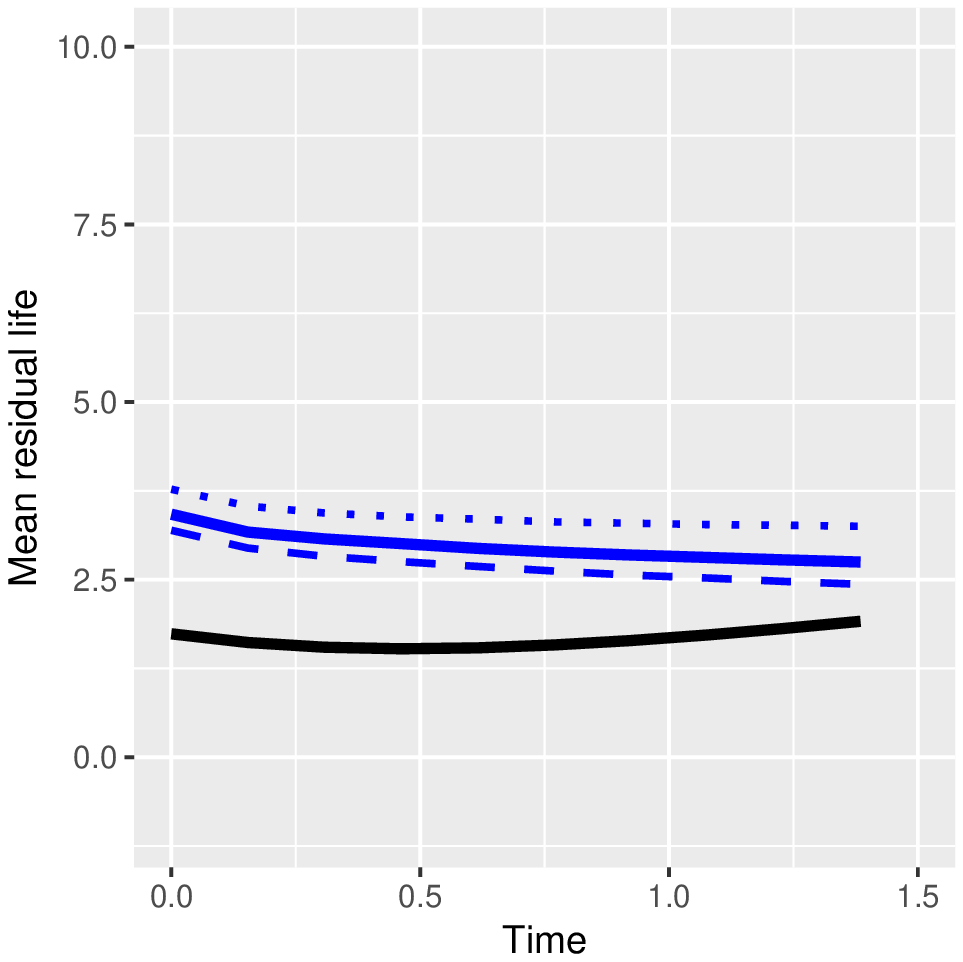}\\
	\caption{Mean residual life function estimation in Study
		1. Row 1: $m(t,\bb\trans\x)$ as a
		function of $\bb\trans\x$
		at $t=1$. Row 2 to Row 5: $m(t,\bb\trans\x)$ as a 
		function of $t$
		at $\bb\trans\x=1.5$ from method ``semiparametric", 
		``PM1", ``PM2", 
		``additive".
		Left to right
		columns: no censoring; 20\% censoring rate; 40\%
		censoring rate. Black
		line: True $m(t,\bb\trans\x)$; Blue line: Median of $\wh
		m(t,\bb\trans\x)$;
		Blue dashed line: 2.5\% empirical percentile curve;
		Blue dotted line: 97.5\% empirical percentile curve.
	}
	\label{fig:simu1curve}
\end{figure}

\newpage

\begin{figure}[H]
	\centering
	\includegraphics[width=5cm]{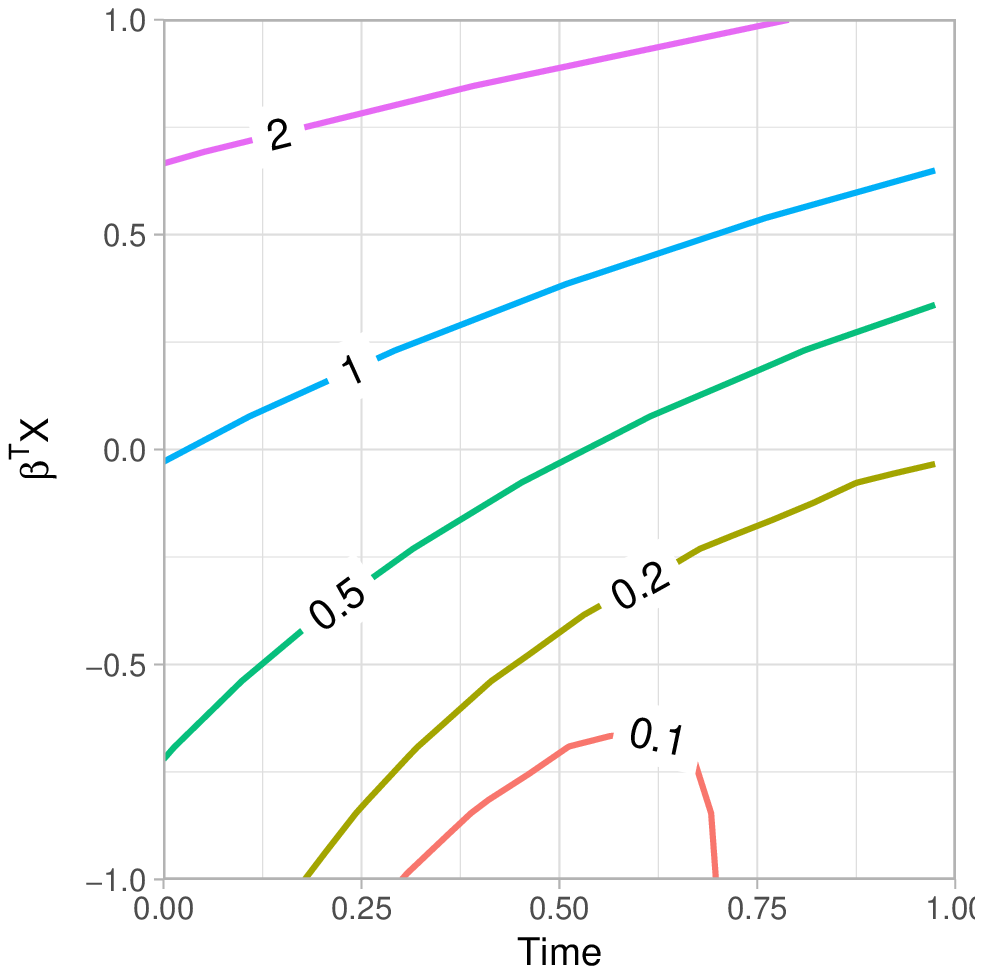}
	\includegraphics[width=5cm]{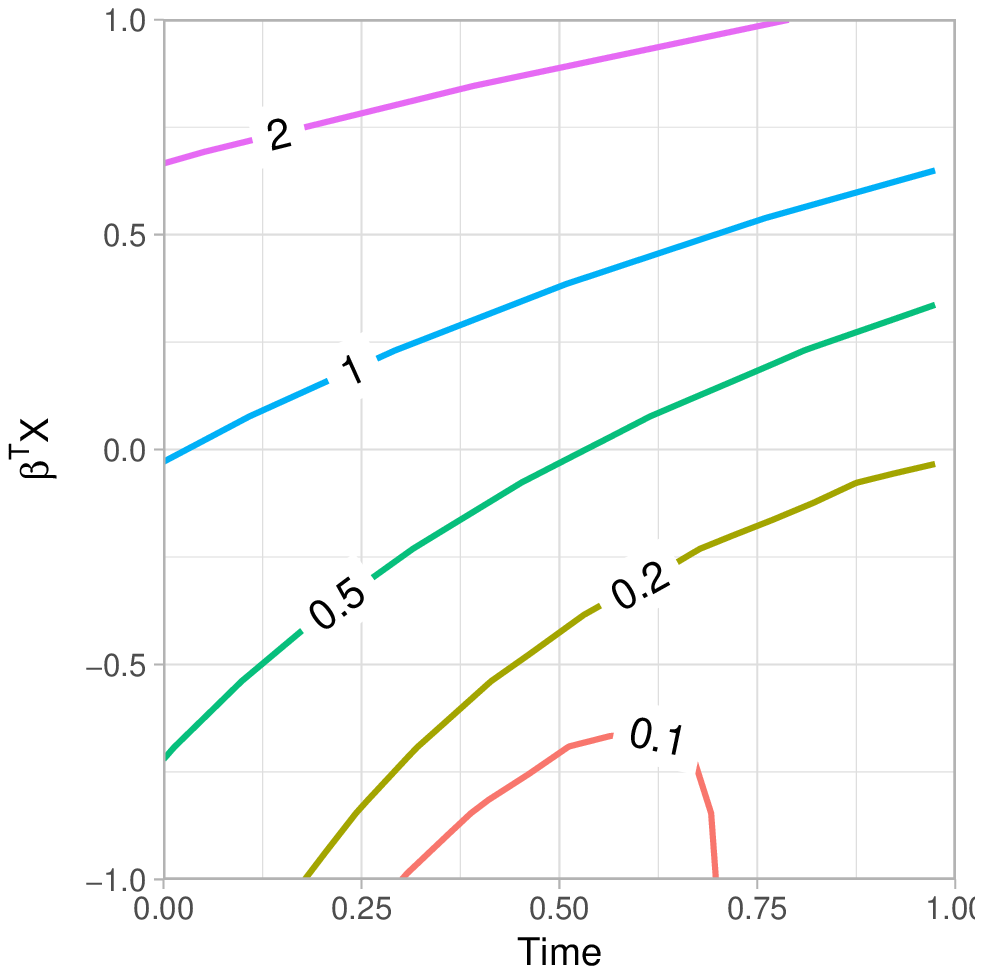}
	\includegraphics[width=5cm]{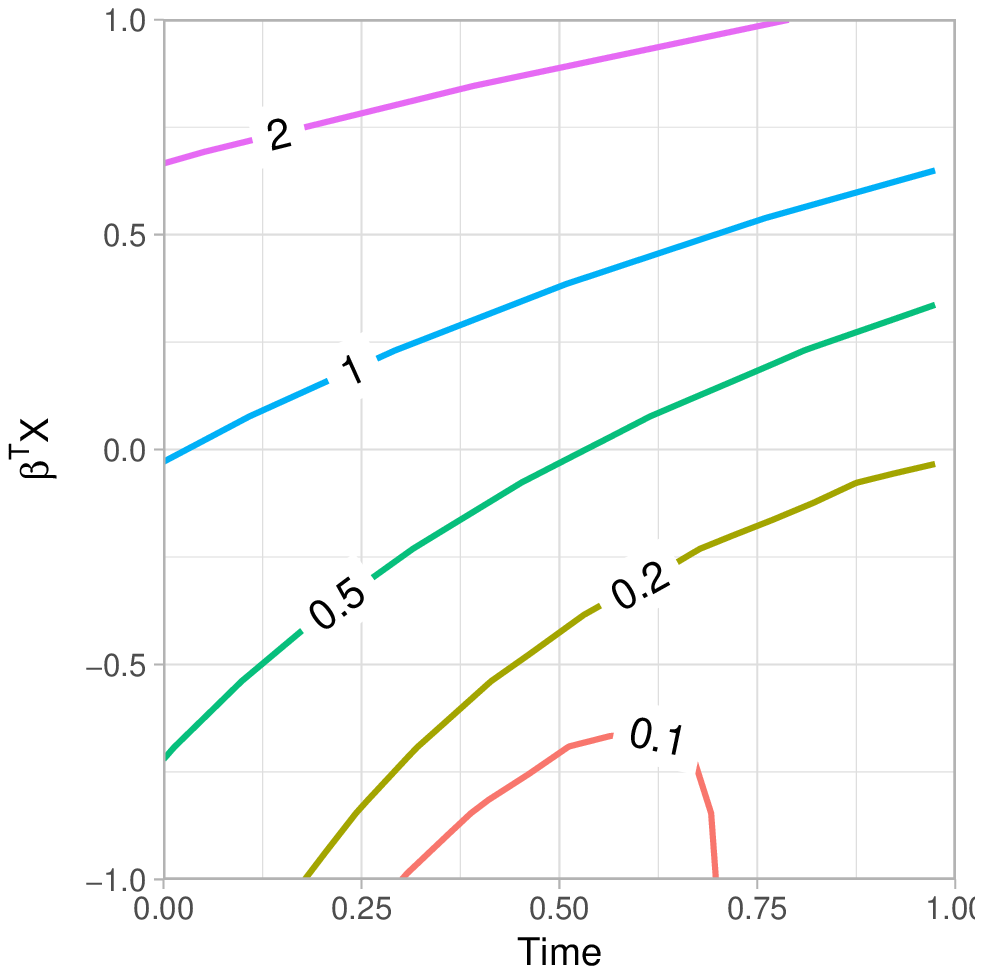}\\
	\includegraphics[width=5cm]{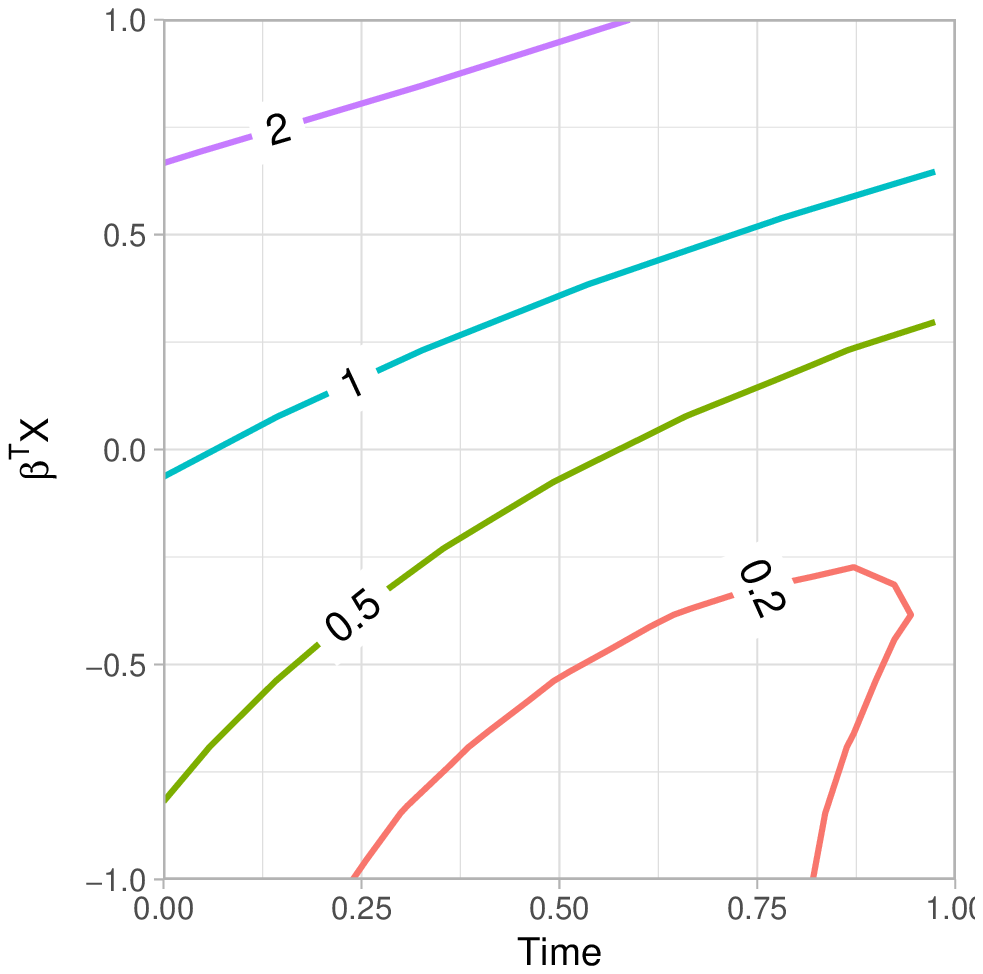}
	\includegraphics[width=5cm]{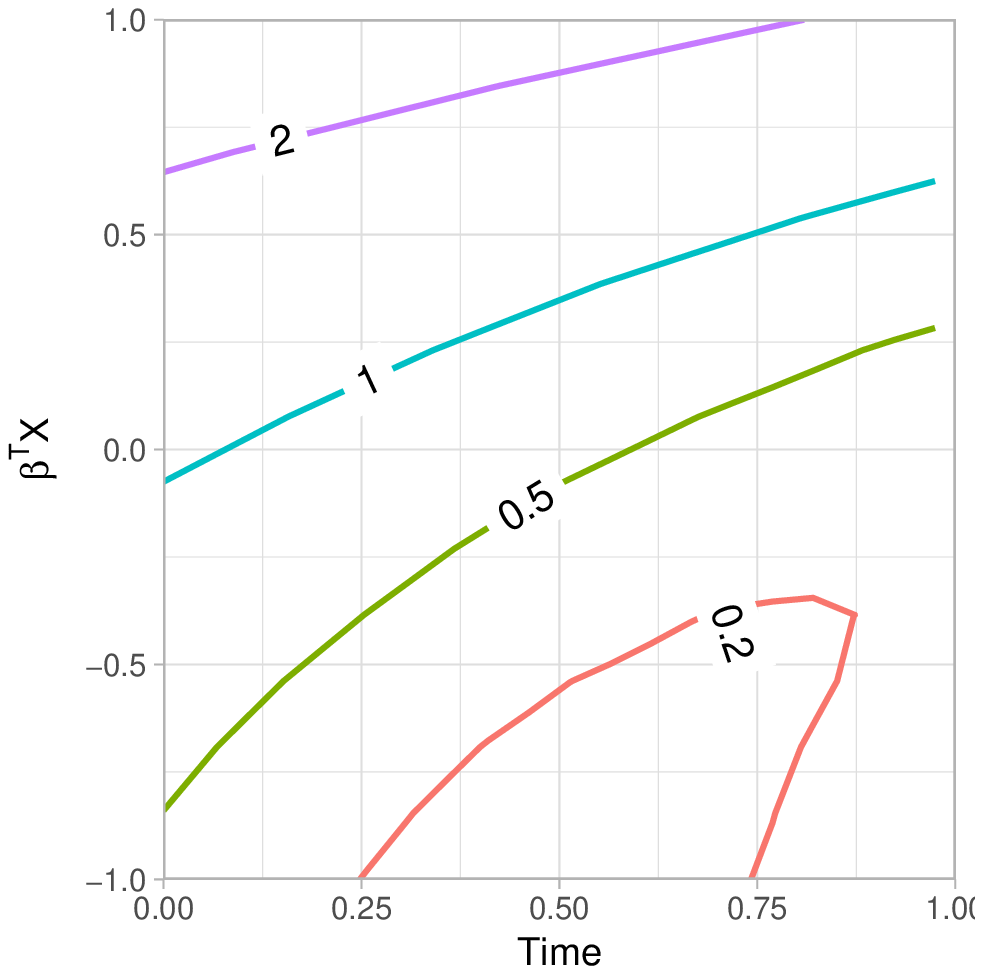}
	\includegraphics[width=5cm]{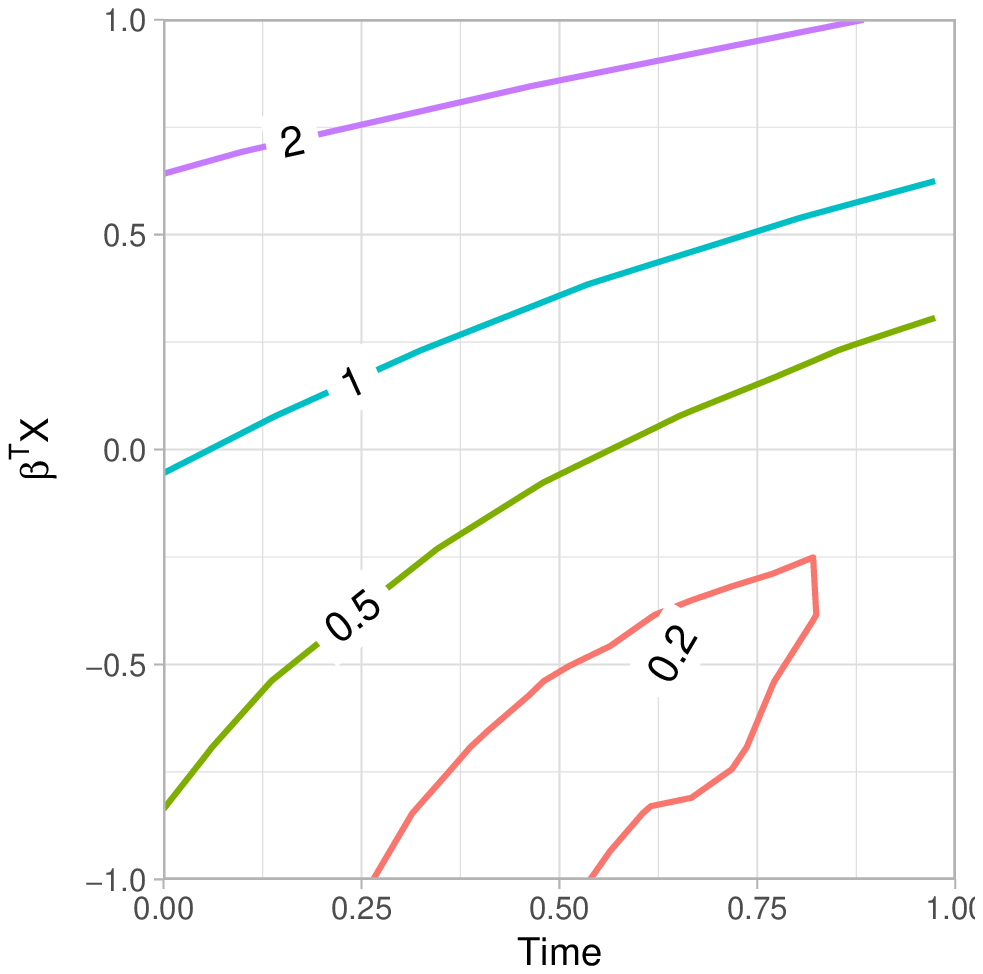}\\
	\includegraphics[width=5cm]{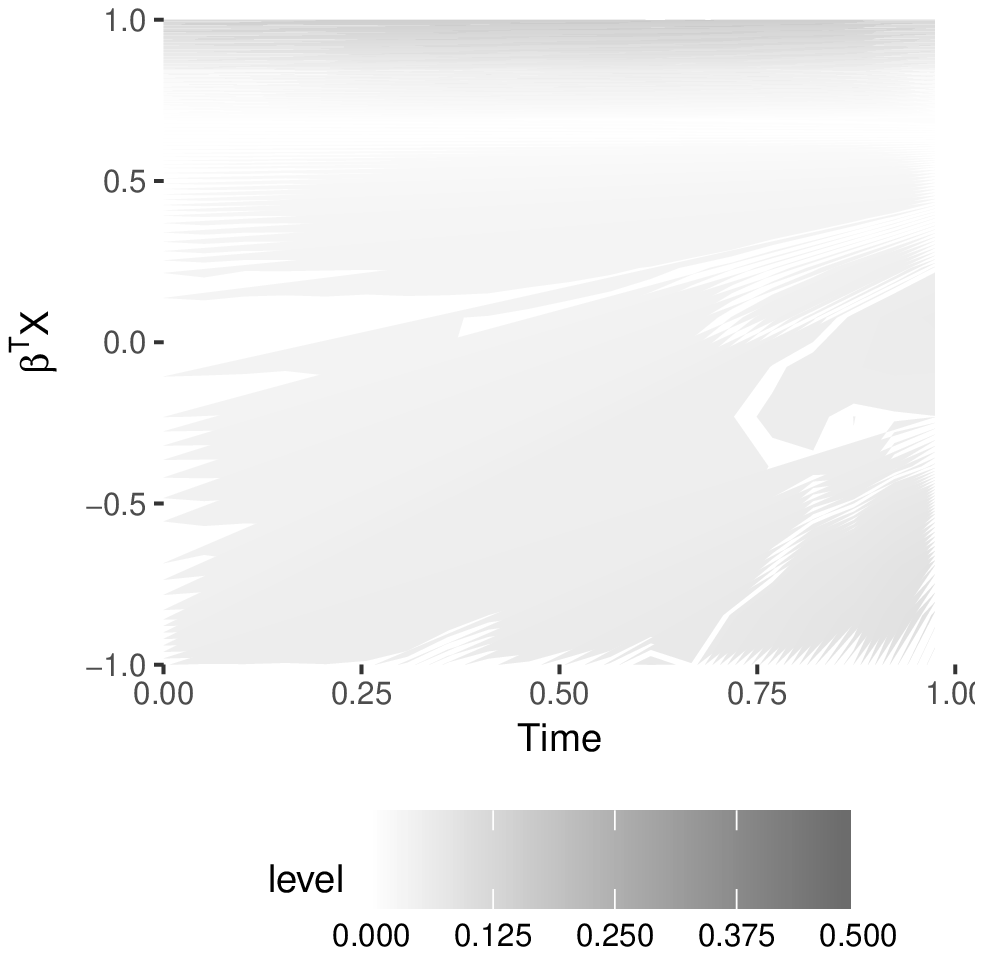}
	\includegraphics[width=5cm]{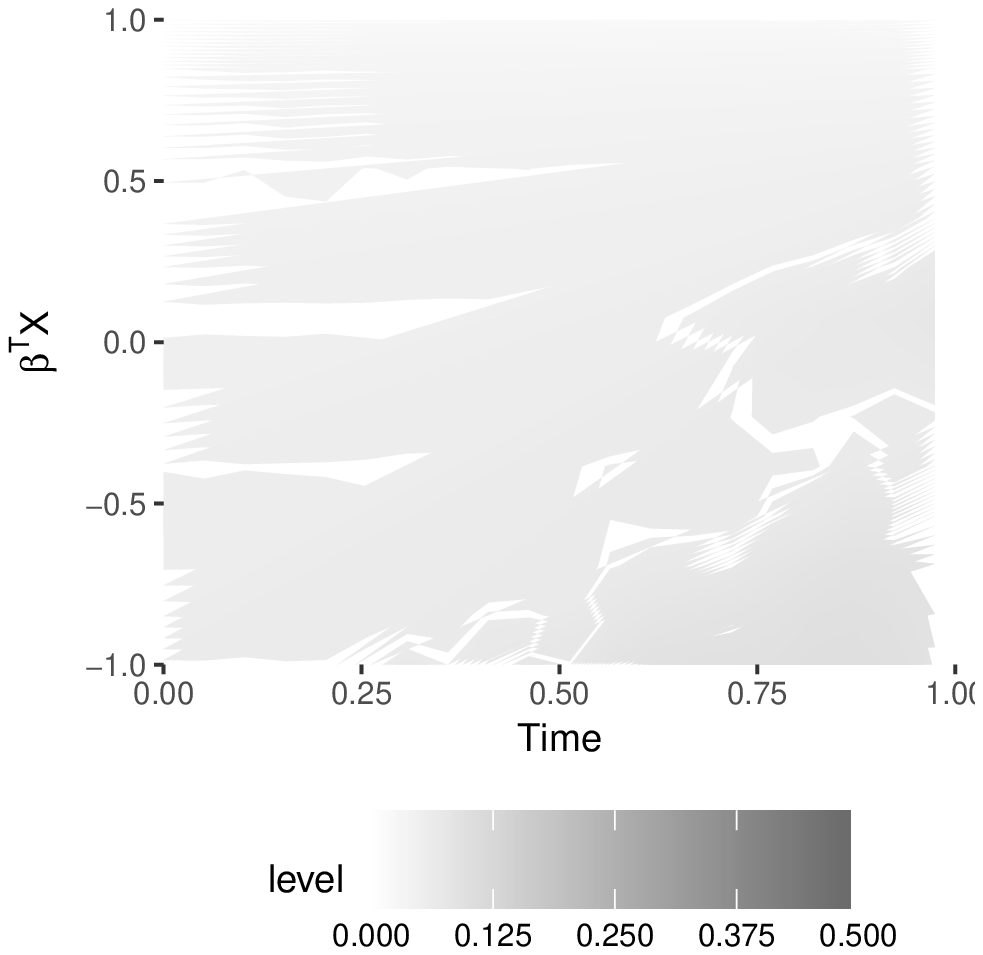}
	\includegraphics[width=5cm]{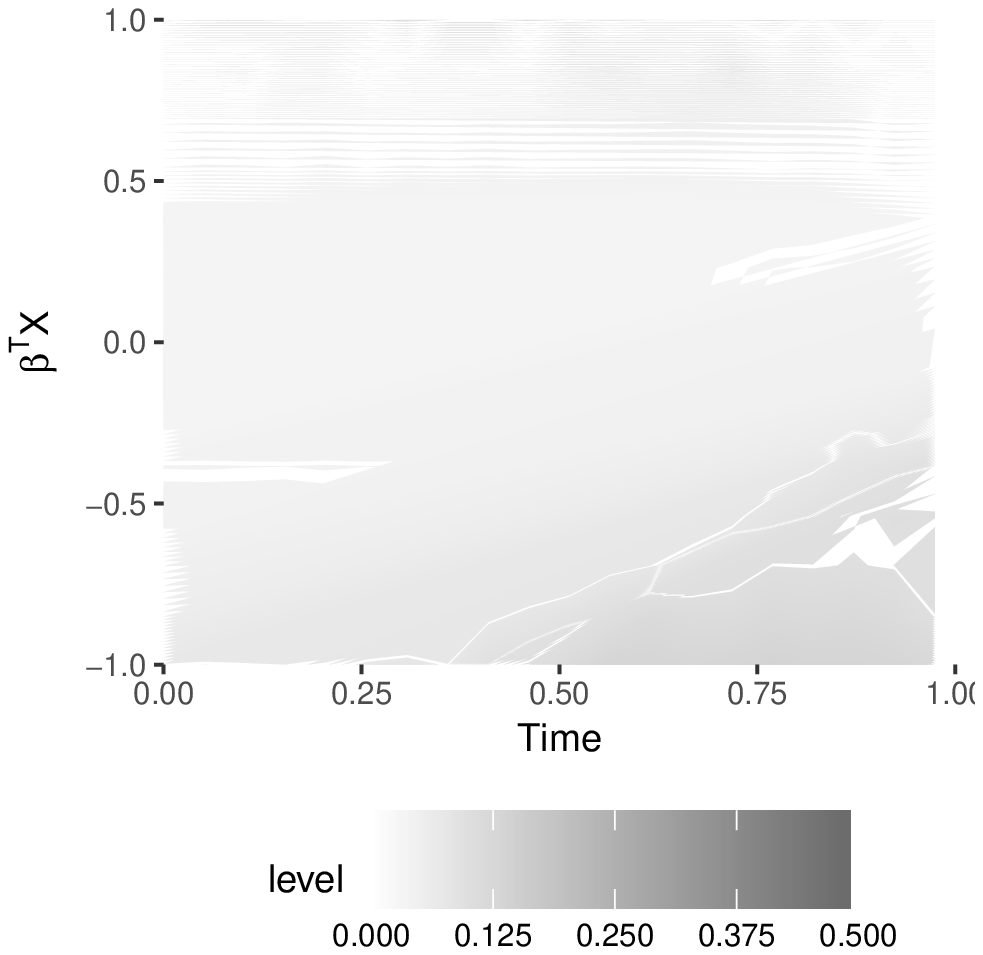}	
	\caption{Performance of the semiparametric method on mean
		residual life
		function of study 2.
		First row: contour plot of true $m(t,\bb\trans\X)$;
		Second
		row: contour plot of $\wh m(t,\bb\trans\X)$;
		Third row: contour plot of $|\wh
		m(t,\bb\trans\X)-m(t,\bb\trans\X)|$.
		Left to right
		columns: no censoring; 20\% censoring rate; 40\%
		censoring rate.
	}
	\label{fig:simu2contour}
\end{figure}
\begin{figure}[H]
	\centering
	\includegraphics[width=4cm]{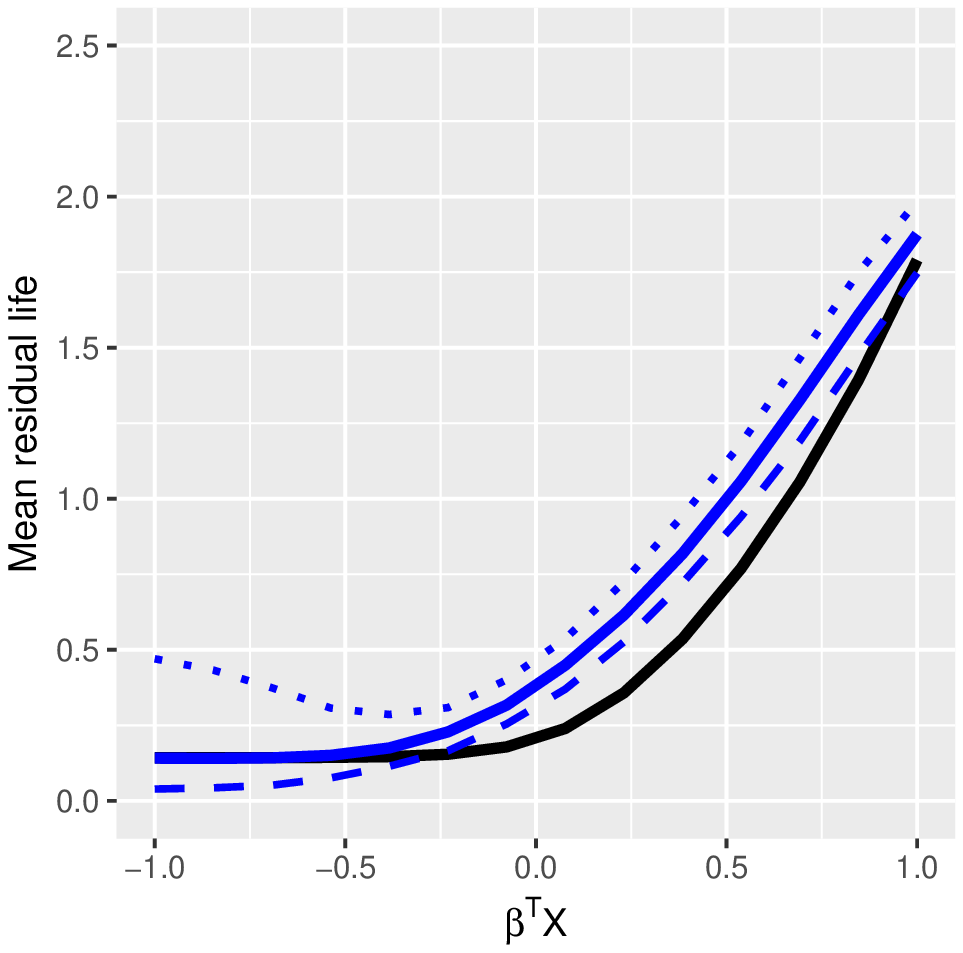}
	\includegraphics[width=4cm]{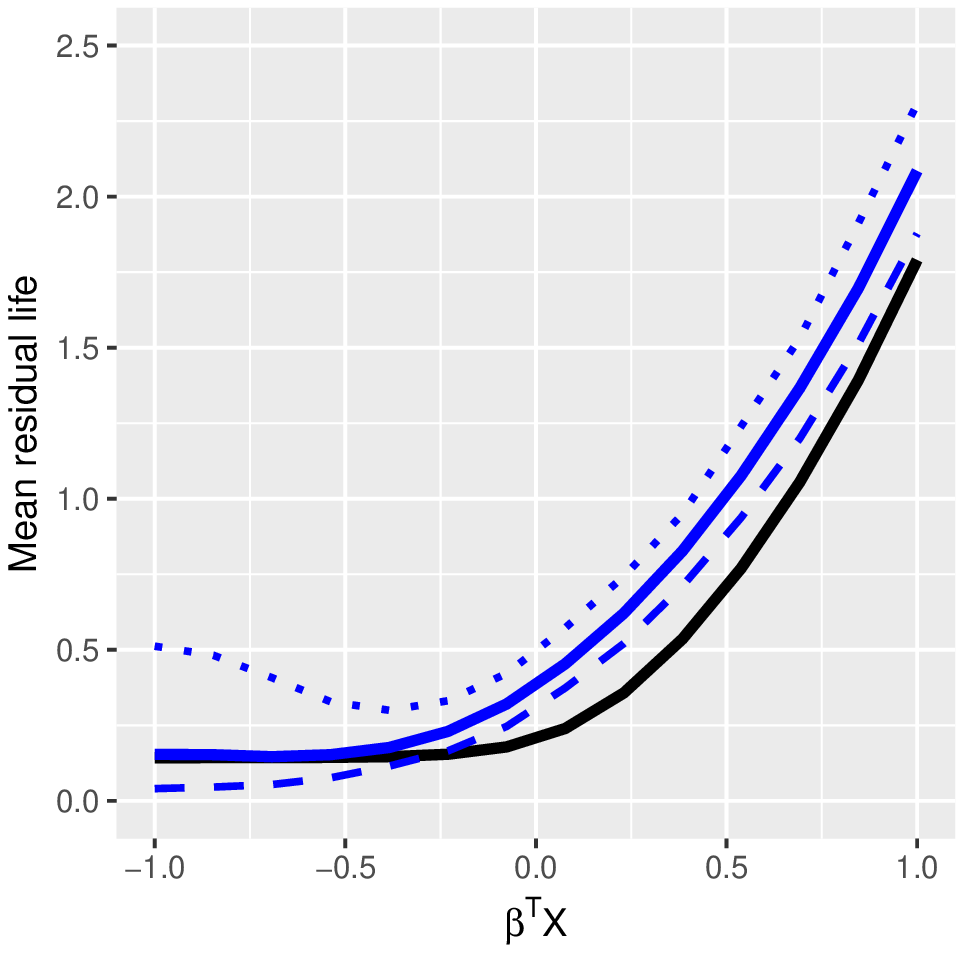}
	\includegraphics[width=4cm]{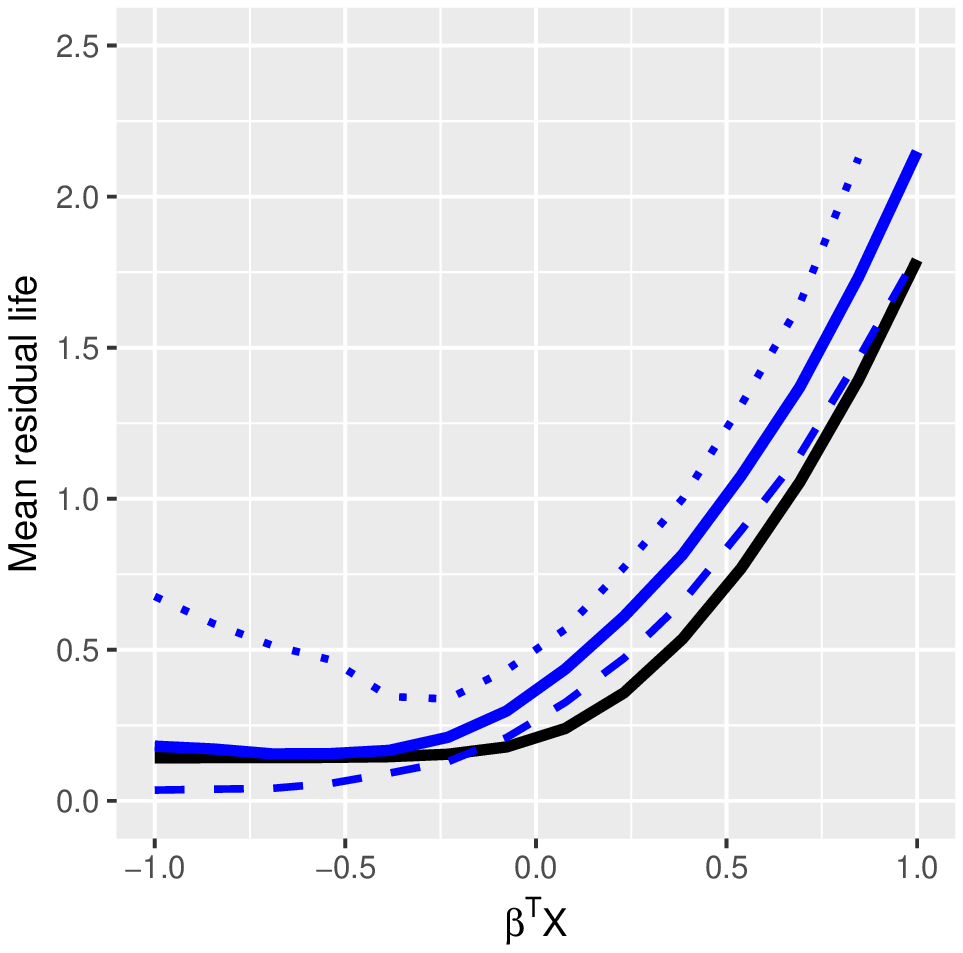}\\
	\includegraphics[width=4cm]{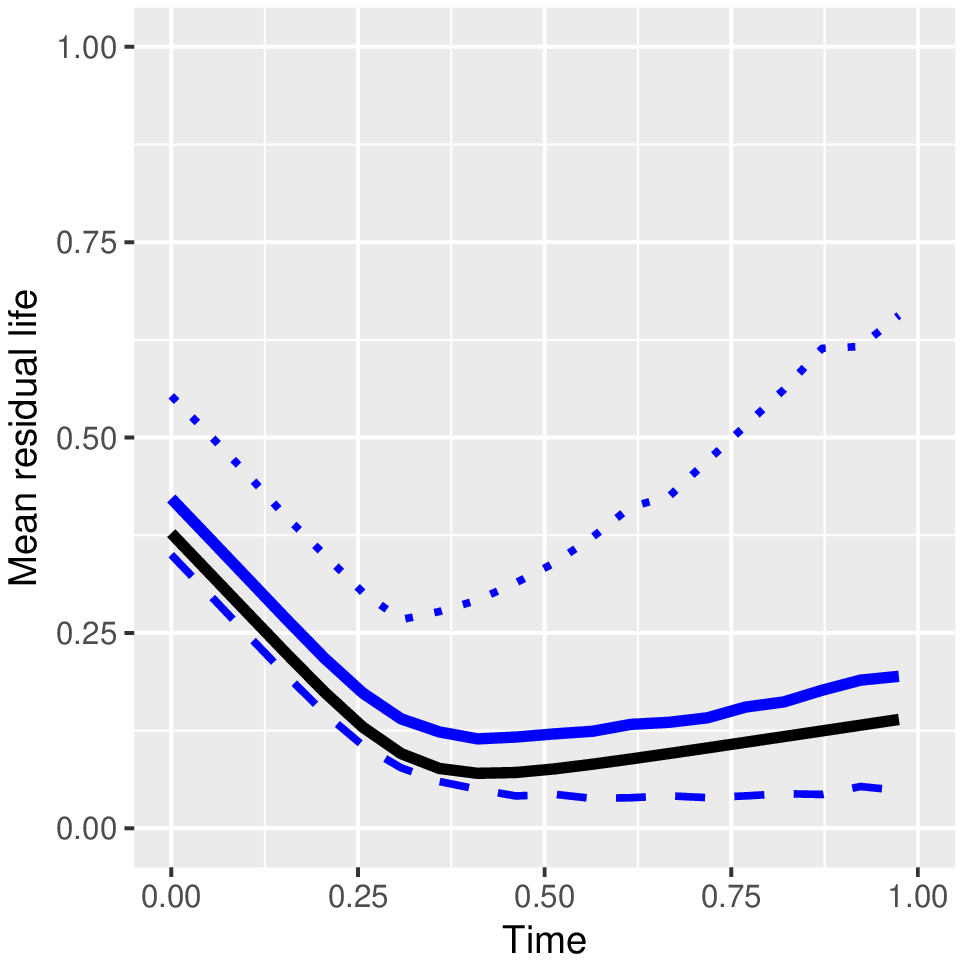}
	\includegraphics[width=4cm]{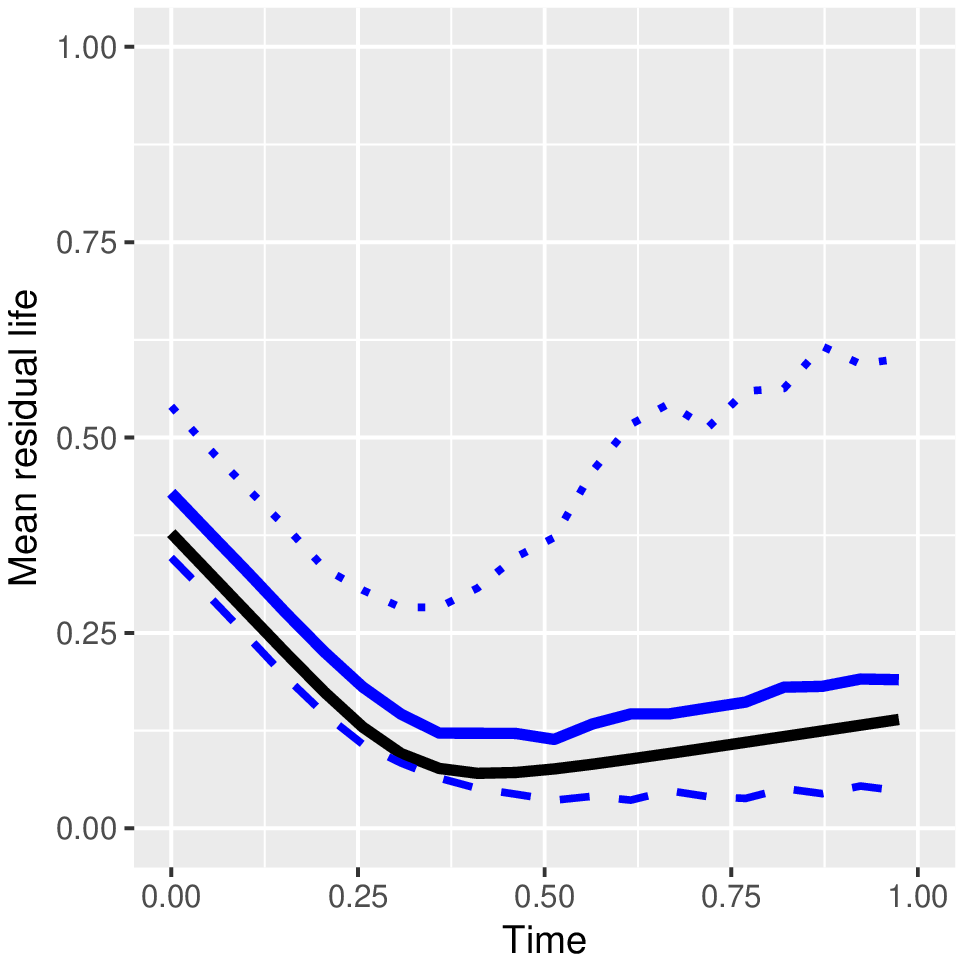}
	\includegraphics[width=4cm]{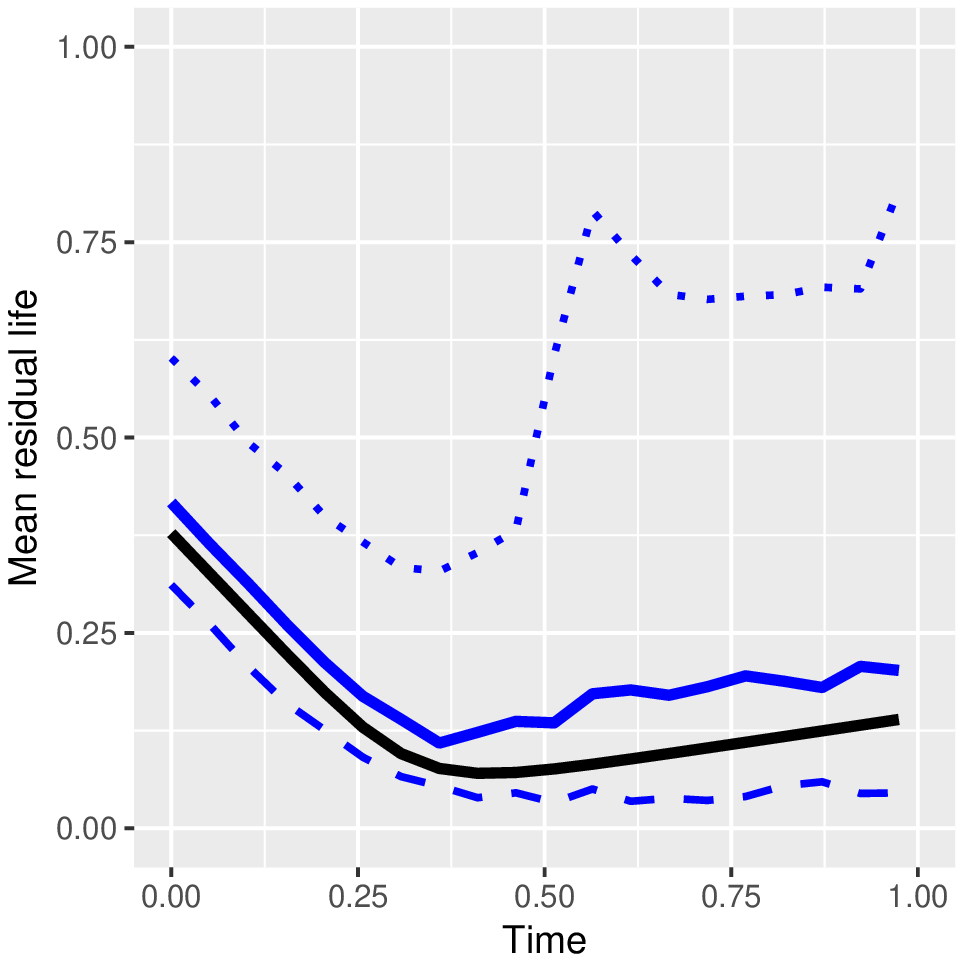}\\
	\includegraphics[width=4cm]{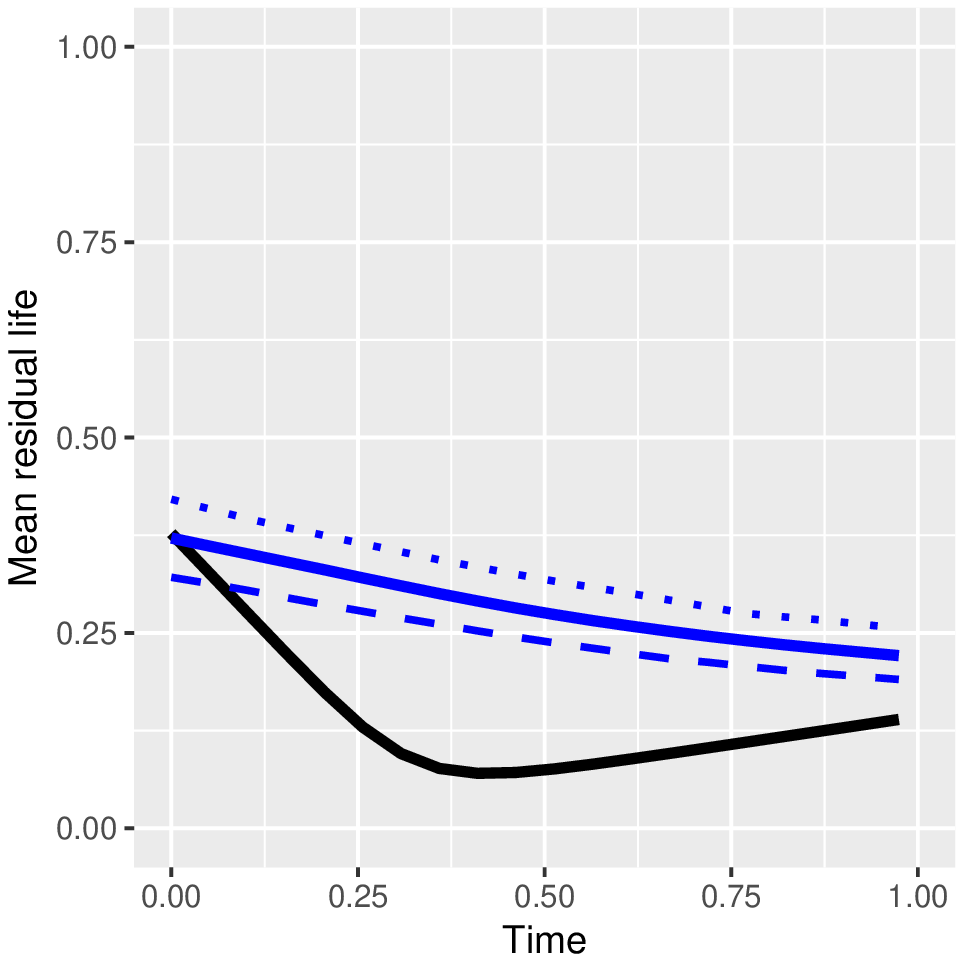}
	\includegraphics[width=4cm]{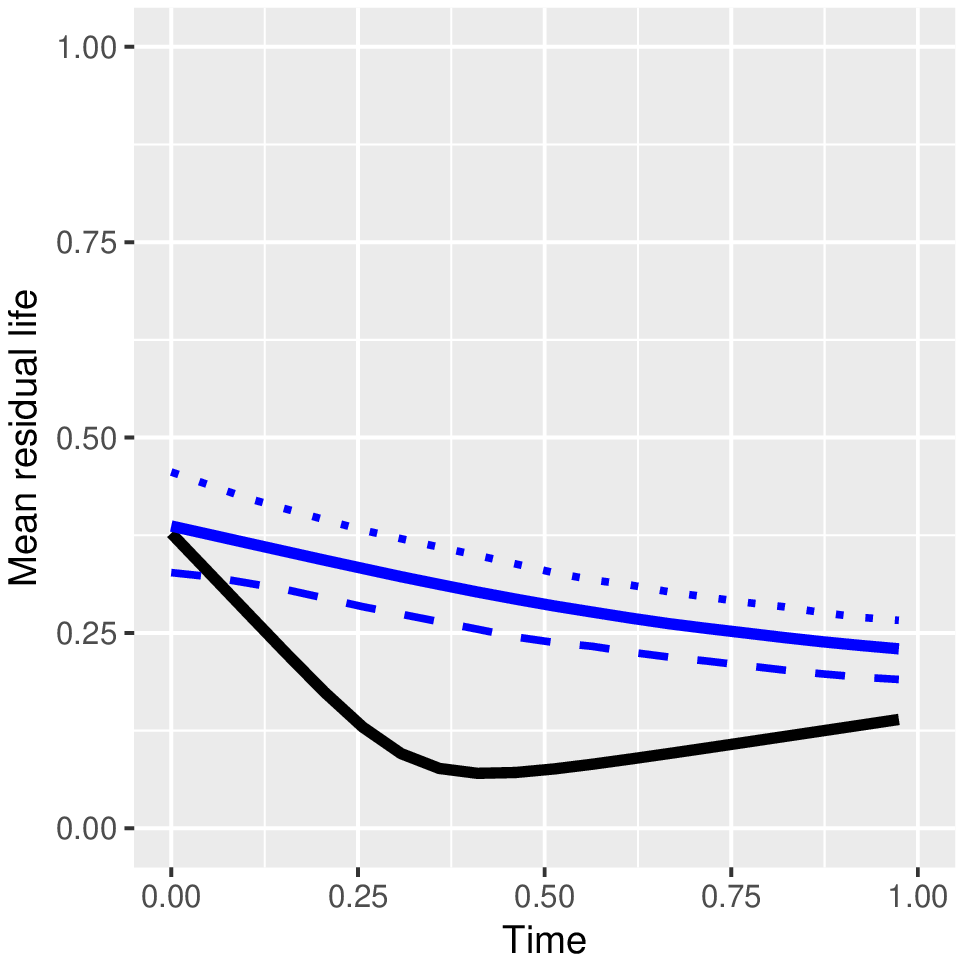}
	\includegraphics[width=4cm]{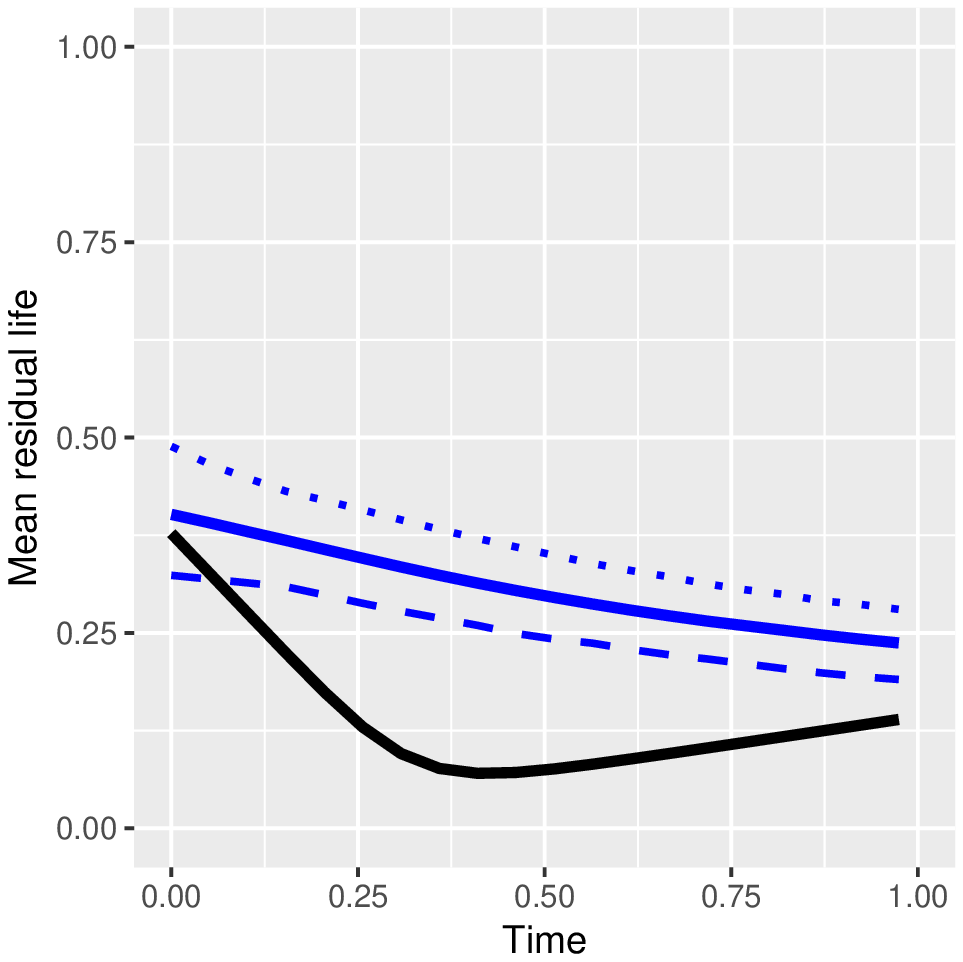}\\
	\includegraphics[width=4cm]{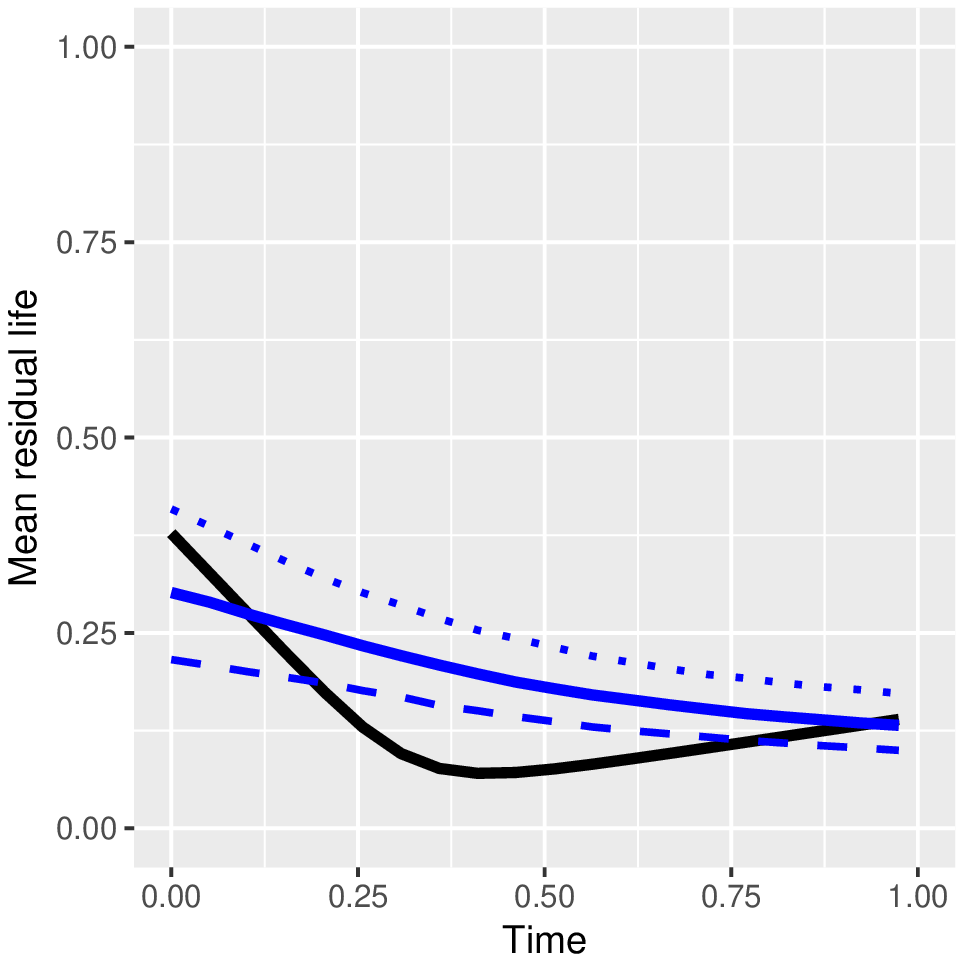}
	\includegraphics[width=4cm]{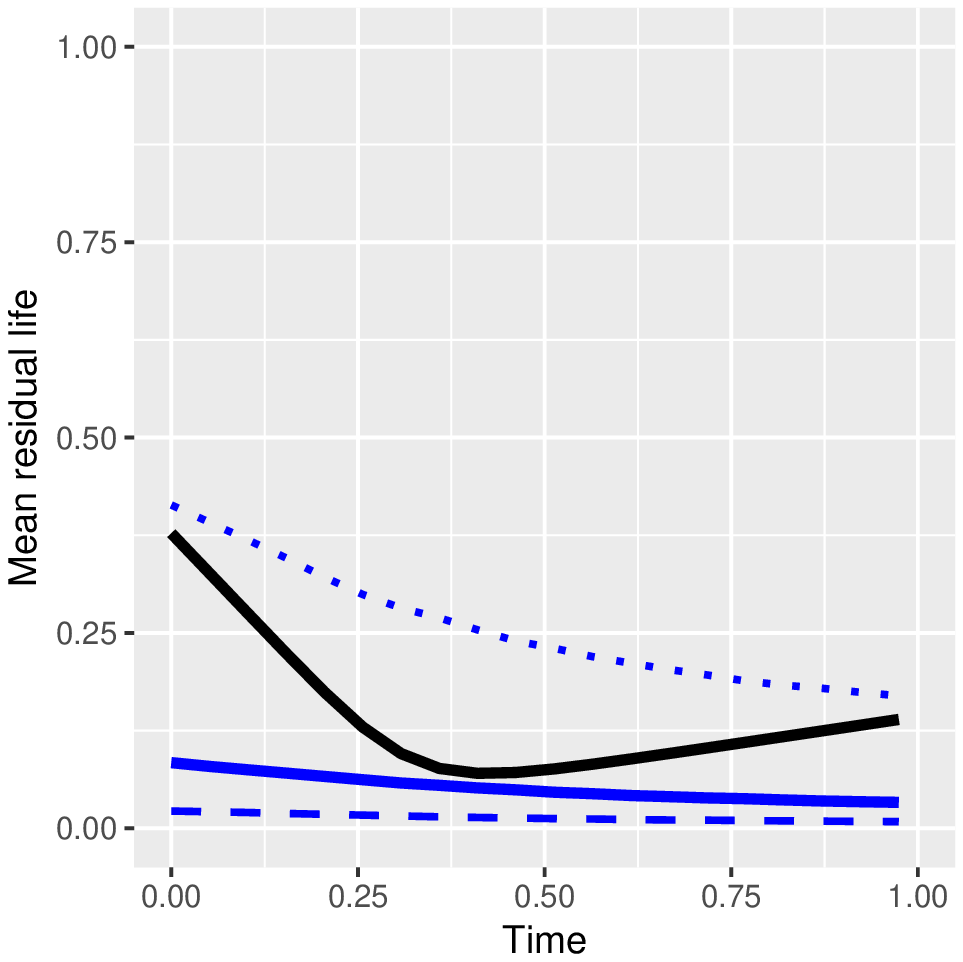}
	\includegraphics[width=4cm]{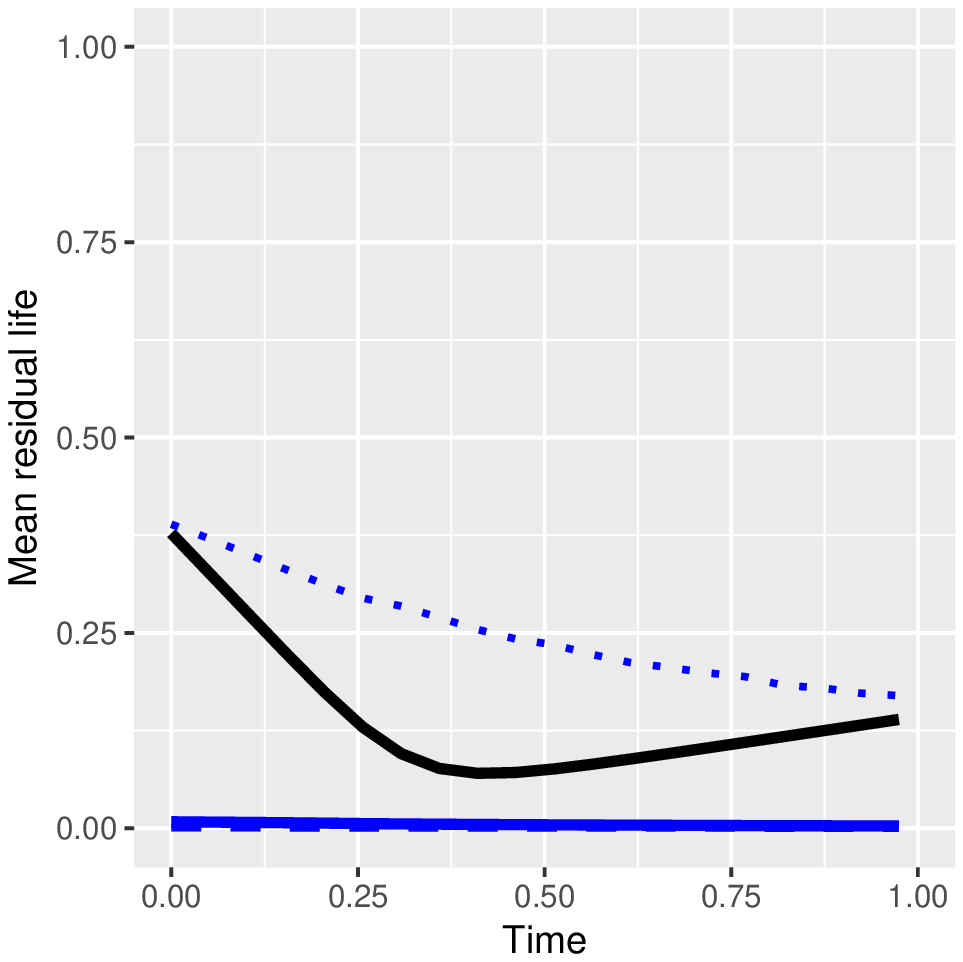}\\
	\includegraphics[width=4cm]{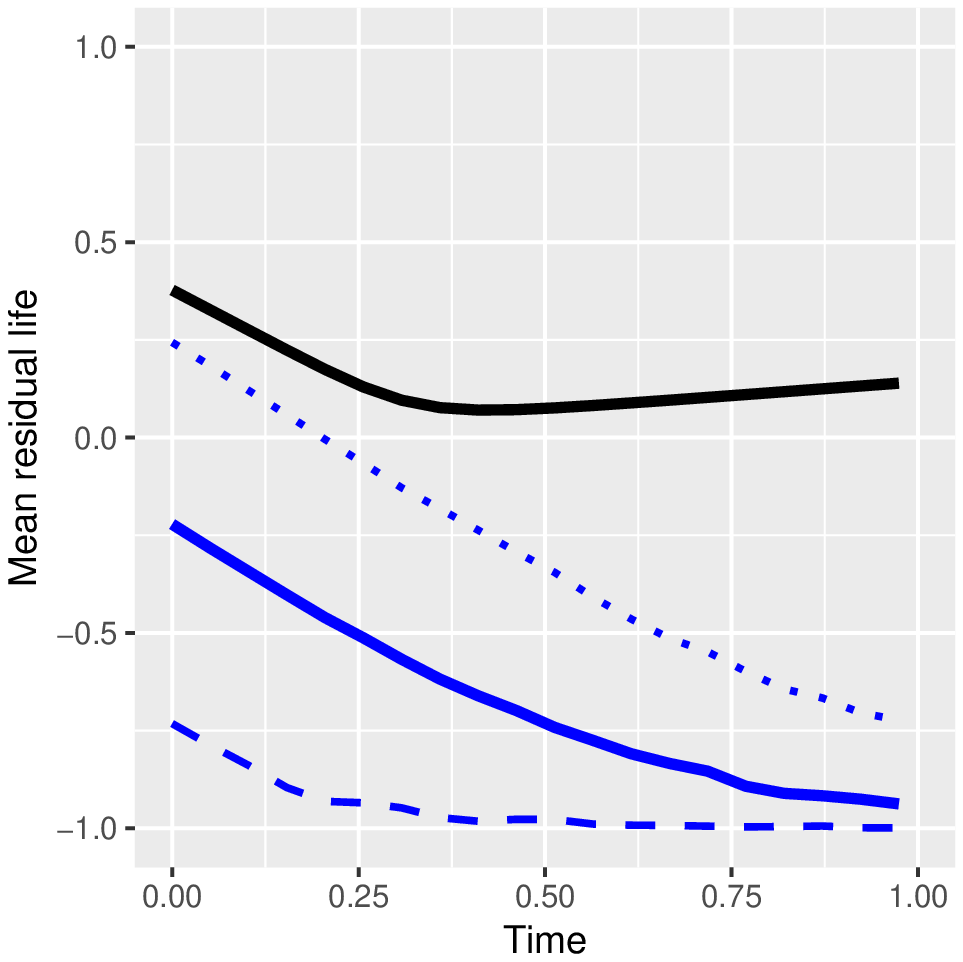}
	\includegraphics[width=4cm]{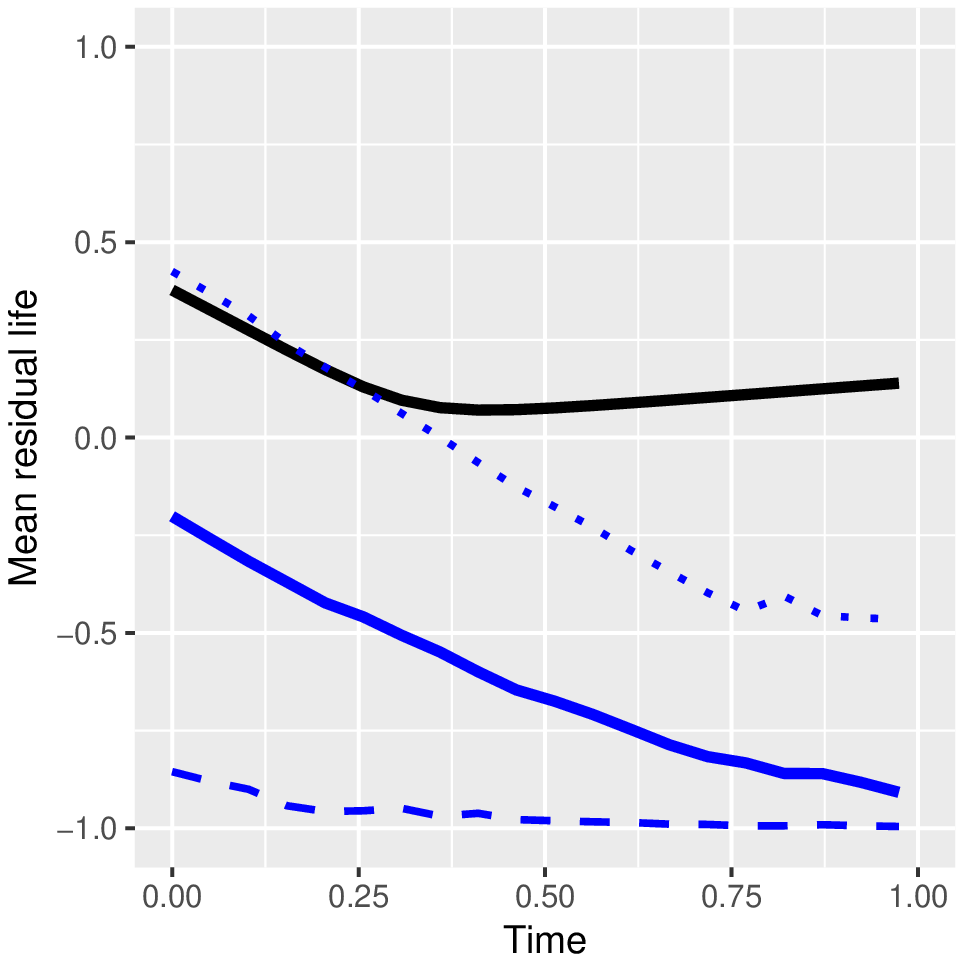}
	\includegraphics[width=4cm]{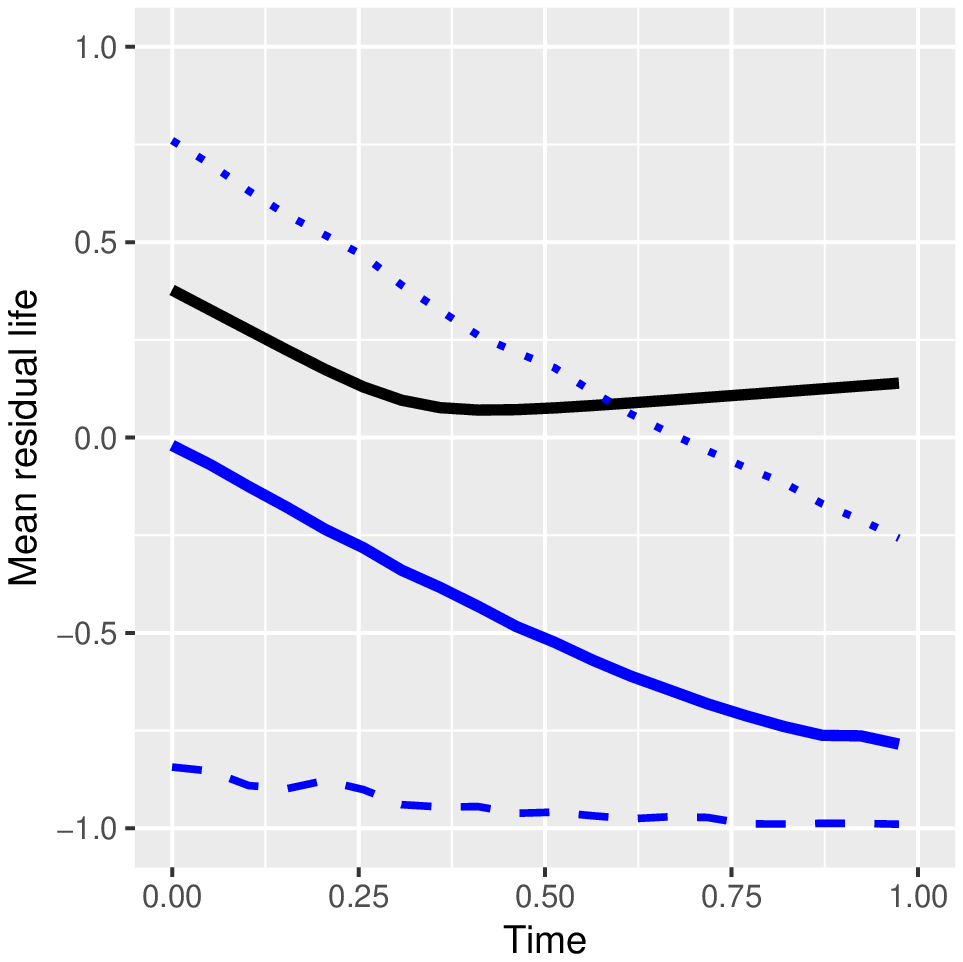}\\
	\caption{Mean residual life function estimation in Study
		2. Row 1: $m(t,\bb\trans\x)$ as a
		function of $\bb\trans\x$
		at $t=0.7$. Row 2 to Row 5: $m(t,\bb\trans\x)$ as a 
		function of $t$
		at $\bb\trans\x=-1$ from method ``semiparametric'', 
		``PM1'', ``PM2", 
		``additive".
		Left to right
		columns: no censoring; 20\% censoring rate; 40\%
		censoring rate. Black
		line: True $m(t,\bb\trans\x)$; Blue line: Median of $\wh
		m(t,\bb\trans\x)$;
		Blue dashed line: 2.5\% empirical percentile curve;
		Blue dotted line: 97.5\% empirical percentile curve.
	}
	\label{fig:simu2curve}
\end{figure}

\newpage

\newpage
\begin{figure}[H]
	\centering
	\includegraphics[width=5cm]{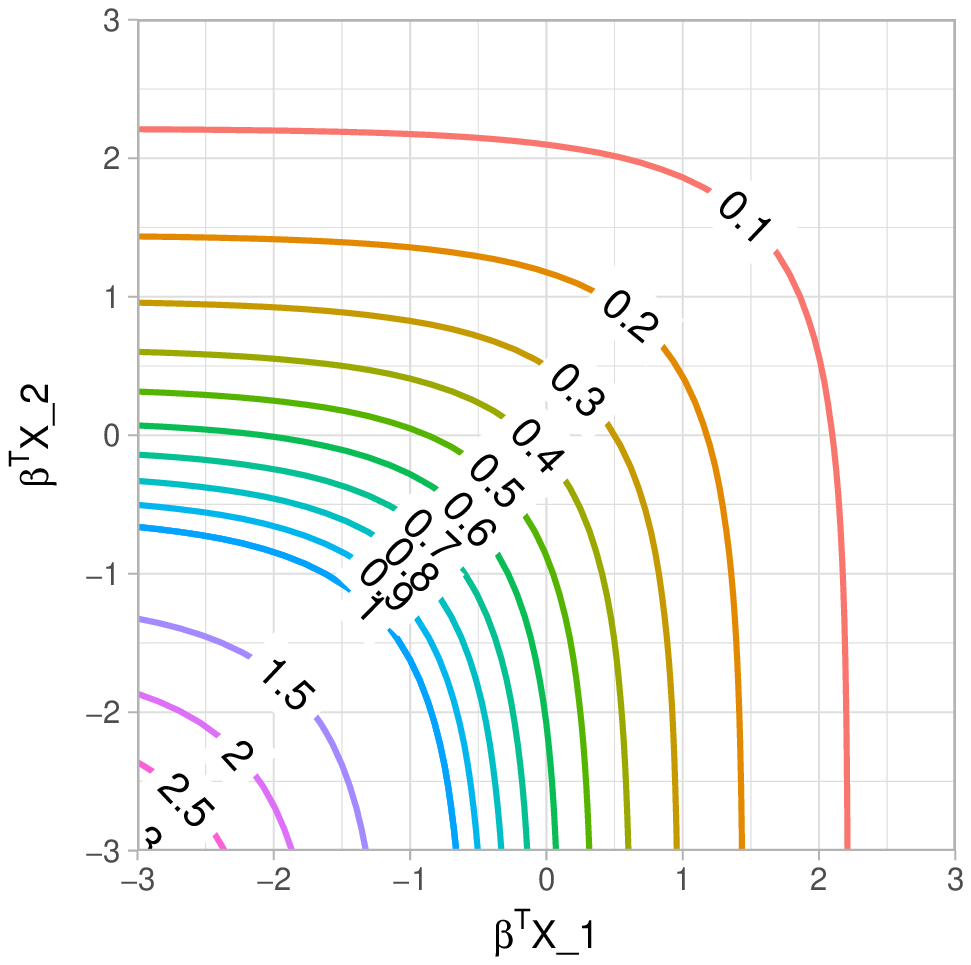}
	\includegraphics[width=5cm]{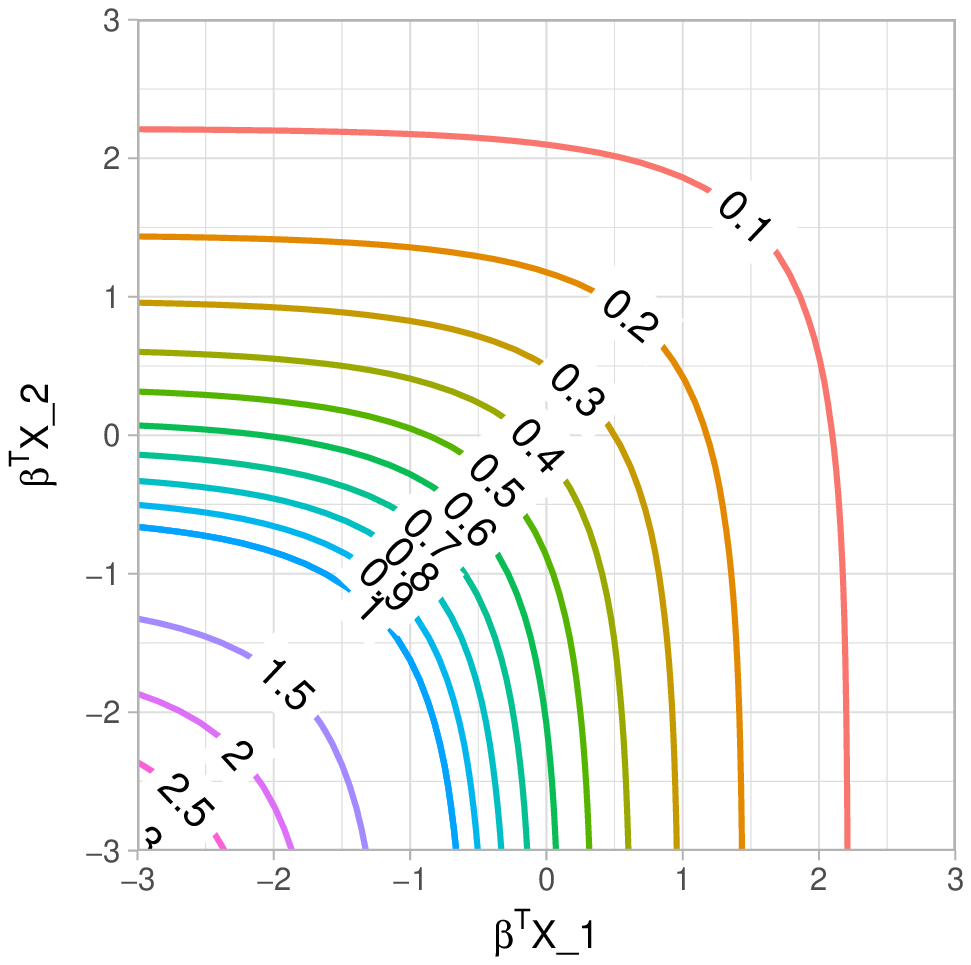}
	\includegraphics[width=5cm]{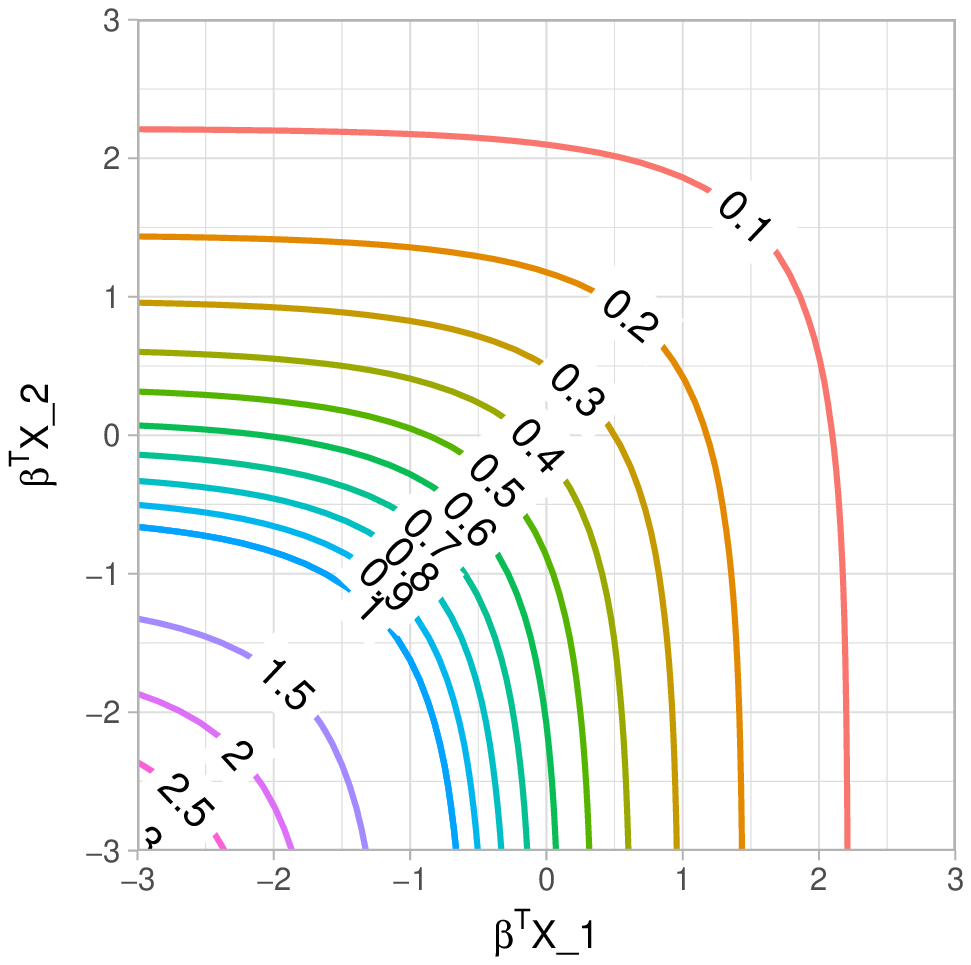}\\
	\includegraphics[width=5cm]{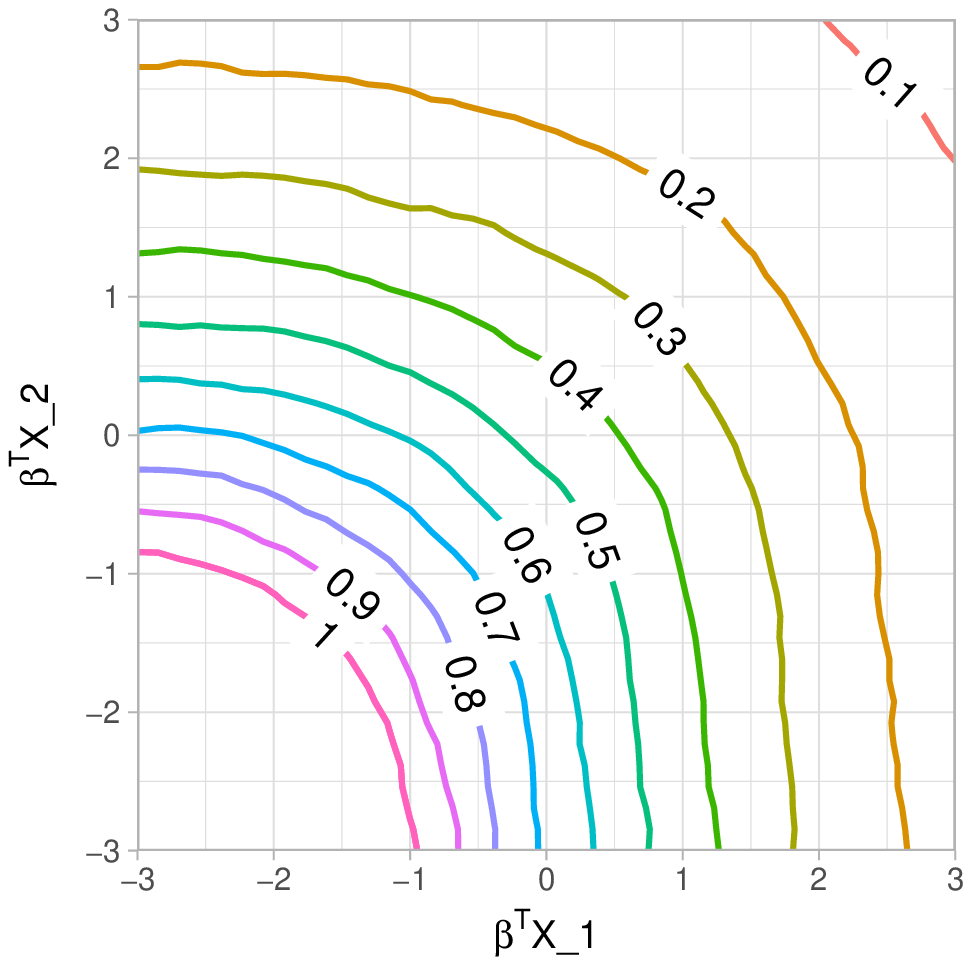}
	\includegraphics[width=5cm]{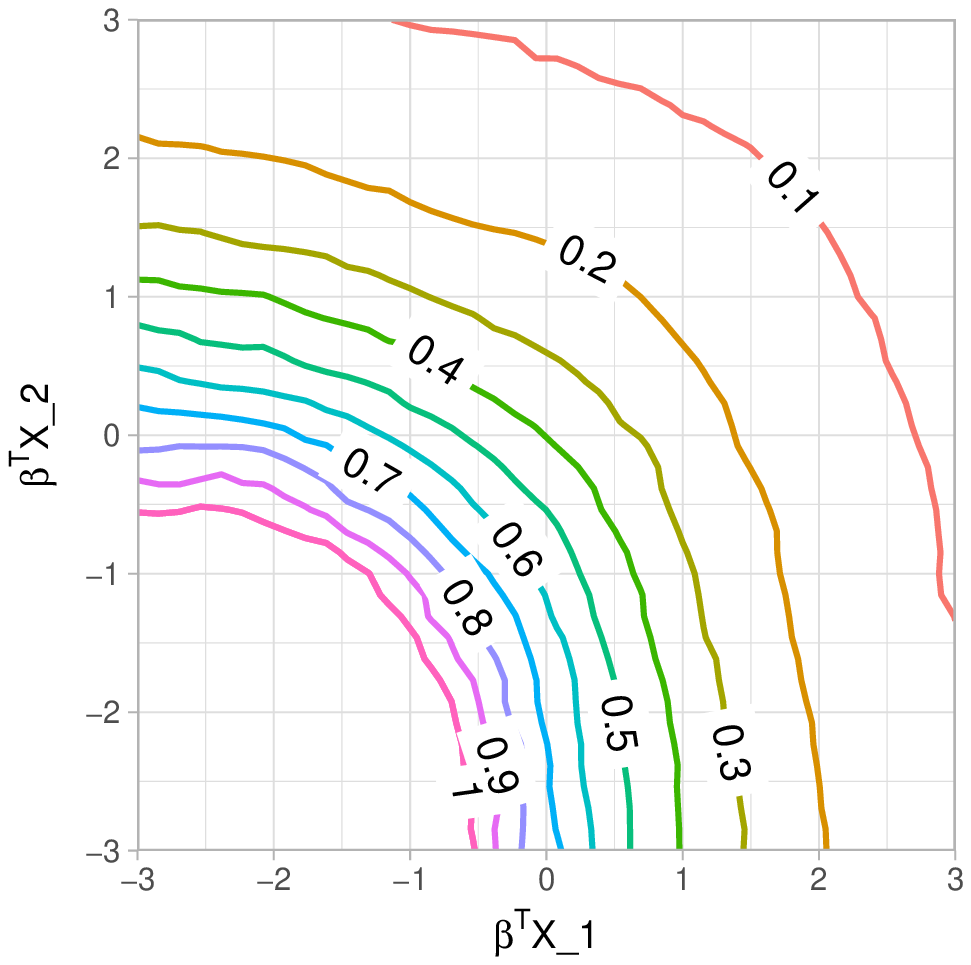}
	\includegraphics[width=5cm]{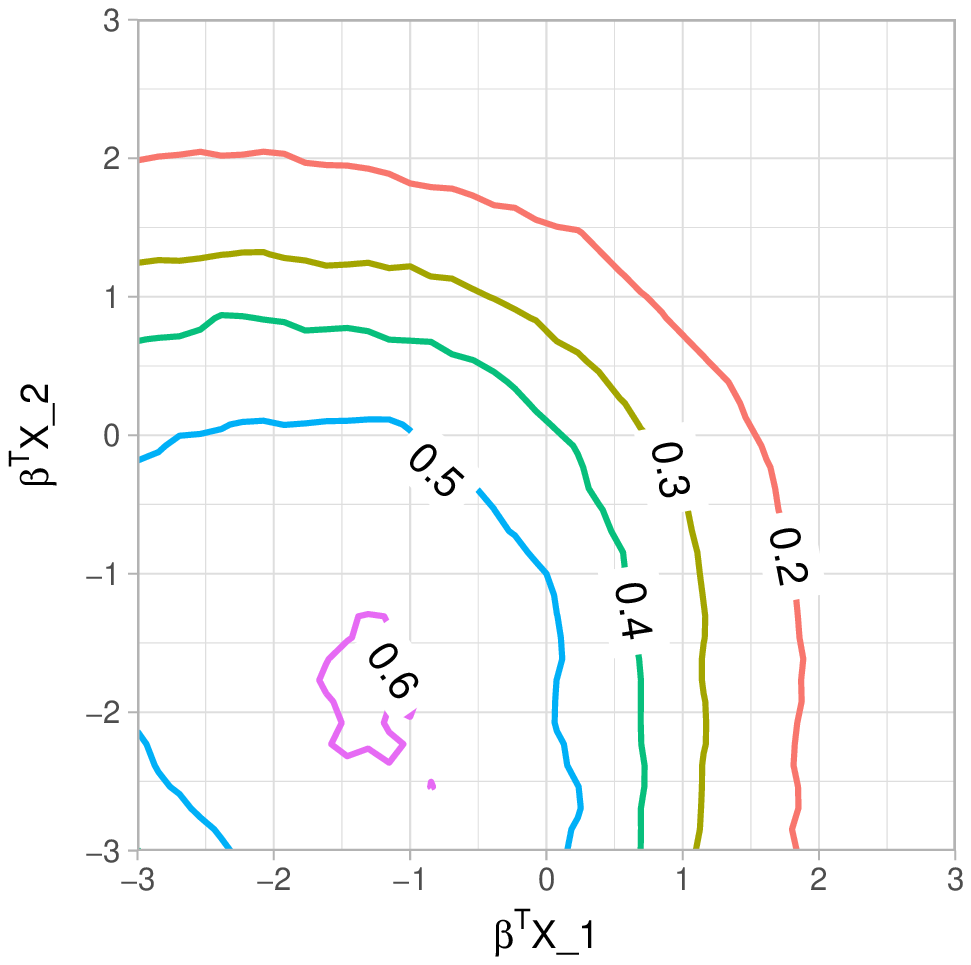}\\
	\includegraphics[width=5cm]{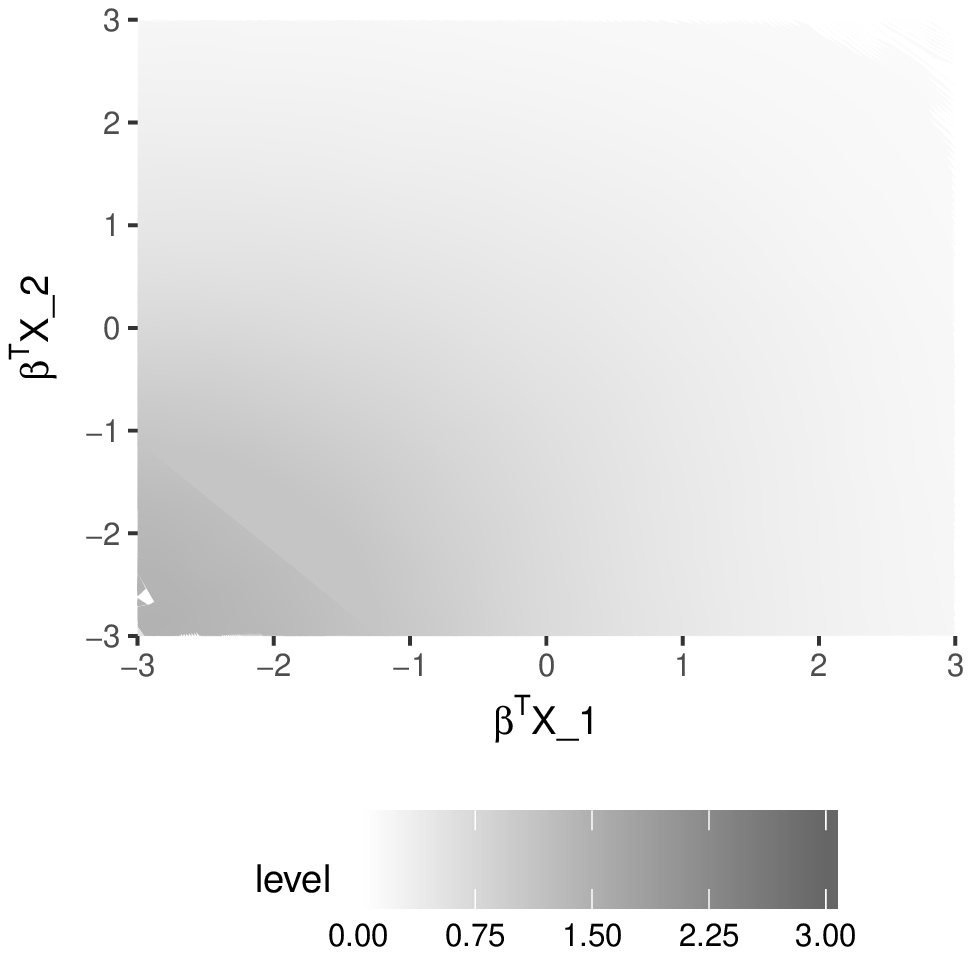}
	\includegraphics[width=5cm]{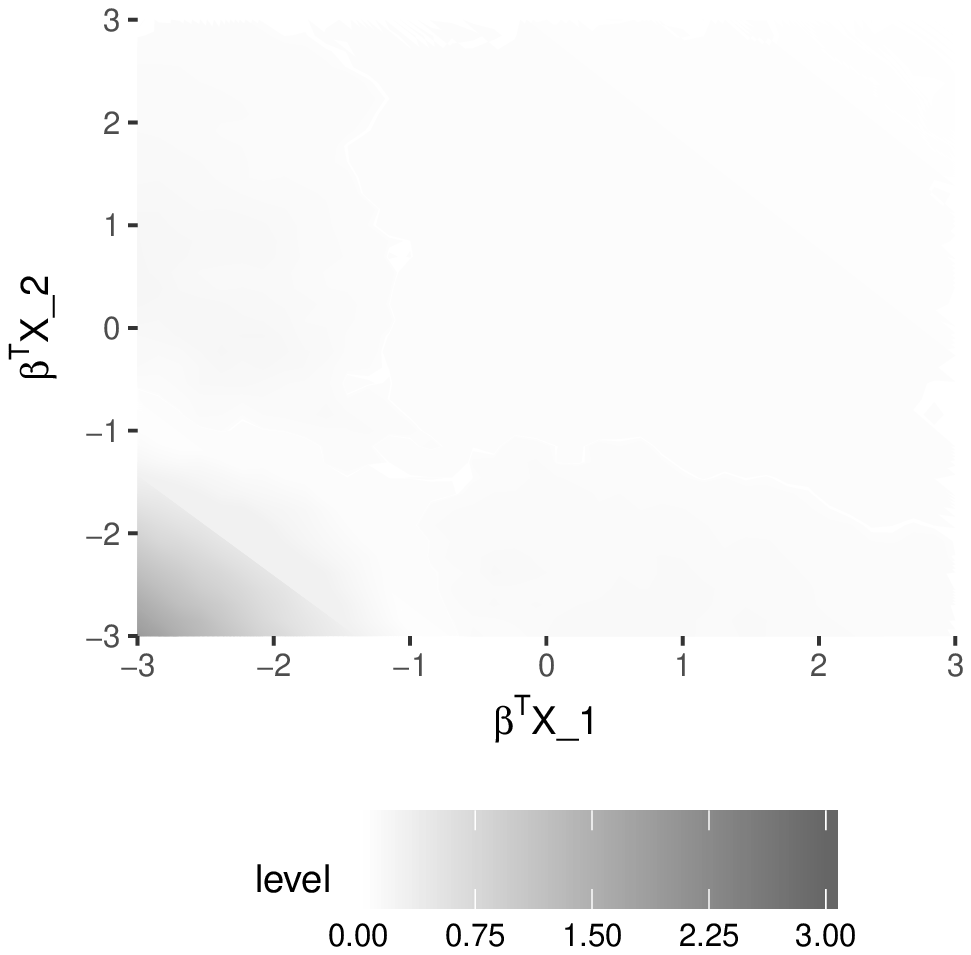}
	\includegraphics[width=5cm]{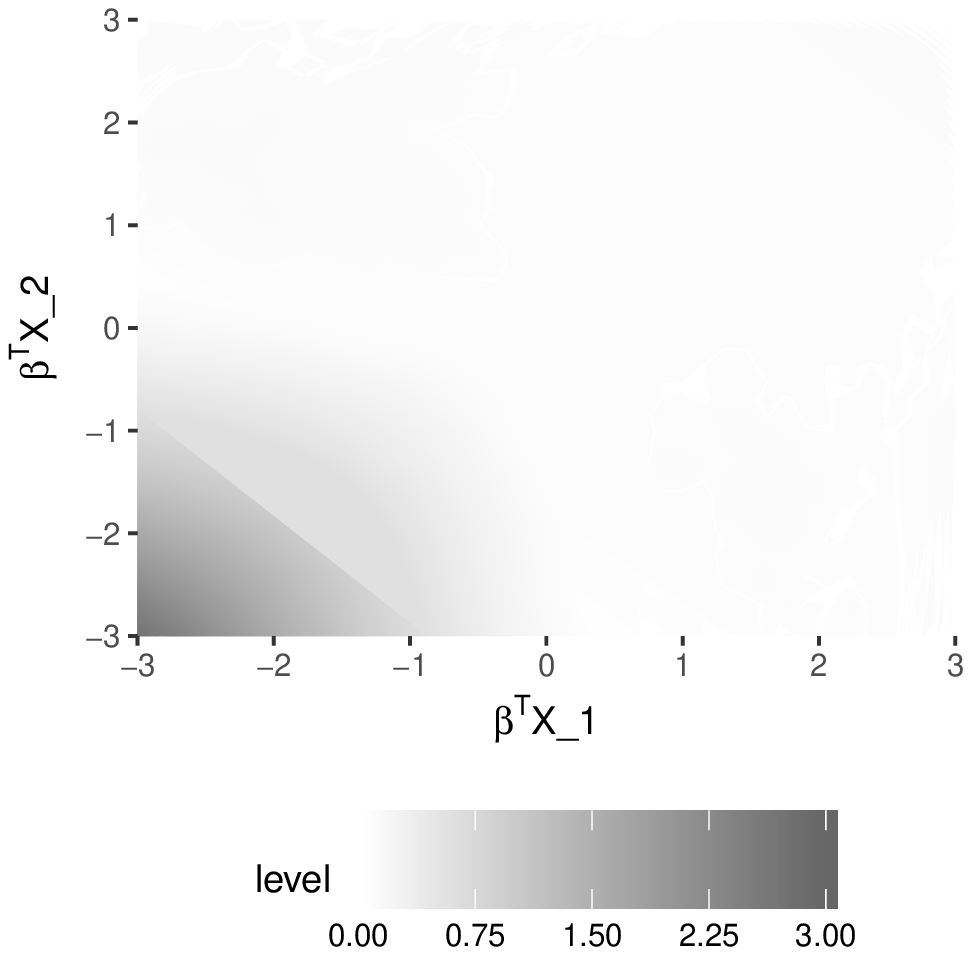}	
	\caption{Performance of the semiparametric method on mean
		residual life
		function at fixed $t$ of study 3.
		First row: contour plot of true $m(t=1,\bb\trans\X)$;
		Second
		row: contour plot of
		averaged $\overline{\wh m}(t=1,\bb\trans\X)$
		over 1000 simulations;
		Third row: contour plot of  $|\overline{\wh
			m}(t=1,\bb\trans\X)-m(t=1,\bb\trans\X)|$. Left to
		right
		columns: no censoring; 20\% censoring rate; 40\%
		censoring rate.
	}
	\label{fig:simu3contour}
\end{figure}
\newpage
\begin{figure}[H]
	\centering
	\includegraphics[width=5cm]{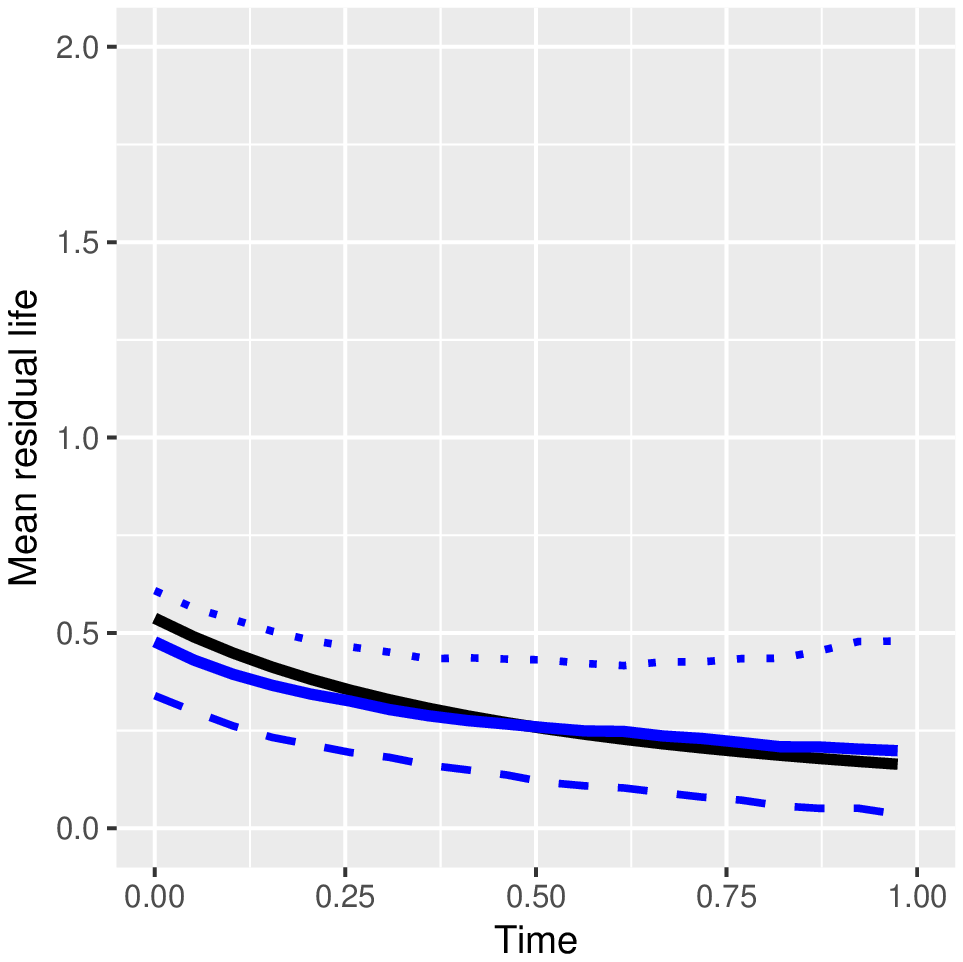}
	\includegraphics[width=5cm]{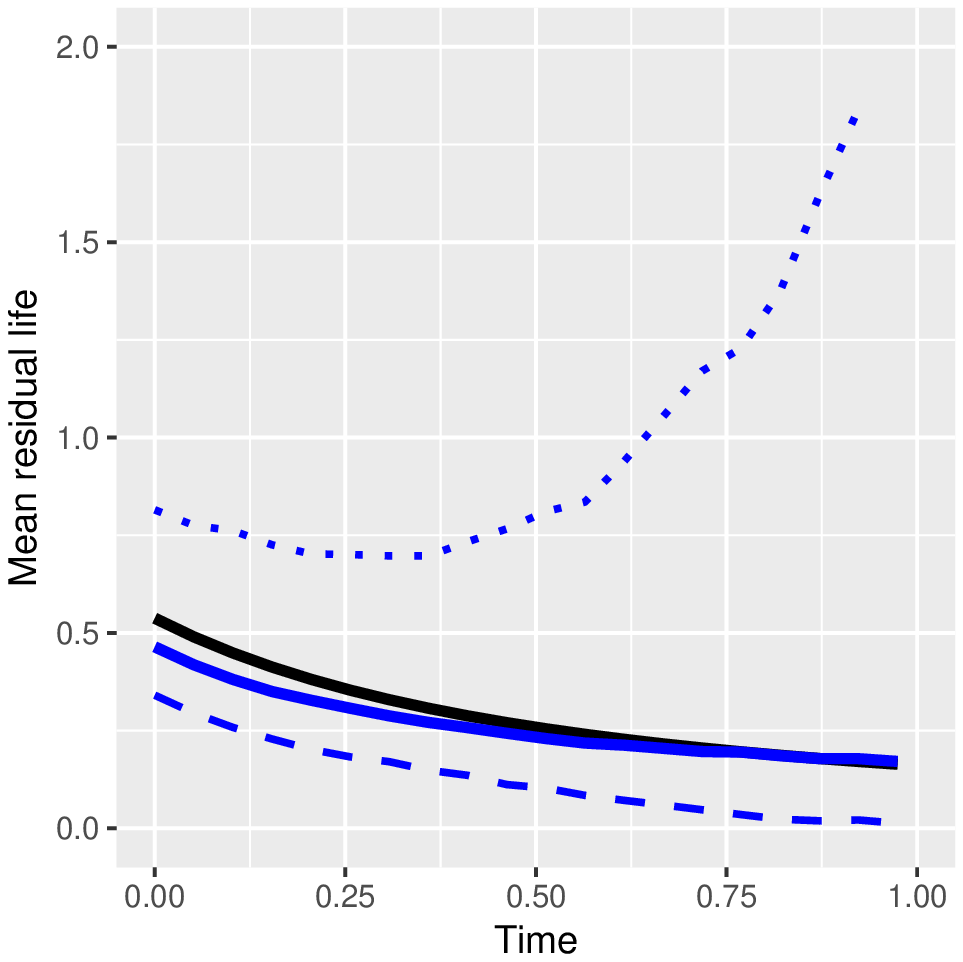}
	\includegraphics[width=5cm]{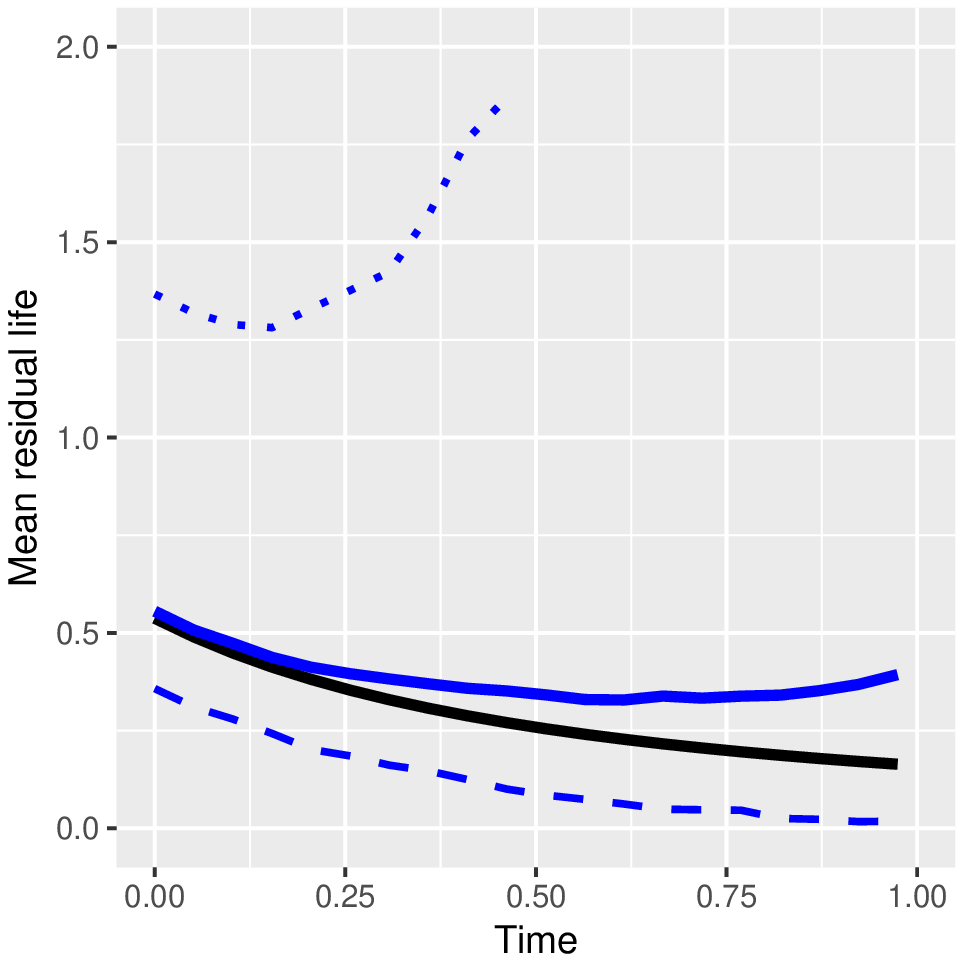}
	\includegraphics[width=5cm]{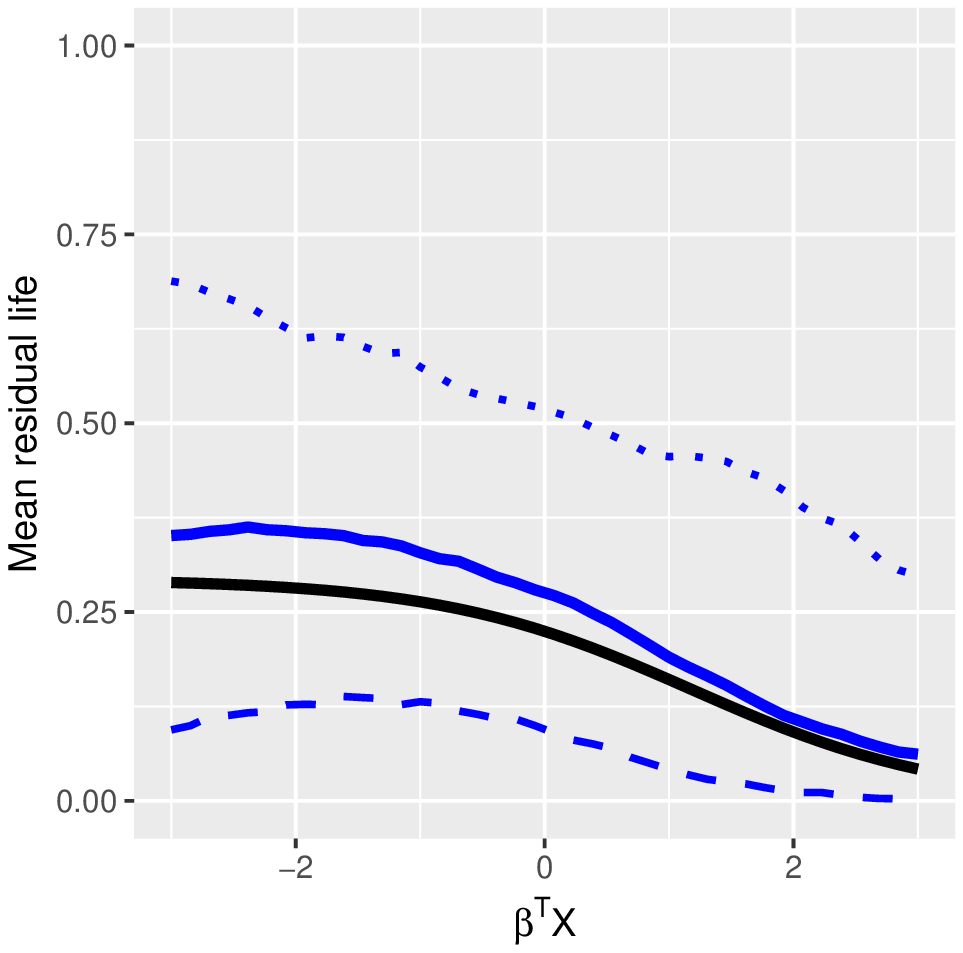}
	\includegraphics[width=5cm]{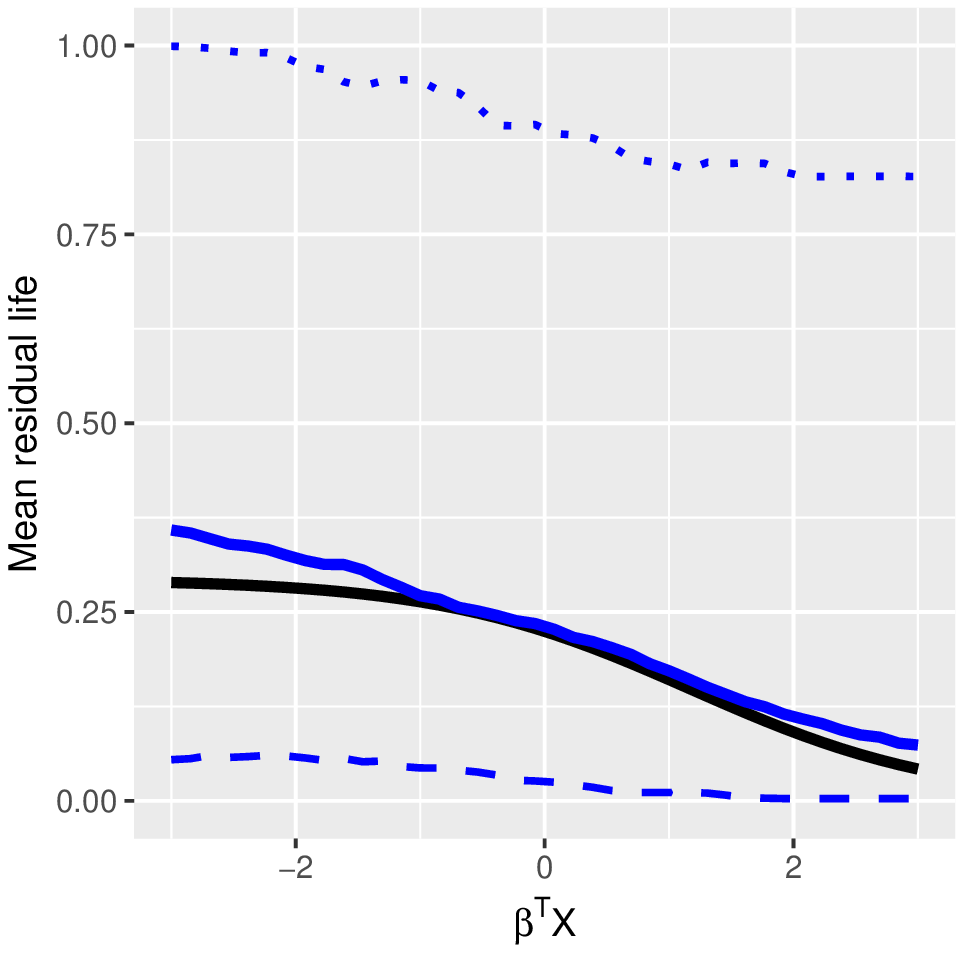}
	\includegraphics[width=5cm]{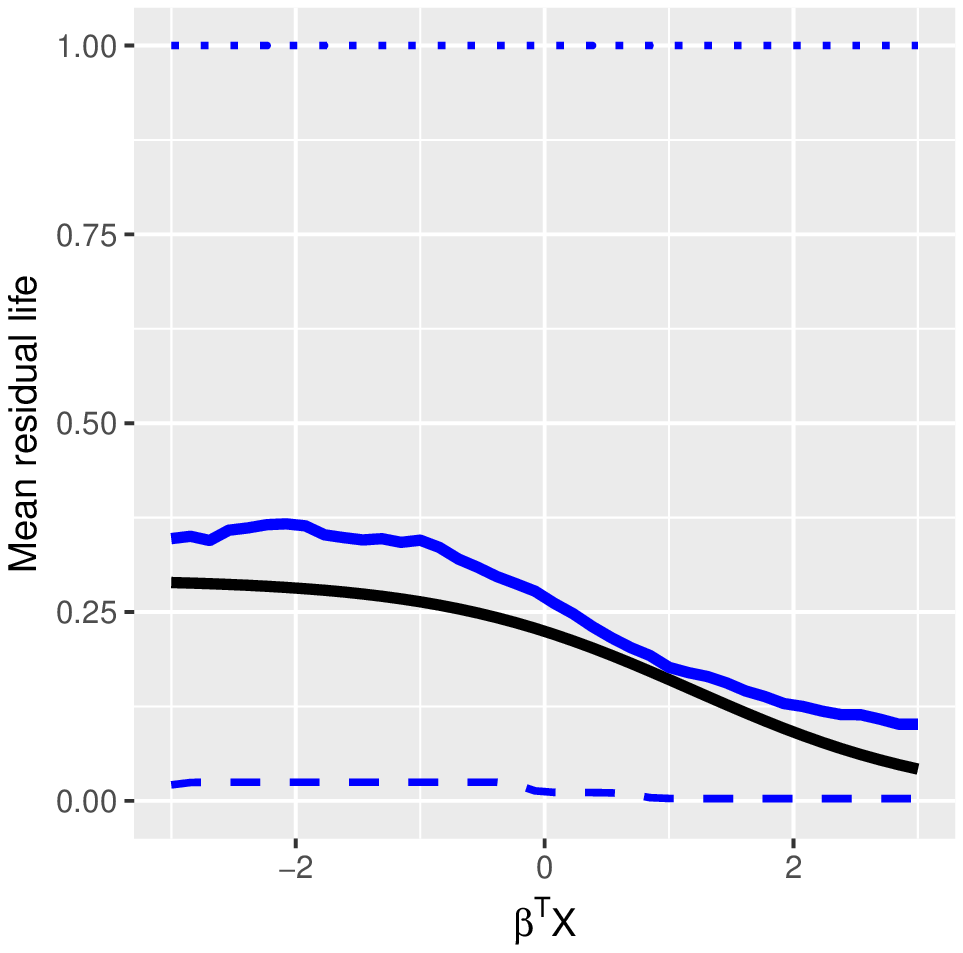}
	\includegraphics[width=5cm]{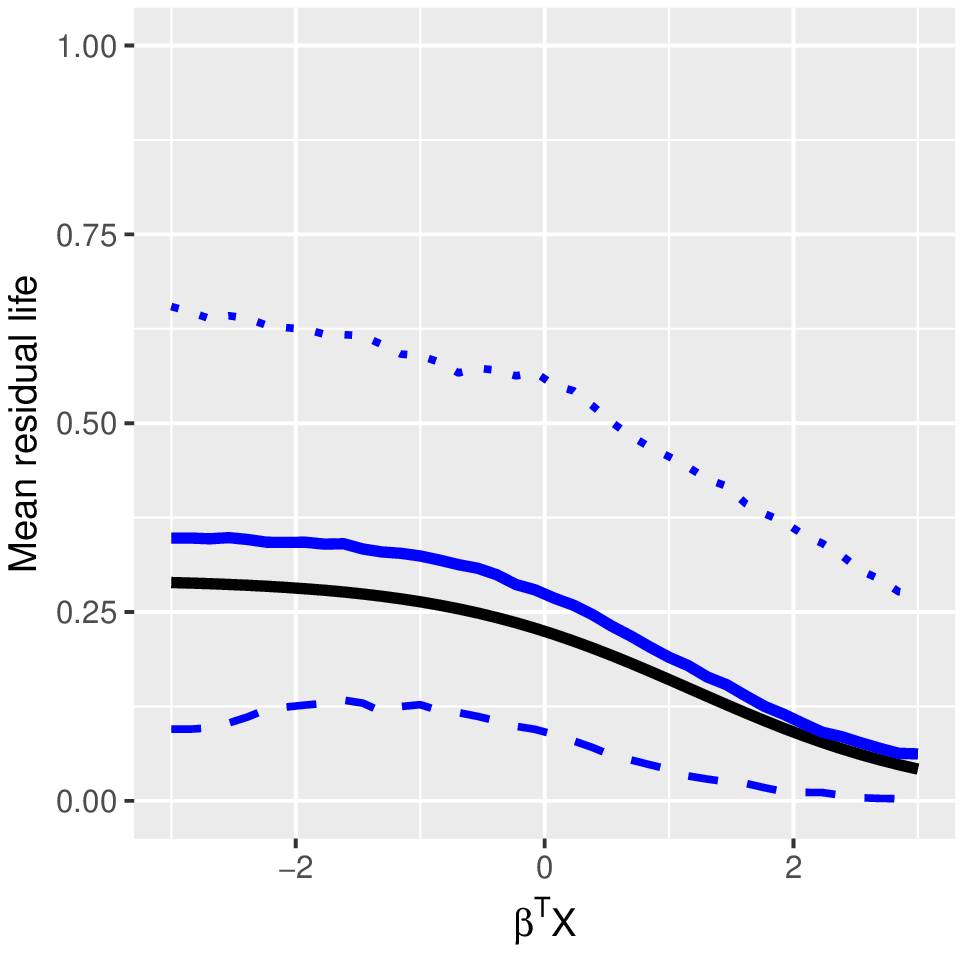}
	\includegraphics[width=5cm]{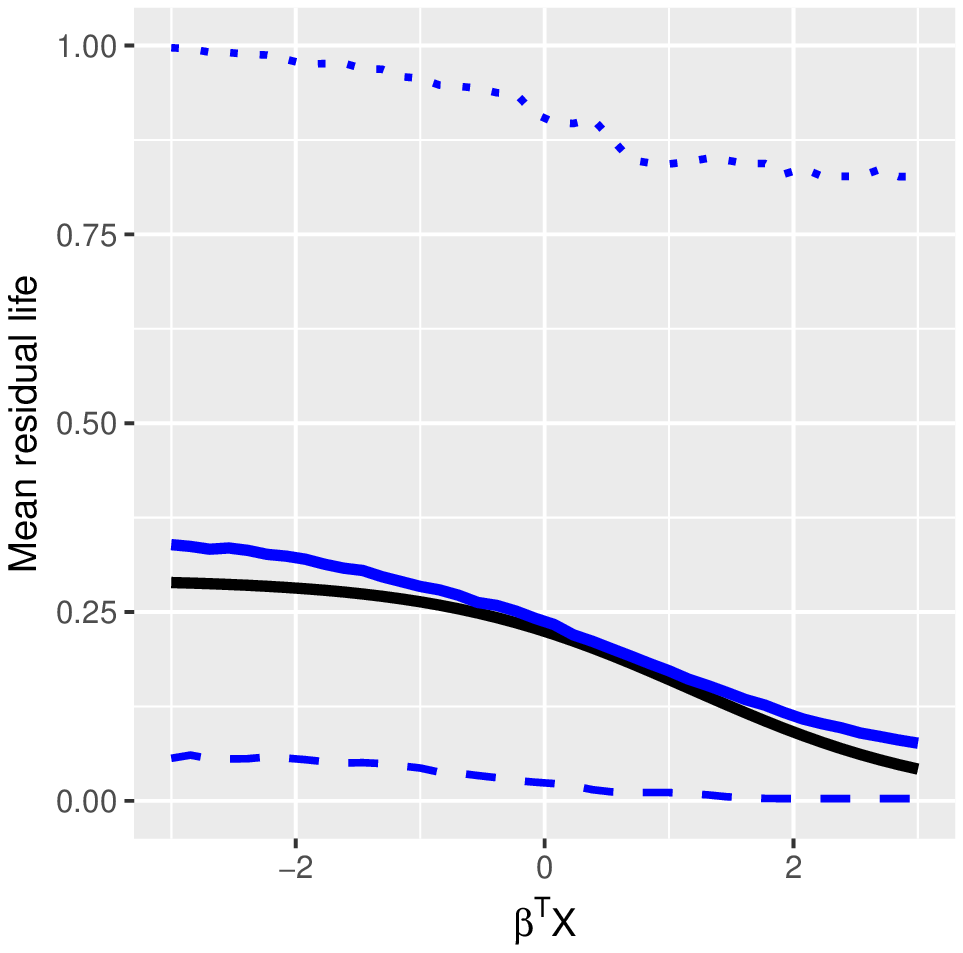}
	\includegraphics[width=5cm]{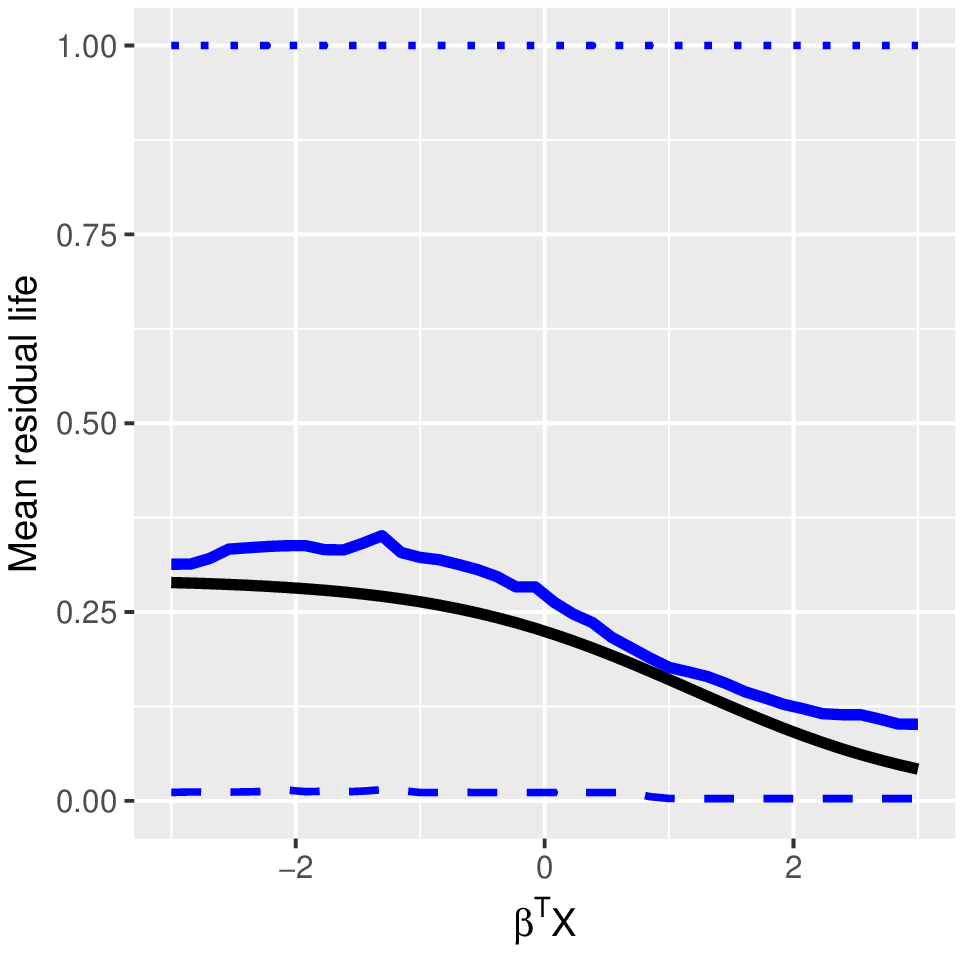}
	\caption{Mean residual life function estimation in Study
		3. Row 1: $m(t,\bb\trans\x)$ as a function of $t$
		at $\bb\trans\x=[1,1]'$. Row 2: $m(t,\bb\trans\x)$ as a
		function of
		$(\bb\trans\x)_1$
		at $t=1$ and $(\bb\trans\x)_2=1$. Row 3:
		$m(t,\bb\trans\x)$ as a function of
		$(\bb\trans\x)_2$
		at $t=1$ and $(\bb\trans\x)_1=1$.
		Left to right
		columns: no censoring; 20\% censoring rate; 40\%
		censoring rate. Black
		line: True $m(t,\bb\trans\x)$; Blue line: Median of $\wh
		m(t,\bb\trans\x)$;
		Blue dashed line: 2.5\% empirical percentile curve;
		Blue dotted line: 97.5\% empirical percentile curve.
	}
	\label{fig:simu3curve}
\end{figure}

\newpage
\begin{figure}[H]
	\centering
	\includegraphics[width=5cm]{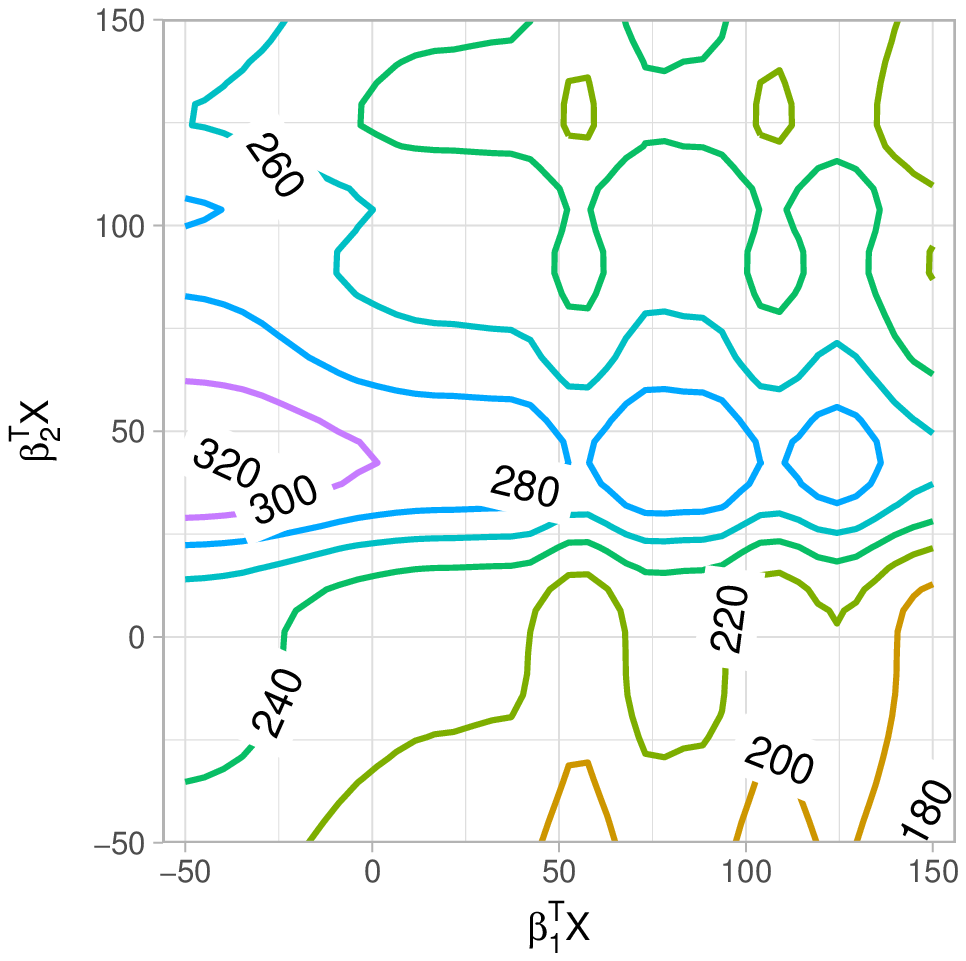}
	\includegraphics[width=5cm]{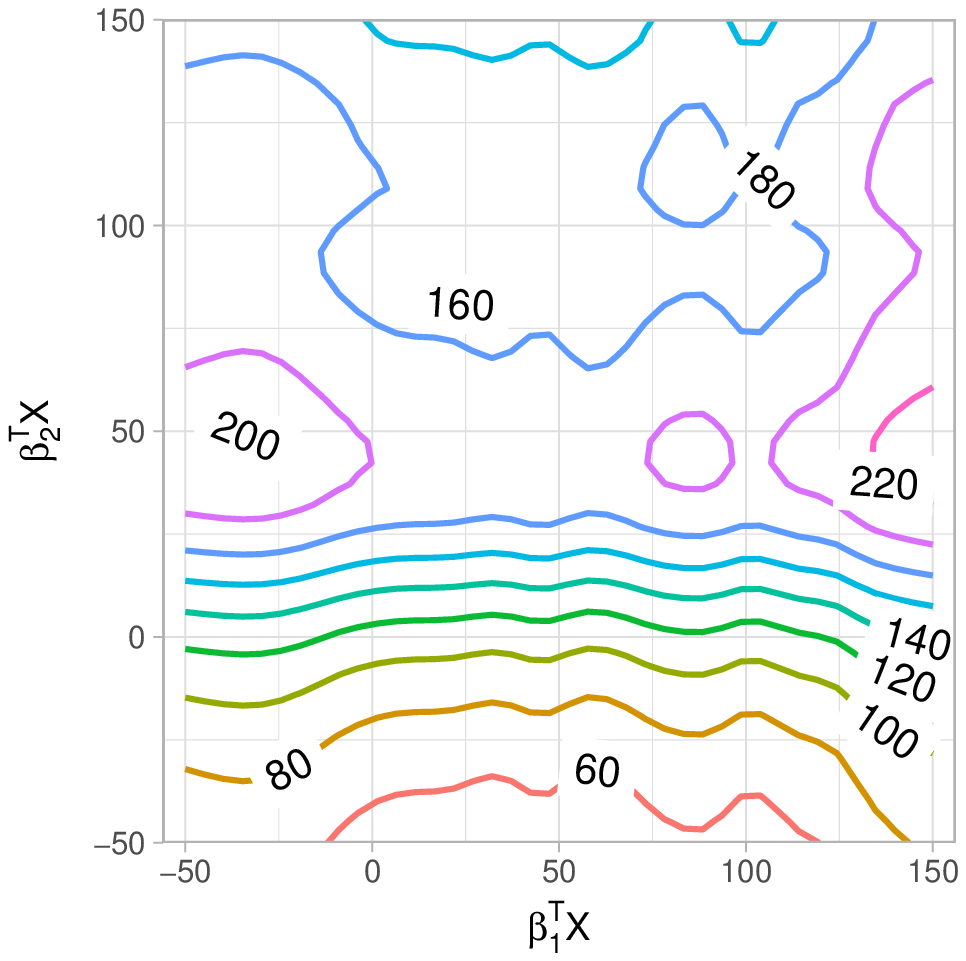}
	\includegraphics[width=5cm]{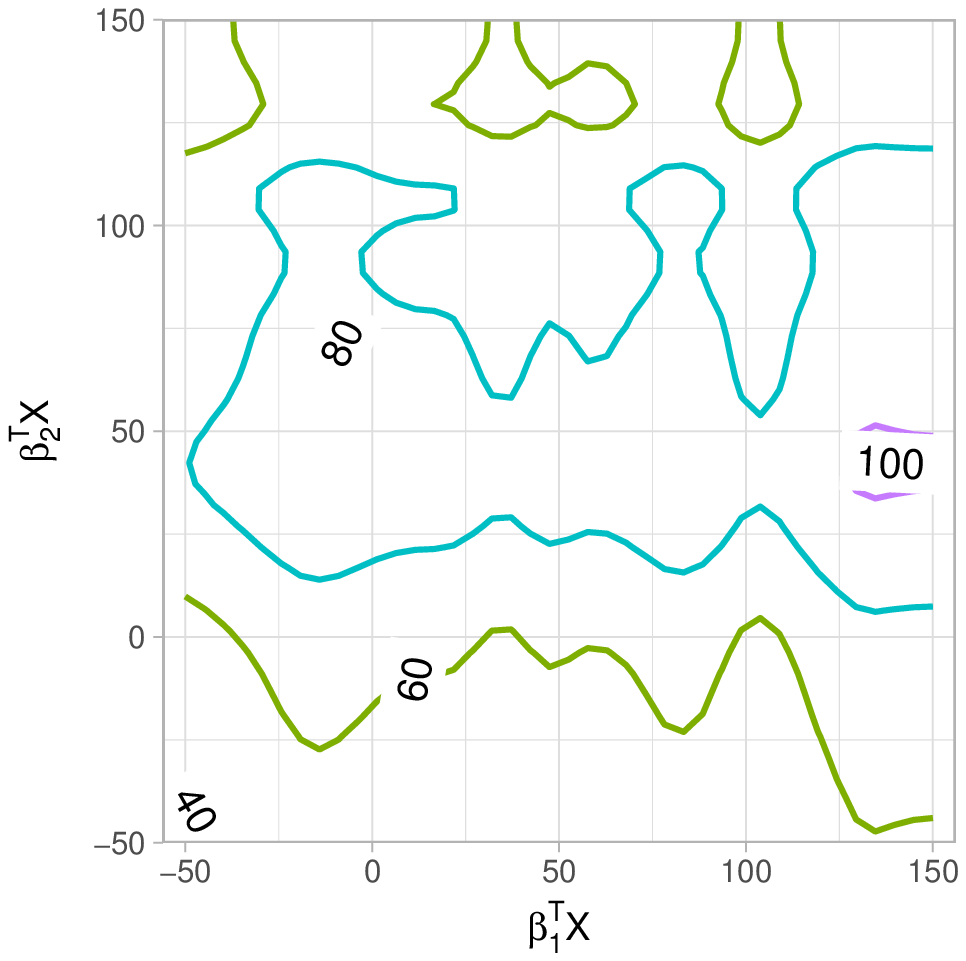}\\
	\includegraphics[width=5cm]{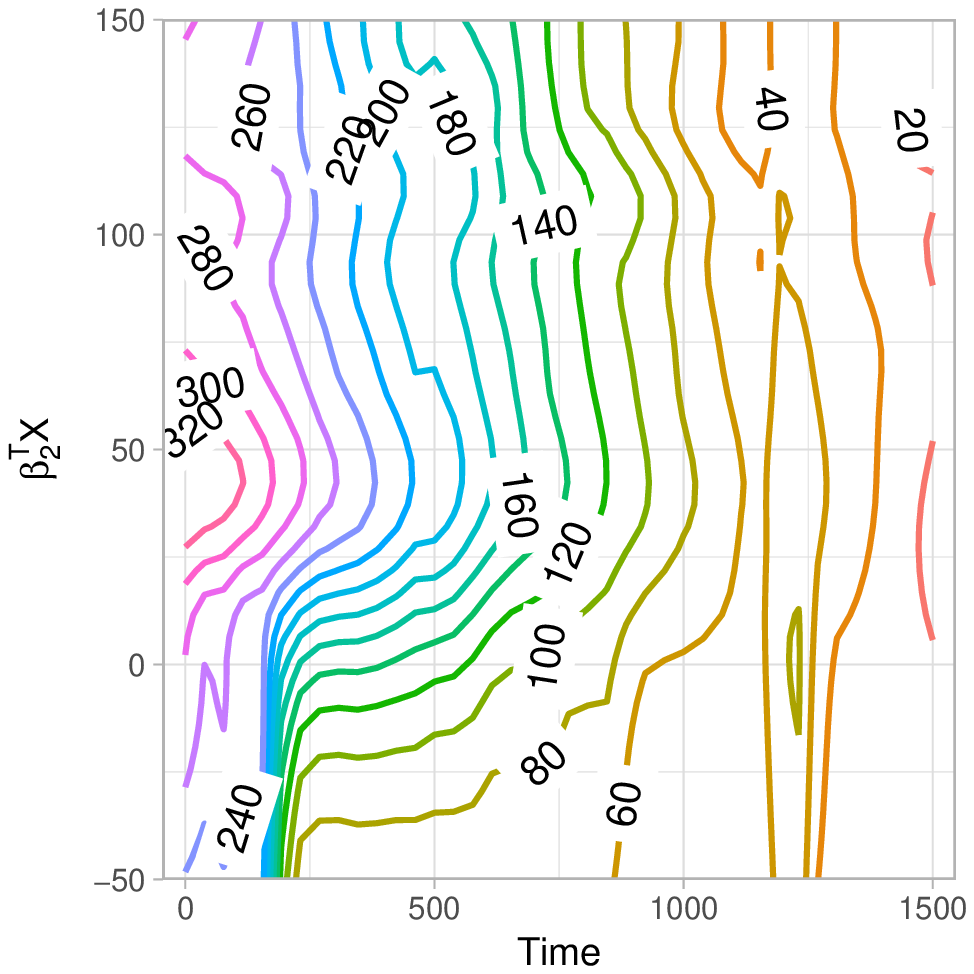}
	\includegraphics[width=5cm]{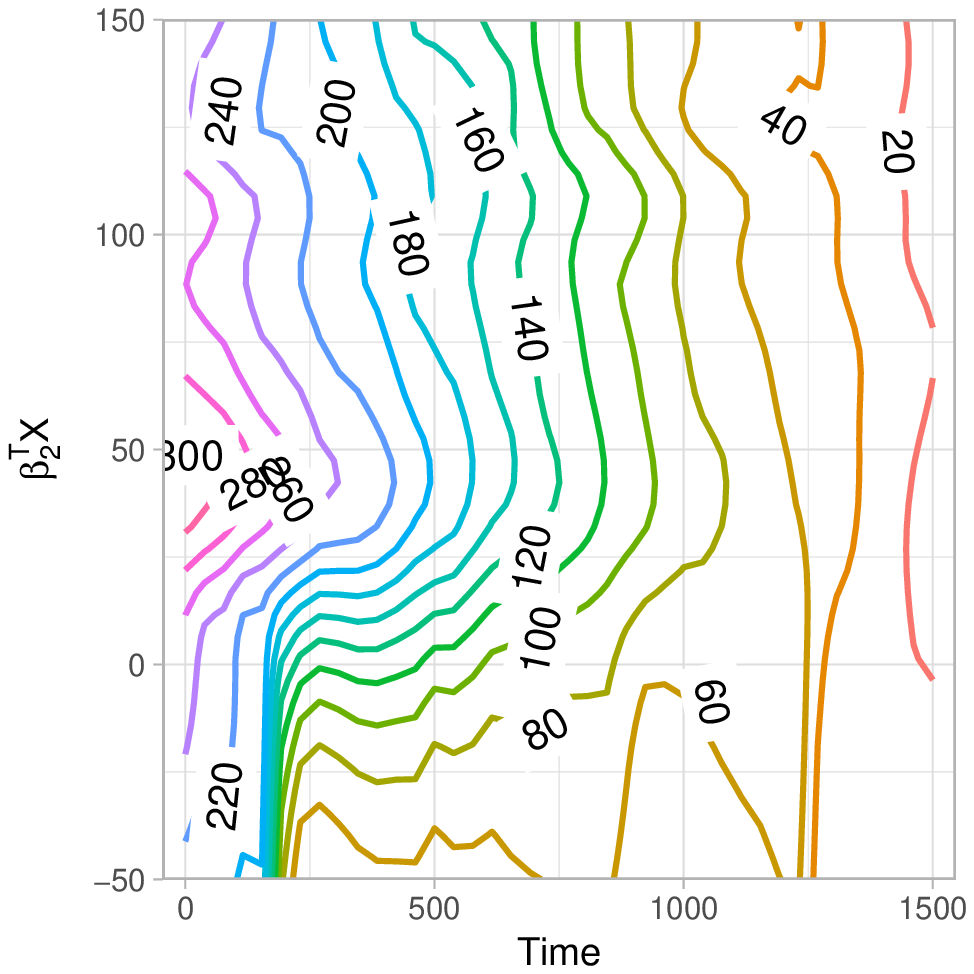}
	\includegraphics[width=5cm]{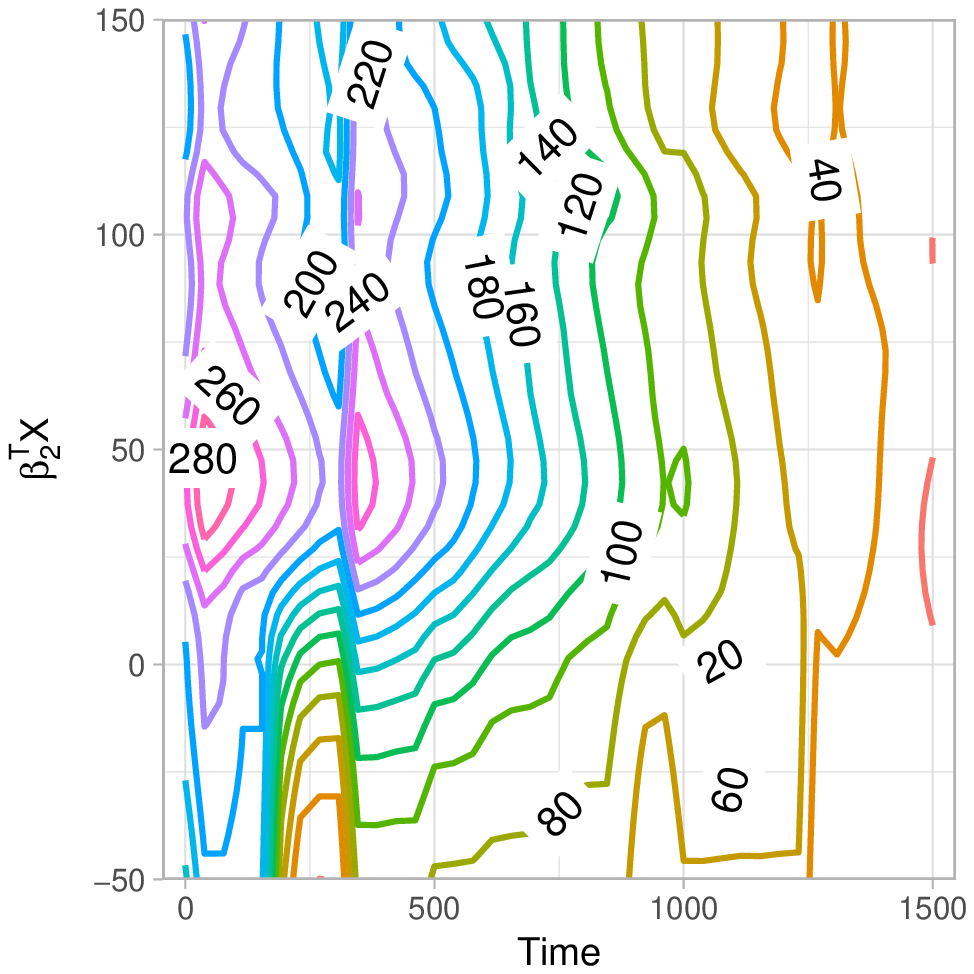}\\
	\includegraphics[width=5cm]{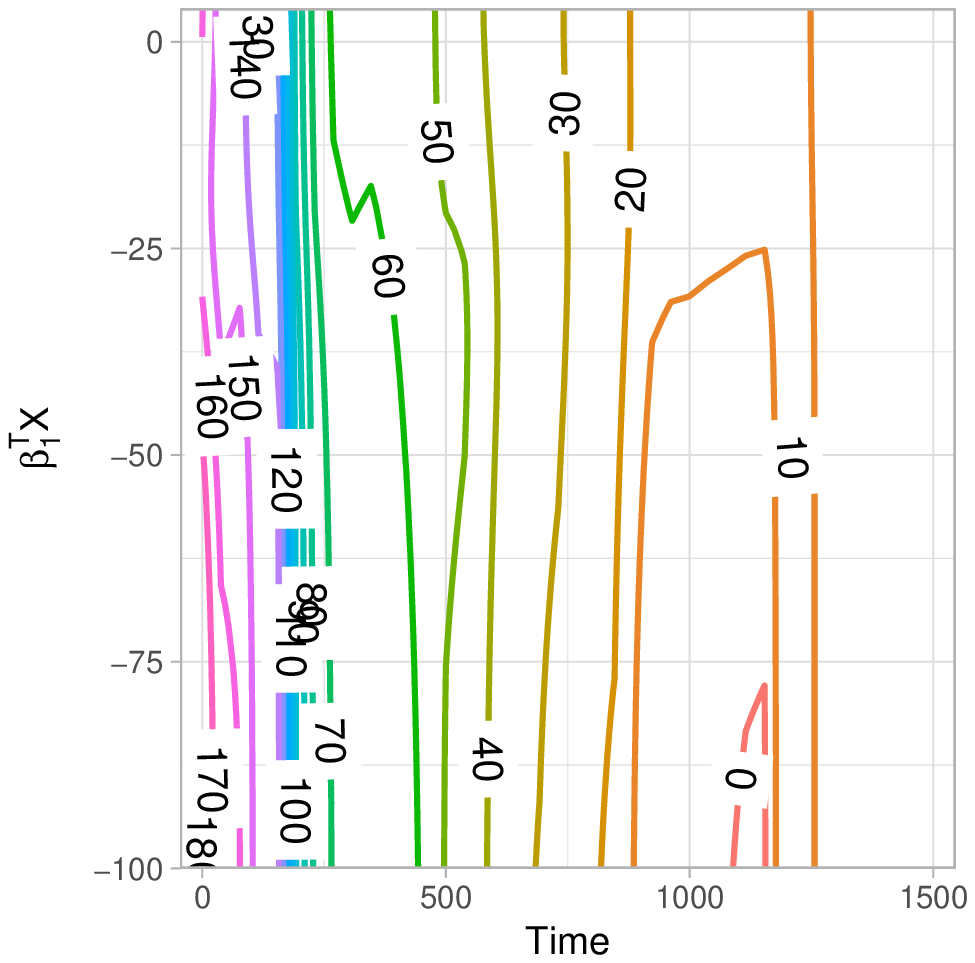}
	\includegraphics[width=5cm]{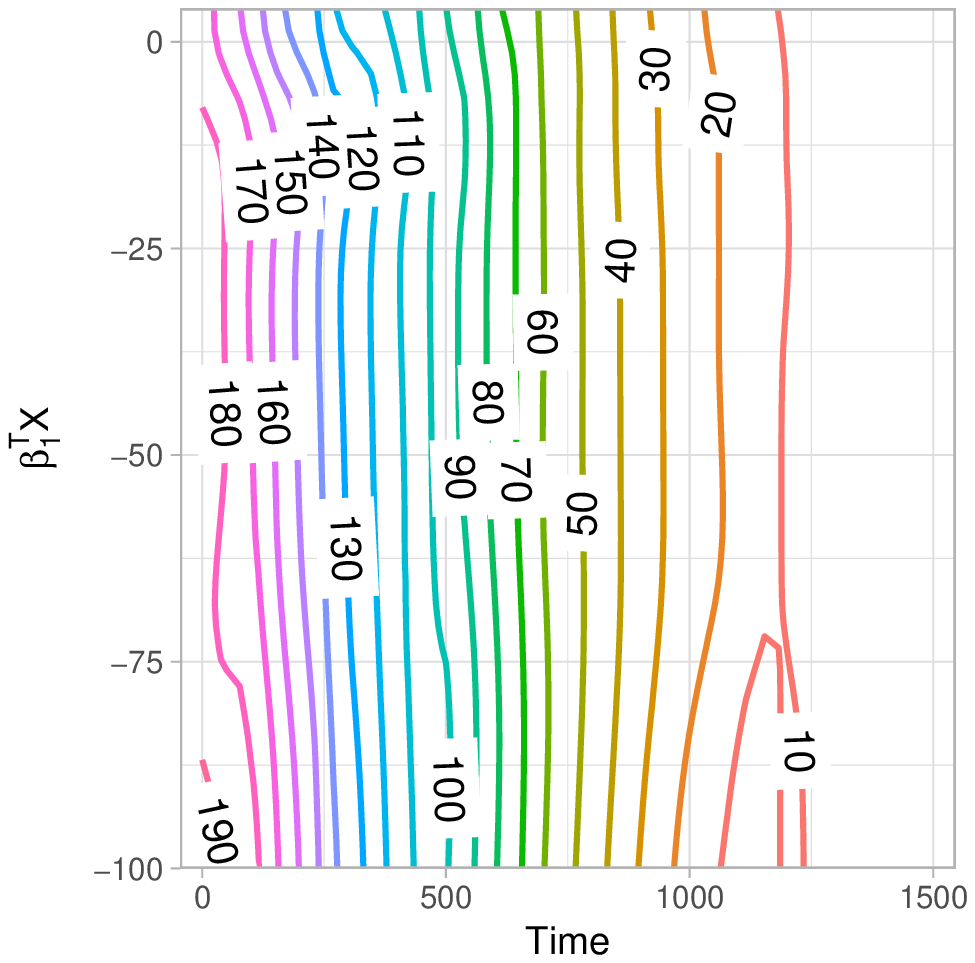}
	\includegraphics[width=5cm]{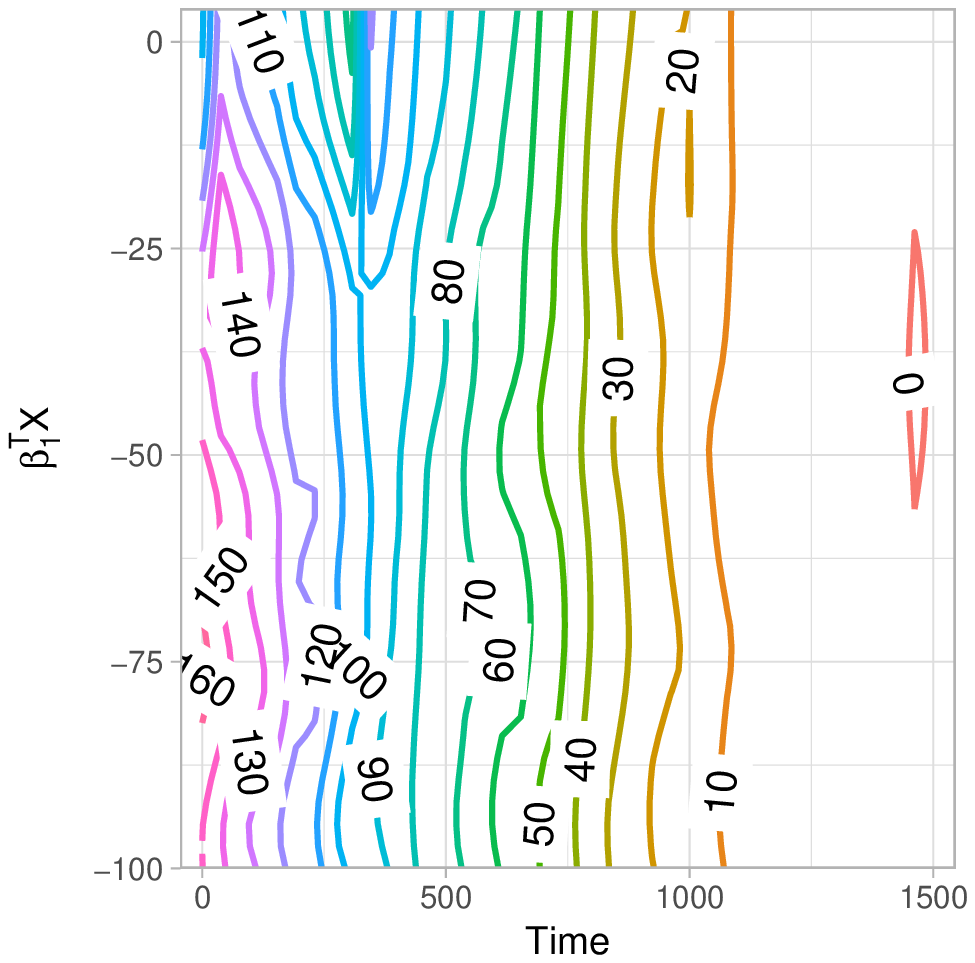}
	\caption{Contour plot of mean residual life difference $\wh
		m_{\rm treat}(t,\bb_{\rm treat}\trans\X)-\wh m_{\rm
			wait}(t,\bb_{\rm
			wait}\trans\X)$. First
		row:
		fix $t$ at 100, 500 and 1000 from left to right. Second
		row:
		fix $\bb_{\rm treat}\trans\X$ at -40, 50 and 100 from
		left to
		right.
		Third row: fix $\bb_{\rm wait}\trans\X$ at -40, 50 and
		100 from
		left
		to right.}
	\label{fig:appcontour2}
\end{figure}
\newpage

\newpage
\bibliographystyle{ECA_jasa}
\bibliography{meanres}

\end{document}